\magnification=1000
\hsize=11.7cm
\vsize=18.9cm
\lineskip2pt \lineskiplimit2pt
\nopagenumbers

%%%%%%%%%%%%%%%%%%%%%%%%%%%
\hoffset=-1truein
\voffset=-1truein
%%%%%%%%%%%%%%%%%%%%%%%%%%%

\advance\voffset by 4truecm
\advance\hoffset by 4.5truecm

\newif\ifentete

%%%% ENTETE COURANTE
\headline{\ifentete\ifodd	\count0 
      \rlap{\head}\hfill\tenrm\llap{\the\count0}\relax
    \else
        \tenrm\rlap{\the\count0}\hfill\llap{\head} \relax
    \fi\else
\global\entetetrue\fi}

\def\entete#1{\entetefalse\gdef\head{#1}}
\entete{}

\input amssym.def
\input amssym.tex

\def\-{\hbox{-}}
\def\.{{\cdot}}
\def\O{{\cal O}}
\def\K{{\cal K}}
\def\F{{\cal F}}
\def\E{{\cal E}}
\def\L{{\cal L}}

\def\P{{\cal P}}

\def\G{{\cal G}}
\def\T{{\cal T}}

\def\U{{\cal U}}

\def\Z{{\cal Z}}
\def\C{{\cal C}}

\def\ch{\frak c\frak h}
\def\ad{\frak a\frak c}

\def\Ab{\frak A\frak b}
\def\st{\frak s\frak t}

\def\Gr{\frak G\frak r}

\def\Fct{\frak F\frak c\frak t}

\def\Ker{\frak K\frak e\frak r}
\def\SC{\frak C\frak C}
\def\CC{\frak C\frak C}

\def\id{\frak i\frak d}
\def\int{\frak i\frak n\frak t}

\def\qq{\quad{\rm and}\quad}

\def\cqq{\quad ,\quad}
\def\ft{\frak f\frak g}

\def\qt{\frak q\frak t}

\def\mod{\frak m\frak o\frak d}
\def\res{\frak r\frak e\frak s}

\def\too{\longrightarrow}
\def\aut{\frak a\frak u\frak t}

\def\Loc{\frak L\frak o\frak c}
\def\loc{\frak l\frak o\frak c}

 3
 2
\font\large=cmr10  scaled \magstep 2
 2
\font\larti=cmti10  scaled \magstep 2
 2
\font\cds=cmr7
\font\cdt=cmti7
\font\cdy=cmsy7

\centerline{\large Ordinary Grothendieck groups}
\smallskip
\centerline{\large  of a Frobenius {\larti P}-category}

\medskip
\centerline{ Lluis Puig }
\medskip
\noindent
{\bf  Abstract:} {\cds In [5] we have introduced  the Frobenius categories {\cdy F} over a finite {\cdt p}-group {\cdt P},
and we have associated to~{\cdy F} --- suitably endowed with some central {\cdt k}*-extensions --- a ``Grothendieck group''
as an inverse limit of Grothendieck groups of categories of modules in characteristic {\cdt p} obtained from  {\cdy F},
determining its rank.  Our purpose here is to introduce an analogous inverse limit of Grothendieck groups of categories of modules in characteristic zero obtained from   {\cdy F}, determining its rank and proving that its extension to a field
is canonically isomorphic to the direct sum of the corresponding extensions of the ``Grothendieck groups'' above associated with the {\cdt centralizers} in {\cdy F} of a suitable set of representatives of the {\cdy F}-classes of elements of {\cdt P}.
}

\bigskip\noindent
{\bf £1\phantom{.} Introduction}
\medskip

£1.1\phantom{.} Let $p$ be a prime number and $\O$ a complete discrete valuation
ring with a {\it field of quotients\/} $\K$ of characteristic zero and a {\it residue field\/}
$k$ of characteristic $p\,;$ we assume that $k$ is algebraically closed and that $\K$
contains ``enough'' roots of unity for the finite family of finite groups we will consider.
Let $G$ be a finite group, $b$ a {\it block\/} of $G$ ---
namely a primitive idempotent in the center $Z(\O G)$ of the group $\O\-$algebra ---
and $(P,e)$ a maximal {\it Brauer $(b,G)\-$pair\/} [5,~1.16]; recall that the {\it Frobenius 
$P\-$category\/} $\F_{\!(b,G)}$ associated with~$b$ is the subcategory of the category of finite groups where the objects are all the subgroups of $P$ and, for any pair of subgroups $Q$ and $R$ of~$P\,,$ the morphisms $\varphi$ from $R$ to $Q$ are the group homomorphisms $\varphi\,\colon R\to Q$ induced by the conjugation of some
element $x\in G$ fulfilling
$$(R,g)\i (Q,f)^x
\eqno £1.1.1\phantom{.}$$
where $(Q,f)$ and $(R,g)$ are the corresponding Brauer $(b,G)\-$pairs contained in 
$(P,e)$ [5, Ch.~3].

\medskip
£1.2\phantom{.} In~[5, Ch. 14] we consider a suitable inverse limit of Grothendieck groups of categories of modules in characteristic $p$ obtained from  $\F_{\!(b,G)}\,,$ which, according to Alperin's Conjecture, should be isomorphic to the 
 Gro-thendieck group of the category of finitely dimensional $kGb\-$modules. As announced in the title,
 our purpose here is to introduce an analogous inverse limit of Grothendieck groups of categories of modules in characteristic zero obtained from  $\F_{\!(b,G)}\,,$ which again, according to Alperin's Conjecture, should be isomorphic to the  Grothendieck group of the category of finitely dimensional $\K Gb\-$modules.

 \medskip
 £1.3\phantom{.} More explicitly, recall that a Brauer $(b,G)\-$pair $(Q,f)$ is called {\it selfcentralizing\/} if, 
 for any $x\in G$ such that $(P,e)^x$ contains $(Q,f)\,,$ we have $C_{P^x}(Q) = Z(Q)$ or, equivalently, 
 if the image $\bar f$ of $f$ in   $Z\big(k \bar C_G(Q)\big)\,,$ where $\bar C_G(Q) = C_G(Q)/Z(Q)\,,$ 
 is a block of {\it defect zero\/} [5,~Corollary~7.3]. As a matter of fact, we restrict ourselves to the full subcategory  
 $\F_{\!(b,G)}^{^{\rm sc}}$   of $\F_{\!(b,G)}$ over the  $\F_{\!(b,G)}\-$objects $Q$ such that the corresponding Brauer $(b,G)\-$pair $(Q,f)$ contained in $(P,e)$ is selfcentralizing; recall that, for such an $\F_{\!(b,G)}\-$object~$Q\,,$ the action of the normalizer $N_G (Q,f)$ on the simple $k\-$algebra $k \bar C_G(Q)\bar f$ determines a central $k^*\-$extension $\hat\F_{\!(b,G)} (Q)$ of the quotient [5, 7.4]
 $$\F_{\!(b,G)} (Q) = \F_{\!(b,G)}(Q,Q)\cong N_G (Q,f)/C_G(Q)
 \eqno £1.3.1.$$

 \medskip
 £1.4\phantom{.}  Then, considering the {\it proper category of 
 $\F_{\!(b,G)}^{^{\rm sc}}\-$chains\/} $\ch^*(\F_{\!(b,G)}^{^{\rm sc}})$ --- namely, the category of functors $\frak q\,\colon \Delta_n \to \F_{\!(b,G)}^{^{\rm sc}}$ from any {\it ordered simplex\/} $\Delta_n$ where the morphisms are defined by the {\it order preserving maps\/} between simplexes and the {\it natural isomorphisms\/} between the corresponding functors over the same simplex [5, A2.8] --- we prove in [5, Ch. 11] that the canonical functor
$$\aut_{\F_{\!(b,G)}^{^{\rm sc}}} : \ch^*(\F_{\!(b,G)}^{^{\rm sc}})\too \Gr
\eqno £1.4.1\phantom{.}$$
mapping any $\F_{\!(b,G)}^{^{\rm sc}}\-$chain $\frak q\,\colon
\Delta_n \to \F_{\!(b,G)}^{^{\rm sc}}$ on the group $\F_{\!(b,G)}^{^{\rm sc}}
(\frak q)$ of {\it natural automorphisms\/} of~$\frak q$  [5, Proposition~A2.10] can lifted to a 
functor
$$\widehat\aut_{\F_{\!(b,G)}^{^{\rm sc}}} : \ch^*(\F_{\!(b,G)}^{^{\rm sc}})\too
 k^*\-\Gr
\eqno £1.4.2\phantom{.}$$
where $\Gr$ and $k^*\-\Gr$ respectively denote the categories of finite groups and of central $k^*\-$exten-sions of finite groups, called {\it finite $k^*\-$groups\/}.

\medskip
£1.5\phantom{.} At this point, denoting by $\frak g_k\,\colon  k^*\-\Gr\to \O\-\mod$
the {\it contravariant\/} functor sending any $k^*\-$group $\hat G$ to the 
$\O\-$extension of the Grothendieck group of the category of finitely dimensional 
$k_*\hat G\-$modules --- noted $\G_k (\hat G)$ --- in [5,~Ch.~4] we introduce the inverse limit
$$\G_k(\F_{\!(b,G)},\widehat\aut_{\F_{\!(b,G)}^{^{\rm sc}}}) = \lim_{\longleftarrow}
(\frak g_k\circ  \widehat\aut_{\F_{\!(b,G)}^{^{\rm sc}}})
\eqno £1.5.1\phantom{.}$$
--- called the {\it (modular) Grothendieck group\/} of $\F_{\!(b,G)}$ ---
which, strictly speaking, depends not only on the Frobenius $P\-$category 
$\F_{\!(b,G)}$ but on the lifting~$\widehat\aut_{\F_{\!(b,G)}^{^{\rm sc}}}\,.$ There, we also prove that 
$${\rm rank}_\O\big(\G_k(\F_{\!(b,G)},\widehat\aut_{\F_{\!(b,G)}^{^{\rm sc}}}) \big) =  \sum_{(\frak
q,\Delta_n)} (-1)^{n}\,{\rm rank}_{\O} \Big(\G_k \big(\hat\F_{\!(b,G)}^{^{\rm sc}}
(\frak q)\big)\Big)
\eqno £1.5.2,$$
where $(\frak q,\Delta_n)$ runs over a set of representatives for the set of isomorphism
classes of {\it regular\/} $\ch^*(\F_{\!(b,G)}^{^{\rm sc}})\-$objects [5, 14.31] which
shows that Alperin's Conjecture actually states [5, I\hskip1pt32]
$$\G_k(\F_{\!(b,G)},\widehat\aut_{\F_{\!(b,G)}^{^{\rm sc}}}) \cong \G_k (G,b)
\eqno £1.5.3.$$

\medskip
£1.6\phantom{.} Since we have the well-known isomorphism [1]
$${}^\K \G_\K (G,b)\cong \bigoplus_{(u,g)\in \U} {}^\K \G_k 
\big(C_G(u),g\big)
\eqno £1.6.1\phantom{.}$$
where $\U$ is a set of representatives for the set of $G\-$conjugacy classes
of  {\it Brauer $(b,G)\-$elements\/} [5,~1.10], where $\G_\K (G,b)$ denotes the $\O\-$extension 
of the Grothendieck group of the category of finitely dimensional $\K Gb\-$modules
and where we set ${}^\K \G_\K (G,b) = \K\otimes_\O \G_\K(G,b)$ for short, as we said in 
[5, I\hskip1pt51 and I\hskip1pt52]  the direct sum
$$\bigoplus_{(u,g)\in \U} {}^\K \G_k \big(\F_{\!(g,C_G(u))},
 \widehat\aut_{\F_{\!(g,C_G(u))}^{^{\rm sc}}}\big)
\eqno £1.6.2\phantom{.}$$
already provides a reasonable definition for the ordinary  Grothendieck group of $\F_{\!(b,G)}\,,$ at least extended to $\K\,.$ But, as we mention there, there is a more reasonable definition as an inverse limit, analogous to definition~£1.5.1.

\medskip
£1.7\phantom{.} Firstly notice that it does not suffice to replace $\frak g_k$ by the {\it contrava-riant\/} functor 
$\frak g_\K\,\colon k^*\-\Gr\to \O\-\mod$  sending any $k^*\-$group $\hat G$ to the $\O\-$extension  of the Grothendieck group of the category of finitely dimensional $\K_*\hat G\-$mo-dules where, considering the canonical group homomorphisms $k^*\to
\O^*\i\K^*$ and $k^*\to \hat G\,,$ we set
$$\K_*\hat G = \K\otimes_{\K k^*} \K \hat G
\eqno £1.7.1;$$
indeed, it suffices to consider the case of the blocks of the $p\-$solvable groups --- discussed in~[6],  for instance --- to understand that the groups 
$\F_{\!(b,G)}^{^{\rm sc}}(\frak q)$ above have to be replaced by the corresponding 
{\it localizers\/} [5, 18.3] and, more generally, that the functor 
$\aut_{\F_{\!(b,G)}^{^{\rm sc}}}$ has to be replaced by the 
{\it $\F_{\!(b,G)}^{^{\rm sc}}\-$localizing functor\/}
introduced in [5, Ch.~18].

\medskip
£1.8\phantom{.} Precisely, recall that for any $\F_{\!(b,G)}^{^{\rm sc}}\-$object $Q$
such that $\F_P (Q)$ is a Sylow $p\-$subgroup of $\F_{\!(b,G)}(Q)\,,$ the group 
$N_P (Q)$ --- viewed as an extension of $\F_P (Q)$ by $Z(Q)$ --- determines an extension
$\L_{(b,G)} (Q)$ of $\F_{\!(b,G)}(Q)$ by $Z(Q)\,,$ containing $N_P (Q)$ ---  called the {\it localizer\/} of $Q$ [5, Theorem~18.6]. Then, considering the category 
$\widetilde\Loc$ where the objects are the pairs $(L,Z)$ formed by a finite group $L$ and a normal $p\-$subgroup $Z$ of $L\,,$ and where the morphisms from $(L,Z)$ to $(L',Z')$  are the $Z'\-$conjugacy classes of group monomorphisms $f\,\colon L\to L'$ fulfilling $f(Z)\i Z'$ [5, 18.12], in [5, Proposition~18.19] we prove the existence, and the uniqueness up to isomorphisms, of a functor 
$$\loc_{\F_{\!(b,G)}^{^{\rm sc}}}: \ch^*(\F_{\!(b,G)}^{^{\rm sc}})\too \widetilde\Loc
\eqno £1.8.1\phantom{.}$$
which lifts $\aut_{\F_{\!(b,G)}^{^{\rm sc}}}$ above {\it via\/} the functor 
$\widetilde\Loc\to \Gr$ sending $(L,Z)$ to $L/Z\,,$ and maps any 
$\F_{\!(b,G)}^{^{\rm sc}}\-$chain $\frak q\,\colon \Delta_n \to 
\F_{\!(b,G)}^{^{\rm sc}}$ such that $\F_P \big(\frak q (n)\big)$ is a Sylow $p\-$subgroup of $\F_{\!(b,G)}\big(\frak q (n)\big)\,,$ on the pair $\big(\L_{(b,G)}
 (\frak q),Z(\frak q (n))\big)$ where $\L_{(b,G)} (\frak q)$ is the converse image of 
 $\F_{(b,G)} (\frak q)$ in $\L_{(b,G)} \big(\frak q (n)\big)\,.$
 \eject

 \medskip
 £1.9\phantom{.} {\it Mutatis mutandis\/}, we can consider the category 
 $k^*\-\widetilde\Loc$ where the objects are the pairs $(\hat L,Z)$ formed by a finite
 $k^*\-$group $\hat L$ and a normal $p\-$subgroup $Z$ of $\hat L\,,$ and then the lifting $\widehat\aut_{\F_{\!(b,G)}^{^{\rm sc}}}$ of 
 $\aut_{\F_{\!(b,G)}^{^{\rm sc}}}$ above determines, {\it via pull-backs\/}, a functor
$$\widehat\loc_{\F_{\!(b,G)}^{^{\rm sc}}} : \ch^*(\F_{\!(b,G)}^{^{\rm sc}})\too
 k^*\-\widetilde\Loc
\eqno £1.9.1\phantom{.}$$
lifting  $\loc_{\F_{\!(b,G)}^{^{\rm sc}}}\,.$ At this point, still denoting by 
$\frak g_\K\,\colon  k^*\-\widetilde\Loc\to \O\-\mod$ the functor sending 
$(\hat L,Z)$ to $\G_\K(\hat L)\,,$ we define the {\it ordinary Grothendieck group\/} of $\F_{\!(b,G)}$ as the following inverse limit
$$\G_\K (\F_{\!(b,G)},\widehat\aut_{\F_{\!(b,G)}^{^{\rm sc}}}) = \lim_{\longleftarrow}
(\frak g_\K\circ  \widehat\loc_{\F_{\!(b,G)}^{^{\rm sc}}})
\eqno £1.9.2\phantom{.}$$
which once again depends on the lifting  $\widehat\aut_{\F_{\!(b,G)}^{^{\rm sc}}}\,.$ Note that, since  we have
$\G_k (\hat L)\cong \G_k (\hat L/Z)\,,$ we also have
$$\frak g_k\circ  \widehat\loc_{\F_{\!(b,G)}^{^{\rm sc}}}\cong \frak g_k\circ  
\widehat\aut_{\F_{\!(b,G)}^{^{\rm sc}}}
\eqno £1.9.3.$$

\medskip
£1.10\phantom{.} Our purpose here is to prove that an equality analogous to 
equality~£1.5.2 holds, namely that we have
$${\rm rank}_\O\big(\G_\K (\F_{\!(b,G)},\widehat\aut_{\F_{\!(b,G)}^{^{\rm sc}}})\big) =  \sum_{(\frak q,\Delta_n)} (-1)^{n}\,{\rm rank}_{\O} 
\Big(\G_\K \big(\hat\L_{\!(b,G)} (\frak q)\big)\Big)
\eqno £1.10.1,$$
where $(\frak q,\Delta_n)$ runs over a set of representatives for the set of isomorphism
classes of {\it regular\/} $\ch^*(\F_{\!(b,G)}^{^{\rm sc}})\-$objects (cf.~£8.3 below) and we set
$$\hat\L_{\!(b,G)} (\frak q) = \hat\F_{\!(b,G)} (\frak q)\times_{\F_{\!(b,G)} (\frak q)}
\L_{\!(b,G)} (\frak q)
\eqno £1.10.2,$$
and that the direct sum~£1.6.2 coincides with the extension to $\K$ of the  ordinary Grothendieck group of $\F_{\!(b,G)}\,,$ namely that we still have
$${}^\K \G_\K (\F_{\!(b,G)},\widehat\aut_{\F_{\!(b,G)}^{^{\rm sc}}}) \cong 
\bigoplus_{(u,g)\in \U} {}^\K \G_k \big(\F_{\!(g,C_G(u))},
 \widehat\aut_{\F_{\!(g,C_G(u))}^{^{\rm sc}}}\big)
\eqno £1.10.3.$$

\medskip
£1.11\phantom{.} A remarkable fact is that neither these statements nor our arguments 
for proving them need to assume that the Frobenius $P\-$category $\F$ we are dealing with
comes from a block of a finite group, but only need the choice of a lifting 
$$\widehat\aut_{\F^{^{\rm sc}}} : \ch^*(\F^{^{\rm sc}})\too k^*\-\Gr
\eqno £1.11.1\phantom{.}$$ of the functor $\aut_{\F^{^{\rm sc}}}$   [5, Proposition~A2.10]. 
Thus, as in [5] for the {\it modular Grothendieck group\/}, we will carry out our purpose over such a triple $(P,\F,\widehat\aut_{\F^{^{\rm sc}}})$ that we call a {\it folded Frobenius $P\-$category\/}. Actually, the reader may ask himself
why the present material has not been included in [5].\break
\eject
\noindent
 The answer is quite simple: because when finishing [5] we had not at all it and our question in [5, I\hskip1pt52] was not yet answered! Naturally, this means that some arguments here are definitely not contained in [5]. On the one
hand, our proof of isomorphism~£1.10.3 needs the rather technical Lemma~£9.4 below which in some sense ``explains''
why the localizers of the $\F\-$selfcentralizing subgroups of $P$ are powerful enough to compute the complete Grothendieck group. On the other hand, even our proof of equality~£1.10.1 needs a more sophisticated machinery
than the proof in [5] of equality~£1.5.2 above.

\medskip
£1.12\phantom{.} This paper is divided on nine sections; notation, terminology and results in [5] are our main reference. In all the paper $P$  is a finite $p\-$group; in section~2  we recall the main facts we need here on Frobenius
$P\-$categories and give a sufficient condition to have a {\it folded Frobenius $P\-$category\/}. In section~3 we give equivalent definitions of the {\it ordinary Grothendieck group\/}
of a folded  Frobenius $P\-$category. Section~4 is devoted to the functoriality of both, the ordinary and the modular Grothendieck groups of folded  Frobenius $P\-$categories, a subject that in [5] has only been partially discussed
in the framework of blocks; on the converse, it is reasonable to hope that the analogous reduction results in [5, Ch. 15]
will admit a translation to the ordinary Grothendieck group  of the  Frobenius $P\-$categories of blocks, but this has not been done here.

\medskip
£1.13\phantom{.} In section~8 we prove equality~£1.10.1, and the previous sections~5, 6 and 7 play an auxiliary role
in that proof. In section~5 we develop a canonical decomposition of the ordinary Grothendieck group analogous to the decomposition of the modular Grothendieck group in~[5,~Ch.~14], except on the fact, pointed out in~£1.7 above,
that we have to replace the functor $\aut_{\F^{^{\rm sc}}}$ by the {\it $\F^{^{\rm sc}}\-$localizing functor\/}.
Section~6 is devoted to some formal transformation of each term of our decomposition, which facilitates the application
 of a vanishing cohomological result. In section~7 we prove this vanishing cohomological result which generalizes 
 [5,~Theorem~6.26]; the existence of such a generalization is a critical point in our paper. Finally, in section~9 we prove isomorphism~£1.10.3; our proof needs  Lemma~£9.4 as metioned above, and a sophisticated counting argument
 to show the equality of dimentions.

\bigskip
\bigskip
\noindent
{\bf £2\phantom{.} The folded Frobenius $P\-$categories}
\bigskip

£2.1\phantom{.} Let $P$ be a finite $p\-$group and denote  by  $\frak i\Gr$ the category formed by the finite groups and by the injective group  homomorphisms, and  by 
$\F_{\!P}$ the subcategory of  $\frak i\Gr$ where the objects are all the  subgroups
 of $P$ and the morphisms are the group homomorphisms induced by conjugation by elements of $P\,.$

 \medskip
£2.2\phantom{.} Recall that a {\it Frobenius  $P\-$category\/} $\F$ is a subcategory 
of $\frak i\Gr$ containing $\F_{\!P}$ where the objects are all the  subgroups of $P$
and the morphisms fulfill the following three conditions [5, 2.8 and Proposition~2.11]
\smallskip
\noindent
£2.2.1\quad {\it For any subgroup $Q$ of $P$ the inclusion functor $(\F)_Q\to 
(\frak i\Gr)_Q$ is full.\/}
\smallskip
\noindent
£2.2.2\quad {\it $\F_P (P)$ is a Sylow $p\-$subgroup of $\F (P)\,.$\/}
\smallskip
\noindent
£2.2.3\quad {\it For any  subgroup $Q$ of $P$ such that we have $\xi \big(C_P (Q)\big)
= C_P\big(\xi (Q)\big)$ for any $\F\-$morphism $\xi\,\colon Q\.C_P(Q)\to P\,,$
any $\F\-$morphism $\varphi\,\colon Q\to P$ and any subgroup $R$ of $N_P\big(\varphi(Q)\big)$ containing $\varphi (Q)$ such that $\F_P(Q)$ contains the action of $\F_R \big(\varphi(Q)\big)$ over $Q$ via $\varphi\,,$ there is an $\F\-$morphism 
$\zeta\,\colon R\to P$ fulfilling $\zeta\big(\varphi (u)\big) = u$ for any $u\in Q\,.$\/}
\smallskip
\noindent
As in [5,~1.2], for any pair of subgroups $Q$ and $R$ of $P\,,$ we denote by $\F (Q,R)$ the set of $\F\-$morphisms from $Q$ to $R$ and set $\F (Q) = \F (Q,Q)\,;$ moreover, recall
that, for any category $\frak C$ and any $\frak C\-$object $C\,,$ $\frak C_C$ 
(or $(\frak C)_C$ to avoid confusion) denotes the category of ``$\frak C\-$morphisms
to $C$'' [5, 1.7].

\medskip
£2.3\phantom{.} Given a Frobenius $P\-$category $\F\,,$ a subgroup $Q$ of $P$ and
a subgroup $K$ of the group ${\rm Aut}(Q)$ of automorphisms of $Q\,,$
we say that $Q$ is {\it fully $K\-$normalized\/} in $\F$ if we have [5, 2.6]
$$\xi \big(N_P^K (Q)\big) = N_P^{\,{}^\xi \!K}\big(\xi (Q)\big)
\eqno £2.3.1$$
for any $\F\-$morphism $\xi\,\colon Q\.N_P^K (Q)\to P\,,$ where $N_P^K (Q)$ is the converse image of $K$ in $N_P (Q)$ {\it via\/} the canonical group homomorphism
$N_P (Q)\to {\rm Aut}(Q)$ and ${}^\xi \! K$ is the image of $K$ in 
${\rm Aut}\big(\xi (Q)\big)$ {\it via\/} $\xi\,.$ Recall that if $Q$ is fully 
$K\-$normalized in $\F$ then we have a new Frobenius $N_P^K(Q)\-$category 
$N_\F^K(Q)$ where, for any pair of subgroups $R$ and $T$ of $N_P^K (Q)\,,$
$\big(N_\F^K(Q)\big)(R,T)$ is the set of group homomorphisms from $T$ to $R$ induced by the $\F\-$morphisms $\psi\,\colon Q\.T\to Q\.R$ which stabilzes $Q$ and induces on it an element of $K$ [5, 2.14 and Proposition~2.16].

\medskip
£2.4\phantom{.} We say that a subgroup $Q$ of $P$ is {\it $\F\-$selfcentralizing\/} if we have
$$C_P\big(\varphi (Q))\i \varphi (Q)
\eqno £2.4.1\phantom{.}$$
 for any $\varphi \in \F (P,Q)\,,$ and we denote by $\F^{^{\rm sc}}$ the full subcategory of $\F$ over the set of $\F\-$selfcentralizing subgroups of $P\,.$ From the case of the Frobenius $P\-$categories associated with a block of a finite group, we know that only makes sense to consider central $k^*\-$extensions of $\F (Q)$ whenever $Q$ is 
$\F\-$selfcentralizing [5,~7.4]; but, if $U$ is a subgroup of $P$ fully $K\-$normalized in~$\F$
for some subgroup $K$ of ${\rm Aut}(U)\,,$ a $N_F^K (U)\-$selfcentralizing subgroup
of $N_P (Q)$ needs not be  $\F\-$selfcentralizing, which is a handicap when comparing
choices of central $k^*\-$extensions in $\F$ and in $N_F^K (U)\,.$
 In order to overcome this difficulty, we consider the {\it $\F\-$radical\/}
subgroups of~$P\,;$ we say that a subgroup $R$ of $P$ is {\it $\F\-$radical\/} if it is
$\F\-$selfcentralizing and we have
$${\bf O}_p\big(\tilde\F (R)\big) = \{1\}
\eqno £2.4.2\phantom{.}$$
where   $\tilde\F (R) = \F (R)/\F_R (R)$ [5, 1.3]; we denote by $\F^{^{\rm rd}}$ the full subcategory of $\F$ over the set of $\F\-$radical subgroups of $P\,.$
\eject

\bigskip
\noindent
{\bf Lemma~£2.5}\phantom{.} {\it Let $\F$ be a Frobenius $P\-$category, $U$
a subgroup of $P$ and $K$ a subgroup of ${\rm Aut}(U)$ containing ${\rm Int}(U)\,.$
If $U$ is fully $K\-$normalized in $\F$ then any $N_\F^K(U)\-$radical subgroup $R$ of $N_P^K(U)$ contains $U$ and, in particular, it is $\F\-$selfcentralizing.\/}
\medskip
\noindent
{\bf Proof:} It is quite clear that the image of $N_{U\.R}(R)$ in $\big(N_\F^K(U)\big)(R)$
is a normal $p\-$subgroup and therefore it is contained in ${\bf O}_p\Big(\big(N_\F^K(U)\big)(R)\Big)\,,$ so that $N_{U\.R}(R) = R$ which forces $U\.R = R\,.$ Moreover,
for any $\F\-$morphism $\psi\,\colon R\to P\,,$ it is clear that $\psi (U)$ is a normal 
subgroup of $\psi (R)\.C_P\big(\psi (R)\big)$ and therefore, since $U$ is also fully centralized in $\F$ [5, Proposition~2.12], it follows from £2.2.3 that there is an 
$\F\-$morphism 
$$\zeta  : \psi (R)\.C_P\big(\psi (R)\big)\too P
\eqno £2.5.1\phantom{.}$$
fulfilling $\zeta \big(\psi (u)\big) = u$ for any $u\in U\,,$ so that the group homomorphism from $R$ to $N_P^K(U)$ mapping $v\in R$ on $\zeta\big(\psi (v)\big)$
is a $N_\F^K(U)\-$morphism; in particular, $\zeta\big(\psi (R)\big)$ is also 
$N_\F^K(U)\-$selfcentralizing and therefore we get
$$\zeta \Big(C_P\big(\psi (R)\big)\Big)\i \zeta\big(\psi (R)\big)
\eqno £2.5.2\phantom{.}$$
which forces $C_P\big(\psi (R)\big)\i \psi (R)\,.$ We are done.

\medskip
£2.6\phantom{.} As a matter of fact, the $\F\-$radical subgroups of $P$ admit a description in terms of the {\it dominant\/} $\F\-$morphisms; let us say that an 
$\F\-$morphism  $\varphi\,\colon Q\to R$ is {\it dominant\/} if it fulfills
$$\varphi\circ\F (Q)\i \F (R)\circ\varphi
\eqno £2.6.1;$$
it is quite clear that the $\F\-$isomorphisms are dominant and that the composition of dominant $\F\-$morphisms is a dominant  $\F\-$morphism. Setting $Q' = \varphi (Q)$ and denoting by $\F(R)_{Q'}$ the stabilizer of $Q'$ in $\F (R)\,,$ this
condition is equivalent to saying that $\varphi$ determines a surjective group homomorphism
$$\F(R)_{Q'}\too \F (Q)
\eqno £2.6.2;$$
moreover, if $Q$ is $\F\-$selfcentralizing then the kernel of this homomorphism coincides with $\F_{Z(Q)}(R)$ [5,~Corollary~4.7]; in this case, the choice of a central 
$k^*\-$extension of $\F (R)$ clearly determines one of $\F (Q)\,.$ On the other hand, since a Sylow $p\-$subgroup of $\F(R)_{Q'}$ maps onto a Sylow $p\-$subgroup of 
$\F(Q)\,,$ if $Q$ is  $\F\-$selfcentralizing and $R$ is fully normalized in $\F$ then we may assume that $\F_P (R)$ contains a Sylow $p\-$subgroup of $\F(R)_{Q'}\,;$ in this case $Q'$ is also fully normalized in $\F$ [5, Proposition~2.12] and it is easily checked that $\varphi$ can be extended to an $\F\-$morphism
$$\hat\varphi : N_P (Q)\too N_P (Q') = N_P(R)_{Q'}\i N_P(R)
\eqno £2.6.3.$$
\eject

\bigskip
\noindent
{\bf Proposition~£2.7}\phantom{.} {\it  Let $\F$ be a Frobenius $P\-$category,
$R$ a subgroup of $P$ fully centralized in $\F$ and $Q$ the converse image of
${\bf O}_p \big(\F (R)\big)$ in $N_P (R)\,.$ Then, the inclusion map $\iota_R^Q
\,\colon R\to Q$ is dominant. In particular, $R$ is $\F\-$radical if and only if any dominant $\F\-$morphism from $R$ is an isomorphism.\/}
\medskip
\noindent
{\bf Proof:} Recall that $R$ is also fully ${\bf O}_p \big(\F (R)\big)\-$normalized [5,~2.10] and therefore $N_P(R)$ contains ${\bf O}_p \big(\F (R)\big)$ [5, Proposition~2.12]; in particular, for any $\tau\in \F (R)\,,$ the composition 
$\iota_R^P\circ \tau\,\colon R\to P$ can be extended to an $\F\-$morphism $\zeta\,\colon QÊ\to P$ [5, statement~2.10.1]; since $\zeta \big(C_P (R)\big) = C_P(R)\,,$
it is quite clear that $\zeta (Q) = Q$ and denoting by $\sigma\,\colon Q\cong Q$ the automorphism determined by~$\zeta\,,$ we have $\iota_R^Q\circ \tau = \sigma\circ
\iota_R^Q$ and, according to condition~£2.2.1, $\sigma$ belongs to~$\F(Q)\,.$
This proves that $\iota_R^Q\,\colon R\to Q$ is dominant.

\smallskip
In particular, if we assume that any dominant $\F\-$morphism from $R$ is an isomorphism, then we get $R = Q$ which proves that $R$ is $\F\-$radical.
Conversely, assume that $R$ is $\F\-$radical and let $\varphi\,\colon R\to Q$ be a dominant $\F\-$morphism; for any $\tau\in \F (R)$ there is $\sigma\in \F (Q)$
fulfilling $\varphi\circ \tau = \sigma\circ \varphi$
and then, setting $R' = \varphi(R)\,,$ it is easily checked that 
$$\sigma\big(N_Q (R')\big) = N_Q (R')
\eqno £2.7.1;$$
in particular, the image $\tau'$ of $\tau$ in $\F (R')$ {\it via\/} the $\F\-$isomorphism
$R\cong R'$ determined by $\varphi$ (cf. condition~£2.2.1) normalizes the image
$\F_Q (R')$ of $N_Q(R')$ in~$\F (R')\,;$ consequently, $\F_Q (R')$ is a normal 
$p\-$subgroup of $\F (R')$ and therefore it is contained in ${\bf O}_p 
\big(\F (R')\big)\,;$ but, since $R'$ is also $\F\-$radical, ${\bf O}_p 
\big(\F (R')\big)$ is just the image of $R'$ in $\F (R')$ and therefore we have 
$N_Q (R') = R'$ which forces $R' = Q\,.$ We are done.

\medskip
£2.8\phantom{.} But, in our general setting, we have to deal with 
{\it $\F^{^{\rm sc}}\-$chains\/} and {\it coherent\/} choices of central 
$k^*\-$extensions for their $\F^{^{\rm sc}}\-$automorphisms group.
Recall that we call {\it $\F^{^{\rm sc}}\-$chain\/} any functor $\frak q\,\colon \Delta_n\to \F^{^{\rm sc}}$ where the $n\-$simplex $\Delta_n$ is considered as a category with
the morphisms defined by the order [5,~A2.2]; then, we consider the category 
$\ch^*(\F^{^{\rm sc}})$ where the objects are all the $\F^{^{\rm sc}}\-$chains 
$(\frak q,\Delta_n)$ and the morphisms from $\frak q\,\colon \Delta_n\to 
\F^{^{\rm sc}}$ to another $\F^{^{\rm sc}}\-$chain  $\frak r\,\colon \Delta_m\to \F^{^{\rm sc}}$ are the pairs $(\nu,\delta)$ formed by an {\it order preserving map\/} 
or, equivalently, a functor  $\delta\,\colon \Delta_m\to \Delta_n$ and by a natural isomorphism $\nu\,\colon \frak q\circ\delta\cong \frak r\,,$ the composition being defined by the composition of maps and of natural isomorphisms [5,~A2.8]; the point is that we have a canonical functor
$$\aut_{\F^{^{\rm sc}}} : \ch^*(\F^{^{\rm sc}})\too \Gr
\eqno £2.8.1$$
mapping  any $\F^{^{\rm sc}}\-$chain $\frak q\,\colon \Delta_n\to \F^{^{\rm sc}}$
to the group of natural automorphisms of $\frak q\,,$ simply noted~$\F (\frak q)$ [5, Proposition~A2.10].
 We define a {\it folded Frobenius\break
 \eject
\noindent
 $P\-$category\/} as the pair formed by a Frobenius $P\-$category $\F$ and 
 by the choice  of a functor 
$$\widehat\aut_{\F^{^{\rm sc}}} : \ch^*(\F^{^{\rm sc}})\too k^*\-\Gr
\eqno £2.8.2\phantom{.}$$
lifting $\aut_{\F^{^{\rm sc}}}\,.$ {\it Mutatis mutandis\/}, we consider the category
$\ch^*(\F^{^{\rm rd}})$ and the canonical functor
$$\aut_{\F^{^{\rm rd}}} : \ch^*(\F^{^{\rm rd}})\too \Gr
\eqno £2.8.3.$$

\bigskip
\noindent
{\bf Theorem~£2.9}\phantom{.} {\it Any functor $\widehat\aut_{\F^{^{\rm rd}}}$
lifting $\aut_{\F^{^{\rm rd}}}$ to the category $k^*\-\Gr$ can be extended to a unique functor
lifting $\aut_{\F^{^{\rm sc}}}\,.$
$$\widehat\aut_{\F^{^{\rm sc}}} : \ch^* (\F^{^{\rm sc}})\too k^*\-\Gr 
\eqno  £2.9.1\phantom{.}$$\/}
\par
\noindent
{\bf Proof:} Let $\frak X$ be a set of $\F\-$selfcentralizing subgroups of $P$ which contains all the $\F\-$radical subgroups of $P$ and is stable by $\F\-$isomorphisms; denoting by $\F^{^\frak X}$
the {\it full\/} subcategory of $\F$ over~$\frak X\,,$ assume that
 $\widehat\aut_{\F^{^{\rm rd}}}$ can be extended
to a unique functor
$$\widehat\aut_{\F^{^{\frak X}}} : \ch^* (\F^{^{\frak X}})\too k^*\-\Gr 
\eqno £2.9.2;$$
assuming that $\frak X$ does not coincide with the set of all the $\F\-$selfcentralizing subgroups  of~$P\,,$ let $V$ be a maximal  $\F\-$selfcentralizing subgroup which is not in $\frak X$ and $\rho \,\colon V\to W$ a dominant $\F\-$morphism such that $\rho (V)$ is a proper normal subgroup of $W$ (cf.~Proposition~£2.7); according to~£2.6, we may assume that $W$ and $\rho (V)$ are fully normalized in $\F\,.$ Then, denoting by $\frak Y$ the union of $\frak X$ with all the subgroups of $P$  $\F\-$isomorphic to~$V\,,$ it is clear that it suffices to prove that $\widehat\aut_{\F^{^{\frak X}}}$ admits a unique extension to $\ch^* (\F^{^{\frak Y}})\,.$

\smallskip
 For any chain  $\frak q\,\colon {\bf \Delta}_n\to \F^{^{\frak Y}}\,,$ denote by $\hat\frak q\colon  \Delta_{n+1}\to \F^{^{\frak Y}}$ the chain extending $\frak q$ 
 such that either  $\frak q (n)\cong V$ and we set $\hat\frak q (n+1) = W$ and 
 $\hat\frak q (n\bullet n\!+\!1)$ is the composition  of $\rho$ with an isomorphism 
   $\frak q (n)\cong V\,,$ or $\frak q (n)\not\cong V$ and we set $\hat\frak q (n\!+\!1) = \frak q (n)$ and  $\hat\frak q (n\bullet n\!+\!1) = {\rm id}_{\frak q (n)}\,;$ in both cases, we have an obvious $\ch^* (\F^{^{\frak Y}})\-$morphism [5, A3.1]
$$({\rm id}_\frak q,\delta^n_{n+1}) : (\hat\frak q,\Delta_{n+1})\too (\frak q,\Delta_n)\eqno £2.9.3\phantom{.}$$
 and the functor $\aut_{\F^{^{\frak Y}}}$ maps $({\rm id}_\frak q,\delta^n_{n+1})$ on a group homomorphism
$$\F(\hat\frak q)\too \F (\frak q)
\eqno £2.9.4\phantom{.}$$
which is surjective since any $\sigma\in \F (\frak q)\i \F\big(\frak q (n)\big)$ can be
extended to an $\F\-$automorphism of $\hat\frak q (n+1)$ (cf.~£2.6.2);
moreover, since $V$ is $\F\-$selfcen-tralizing, the kernel of this homomorphism is a 
$p\-$group (cf.~£2.6); then, the functor $\widehat\aut_{\F^{^\frak X}}$ and the structural inclusion $\F(\hat\frak q)\i \F\big(\hat\frak q (n+1)\big)$ determine
 a $k^*\-$group $\hat\F(\hat\frak q)$ and this $k^*\-$group induces  a unique central $k^*\-$extension 
 $\hat\F (\frak q)$ of~$\F (\frak q)$ such that we have a $k^*\-$group homomorphism
$$\hat\F(\hat\frak q)\too \hat\F (\frak q)
\eqno £2.9.5\phantom{.}$$
lifting homomorphism~£2.9.4.
\eject

\smallskip
 Now, for any $\ch^*(\F^{^\frak Y})\-$morphism $(\nu,\delta)\,\colon (\frak r,\Delta_m)\to (\frak q,\Delta_n)\,,$  we have to exhibit a $k^*\-$group homomorphism 
 $\hat\F (\frak r)\to \hat\F (\frak q)$ lifting $\aut_{\F^{^\frak Y}} (\nu,\delta)\,.$ Firstly, denoting by $Q$ and $R$ the respective images of $\frak r\big(\delta (n)\big)$ and  $\frak r (m)$  in $\hat\frak r (m\!+\!1)$ by the group homomorphisms $\hat\frak r 
 (\delta (n)\bullet m\!+\!1)$ and $\hat\frak r (m\bullet m\!+\!1)\,,$ assume  that $Q$ is normal in $R\,;$ in this case, either $Q\cong V$ and, since $\rho (V)$ is fully normalized in $\F\,,$ according to 2.6.3 there is an  $\F\-$morphism $\hat\rho_Q\,\colon N_P (Q) \to N_P(W)$ extending the composition with  $\rho$ of such an isomorphism, which allows us to set  $U = \hat\rho_Q (R)\.W\,,$ or  $Q\not\cong V$ and we set $U = R\,;$ then, we consider the chains
$$\hat\frak q^R : \Delta_{n+2}\too \F^{^\frak Y}\qq 
\hat\frak r^Q : \Delta_{m+2}\too \F^{^\frak Y}
\eqno £2.9.6\phantom{.}$$
respectively extending $\hat\frak q$ and $\hat\frak r\,,$ fulfilling 
$$\hat\frak q^R(n+2) = U = \hat\frak r^Q(m+2)
\eqno £2.9.7\phantom{.}$$ 
and respectively mapping $(n\!+\!1\bullet n\!+\!2)$ and 
$(m\!+\!1\bullet m\!+\!2)$ either on the inclusion map $W\to \hat\rho_Q(R)\.W$ and, whenever $Q\not= R\,,$
 on the $\F\-$morphism $R\to \hat\rho_Q (R)\.W$ induced  by $\hat\rho_Q$ if $Q\cong V\,,$ or on the group homomorphism $\frak q (n)\to R$ determined by $(\nu_n)^{-1}$ and on ${\rm id}_R$ whenever $Q\not\cong V\,.$

\smallskip
Then, we have evident $\ch^*(\F^{^\frak Y})\-$morphisms
$$(\hat\frak q^R,\Delta_{n+2})\too (\hat\frak q,\Delta_{n+1}) \qq 
(\hat\frak r^Q,\Delta_{m+2})\too (\hat\frak r,\Delta_{m+1})
\eqno £2.9.8\phantom{.}$$
and, considering  the maps
$$\Delta_{n+2}\buildrel \sigma_n\over\longleftarrow\Delta_1\buildrel \sigma_m\over\too \Delta_{m+2}\qq \Delta_{n+1}\buildrel \tau_n\over
\longleftarrow\Delta_0\buildrel \tau_m\over\too \Delta_{m+1}
\eqno £2.9.9\phantom{.}$$
 induced by the sum with $n+1$ and $m+1$ respectively, the  
 $\ch^*(\F^{^\frak Y})\-$morphisms above  determine the following 
 $\ch^*(\F^{^\frak X})\-$morphisms
$$(\hat\frak q^R\!\circ \sigma_n,\Delta_{1})\too (\hat\frak q\circ \tau_n,
\Delta_{0}) \!\!\qq\!\!  (\hat\frak r^Q\!\circ \sigma_m,\Delta_{1})\too (\hat\frak r \circ
\tau_m,\Delta_{0})
\eqno £2.9.10.$$
Thus, the functor $\widehat\aut_{\F^{^\frak X}}$ maps these morphisms on 
$k^*\-$group homomorphisms
$$\hat\F (\hat\frak q^R\!\circ \sigma_n)\too 
\hat\F (\hat\frak q \circ  \tau_n) \qq 
\hat\F (\hat\frak r^Q\!\circ \sigma_m)\too 
\hat\F (\hat\frak r \circ \tau_m) 
\eqno £2.9.11.$$

\smallskip
But note that $\F(\hat\frak q^R)\,,$ $\F(\hat\frak q)\,,$   $\F(\hat\frak r^Q)$ 
and $\F (\hat r)$ are respectively contained in $\F (\hat\frak q^R\!\circ \sigma_n)\,,$ 
$\F (\hat\frak q \circ \tau_n)\,,$ $\F (\hat\frak r^Q\!\circ \sigma_m)$
and $\F (\hat\frak r \circ \tau_m) \,,$ and that,  considering the corresponding
$k^*\-$subgroups, the homomorphisms~£2.9.11 induce $k^*\-$group homomorphisms
$$\hat\F (\hat\frak q^R)\too \hat\F (\hat\frak q ) \qq \hat\F (\hat\frak r^Q) \too 
\hat\F (\hat\frak r) 
\eqno £2.9.12;$$
moreover, the right-hand one is surjective  and it is quite clear that
$$\hat\F (\hat\frak r^Q)\i \hat\F (\hat\frak q^R)\i \hat\F(U)
\eqno £2.9.13.$$
 Consequently, we get a {\it unique\/} $k^*\-$group homomorphism $\widehat\aut_{\F^{^\frak Y}}(\nu,\delta)$ from $\hat\F (\frak r)$
 to $\hat\F (\frak q)$ such that the following diagram is commutative
$$\matrix{\hat\F (\hat\frak r)\hskip-8pt&\leftarrow\hat\F (\hat\frak r^Q)\i 
&\hskip-8pt\hat\F (\hat\frak q^R)\cr
\downarrow\hskip-8pt&\phantom{\bigg\downarrow} &\hskip-8pt\downarrow\cr
 \hat\F (\frak r)\hskip-8pt&&\hskip-8pt\hat\F (\hat\frak q) \cr
\phantom{\bigg\downarrow}& \searrow&\hskip-8pt\downarrow\cr
& &\hskip-8pt\hat\F (\frak q)\cr}
\eqno £2.9.14.$$

\smallskip
Indeed, if $Q\not\cong V$ or $Q = R\,,$ these statements are clear; if $Q\cong V$ 
and $Q\not= R$ then, since $\F (\hat\frak r)$ stabilizes $Q\,,$ the action on $Q$ of any $\sigma\in \F (\hat\frak r)$ can be extended to an $\F\-$automorphism $\hat\sigma$ of $W$ (cf.~£2.6.2); moreover, since $W$ is normal in $ \hat\rho_Q(R)\.W$ and $W$ is fully normalized in $\F\,,$ $\iota_W^P\circ \hat\sigma$ can still be extended to 
an $\F\-$morphism 
$$\hat\zeta : \hat\rho_Q(R)\.W\too P
\eqno £2.9.15\phantom{.}$$
which, {\it via\/} $\hat\rho_Q\,,$ induces an $\F\-$morphism $\psi\,\colon R\to P\,;$ then, the restriction of this homomorphism to $Q$ coincides with the restriction of 
$\sigma$ and therefore, since $Q$ is $\F\-$selfcentralizing, it follows from Proposition~4.6 in [5] that $\psi$ and $\iota_R^P\circ\sigma$ are $Z(Q)\-$conjugate. In particular, we get
$$\hat\zeta \big(\hat\rho_Q(R)\.W\big) = \hat\rho_Q(R)\.W
\eqno £2.9.16\phantom{.}$$
and thus, up to modifying our choices, we may assume that the corresponding 
$\F\-$automorphism $\hat{\hat\sigma}$ of $\hat\rho_Q(R)\.W$ extends $\sigma$
and $\hat\sigma\,;$ then, $\hat{\hat\sigma}$ belongs to~$\F (\hat\frak r^Q)$ and
to $\hat\F (\hat\frak q^R)\,,$ which proves the surjectivity of the right-hand homomorphism in~£2.9.12 and easily implies the left-hand inclusion in~£2.9.13.

\smallskip
 Consider another $\ch^*(\F^{^\frak Y})\-$morphism $(\mu,\varepsilon)\,\colon 
 (\frak t,\Delta_\ell)\to (\frak r,\Delta_m)\,,$ and  respectively denote by $T\,,$ $R$ and $Q$ the images of $\frak t(\ell)\,,$ $\frak t \big(\varepsilon (m)\big)$ and $\frak t \big((\varepsilon\circ \delta)(n)\big)$  in~$\hat \frak t (\ell +1)\,;$ assume that $R$ and $Q$ are both normal in  $T\,;$ as above, we consider the chains $\hat\frak q^R\,,$
$\hat \frak q^T\,,$ $\hat\frak r^Q\,,$ $\hat\frak r^T\,,$ $\hat\frak t^Q$ and $\hat \frak t^R\,,$
and moreover we need the chains
$$\eqalign{\hat\frak q^{R,T} &: \Delta_{n+3}\too \F^{^\frak Y}\cr 
\hat\frak r^{T,Q} &: \Delta_{m+3}\too \F^{^\frak Y}\cr
 \hat\frak t^{R,Q} &: \Delta_{\ell+3}\too \F^{^\frak Y}\cr}
\eqno £2.9.17\phantom{.}$$
respectively extending $\hat\frak q^R\,,$ $\hat\frak r^T$ and $\hat \frak t^R\,,$ fulfilling
$$\hat\frak q^{R,T}(n+3) = \hat\frak r^{T,Q}(m+3) = \hat\frak t^{R,Q}(\ell +3) 
= \hat\frak q^T (n +2)
\eqno £2.9.18,$$ 
and respectively mapping  $(n\!+\!2\bullet n\!+\!3)\,,$  
$(m\!+\!2\bullet m\!+\!3)$ and $(\ell\!+\!2\bullet \ell\!+\!3)$ on the inclusion 
$\hat\frak q^R (n\!+\!2)\i \hat\frak q^T (n\!+\!2)$ in both cases,  either
  on the homomorphism $T\to \hat\rho_Q(T)\.W$ induced by $\hat\rho_Q$ 
  if $Q\cong V$ and $Q\not= R\,,$ or on ${\rm id}_{\hat\frak q^T (n +2)}$ if 
$Q\not\cong V$ or $Q = R\,,$ and once again,  either  on the homomorphism $T\to \hat\rho_Q(T)\.W$ induced by $\hat\rho_Q$ if $Q\cong V$ and $Q\not= R\,,$ or on
${\rm id}_{\hat\frak q^T (n +2)}$ if  $Q\not\cong V$ or $Q = R\,.$

 As above, it is easily checked that applying the functor $\widehat\aut_{\F^{^\frak X}}$ 
 to the evident $\ch^*(\F^{^\frak X})\-$morphisms, we get
$k^*\-$group homomorphisms
$$\eqalign{\widehat\aut_{\F^{^\frak X}} &: \hat\F (\hat\frak q^{R,T})\too \hat\F (\hat\frak q^R)\cr
\widehat\aut_{\F^{^\frak X}} &: \hat\F (\hat\frak r^{T,Q})\too \hat\F (\hat\frak r^{Q})\cr
\widehat\aut_{\F^{^\frak X}} &:  \hat\F (\hat\frak t^{R,Q})\too \hat\F (\hat\frak t^{R})\cr}
\eqno £2.9.19\phantom{.}$$
and moreover it is quite clear that $\hat\F (\hat\frak t^{R,Q}) = \hat\F (\hat\frak t^{Q})\,.$
Consequently, the functoriality of $\widehat\aut_{\F^{^\frak X}}$ guarantees the commutativity of
the following diagram
$$\matrix{\hat\F (\hat\frak t)&\hskip-5pt\leftarrow&\hskip-5pt
\hat\F (\hat\frak t^{R,Q})&\hskip-5pt  = &\hskip-5pt\hat\F (\hat\frak t^Q)&\hskip-5pt&\hskip-5pt\i &\hskip-5pt&\hskip-5pt\hat\F (\hat\frak q^T)&\hskip-5pt = &\hskip-5pt\hat\F (\hat\frak q^T)\cr
\Vert&\hskip-5pt&\hskip-5pt\downarrow&\hskip-5pt\phantom{\big\downarrow}&\hskip-5pt&\hskip-5pt&\hskip-5pt\phantom{\Big\downarrow}&\hskip-5pt&\hskip-5pt\cup\cr
\hat\F (\hat\frak t)&\hskip-5pt\leftarrow&\hskip-5pt\hat\F (\hat\frak t^R)&\hskip-5pt\i
&\hskip-5pt\hat\F (\hat\frak r^T)&\hskip-5pt\leftarrow&\hskip-5pt\hat\F 
(\hat\frak r^{T,Q})&\hskip-5pt\i
&\hskip-5pt\hat\F (\hat\frak q^{R,T})\cr
\downarrow&\hskip-5pt&\hskip-5pt&\hskip-5pt\phantom{\Big\downarrow}&\hskip-5pt\downarrow&\hskip-5pt&\hskip-5pt\downarrow&\hskip-5pt&\hskip-5pt\downarrow&\hskip-5pt&\hskip-5pt\downarrow\cr
\hat\F (\frak t)&\hskip-5pt&\hskip-5pt&\hskip-5pt&\hskip-5pt \hat\F (\hat\frak r)&\hskip-5pt\leftarrow&\hskip-5pt
\hat\F (\hat\frak r^Q)&\hskip-5pt\i  &\hskip-5pt\hat\F (\hat\frak q^R)\cr 
&\hskip-5pt&\hskip-5pt\searrow&\hskip-5pt&\hskip-5pt\downarrow&\hskip-5pt\phantom{\Bigg\downarrow}&\hskip-5pt&\hskip-5pt &\hskip-5pt\downarrow\cr
 &\hskip-5pt&\hskip-5pt&\hskip-5pt&\hskip-5pt\hat\F (\frak r)&\hskip-5pt&\hskip-5pt&\hskip-5pt&\hskip-5pt\hat\F (\hat\frak q)&\hskip-5pt =&\hskip-5pt\hat\F (\hat\frak q)\cr 
&\hskip-5pt&\hskip-5pt&\hskip-5pt&\hskip-5pt\phantom{\Bigg\downarrow}&\hskip-5pt&\hskip-5pt \searrow&\hskip-5pt&\hskip-5pt\downarrow\cr
&\hskip-5pt&\hskip-5pt&\hskip-5pt &\hskip-5pt&\hskip-5pt&\hskip-5pt&\hskip-5pt&\hskip-5pt\hat\F (\frak q)\cr}
\eqno £2.9.20;$$
thus, by uniqueness, in this case we obtain 
$$\widehat\aut_{\F^{^\frak Y}}(\nu,\delta)\circ\widehat\aut_{\F^{^\frak Y}}(\mu,\varepsilon) =\widehat\aut_{\F^{^\frak Y}}\big((\nu,\delta)\circ (\mu,\varepsilon)\big)
\eqno £2.9.21.$$

\smallskip
 Secondly, assume that the image of $\frak r\big(\delta (n)\big)$ by 
 $\frak r (\delta (n)\!\bullet\! m)$ is not normal in $\frak r (m)\,;$ let $m'$ be the maximal element in $\Delta_m - \Delta_{\delta (n) -1}$ such that  the image of 
 $\frak r\big(\delta (n)\big)$  by $\frak r (\delta (n)\!\bullet\! m')$ is normal in~$\frak r (m')$ and denote by $R_{(\nu,\delta)}$ the normalizer of the image of $\frak r\big(\delta (n)\big)$
 in $\frak r (m'+1)\,,$ by $\frak r_{(\nu,\delta)}
\colon \Delta_{m+1}\to \F^{^\frak Y}$ the functor fulfilling 
$$\frak r_{(\nu,\delta)}\circ \delta^m_{m'+1} = \frak r\qq
\frak r_{(\nu,\delta)}(m'+1) = R_{(\nu,\delta)}
\eqno £2.9.22\phantom{.}$$
 and mapping $(m'\!+\!1\bullet m'\!+\!2)$ on the inclusion map $R_{(\nu,\delta)}\to\frak r (m'+1)\,,$ and by $\frak r'_{(\nu,\delta)}$ the restriction of $\frak r_{(\nu,\delta)}$
to $\Delta_{m'+1}\,;$ then, it is quite clear that $\F (\frak r_{(\nu,\delta)}) 
= \F (\frak r)$ and it is easily checked that $\hat\F (\frak r_{(\nu,\delta)}) 
= \hat\F (\frak r)\,;$ moreover, we have an evident $\ch^*(\F^{^\frak Y})\-$morphism  
$$(\nu',\delta') : (\frak r'_{(\nu,\delta)},\Delta_{m'+1})\too (\frak q,
\Delta_{n})
\eqno £2.9.23\phantom{.}$$
such that 
$$(\nu',\delta')\circ ({\rm id}_{\frak r'_{(\nu,\delta)}},\iota^m_{m'}) =
(\nu,\delta)\circ  ({\rm id}_\frak r,\delta^m_{m'+1})
\eqno £2.9.24\phantom{.}$$
 where $\iota^m_{m'}\colon\Delta_{m'+1} \to \Delta_{m+1}$ denotes the natural inclusion, and in~£2.9.14 we already have defined 
$\widehat\aut_{\F^{^\frak Y}}({\rm id}_\frak r,\delta^m_{m'+1}) = 
{\rm id}_{\hat\F (\frak r)}$ and $\widehat\aut_{\F^{^\frak Y}} (\nu',\delta')\,;$ on the other hand, arguing by induction on $\vert\frak r (m)\vert/\vert\frak q (n)\vert$, we may assume that $\widehat\aut_{\F^{^\frak Y}}
({\rm id}_{\frak r'_{(\nu,\delta)}},\iota^m_{m'})$ is already defined and then we set
$$\widehat\aut_{\F^{^\frak Y}}(\nu,\delta) =
\widehat\aut_{\F^{^\frak Y}} (\nu',\delta') \circ\widehat\aut_{\F^{^\frak Y}}
({\rm id}_{\frak r'_{(\nu,\delta)}},\iota^m_{m'})
\eqno £2.9.25.$$

\smallskip
  For another $\ch^*(\F^{^\frak Y})\-$morphism $(\mu,\varepsilon)\,\colon (\frak t,\Delta_\ell)\to (\frak r, \Delta_m)\,,$  we claim that
$$\widehat\aut_{\F^{^\frak Y}}(\nu,\delta)\circ\widehat\aut_{\F^{^\frak Y}}(\mu,\varepsilon) = \widehat\aut_{\F^{^\frak Y}}\big((\nu,\delta)\circ (\mu,\varepsilon)\big)
\eqno £2.9.26;$$
we argue by induction firstly on $\vert\frak t (\ell)\vert/\vert\frak q (n)\vert$ and after
on $\vert\frak t (\ell)\vert/\vert\frak r (m)\vert\,.$ First of all, we assume that the image of $\frak r\big(\delta (n)\big)$ in~$\frak r (m)$ by $\frak r\big(\delta (n)\bullet m\big)$ is not normal; with the notation above, denote by $\ell'$ the maximal element in 
$\Delta_\ell - \Delta_{(\varepsilon\circ\delta) (n) -1}$ such that  the image of $\frak t\big((\varepsilon\circ\delta) (n)\big)$ 
by  $\frak t \big((\varepsilon\circ\delta) (n)\bullet \ell'\big)$ is normal in
$\frak t (\ell')\,;$ then, it is clear that
$\varepsilon (m')\le \ell' < \varepsilon (m)$ and easily checked that we have a $\ch^*(\F^{^\frak Y})\-$morphism
$$(\mu_{(\nu,\delta)},\varepsilon_{(\nu,\delta)}) : (\frak t_{(\nu,\delta) \circ
(\mu,\varepsilon)},\Delta_{\ell +1})\too (\frak r_{(\nu,\delta)} ,\Delta_{m +1})
\eqno £2.9.27\phantom{.}$$
such that 
$$({\rm id}_\frak r,\delta^m_{m'+1})\circ (\mu_{(\nu,\delta)}, \varepsilon_{(\nu,\delta)}) = (\mu,\varepsilon)\circ ({\rm id}_\frak t,\delta^\ell_{\ell'+1})
\eqno £2.9.28,$$
 that $\varepsilon_{(\nu,\delta)} (m' +1) =\ell' +1$ and that 
 $(\mu_{(\nu,\delta)})_{m'+1}$ is determined by~$\mu_{m'+1}$ and
  $\frak t \big(\ell'\! +\!1 \bullet \varepsilon (m' \!+\!1)\big)\,;$ moreover, we consider the corresponding restriction
$$(\mu'_{(\nu,\delta)},\varepsilon'_{(\nu,\delta)}) : (\frak t'_{(\nu,\delta) \circ
(\mu,\varepsilon)},\Delta_{\ell' +1})\too (\frak r'_{(\nu,\delta)} ,\Delta_{m' +1})
\eqno £2.9.29\phantom{.}$$
which obviously fulfills
$$({\rm id}_{\frak r'_{(\nu,\delta)}},\iota^m_{m'})\circ
(\mu_{(\nu,\delta)},\varepsilon_{(\nu,\delta)}) =  (\mu'_{(\nu,\delta)},\varepsilon'_{(\nu,\delta)})
\circ ({\rm id}_{\frak t'_{(\nu,\delta)\circ (\mu,\varepsilon)}},\iota^\ell_{\ell'})
\eqno £2.9.30.$$

\smallskip
Now, it is easily checked that the composition $(\nu',\delta')\circ
(\mu'_{(\nu,\delta)},\varepsilon'_{(\nu,\delta)})$ coincides with the corresponding morphism~£2.9.23 for the $\ch^*(\F^{^\frak Y})\-$morphism $(\nu,\delta)\circ (\mu,\varepsilon)$ and therefore, by the very definition~£2.9.25, we have
$$\eqalign{\widehat\aut_{\F^{^\frak Y}}&\big((\nu,\delta)\circ (\mu,\varepsilon)\big)\cr
&=  \widehat\aut_{\F^{^\frak Y}} \big((\nu',\delta')\circ (\mu'_{(\nu,\delta)}, \varepsilon'_{(\nu,\delta)})\big) \circ
\widehat\aut_{\F^{^\frak Y}} \big({\rm id}_{\frak t'_{(\nu,\delta) \circ (\mu,\varepsilon)}},\iota^\ell_{\ell'}\big)\cr}
\eqno £2.9.31;$$
but, since $\vert R_{(\nu,\delta)}\vert/\vert \frak q (n)\vert < \vert\frak t (\ell)\vert/
\vert\frak q (n)\vert\,,$ it follows from the induction hypo-thesis that
$$\widehat\aut_{\F_{\!\rm nc}} \big((\nu',\delta') \circ
(\mu'_{(\nu,\delta)},\varepsilon'_{(\nu,\delta)})\big) = 
\widehat\aut_{\F^{^\frak Y}}(\nu',\delta')\circ\widehat\aut_{\F^{^\frak Y}}(\mu'_{(\nu,\delta)},\varepsilon'_{(\nu,\delta)})
\eqno £2.9.32;$$
similarly, since we have $\vert \frak t (\ell)\vert/\vert R_{(\nu,\delta)}\vert < 
\vert\frak t (\ell)\vert/\vert\frak q (n)\vert$ and
$$\widehat\aut_{\F^{^\frak Y}}(\mu_{(\nu,\delta)},\varepsilon_{(\nu,\delta)}) = \widehat\aut_{\F^{^\frak Y}} (\mu,\varepsilon)
\eqno £2.9.33,$$ 
we still get
$$\eqalign{&\widehat\aut_{\F^{^\frak Y}}\big((\nu,\delta)\circ (\mu,\varepsilon)\big)\cr 
& = \widehat\aut_{\F^{^\frak Y}}(\nu',\delta') \circ \widehat\aut_{\F^{^\frak Y}}
(\mu'_{(\nu,\delta)},\varepsilon'_{(\nu,\delta)}) \circ \widehat\aut_{\F^{^\frak Y}} ({\rm id}_{\frak t'_{(\nu,\delta)\circ (\mu,\varepsilon)}},\iota^\ell_{\ell'})\cr
&= \widehat\aut_{\F^{^\frak Y}}(\nu',\delta') \circ \widehat\aut_{\F^{^\frak Y}}
\big(({\rm id}_{\frak r'_{(\nu,\delta)}},\iota^m_{m'})\circ
(\mu_{(\nu,\delta)},\varepsilon_{(\nu,\delta)})\big)\cr  
&= \widehat\aut_{\F^{^\frak Y}}(\nu,\delta)\circ\widehat\aut_{\F^{^\frak Y}}
(\mu,\varepsilon)\,.\cr}
\eqno £2.9.34.$$

\smallskip
 Finally, we may assume that the image of $\frak r\big(\delta (n)\big)$  by 
 $\frak r\big(\delta (n)\bullet m\big)$ is normal in~$\frak r (m)\,,$ so that the image of $\frak t\big((\varepsilon\circ\delta)(n)\big)$ by
$\frak t\big((\varepsilon\circ\delta)(n)\!\bullet\!\varepsilon(m)\big)$ is normal 
in~$\frak t\big(\varepsilon(m)\big)\,;$ in particular, denoting by $\ell'$ the maximal element in~$\Delta_\ell - \Delta_{(\varepsilon\circ\delta) (n) -1}$ such that  the image of $\frak t\big((\varepsilon\circ\delta) (n)\big)$ 
by  $\frak t \big((\varepsilon\circ\delta) (n)\bullet \ell'\big)$ is normal 
in~$\frak t (\ell')\,,$ we have $\varepsilon (m)\le \ell'\,.$ If $\ell' = \ell$ then we may assume that the image of $\frak t\big(\varepsilon (m)\big)$ is not normal 
in~$\frak t (\ell)$  and, denoting by~$\ell''\ge \varepsilon (m)$ the maximal element
 in~$\Delta_\ell$ such that  the image of 
$\frak t\big(\varepsilon (m)\big)$ by $\frak t (\varepsilon (m)\!\bullet\! \ell'')$ is normal
in~$\frak r (\ell'')\,,$ by our very definition (cf.~£2.9.25) we have
$$\widehat\aut_{\F^{^\frak Y}}(\mu,\varepsilon) =
\widehat\aut_{\F^{^\frak Y}} (\mu',\varepsilon') \circ\widehat\aut_{\F^{^\frak Y}}
({\rm id}_{\frak t'_{(\mu,\varepsilon)}},\iota^\ell_{\ell''})
\eqno £2.9.35;$$
but, according to equality~£2.9.21, we have
$$\widehat\aut_{\F^{^\frak Y}}(\nu,\delta)\circ \widehat\aut_{\F^{^\frak Y}} (\mu',\varepsilon') =
\widehat\aut_{\F^{^\frak Y}}\big((\nu,\delta)\circ (\mu',\varepsilon')\big)
\eqno £2.9.36;$$
hence, since in the compositions of $(\nu,\delta)$ with $(\mu,\varepsilon)$ and
of $\big((\nu,\delta)\circ (\mu',\varepsilon')\big)$ with $ ({\rm id}_{\frak t'_{(\mu,\varepsilon)}},\iota^\ell_{\ell''})$ the first induction indices coincide with each other and the second ones strictly decreasse, it follows from the induction hypothesis that
$$\eqalign{\widehat\aut_{\F^{^\frak Y}}(\nu,\delta)&\circ
 \widehat\aut_{\F^{^\frak Y}} (\mu,\varepsilon)\cr
&= \widehat\aut_{\F^{^\frak Y}}(\nu,\delta)\circ \widehat\aut_{\F^{^\frak Y}}(\mu',\varepsilon')
\circ\widehat\aut_{\F^{^\frak Y}} ({\rm id}_{\frak t'_{(\mu,\varepsilon)}},\iota^\ell_{\ell''})
\cr
&= \widehat\aut_{\F^{^\frak Y}}\big((\nu,\delta)\circ (\mu',\varepsilon')\big)
\circ\widehat\aut_{\F^{^\frak Y}} ({\rm id}_{\frak
t'_{(\mu,\varepsilon)}},\iota^\ell_{\ell''})\cr
&= \widehat\aut_{\F^{^\frak Y}}\big((\nu,\delta)\circ (\mu,\varepsilon)\big)\cr}
\eqno £2.9.37.$$

\smallskip
Otherwise, we have a $\ch^*(\F^{^\frak Y})\-$morphism
$$(\mu'_{(\nu,\delta)},\varepsilon'_{(\nu,\delta)}) : (\frak t'_{(\nu,\delta) \circ
(\mu,\varepsilon)},\Delta_{\ell' +1})\too (\frak r ,\Delta_m)
\eqno £2.9.38\phantom{.}$$
fulfilling
$$(\mu'_{(\nu,\delta)},\varepsilon'_{(\nu,\delta)})\circ 
({\rm id}_{\frak t'_{(\nu,\delta) \circ (\mu,\varepsilon)}},\iota^\ell_{\ell'}) =
(\mu,\varepsilon) \circ ({\rm id}_\frak t, \delta^\ell_{\ell'+1})
\eqno £2.9.39;$$
 as above, it is easily checked that the composition $(\nu,\delta)\circ
(\mu'_{(\nu,\delta)},\varepsilon'_{(\nu,\delta)})$ coincides with the corresponding morphism~£2.9.23 for the $\ch^*(\F^{^\frak Y})\-$morphism $(\nu,\delta)\circ (\mu,\varepsilon)$ and therefore, by the very definition~£2.9.25, we have
$$\eqalign{\widehat\aut_{\F^{^\frak Y}}&\big((\nu,\delta)\circ (\mu,\varepsilon)\big)\cr
&= \widehat\aut_{\F^{^\frak Y}} \big((\nu,\delta)\circ (\mu'_{(\nu,\delta)},\varepsilon'_{(\nu,\delta)})\big) \circ
\widehat\aut_{\F^{^\frak Y}} \big({\rm id}_{\frak t'_{(\nu,\delta) \circ (\mu,\varepsilon)}},\iota^\ell_{\ell'}\big)\cr}
\eqno £2.9.40;$$
since $\ell'\not= \ell$ and $\widehat\aut_{\F_{\!\rm nc}}({\rm id}_\frak t,\delta^\ell_{\ell'+1}) = {\rm id}_{\hat\F (\frak t)}\,,$
it follows from the induction hypothesis applied to the composition of $(\nu,\delta)$
with $(\mu'_{(\nu,\delta)},\varepsilon'_{(\nu,\delta)})$ that
$$\widehat\aut_{\F^{^\frak Y}} \big((\nu,\delta)\circ
(\mu'_{(\nu,\delta)},\varepsilon'_{(\nu,\delta)})\big) = \widehat\aut_{\F^{^\frak Y}}
(\nu,\delta)\circ \widehat\aut_{\F^{^\frak Y}}(\mu'_{(\nu,\delta)},\varepsilon'_{(\nu,\delta)})
\eqno £2.9.41;$$
moreover, if  $\vert \frak q (n)\vert < \vert \frak r (m)\vert\,,$ we can apply the induction hypothesis to both members of equality~£2.9.39 and then  we get
$$\widehat\aut_{\F^{^\frak Y}}(\mu'_{(\nu,\delta)},\varepsilon'_{(\nu,\delta)})\circ 
\widehat\aut_{\F^{^\frak Y}} \big({\rm id}_{\frak t'_{(\nu,\delta) \circ
(\mu,\varepsilon)}},\iota^\ell_{\ell'}\big)
 = \widehat\aut_{\F^{^\frak Y}}(\mu,\varepsilon)
 \eqno £2.9.42;$$
consequently, once again we have
$$\widehat\aut_{\F_{\!\rm nc}}\big((\nu,\delta)\circ (\mu,\varepsilon)\big) = 
\widehat\aut_{\F_{\!\rm nc}} (\nu,\delta)\circ \widehat\aut_{\F_{\!\rm nc}}(\mu,\varepsilon)
 \eqno £2.9.43.$$

 \smallskip
 If $\vert \frak q (n)\vert = \vert \frak r (m)\vert$ then it follows from the  definition of $\widehat\aut_{\F^{^\frak Y}}(\mu,\varepsilon)$ and $\widehat\aut_{\F^{^\frak Y}}\big((\nu,\delta)\circ (\mu,\varepsilon)\big)$ (cf.~£2.9.25) that $\ell'$ coincides with both indices,
 that we get $\frak t'_{(\mu,\varepsilon)} = \frak t'_{(\nu,\delta)\circ  (\mu,\varepsilon)}$ and that the homomorphism~£2.9.23
 $$(\frak t'_{(\nu,\delta)\circ  (\mu,\varepsilon)},\Delta_{\ell' +1})
 \too (\frak q,\Delta_n)
 \eqno £2.9.44\phantom{.}$$
 corresponding to the  composition $(\nu,\delta)\circ (\mu,\varepsilon)$ coincides
 with $(\nu,\delta)\circ (\mu',\varepsilon')\,;$ at this point, we can apply 
 equality~£2.9.21 to obtain
 $$\widehat\aut_{\F^{^\frak Y}}(\nu,\delta)\circ\widehat\aut_{\F^{^\frak Y}}
 (\mu',\varepsilon') =\widehat\aut_{\F^{^\frak Y}}\big((\nu,\delta)\circ (\mu',\varepsilon')\big)
\eqno £2.9.45;$$
then, composing this equality with $\widehat\aut_{\F^{^\frak Y}}
({\rm id}_{\frak t'_{(\mu,\varepsilon)}},\iota^\ell_{\ell'})\,,$ from definition 
£2.9.25 we get
$$\widehat\aut_{\F^{^\frak Y}}(\nu,\delta)\circ\widehat\aut_{\F^{^\frak Y}}
 (\mu,\varepsilon) =\widehat\aut_{\F^{^\frak Y}}\big((\nu,\delta)\circ (\mu,\varepsilon)\big)
\eqno £2.9.46.$$
We are done.

\bigskip
\bigskip
\noindent
{\bf £3\phantom{.} Ordinary Grothendieck group of a folded Frobenius $P\-$category}
\bigskip

£3.1\phantom{.} Our settting is a finite $p\-$group $P\,,$ a Frobenius 
$P\-$category~$\F$ and a functor
$$\widehat\aut_{\F^{^{\rm sc}}} : \ch^*(\F^{^{\rm sc}})\too k^*\-\Gr
\eqno £3.1.1\phantom{.}$$
lifting the functor $\aut_{\F^{^{\rm sc}}}$ (cf.~£2.8); if $\frak q\,\colon \Delta_n\to \F^{^{\rm sc}}$ is an  $\F^{^{\rm sc}}\-$chain, we simply denote by $\hat\F (\frak q)$
the image of $(\frak q,\Delta_n)$ {\it via\/} $\widehat\aut_{\F^{^{\rm sc}}}\,.$ Note that, if $P'$ is a second 
finite $p\-$group, $\F'$ a Frobenius $P'\-$category and $\alpha\,\colon P'\to P$ an 
{\it $(\F',\F)\-$functorial \/} group homomorphism [5, 12.1] mapping any $\F'\-$radical subgroup of $P'$ (cf.~£2.4) on a  $\F\-$selfcentralizing subgroup of $P\,,$ the {\it Frobenius functor\/} $\frak f_\alpha\,\colon \F'\to \F$ [5, 12.1] 
and the lifting  $\widehat\aut_{\F^{^{\rm sc}}}$ determine a functor
$$\widehat\aut_{\F'^{^{\rm rd}}} : \ch^*(\F'^{^{\rm rd}})\too 
\ch^*(\F^{^{\rm sc}})\too k^*\-\Gr
\eqno £3.1.2.$$
lifting $\aut_{\F'^{^{\rm rd}}}$ and then, according to Theorem~£2.9, this lifting can be
 uniquely extended to a functor  $\widehat\aut_{\F'^{^{\rm sc}}}$ lifting 
 $\aut_{\F'^{^{\rm sc}}}\,; $ in particular, for any  subgroup $U$ of $P$ fully $K\-$normalized in $\F$ for some subgroup $K$ of ${\rm Aut}(U)\,,$ it follows from 
 Lemma~£2.5 and Theeorem~£2.9 that our {\it folded Frobenius $P\-$category\/} induces a  {\it folded Frobenius 
 $N_P^K(U)\-$category\/} formed by $N_\F^K (U)$ and by~$\widehat\aut_{N_\F^K (U)^{^{\rm sc}}}\,.$

\medskip
£3.2\phantom{.} For any $\F\-$selfcentralizing subgroup $Q$ of~$P$ fully normalized in $\F\,,$ in [5, Theorem~18.6] we prove the existence, and the uniqueness up to ismorphisms, of a finite group $\L (Q)$ --- the  {\it $\F\-$localizer\/} of~$Q$ --- such that $N_P (Q)$ is a Sylow $p\-$subgroup of 
$\L (Q)$ and $Z(Q)$ a normal subgroup fulfilling $\L(Q)/Z(Q)\cong \F (Q)\,.$
More generally,  in [5, Proposition~18.19] we prove the existence of the 
{\it $\F\-$localizing functor\/}
$$\loc_{\F^{^{\rm sc}}} : \ch^*(\F^{^{\rm sc}})\too \widetilde{\Loc}
\eqno £3.2.1\phantom{.}$$
which lifts $\aut_{\F^{^{\rm sc}}}$ {\it via\/} the canonical functor 
$\widetilde{\Loc}\to \Gr$ sending any $\widetilde{\Loc}\-$object $(L,Z)$ to $L/Z\,,$ and maps any  $\F\-$chain $\frak q\,\colon \Delta_n\to \F$ such that $\frak q(n)$ is fully normalized in $\F\,,$ on the pair $\big(\L (\frak q),Z(\frak q(n))\big)$ where
$\L (\frak q)$ denotes the converse image of $\F (\frak q)$ in 
$\L\big(\frak q (n)\big)\,;$ actually, for any  $\F\-$chain $\frak q\,\colon 
\Delta_n\to \F$ we set $\loc_{\F^{^{\rm sc}}}(\frak q,\Delta_n) = \big(\L (\frak q), Z(\frak q(n))\big)$ and identify $N_P (\frak q)$ with the converse image of 
$\F_P (\frak q)$ in $\L(\frak q)$ [5, Proposition~18.16].
\eject

\medskip
£3.3\phantom{.} As in~£1.9 above, we consider the category $k^*\-\widetilde\Loc$ where the objects are the pairs $(\hat L,Z)$ formed by a {\it finite\/} $k^*\-$group $\hat L$ and a normal $p\-$subgroup $Z$ of $\hat L$ and where, coherently, the morphisms from $(\hat L,Z)$ to $(\hat L',Z')$ are the $Z'\-$conjugacy classes of $k^*\-$group
homomorphisms $\hat f\,\colon \hat L\to \hat L'$ fulfilling $\hat f(Z)\i Z'\,.$ Then, consider the functor
$$\widehat{\loc}_{\F^{^{\rm sc}}} : \ch^* (\F^{^{\rm sc}})\too k^*\-\widetilde\Loc
\eqno £3.3.1\phantom{.}$$
mapping  any $\F^{^{\rm sc}}\-$chain $\frak q\,\colon \Delta_n\to \F^{^{\rm sc}}$
such that $\frak q (n)$ is fully normalized in $\F$ on the pair formed by the $k^*\-$group
$$\hat\L (\frak q) = \L (\frak q) \times_{\F (\frak q)} \hat\F(\frak q)
\eqno £3.3.2\phantom{.}$$
and by the $p\-$subgroup lifting $Z\big(\frak q (n)\big)\,,$ and any 
$ \ch^* (\F^{^{\rm sc}})\-$morphism 
$$(\nu,\delta) : (\frak r,\Delta_m)\too  (\frak q,\Delta_n)
\eqno £3.3.3\phantom{.}$$
 where $\frak q (n)$ and $\frak r (m)$ are fully normalized in 
$\F\,,$ on the $Z\big(\frak q (n)\big)\-$conjugacy class of
$$\widehat{\loc}_{(\nu,\delta)} = \loc_\F (\nu,\delta) \times_{\aut_{\F^{^{\rm sc}}}
(\nu,\delta)} \widehat\aut_{\F^{^{\rm sc}}}(\nu,\delta)
\eqno £3.3.4.$$

\medskip
£3.4\phantom{.}  For any central $k^*\-$extension~$\hat G$ of a finite group $G\,,$
recall that we respectively denote by $\G_\K   (\hat G)$ and $\G_k (\hat G)$ the 
{\it scalar extensions\/} from $\Bbb Z$ to~$\O$ of the Grothendieck groups  of~the categories of finitely dimensional $\K_*\hat G\-$ and $k_*\hat G\-$modules; it is well-known that we have {\it contravariant\/} functors
$$\frak g_\K : k^*\-{\Gr}\too \O\-\mod\qq \frak g_k : k^*\-{\Gr}\too \O\-\mod
\eqno  £3.4.1\phantom{.}$$
respectively mapping $\hat G$ on $\G_\K (\hat G)$ and $\G_k (\hat G)\,,$ and 
any $k^*\-$group homomorphism $\hat\varphi\, \,\colon\hat G\to\hat G'$ on the corresponding {\it restriction\/} maps;  it is clear that these restriction maps only depend on  the $\hat G'\-$conjugacy class of $\varphi\,;$
in particular, these functors have an obvious extension to the category 
$k^*\-\widetilde\Loc$ respectively mapping any $k^*\-\widetilde\Loc\-$object 
$(\hat L,Z)$ on $\G_\K (\hat L)$ and $\G_k (\hat L)\,.$ 
Moreover,  recall that the  {\it Brauer decomposition maps\/} define
 a natural map
$$\partial :  \frak g_\K\too \frak g_k
\eqno  £3.4.2\phantom{.}$$
admitting a natural section.

\medskip
£3.5\phantom{.} Now, we consider the composed functor
$$\ch^* (\F^{^{\rm sc}})\buildrel \widehat\loc_{\F^{^{\rm sc}}} \over{\hbox to
30pt{\rightarrowfill}}  k^*\-\widetilde\Loc \buildrel \frak g_\K 
\over{\hbox to 20pt{\rightarrowfill}}  \O\-\mod
\eqno  £3.5.1\phantom{.}$$
and we define the {\it ordinary Grothendieck group of the folded Frobenius $P\-$cate-gory $(\F,
\widehat\aut_{\F^{^{\rm sc}}})$\/} as the inverse limit
$$\G_\K (\F,\widehat\aut_{\F^{^{\rm sc}}}) = \lim_{\longleftarrow}\,
(\frak g_\K\circ \widehat\loc_{\F^{^{\rm sc}}})
\eqno  £3.5.2;$$\break
\eject
\noindent
 it is a finitely generated free $\O\-$module. Since the natural map
from $\widehat\loc_{\F^{^{\rm sc}}}$ to $\widehat\aut_{\F^{^{\rm sc}}}$ induces 
a natural isomorphism $\frak g_k\circ \widehat\loc_{\F^{^{\rm sc}}}\cong \frak g_k
\circ \widehat\aut_{\F^{^{\rm sc}}}\,,$ the natural map~£3.4.2 determines a {\it decomposition map\/}
$$\partial_{(\F,\widehat\aut_{\F^{^{\rm sc}}})} : \G_\K (\F,\widehat\aut_{\F^{^{\rm sc}}})\too \G_k (\F,\widehat\aut_{\F^{^{\rm sc}}})
\eqno £3.5.3\phantom{.}$$
admitting a section. Note that, if $\F$ admits an $\F\-$normal $\F\-$selfcentralizing subgroup $Q$ of $P$ then we have evident isomorphisms [5, Proposition~19.5]
$$\G_\K (\F,\widehat\aut_{\F^{^{\rm sc}}}) \cong \G_\K\big(\hat\L (Q)\big)\qq 
\G_k (\F,\widehat\aut_{\F^{^{\rm sc}}}) \cong \G_k\big(\hat\L (Q)\big)
\eqno £3.5.4\phantom{.}$$ 
and $\partial_{(\F,\widehat\aut_{\F^{^{\rm sc}}})}$ is the corresponding {\it Brauer decomposition map\/}.

\medskip
£3.6\phantom{.} With the notation in~£3.1 above, in order to relate the ordinary Grothendieck
groups $\G_\K (\F,\widehat\aut_{\F^{^{\rm sc}}})$ and $\G_\K (\F',\widehat\aut_{\F'^{^{\rm sc}}})$ we prove below that, denoting by $\frak i\,\colon 
\F^{^{\rm rd}}\to \F^{^{\rm sc}}$ the inclusion functor, the ordinary Grothendieck
group $\G_\K (\F,\widehat\aut_{\F^{^{\rm sc}}})$ can also be defined from the functor
$$\widehat\loc_{\F^{^{\rm rd}}} = \widehat\loc_{\F^{^{\rm sc}}}\circ \ch^*(\frak i)
\eqno £3.6.1,$$
as in the case of the modular Grothendieck group [5, Proposition~14.6]. Recall that, for any category $\frak C$ and any {\it contravariant\/} functor $\frak m\,\colon \ch^*(\frak C)\to \O\-\mod\,,$ the inverse limit 
${\displaystyle \lim_{\longleftarrow}}\,\frak m$ coincides
with the {\it $0\-$cohomology group\/} of the {\it dif-ferential complex\/} given by 
the {\it differential maps\/} [5,~A3.11.2]
$$d_n : \prod_{\frak q} \frak m (\frak q,\Delta_n)\too
\prod_{\frak r}\frak m (\frak r,\Delta_{n+\!1})
\eqno £3.6.2,$$
where $\frak q$ and~$\frak r$ respectively run over the sets of functors 
$\Fct (\Delta_n,\frak C)$ and $\Fct (\Delta_{n+1},\frak C)\,,$ sending any family 
$a = (a_\frak q)_\frak q$ to the family $d_n (a) =\big(d_n (a)_\frak r \big)_\frak r$ defined by
$$d_n (a)_\frak r = \sum_{i = 0}^{n+\!1} (-1)^i \frak m ({\rm id}_{\frak r\circ
\delta^n_i}, \delta^n_i)(a_{\frak r\circ \delta^n_i})
\eqno £3.6.3.$$

\bigskip
\noindent
{\bf Proposition £3.7}\phantom{.} {\it  With the notation above, the functor 
$\ch^* (\frak i)$ induces an $\O \-$module 
 isomorphism\/}
$$\G_\K (\F,\widehat\aut_{\F^{^{\rm sc}}}) \cong \lim_{\longleftarrow}\,\big(\frak g_\K\circ\,\widehat{\!\loc\!}\,_{\F^{^{\rm rd}}}\big)
\eqno £3.7.1.$$
\par
\noindent
{\bf Proof:} First of all, we prove that the homomorphism from 
${\displaystyle \lim_{\longleftarrow}}\,(\frak g_\K\circ\,\widehat{\!\loc\!}
\,_{\F^{^{\rm sc}}})$ to ${\displaystyle \lim_{\longleftarrow}}\,\big(\frak g_\K\circ
\,\widehat{\!\loc\!}\,_{\F^{^{\rm rd}}}\big)$ induced by 
$\ch^* (\frak i)$ is injective;  consider a family~$X = (X_Q)_{Q}\,,$ 
where $Q$ runs over the set of all the $\F\-$selfcentralizing subgroups  of~$P$ and $X_Q$ belongs to $\G_\K\big(\skew5\hat\L (Q)\big)\,,$ and assume that  
$d^{^{\rm sc}}_0 (X) = 0$  where
$d^{^{\rm sc}}_0$ denotes the corresponding differential map (cf.~£3.6.3).
\eject

\smallskip
For such a $Q\,,$ consider a dominant $\F\-$morphism $\rho\,\colon Q\to R$ such that 
$\vert R\vert/\vert Q\vert$ is maximal and the {\it chain\/} $\frak r\,\colon  {\Delta}_1\to \F^{^{\rm sc}}$ mapping $0$ on $Q\,,$ $1$ on~$R$ and 
$0\!\bullet\! 1$ on $\rho\,;$ then, we get (cf.~£3.6.3)
$$0 = {\rm res}_{\,\widehat{\!\loc\!}\,_{\F^{^{\rm sc}}} ({\rm id}_{\frak r \circ
\delta^0_0},\delta^0_0)} (X_R) -  {\rm res}_{\,\widehat{\!\loc\!}
\,_{\F^{^{\rm sc}}} ({\rm id}_{\frak r\circ  \delta^0_1},\delta^0_1)} (X_{Q}) 
\eqno £3.7.2;$$
but, it follows from Proposition~£2.7 that $R$ is $\F\-$radical and we know that the homomorphism $\F(R)_{\rho (Q)}\to \F (Q)$ induced by $\rho$  is surjective (cf.~£2.6.2); in particular, the functor $\widehat{\loc}_{\F^{^{\rm sc}}}$ provides a $k^*\-$group isomorphism $\skew5\hat\L (\frak r)\cong \skew5\hat\L (Q)\,,$ so that we also get the isomorphism 
$${\rm res}_{\,\widehat{\!\loc\!}\,_{\F^{^{\rm sc}}}({\rm id}_{\frak r\circ  \delta^0_1},\delta^0_1)} :
\G_\K\big(\skew5\hat{\L} (Q)\big)\cong \G_\K\big(\skew5\hat{\L}  (\frak r)\big)
\eqno £3.7.3;$$
thus, if $X$ belongs to the kernel of the homomorphism induced by~$\ch^* (\frak i)\,,$ 
we get~$X_Q = 0\,.$

\smallskip
In order to prove the surjectivity, assume that a family $Y = (Y_R)_{R}\,,$ where
$R$ runs over all the $\F\-$radical subgroups of $P$ and $Y_R\in \G_\K\big(\skew5\hat{\L} (R)\big)\,,$ belongs to the kernel of
the corresponding differential map $d^{^{\rm rd}}_0\,;$ we have to extend $Y$
to a family $X =  (X_Q)_{Q}\,,$ where $Q$ runs over the set of all the $\F\-$selfcentralizing subgroups of~$P$ and $X_Q\in \G_\K\big(\skew5\hat{\L} (Q)\big)\,,$  fulfilling  $d^{^{\rm sc}}_0 (X) = 0\,.$ For any $\F\-$selfcentralizing subgroup $Q$ of~$P\,,$ with the notation above we can define 
$$X_{Q} = \big({\rm res}_{\,\widehat{\!\loc\!}\,_{\F^{^{\rm sc}}} 
({\rm id}_{\frak r \circ  \!\delta^0_1}, \delta^0_1)}\big)^{-1} \big({\rm res}_{\,\widehat{\!\loc\!}\, _{\F^{^{\rm sc}}}({\rm id}_{\frak r \circ  \!\delta^0_0},\delta^0_0)} (Y_R)\big)
\eqno £3.7.4;$$
then,  we claim that  $d^{^{\rm sc}}_0 (X) = 0\,.$

\smallskip
 Let $\frak q\,\colon {\Delta}_1\to \F^{^{\rm sc}}$ be a {\it chain\/}, set $Q = \frak q
(0)\,,$  $Q' = \frak q (1)$ and $\varphi = \frak q (0\bullet 1)\,;$ arguing by
induction on~$\vert P\,\colon Q\vert$ and then on $\vert Q'\vert/\vert Q\vert\,,$ we will prove that we have $d^{^{\rm sc}}_0 (X)_{\frak q} = 0\,.$ If $Q= P$ then $\varphi$
is an $\F^{^{\rm rd}}\-$morphism, so that $\frak q$ is an $\F^{^{\rm rd}}\-$chain, and therefore we have
$$d^{^{\rm sc}}_0 (X)_{\frak q} = d^{^{\rm sc}}_0 (Y)_{\frak q} = 0
\eqno £3.7.5;$$
moreover, if we have $Q\not= P$ and $\varphi$ is an isomorphism, we may assume that we have chosen $\rho\circ\varphi^{-1}\,\colon Q'\to R$ as a dominant homomorphism 
for $Q'$ and considered the corresponding chain $\frak r'\,\colon \Delta_1
\to \F^{^{\rm sc}}\,,$
 so that we have defined
$$X_{Q'} = \big({\rm res}_{\,\widehat{\!\loc\!}\,_{\F^{^{\rm sc}}} 
({\rm id}_{\frak r' \circ \!\delta^0_1}, \delta^0_1)}\big)^{-1} \big({\rm res}_{\,\widehat{\!\loc\!}\,_{\F^{^{\rm sc}}}  
({\rm id}_{\frak r' \circ  \!\delta^0_0},\delta^0_0)} (Y_R)\big)
\eqno £3.7.6;$$
then, the functoriality of $\,\widehat{\!\loc\!}\,_{\F^{^{\rm sc}}}$ forces $X_Q = 
{\rm res}_{\varphi}(X_{Q'})$ and therefore we get (cf.~£3.6.3)
$$d^{^{\rm sc}}_0 (X)_{\frak q} = {\rm res}_{\,\widehat{\!\loc\!}\,_{\F^{^{\rm sc}}} ({\rm id}_{\frak q \circ\delta^0_0},\delta^0_0)} (X_{Q'}) -  {\rm res}_{\,\widehat{\!\loc\!}\,_{\F^{^{\rm sc}}} ({\rm id}_{\frak q\circ  \delta^0_1},\delta^0_1)} (X_{Q})= 0  \eqno £3.7.7.$$

\smallskip
From now on, we assume that $\vert P\,\colon Q\vert\not= 1\not=\vert Q'\vert/\vert Q\vert\,;$ firstly note that, for another $\F\-$selfcentralizing subgroup $Q''$ of~$P$ and an $\F\-$morphism $\varphi'\,\colon Q'\to Q''\,,$ considering the  chains
$$\frak q' : {\Delta}_1 \too \F^{^{\rm sc}}\quad ,
\quad\frak q'' : {\Delta}_1\too \F^{^{\rm sc}}\qq \frak q^{Q''}
: {\Delta}_2\too \F^{^{\rm sc}}
\eqno £3.7.8\phantom{.}$$ 
respectively mapping $0$ on $Q'\,,$ $Q$ and $Q\,,$ $1$ on $Q''\,,$ $Q''$ and $Q'\,,$  
$0\!\bullet\! 1$ on $\varphi'\,,$ $\varphi'\circ\varphi$ and $\varphi\,,$ $2$ on~$Q''\,,$ and $1\!\bullet\! 2$ 
on~$\varphi'\,,$ the equality $d^{^{\rm sc}}_1\circ d^{^{\rm sc}}_0 = 0$ implies that (cf.~£3.6.3)
$$\eqalign{&{\rm res}_{\,\widehat{\!\loc\!}\,_{\F^{^{\rm sc}}} ({\rm id}_{\frak q^{Q''} \circ \delta^1_1}, \delta^1_1)} \big(d^{^{\rm sc}}_0 (X)_{\frak q''}\big)\cr 
 &= {\rm res}_{\,\widehat{\!\loc\!}\,_{\F^{^{\rm sc}}}
({\rm id}_{\frak q^{Q''} \circ \delta^1_2},\delta^1_2)} \big(d^{^{\rm sc}}_0
(X)_{\frak q}\big)  +   {\rm res}_{\,\widehat{\!\loc\!}\,_{\F^{^{\rm sc}}}
 ({\rm id}_{\frak q^{Q''} \circ  \delta^1_0},\delta^1_0)} \big(d^{^{\rm sc}}_0 (X)_{\frak q'}\big)\cr}
\eqno £3.7.9;$$
but, by our induction hypothesis, we already know that  $d^{^{\rm sc}}_0 (X)_{\frak q'}  = 0\,.;$ hence, we get
$${\rm res}_{\,\widehat{\!\loc\!}\,_{\F^{^{\rm sc}}} ({\rm id}_{\frak q^{Q''} \circ \delta^1_1}, \delta^1_1)} \big(d^{^{\rm sc}}_0 (X)_{\frak q''}\big)
 = {\rm res}_{\,\widehat{\!\loc\!}\,_{\F^{^{\rm sc}}}
({\rm id}_{\frak q^{Q''} \circ \delta^1_2},\delta^1_2)} \big(d^{^{\rm sc}}_0
(X)_{\frak q}\big)
\eqno £3.7.10.$$

\smallskip
On the one hand, set $N = N_R\big(\rho (Q)\big)$ and consider 
the {\it chains\/}
$$\frak n : {\Delta}_1\too \F^{^{\rm sc}}\quad ,\quad \frak i_N^R : {\Delta}_1\too  \F^{^{\rm sc}}\quad{\rm and} \quad \frak r : {\Delta}_1
\too \F^{^{\rm sc}}
\eqno £3.7.11\phantom{.}$$
where $\frak r$ is defined as above, $\frak i_{N}^R$ by the inclusion
$N\i R$  and $\frak n$ by the restriction of $\rho$ from $Q$ to $N\,;$ then, we have $d^{^{\rm sc}}_0 (X)_{\frak r} = 0$ by definition~£3.7.4. If $\rho$ is an isomorphism
then $\frak n = \frak r$ and therefore $d^{^{\rm sc}}_0 (X)_{\frak n} = 0\,.$ 
Otherwise, it follows from the induction hypothesis that we still have $d^{^{\rm sc}}_0 (X)_{\frak i_{N}^R} = 0\,;$ but, since the image in $\F (Q)$ of any element~$\sigma$ of~$\F (N)_{\rho (Q)}$ can be~lifted to and element $\hat\sigma$ of $\F (R)$ which stabilizes~$\rho (Q)$ and therefore it stabilizes $N\,,$ the difference  between $\sigma$ and the restriction of $\hat\sigma$ to $N$ belongs to $\F_{Z (\rho(Q))} (N)$
[5,~Corollary~£4.7]; consequently, any element of~$\F (N)_{\rho (Q)}$ can be lifted to and element of  $\F (R)\,.$ Thus, {\it mutatis mutandis,\/} considering the 
chain $\frak n^R\,\colon  \Delta_2
\too \F^{^{\rm sc}}$ extending $\frak n$ and mapping $2$ on $R$ and
$1\!\bullet\! 2$ on the inclusion map $N\to R\,,$ we have a
$k^*\-\widetilde\Loc\-$isomorphism
$$\widehat{\!\loc\!}\,_{\F^{^{\rm sc}}} ({\rm id}_{\frak n^R\circ \delta^1_2},\delta^1_2) : 
\skew5\hat{\L} (\frak n^R)\cong \skew5\hat{\L} (\frak n)
\eqno £3.7.12;$$
hence, since $\frak n^R\circ \delta^1_1 = \frak r\,,$ by
equality~£3.7.10 we also have $d^{^{\rm sc}}_0 (X)_{\frak n} = 0\,.$

\smallskip
On the other hand, set $N' = N_{Q'}\big(\varphi (Q)\big)$ and consider the chains
$$\frak n' : {\Delta}_1\too \F^{^{\rm sc}}\quad ,\quad \frak i_{N'}^{Q'} : 
{\Delta}_1 \too \F^{^{\rm sc}}\qq  \frak q : {\Delta}_1\too \F^{^{\rm sc}}
\eqno £3.7.13\phantom{.}$$
where $\frak n'$ is defined by the restriction of $\varphi$ from $Q$ to $N'$ and 
$\frak i_{N'}^{Q'}$ by the inclusion $N'\i Q'\,;$ by our induction
hypothesis we have~$d^{^{\rm sc}}_0 (X)_{\frak i_{N'}^{Q'}} = 0\,,$ and  it is 
clear that any element of~$\F (Q')_{\varphi (Q)}$ stabilizes~$N'\,.$ As above,
considering the chain $\frak n'^{Q'}\,\colon  \Delta_2\to \F^{^{\rm sc}}$
extending~$\frak n'\,,$ and mapping $2$ on $Q'$ and $1\bullet 2$ on the
inclusion map $N'\to Q'\,,$ we still have a
$k^*\-$group exoisomorphism 
$$\widehat{\!\loc\!}\,_{\F^{^{\rm sc}}} ({\rm id}_{\frak n'^{Q'}\circ \delta^1_1},
\delta^1_1) : \skew5\hat{\L} (\frak n'^{Q'})\cong
\skew5\hat{\L} (\frak q)
\eqno £3.7.14;$$
hence, according to equality~£3.7.10, in order to prove that $d^{^{\rm sc}}_0 (X)_{\frak q} = 0\,,$ it suffices to prove that we have   $d^{^{\rm sc}}_0 (X)_{\frak n'} = 0\,.$

\smallskip
 That is to say, we may assume that $Q'= N'$ normalizes $\varphi (Q)\,;$ in this case,
it follows from [5, Proposition~2.7] that there is an $\F\-$morphism 
$\zeta'\,\colon Q'\to P$ such that $\zeta'\big(\varphi (Q)\big)$ is fully normalized  
in $\F$ and then, from [5,~condition~£2.8.2] that there are an $\F\-$morphism $\eta\,\colon N\to P$ and an element $\sigma\in \F (Q)$ fulfilling  $\eta \Big(\rho\big(\sigma(u)\big)\Big) =
\zeta'\big(\varphi (u)\big)$ for any $u\in Q\,;$ actually, up to modifying  our
choice of~$\zeta'\,,$ we may assume that $\sigma = {\rm id}_Q\,.$ Now, $\eta (N)$ and
$\zeta' (Q')$ normalize $\zeta'\big(\varphi (Q)\big)$ and we consider the group $N'' = 
\langle \eta (N), \zeta' (Q')\rangle$ and the  chains
$$\frak e : {\Delta}_1\too  \F^{^{\rm sc}}\qq \frak n'' : {\Delta}_1\too \F^{^{\rm sc}}
\qq \frak n^{N''} : {\Delta}_2\too \F^{^{\rm sc}}
\eqno £3.7.15\phantom{.}$$
where $\frak e$ and $\frak n''$ are respectively defined by the homomorphisms
$N\to N''$ and $Q\to N''$ determined by~$\eta$ and~$\zeta'\circ\varphi\,,$
and $\frak n^{N''}$ extends $\frak n$ mapping $2$ on $N''$ and $1\!\bullet 2$
on~the~homomorphism $N\to N''$ determined by $\eta\,.$

\smallskip
 We already know that  $d^{^{\rm sc}}_0 (X)_{\frak n} = 0$ and it follows from our induction hypothesis that $d^{^{\rm sc}}_0 (X)_{\frak e} = 0\,;$ but, any element 
 of $\F (Q)$ can be lifted {\it via\/} $\rho$ to an element of $\F (R)_{\rho (Q)}$ (cf.~£2.6.2) which stabilizes $N\,.$ In particular, for any element $\sigma$ in
$$\F(N'')_{\zeta'(\varphi (Q))} = \F (\frak n'')
\eqno £3.7.16\,,$$
 we have an element $\hat\sigma$
 in $\F (N)$ such that $\eta\circ\hat\sigma$ coincides with $\sigma\circ\eta$ over 
 $\rho (Q)\,,$ and therefore it follows from [5,~Proposition~£4.6] that it suffices to modify our choice of~$\hat\sigma$ by composing it with the conjugation by a suitable element of $Z\big(\rho (Q)\big)\i N$ to get the equality  $\eta\circ\hat\sigma = \sigma\circ\eta\,.$ In conclusion, we get the $k^*\-\widetilde\Loc\-$isomorphism
$$\widehat{\!\loc\!}\,_{\F^{^{\rm sc}}} ({\rm id}_{\frak n^{N''}\circ \delta^1_1},
\delta^1_1) : \skew5\hat{\L}  (\frak n^{N''})\cong \skew5\hat{\L}  (\frak n'')
\eqno £3.7.17;$$
consequently, by equality~£3.7.10, we still get $d^{^{\rm sc}}_0 (X)_{\frak n''}  = 0\,.$

\smallskip
Similarly, consider the chains
$$\frak z' : {\Delta}_1\too \F^{^{\rm sc}}\cqq  \frak n'' : {\Delta}_1\too 
\F^{^{\rm sc}}\qq \frak q^{N''} : {\Delta}_2\too \F^{^{\rm sc}}
\eqno £3.7.18\phantom{.}$$
where $\frak z'$ is defined by the homomorphism $Q'\to N''$ determined by
$\zeta'\,,$ and~$\frak q^{N''}$ extends $\frak q = \frak n'$ mapping $2$ on $N''$ and 
$1\bullet 2$ on the homomorphism $Q'\to N''$ determined by~$\zeta'\,;$
it follows from our induction hypothesis that $d^{^{\rm sc}}_0 (X)_{\frak z'} = 0$ and we already  know that $d^{^{\rm sc}}_0 (X)_{\frak n''} = 0\,.$  Since $N''$ contains 
$$\eta \Big(Z\big(\rho (Q)\big)\Big) = Z\Big(\zeta'\big(\varphi
(Q)\big)\Big)
\eqno £3.7.19,$$
 it follows from [5, statement~2.10.1] that the automorphism of~$\zeta'\big(\varphi (Q)\big)$ determined by any $\sigma'\in \F(Q')_{\varphi (Q)}$ {\it via\/} $\zeta'$ can be extended to an $\F^{^{\rm sc}}\-$automorphism $\sigma''$ of $N''$ and then,
arguing as above, it follows from [5, Proposition~4.6] that we may choose
$\sigma''$ fulfilling 
$$\sigma''\big(\zeta' (u')\big) = \zeta'\big(\sigma' (u')\big)
\eqno £3.7.20\phantom{.}$$
 for any
$u'\in Q'\,.$ In conclusion, we get the $k^*\-\widetilde\Loc\-$isomomorphism
$$\widehat{\!\loc\!}\,_{\F^{^{\rm sc}}} ({\rm id}_{\frak q^{N''}\circ \delta^1_2},\delta^1_2) :
\skew5\hat{\L} (\frak q^{N''})\cong  \skew5\hat{\L} (\frak q)
\eqno £3.7.21;$$
consequently, again by equality~£3.7.10, we get $d^{^{\rm sc}}_0 (X)_{\frak q}  = 0\,.$ We are done.

\bigskip
\bigskip
\noindent
{\bf£ 4\phantom{.} Functoriality of the Grothendieck groups of  folded Frobenius $P\-$categories}
\bigskip

£4.1\phantom{.} With the notation in \S3, let $P'$ be a second finite $p\-$group, $\F'$ a Frobenius $P'\-$category and $\alpha\,\colon P'\to P$ an  
{\it $(\F',\F)\-$functorial \/} group homomorphism [5, 12.1] mapping any $\F'\-$radical subgroup of $P'$ on a  $\F\-$selfcentralizing subgroup of $P\,.$ More generally, let 
$\frak X'$ be a set of subgroups $Q'$ of~$P'$ such that $\alpha (Q')$ is $\F\-$selfcentralizing, which contains the $\F'\-$radical subgroups of $P'$ and is stable by 
$\F'\-$isomorphisms, and denote by $\F'^{^{\frak X'}}$ the {\it full\/} subcategory of 
$\F'$ over $\frak X'\,.$ Then the restriction $\frak f^{^{\rm \frak X'}}_\alpha\,\colon 
\F'^{^{\frak X'}}\to \F^{^{\rm sc}}$ of the {\it Frobenius functor\/}
$\frak f_\alpha :  \F'\too \F$  
induces a natural map
$$\aut_{\frak f^{^{\frak X'}}_\alpha} : \aut_{\F'^{^{\frak X'}}}
\too \aut_{\F^{^{\rm sc}}}\circ \ch^*(\frak f^{^{\frak X'}}_\alpha)
\eqno  £4.1.1\phantom{.}$$
and therefore the {\it pull-back}
$$\matrix{\aut_{\F'^{^{\frak X'}}} &\buildrel \aut_{\frak f^{^{\frak X'}}_\alpha}\over {\hbox to 30pt{\rightarrowfill}} &\aut_{\F^{^{\rm sc}}}\circ 
\ch^*(\frak f^{^{\frak X'}}_\alpha)\cr
\big\uparrow&\phantom{\Big\uparrow}&\big\uparrow\cr
\widehat\aut_{\F'^{^{\frak X'}}} &\buildrel \widehat\aut_{\frak f^{^{\frak X'}}_\alpha}\over{\hbox to 30pt{\rightarrowfill}} &\widehat\aut_{\F^{^{\rm sc}}}
\circ \ch^*(\frak f^{^{\frak X'}}_\alpha)\cr}
\eqno £4.1.2\phantom{.}$$
lifts the functor $\aut_{\F'^{^{\frak X'}}}$ to the category $k^*\-\Gr$ (cf.~£3.1). Moreover,  we already know from Theorem £2.9 that this lifting can be extended to a unique functor $\widehat\aut_{\F'^{^{\rm sc}}}$ lifting~$\aut_{\F'^{^{\rm sc}}}$
and therefore we have an ordinary and a modular Grothendieck groups
$$\eqalign{\G_\K (\F',\widehat\aut_{\F'^{^{\rm sc}}}) &= {\displaystyle \lim_{\longleftarrow}}\,(\frak g_\K\circ \widehat\loc_{\F'^{^{\rm sc}}})\cr 
\G_k (\F',\widehat\aut_{\F'^{^{\rm sc}}}) &= {\displaystyle \lim_{\longleftarrow}}
\,(\frak g_k\circ \widehat\loc_{\F'^{^{\rm sc}}})\cr}
\eqno £4.1.3.$$

\medskip
£4.2\phantom{.} Now, according to Proposition~£3.7, in order to define the {\it restriction\/} $\O\-$module homomorphisms from the Grothendieck groups of the folded Frobenius $P\-$category 
$(\F,\widehat\aut_{\F^{^{\rm sc}}})$ to these ones, we have to exhibit a natural map
$$\loc_{\frak f^{^{\frak X'}}_\alpha} : \loc_{\F'^{^{\frak X'}}}
\too \loc_{\F^{^{\rm sc}}}\circ \ch^*(\frak f^{^{\frak X'}}_\alpha)
\eqno  £4.2.1\phantom{.}$$
lifting the natural map $\aut_{\frak f^{^{\frak X'}}_\alpha}\,;$ explicitly, setting~$\L' (\frak q') = 
\loc_{\F'^{^{\frak X'}}}(\frak q',\Delta_n)$ for any  
$\F'^{^{\frak X'}}\-$chain $\frak q'\,\colon \Delta_n\to \F'^{^{\frak X'}}\,,$ we have to exhibit a suitable  group homomorphism
$$\lambda_{\frak q'} : \L'(\frak q')\too \L (\frak f^{^{\frak X'}}_\alpha\circ \frak q')
\eqno  £4.2.2\phantom{.}$$
sending $Z\big(\frak q'(n)\big)$ to a subgroup of $Z\big(\alpha
(\frak q'(n))\big)$ and lifting the group homomorphism 
$$(\aut_{\frak f^{^{\frak X'}}_\alpha})_{(\frak q',\Delta_n)} : \F'(\frak q')\too 
\F (\frak f^{^{\frak X'}}_\alpha\circ \frak q')
\eqno  £4.2.3\phantom{.}$$
which maps $\sigma'\in \F'(\frak q')$ on $\frak f^{^{\frak X'}}_\alpha * \sigma'\,.$

\bigskip
\noindent
{\bf Proposition~£4.3}\phantom{.} {\it  With the notation above,
 for any   $\F'^{^{\frak X'}}\!\!\-$chain $\frak q'\,\colon \Delta_n\to \F'^{^{\frak X'}}$
 such that $\frak q' (n)$ is fully normalized in $\F'$ and that $\F'_{\!P'}\big(\frak q'(n)\big)$ contains a Sylow $p\-$subgroup of $\F' (\frak q')\,,$ there exists a group homomorphism
$$\lambda_{\frak q'} : \L'(\frak q')\too \L (\frak f^{^{\frak X'}}_\alpha\circ \frak q')
\eqno  £4.3.1\phantom{.}$$
 lifting $(\aut_{\frak f^{^{\frak X'}}_\alpha})_{(\frak q',\Delta_n)}\,\colon \F'(\frak q')\to \F (\frak f^{^{\frak X'}}_\alpha\circ \frak q')$ and mapping $u\in N_{P'}(\frak q')$
on $\alpha (u)\,.$ Moreover, the group $Z\big(\alpha(\frak q'(n))\big)^{\alpha (N_{P'}(\frak q'))}$ acts transitively on the set of 
such homomorphisms.\/}
\medskip
\noindent
{\bf Proof:}  In order to apply [5,~Lemma~18.8], let us consider the groups 
$\L'(\frak q')$ and $\L(\frak f^{^{\frak X'}}_\alpha\circ\frak q')\,,$ and the group homomorphisms
$$\eqalign{(\aut_{\frak f^{^{\frak X'}}_\alpha})_{(\frak q',\Delta_n)}\circ \pi_{\frak q'} : 
\L'(\frak q')&\too  \F (\frak f^{^{\frak X'}}_\alpha\circ \frak q')\cr
 \alpha_{\frak q'} : N_{P'}(\frak q')&\too  N_P(\frak f^{^{\frak X'}}_\alpha\circ\frak q')\i \L(\frak f^{^{\frak X'}}_\alpha\circ\frak q')\cr}
 \eqno £4.3.2$$
 where $\pi_{\frak q'}\,\colon \L'(\frak q')\to \F'(\frak q')$ denotes the structural homomorphism and $\alpha_{\frak q'}$ the restriction of~$\alpha$ from 
 $ N_{P'}(\frak q')$ to $N_P(\frak f^{^{\frak X'}}_\alpha\circ\frak q')\,;$ let $R'$ be
 a subgroup of~$P'$  and, setting $Q' = \frak q' (n)\,,$ $x'$ an element of $\L' (\frak q')\i \L'(Q')$ such that $R'^{x'}\i P'\,;$  then, $\frak f^{^{\frak X'}}_\alpha * \pi_{\frak q'}(x')$ belongs to  
 $\F (\frak f^{^{\frak X'}}_\alpha\circ \frak q')$ which is contained in $\F\big(\alpha(Q')\big)\,,$
 and $x'$ actually determines an $\F_{\L' (Q'),Q'}\-$morphism from $R'$ 
 to $N_{P'}(Q')$ [5,~17.2]; in particular, according to [5, Theorem~18.6], $x'$ also determines an $N_{\F',Q'} (Q')\-$morphism $\xi'$ from $R'$  to~$N_{P'}(Q')$
 [5, 17.6].

 \smallskip
Moreover,  it follows from [5, Proposition~2.7] that there is an $\F\-$mor-phism 
$$\zeta : N_P\big(\alpha(Q')\big) \too P
\eqno  £4.3.3\phantom{.}$$
such that $Q = \zeta\big(\alpha(Q')\big)$ is fully normalized in $\F\,;$ in particular, we can consider the {\it $\F\-$localizer\/} $\L(Q)$ of $Q\,,$ the Frobenius $N_P(Q)\-$category $N_\F(Q)$ and the {\it $N_\F(Q)\-$locality\/} $N_{\F,Q}(Q)$ [5,~17.6]. Then, 
 since $\alpha$ is {\it $(\F',\F)\-$functorial\/}, denoting by $\alpha_{Q'}$ the restriction of $\alpha$ from  $ N_{P'}(Q')$ to $N_P\big(\alpha(Q')\big)\,,$ it is quite clear that the composition
 $$\zeta\circ \alpha_{Q'} : N_{P'}(Q')\too N_P(Q)
 \eqno £4.3.4\phantom{.}$$
 is an {\it $(N_{\F'}(Q'),N_\F (Q))\-$functorial\/} group homomorphism [5, 12.1] which induces a functor
$$\frak l_{\zeta\circ \alpha_{R'} } : N_{\F',Q'}(Q')\too N_{\F,Q}(Q)
\eqno £4.3.5.$$
Thus, $\xi'$ determines an $N_{\F,Q} (Q)\-$morphism $\xi$ from $R = \zeta\big(\alpha (R')\big)$  to~$N_{P}(Q)\,.$

\smallskip
But, according to [5, Theorem~18.6] again, the categories $N_{\F,Q}(Q)$ and
$\F_{\L (Q),Q}$ [5,~17.2] are equivalent to each other and therefore there is $x\in \L (Q)$ such that $\xi (u) = u^x$ for any $u\in R$ and that we have $\pi_Q (x) = {}^{\zeta\circ \alpha_{Q'}} 
\big(\pi_{Q'}(x')\big)$ where $\pi_Q\,\colon \L (Q)\to \F (Q)$ and $\pi_{Q'}
\,\colon \L' (Q')\to \F' (Q')$ denote the structural homomorphisms.  Finally, considering 
the $\F\-$chain $\frak q\,\colon \Delta_n\to \F^{^{\rm sc}}$ mapping 
$i\in \Delta_{n-1}$ on $\alpha (\frak q' (i))\,,$  $n$ on $Q\,,$ any 
$\Delta_{n-1}\-$morphism $j\bullet i$ on $(\frak f^{^{\frak X'}}_\alpha\circ \frak q')
(j\bullet i)$ and $n\!-\!1\bullet n$ on the composition of $(\frak f^{^{\frak X'}}_\alpha\circ \frak q')(n\!-\!1\bullet n)$ with the isomorphism $\alpha (Q')\cong Q$ determined by $\zeta\,,$ we have an obvious $\ch^*(\F^{^{\rm sc}})\-$isomorphism 
$$(\nu_\zeta,{\rm id}_{\Delta_n}) : (\frak q,\Delta_n)\too (\frak f^{^{\frak X'}}_\alpha\circ \frak q', \Delta_n)
\eqno £4.3.6\phantom{.}$$ 
which the functor $\loc_{\F^{^{\rm sc}}}$ sends to a group isomorphism
$$\lambda_\zeta : \L (\frak q)\cong \L (\frak f^{^{\frak X'}}_\alpha\circ \frak q')
\eqno £4.3.7.$$

\smallskip
At this point, it is easily checked that $\alpha (u'^{x'}) = 
\alpha (u')^{\lambda_\zeta (x)}$ for any $u'\in R'$ and that we have 
$$\pi_{\alpha (Q')}\big(\lambda_\zeta (x)\big) = {}^{\alpha_{Q'}}\big(\pi_{Q'}(x')\big)
\eqno £4.3.8\,;$$
that is to say, in our setting of [5, Lemma~18.8], condition~18.8.1 holds
and therefore there is a group homomorphism
$$\lambda_{\frak q'} : \L'(\frak q')\too \L (\frak f^{^{\frak X'}}_\alpha\circ \frak q')
\eqno  £4.3.9\phantom{.}$$
fulfilling the announced conditions. Moreover, always by the same lemma, the Abelian group 
$Z\big((\frak f^{^{\frak X'}}_\alpha\circ\frak q')(n)\big)^{\alpha (N_{P'}
(\frak q'))}$ acts transitively over the set of such group
homomorphisms.

\bigskip
\noindent
{\bf Theorem~£4.4}\phantom{.} {\it With the notation and the hypothesis above,
there is a unique natural map
$$\loc_{\frak f^{^{\frak X'}}_\alpha} : \loc_{\F'^{^{\frak X'}}}
\too \loc_{\F^{^{\rm sc}}}\circ \ch^*(\frak f^{^{\frak X'}}_\alpha)
\eqno  £4.4.1\phantom{.}$$
sending any $\F'^{^{\frak X'}}\!\-$chain $\frak q'\,\colon \Delta_n\to \F'^{^{\frak X'}}$
 such that $\frak q' (n)$ is fully normalized in~$\F'$ and that $\F'_{\!P'}\big(\frak q'(n)\big)$ contains a Sylow $p\-$subgroup of $\F' (\frak q')\,,$ to the $Z\big(\alpha (\frak q'(n))\big)\-$conjugacy class of
 $$\lambda_{\frak q'} : \L'(\frak q')\too \L (\frak f^{^{\frak X'}}_\alpha\circ \frak q')
\eqno  £4.4.2.$$
In particular, this natural map lifts $\aut_{\frak f^{^{\frak X'}}_\alpha}$ and, if 
$\frak r'\,\colon \Delta_m\to \F'^{^{\frak X'}}\!\!$ is a chain,
then we have $\lambda_{\frak r'}(u')  = \alpha (u')$ for a suitable representative $\lambda_{\frak r'}$ of 
$\,(\loc_{\frak f^{^{\frak X'}}_\alpha})_{(\frak r',\Delta_m)}$
and any $u'\in N_{P'}(\frak r')\,.$\/}

\medskip
\noindent
{\bf Proof:} For any  $\F'^{^{\frak X'}}\!\-$chain $\frak q'\,\colon \Delta_n\to \F'^{^{\frak X'}}\,,$  it follows from [5, Proposition~2.7] that there is an 
$\F'\-$morphism $\zeta'\,\colon N_{P'}\big(\frak q'(n)\big)\to P'$ such that 
$\zeta'\big(\frak q' (n)\big)$ is fully normalized in $\F'$ and we may assume that 
$\F'_{\!P'}\big(\frak q'(n)\big)$ contains a Sylow $p\-$subgroup of $\F' (\frak q')\,;$ consider the $\F'^{^{\frak X'}}\!\-$chain $\hat\frak q'\,\colon \Delta_n\to 
\F'^{^{\frak X'}}$ which coincides with $\frak q'$ over $\Delta_{n-1}$ and maps $n$ on 
 $\zeta'\big(\frak q' (n)\big)$ and $(n\!-\!1\bullet n)$ on the composition of 
 $\frak q'(n\!-\!1\bullet n)$ with the $\F'\-$isomorphism $\zeta'_*\,\colon\frak q'(n)\cong  \zeta'\big(\frak q' (n)\big)$ determined by~$\zeta'\,;$ we have  an obvious 
 $\ch^*(\F'^{^{\frak X'}})\-$isomorphism 
$$(\nu'_{\zeta'},{\rm id}_{\Delta_n}) : (\frak q',\Delta_n)\cong (\hat\frak q', \Delta_n)
\eqno £4.4.3\phantom{.}$$ 
which the functor $\loc_{\F'^{^{\frak X'}}}$ sends to a class of group isomorphisms
$$\lambda'_{\zeta'} : \L' (\frak q')\cong \L'(\hat\frak q')
\eqno £4.4.4\phantom{.}$$
and we may assume that $\lambda'_{\zeta'}(u') = \zeta' (u')$ for any $u'\in N_{P'}(\frak q')\,.$
\eject

\smallskip
On the other hand, we have the $\F\-$morphism
$$\frak f^{^{\frak X'}}_\alpha (\zeta') : \alpha \Big( N_{P'}\big(\frak q'(n)\big)\Big)
\too \alpha (P')
\eqno £4.4.5\phantom{.}$$
which can be restricted to the 
$\F\-$isomorphism 
$$\frak f^{^{\frak X'}}_\alpha (\zeta'_*) :  \alpha\big(\frak q'(n)\big)\cong  
\alpha\Big(\zeta'\big(\frak q' (n)\big)\Big)
\eqno £4.4.6;$$
then, considering the $\F^{^{\rm sc}}\-$chains $\frak f^{^{\frak X'}}_\alpha\circ\frak q'$ and $\frak f^{^{\frak X'}}_\alpha\circ\hat\frak q'\,,$ this $\F\-$isomorphism
induces  an obvious  $\ch^*(\F^{^{\rm sc}})\-$isomorphism 
$$(\nu_{\frak f^{^{\frak X'}}_\alpha (\zeta'_*) },{\rm id}_{\Delta_n}) : (\frak f^{^{\frak X'}}_\alpha\circ\frak q',\Delta_n)\cong (\frak f^{^{\frak X'}}_\alpha\circ\hat\frak q', \Delta_n)
\eqno £4.4.7\phantom{.}$$ 
which the functor $\loc_{\F^{^{\rm sc}}}$ sends to a class of group isomorphisms
$$\lambda_{\frak f^{^{\frak X'}}_\alpha (\zeta'_*) } : \L (\frak f^{^{\frak X'}}_\alpha
\circ\frak q')\cong \L(\frak f^{^{\frak X'}}_\alpha\circ\hat\frak q')
\eqno £4.4.8\phantom{.}$$
and we may assume that $\lambda_{\frak f^{^{\frak X'}}_\alpha (\zeta'_*) }(u) = 
\big(\frak f^{^{\frak X'}}_\alpha (\zeta')\big) (u)$ for any $u\in N_P 
(\frak f^{^{\frak X'}}_\alpha \circ\frak q')\,.$

\smallskip
Moreover, Proposition~£4.3 provides a group homomorphism
$$\lambda_{\hat\frak q'} : \L'(\hat\frak q')\too 
\L (\frak f^{^{\frak X'}}_\alpha\circ\hat\frak q')
\eqno £4.4.9\phantom{.}$$ 
fulfilling $\lambda_{\hat\frak q'}(u') = \alpha (u')$ for any $u'\in N_{P'}
(\hat\frak q')\,,$ and we define $\loc_{\frak f^{^{\frak X'}}_\alpha}$ as the map 
sending $(\frak q',\Delta_n)$ to the $Z\big(\alpha(\frak q'(n))\big)\-$conjugacy class of 
$$\lambda_{\frak q'} = (\lambda_{\frak f^{^{\frak X'}}_\alpha (\zeta'_*) })^{-1}
\circ \lambda_{\hat\frak q'}\circ \lambda'_{\zeta'} 
\eqno £4.4.10;$$
first of all, note that for any $u'\in N_{P'} (\frak q')$ we get
$$\eqalign{\lambda_{\frak q'} (u') &= \big((\lambda_{\frak f^{^{\frak X'}}_\alpha (\zeta'_*) })^{-1}\circ \lambda_{\hat\frak q'}\big) \big( \lambda'_{\zeta'} (u')\big) = (\lambda_{\frak f^{^{\frak X'}}_\alpha (\zeta'_*) })^{-1}\Big(\lambda_{\hat\frak q'} \big( \zeta' (u')\big)\Big)\cr
&= (\lambda_{\frak f^{^{\frak X'}}_\alpha (\zeta'_*) })^{-1}
\big((\alpha\circ\zeta')(u')\big)\cr
&= (\lambda_{\frak f^{^{\frak X'}}_\alpha (\zeta'_*) })^{-1}
\big((\frak f^{^{\frak X'}}_\alpha (\zeta')\circ \alpha)(u')\big) = \alpha (u')\cr}
\eqno £4.4.11.$$

\smallskip
In order to prove the naturality of this correspondence, let
$$(\mu',\delta) : (\frak r',\Delta_m)\too (\frak q',\Delta_n)
\eqno £4.4.12$$
be a $\ch^*(\F'^{^{\frak X'}})\-$morphism and consider an 
$\F'\-$morphism $\xi'\,\colon \frak r'(m)\to P'$ such that 
$\xi'\big(\frak r' (m)\big)$ is fully normalized in $\F'$ and  that 
$\F'_{\!P'}\big(\frak q'(n)\big)$ contains a Sylow $p\-$subgroup of $\F' (\frak q')\,,$  and the
 $\F'^{^{\frak X'}}\!\-$chain $\hat\frak r'\,\colon \Delta_m\to \F'^{^{\frak X'}}$
 which coincides with $\frak r'$ over $\Delta_{m-1}$ and maps $m$ on 
 $\xi'\big(\frak r' (m)\big)$ and $(m\!-\!1\bullet m)$ on the composition of  $\frak r'(m\!-\!1\bullet m)$ with the $\F'\-$isomorphism $\xi'_*\,\colon
 \frak r'(m)\cong  \xi'\big(\frak r' (m)\big)$ determined by~$\xi'\,.$ Then, 
 the $\ch^*(\F'^{^{\frak X'}})\-$isomorphisms
$(\nu'_{\zeta'},{\rm id}_{\Delta_n})$ and $(\nu'_{\xi'},{\rm id}_{\Delta_m})$ induce 
the commutative $\ch^*(\F'^{^{\frak X'}})\-$diagram
$$\matrix{(\frak r',\Delta_m)&\buildrel (\mu',\delta)\over
{\hbox to 30pt{\rightarrowfill}} &(\frak q',\Delta_n)\cr
\wr\Vert&\phantom{\big\uparrow}&\wr\Vert\cr
(\hat\frak r',\Delta_m)&\buildrel (\hat\mu',\delta)\over
{\hbox to 30pt{\rightarrowfill}} &(\hat\frak q',\Delta_n)\cr}
\eqno £4.4.13.$$
Moreover, $\frak f^{^{\frak X'}}_\alpha$ maps this $\ch^*(\F'^{^{\frak X'}})\-$diagram on the commutative $\ch^*(\F^{^{\rm sc}})\-$dia-gram
$$\matrix{(\frak f^{^{\frak X'}}_\alpha\circ\frak r',\Delta_m)
&\buildrel (\frak f^{^{\frak X'}}_\alpha * \mu',\delta)\over
{\hbox to 40pt{\rightarrowfill}} &(\frak f^{^{\frak X'}}_\alpha\circ\frak q',\Delta_n)\cr
\wr\Vert&\phantom{\big\uparrow}&\wr\Vert\cr
(\frak f^{^{\frak X'}}_\alpha\circ\hat\frak r',\Delta_m)&\buildrel (\frak f^{^{\frak X'}}_\alpha * \hat\mu',\delta)\over {\hbox to 40pt{\rightarrowfill}} 
&(\frak f^{^{\frak X'}}_\alpha\circ\hat\frak q',\Delta_n)\cr}
\eqno £4.4.14.$$

\smallskip
Now, it is quite clear that it suffices to prove the commutativity up to $Z\big(\alpha(\hat\frak q'(n))\big)\-$conjugation of the following diagram
of group homomorphisms
$$\matrix{\L(\frak f^{^{\frak X'}}_\alpha\circ\hat\frak r')&\buildrel 
 \lambda_{\frak f^{^{\frak X'}}_\alpha * \hat\mu'}\over {\hbox to 40pt{\rightarrowfill}} 
&\L(\frak f^{^{\frak X'}}_\alpha\circ\hat\frak q')\cr
\hskip-15pt \lambda_{\hat\frak r'}\big\uparrow&\phantom{\Big\uparrow}&\big\uparrow
\hskip2pt \lambda_{\hat\frak q'}\hskip-15pt\cr
\L'(\hat\frak r')&\buildrel \lambda'_{\hat\mu'}\over
{\hbox to 40pt{\rightarrowfill}} &\L'(\hat\frak q')\cr}
\eqno £4.4.15\phantom{.}$$
where $\lambda_{\frak f^{^{\frak X'}}_\alpha * \hat\mu'}$ and $\lambda'_{\hat\mu'}$ are respective representatives of  $ \loc_{\F^{^{\rm sc}}}(\frak f^{^{\frak X'}}_\alpha * \hat\mu',\delta)$ and
$ \loc_{\F'^{^{\frak X'}}}(\hat\mu',\delta)\,;$ this commutativity up to 
$Z\big(\alpha(\hat\frak q'(n))\big)\-$conjugation is a consequence of the uniqueness part of [5,~Lemma~18.8] and of [5,~Remark~18.9] applied to the groups $\L'(\hat\frak r')$ and $\L(\frak f^{^{\frak X'}}_\alpha\circ\hat\frak q')\,,$ and to the group homomorphisms $\lambda_{\frak f^{^{\frak X'}}_\alpha * \hat\mu'}\circ 
 \lambda_{\hat\frak r'}$ and $ \lambda_{\hat\frak q'}\circ  \lambda'_{\hat\mu'}\,;$ indeed, denoting by 
$$\eqalign{\pi_{\frak f^{^{\frak X'}}_\alpha\circ\hat\frak q'} : \L(\frak f^{^{\frak X'}}_\alpha \circ \hat\frak q')&\too \F(\frak f^{^{\frak X'}}_\alpha\circ\hat\frak q')\cr 
\pi_{\frak f^{^{\frak X'}}_\alpha\circ\hat\frak r'} : \L(\frak f^{^{\frak X'}}_\alpha \circ \hat\frak r')&\too \F(\frak f^{^{\frak X'}}_\alpha\circ\hat\frak r')\cr}
\eqno £4.4.16\phantom{.}$$
the structural homomorphisms, it follows from [5, Proposition~18.16] and from Proposition~£4.3 that we have (cf.~£4.1.1)
$$\eqalign{\pi_{\frak f^{^{\frak X'}}_\alpha\circ\hat\frak q'}\circ 
\lambda_{\frak f^{^{\frak X'}}_\alpha * \hat\mu'}\circ  \lambda_{\hat\frak r'} &= \aut_{\F^{^{\rm sc}}}(\frak f^{^{\frak X'}}_\alpha * \hat\mu',\delta)\circ
\pi_{\frak f^{^{\frak X'}}_\alpha\circ\hat\frak r'}\circ  \lambda_{\hat\frak r'}\cr
&= \aut_{\F^{^{\rm sc}}}(\frak f^{^{\frak X'}}_\alpha * \hat\mu',\delta)\circ
(\aut_{\frak f^{^{\frak X'}}_\alpha })_{\frak r'}\circ \pi_{\hat\frak r'}\cr
&= (\aut_{\frak f^{^{\frak X'}}_\alpha })_{\frak q'}\circ \aut_{\F'^{^{\rm sc}}}(\hat\mu',\delta)\circ \pi_{\hat\frak r'}\cr
&= (\aut_{\frak f^{^{\frak X'}}_\alpha })_{\frak q'}\circ \pi_{\frak f^{^{\frak X'}}_\alpha\circ\hat\frak r'} \circ  \lambda'_{\hat\mu'}\cr
 &= \pi_{\frak f^{^{\frak X'}}_\alpha\circ\hat\frak q'}\circ \lambda_{\hat\frak q'}\circ  \lambda'_{\hat\mu'}\cr}
 \eqno £4.4.17.$$

 \smallskip
 Finally, set $\hat Q' = \hat\frak q' (n)$ and consider  an $\F\-$morphism 
$$\zeta : N_P\big(\alpha(\hat Q')\big) \too P
\eqno  £4.4.18\phantom{.}$$
such that $\hat Q = \zeta\big(\alpha(\hat Q')\big)$ is fully normalized in $\F$
[5,~Proposition~2.7] and that $\zeta\big(N_P(\frak f^{^{\frak X'}}_\alpha
\circ \hat\frak q')\big)$ is contained in  $N_P (\hat Q)\,;$ denote by
 $\zeta_*\,\colon \alpha (\hat Q')\cong \hat Q$ the $\F\-$isomorphism determined by $\zeta\,,$ and
  by~$\lambda_{\zeta_*}$ a representative of $\loc_{\F^{^{\rm sc}}}(\nu_{\zeta_*},\iota_n)$ where 
 $$(\nu_{\zeta_*},\iota_n) : (\frak f^{^{\frak X'}}_\alpha\circ \hat\frak q',\Delta_n)\too (\hat Q,\Delta_0)
\eqno £4.4.19\phantom{.}$$
is the $\ch^*(\F^{^{\rm sc}})\-$morphism  formed by the map $\iota_n\,\colon \Delta_0\to \Delta_n$ sending $0$ 
to~$n$ and by the natural map $\nu_{\zeta_*}$ determined by~$\zeta_*\,;$ moreover, according to 
  [5,~Proposition~18.16], we may assume that $\lambda_{\zeta_*}(u) = \zeta (u)$ for any $u\in N_P(\frak f^{^{\frak X'}}_\alpha\circ \hat\frak q')\,.$

 \smallskip  
 On the other hand, denoting by $\hat R'$ the image by 
 $\hat\frak r'(\delta (n)\bullet m)$ of $\hat\frak r' (\delta (n))$ in~$\hat\frak r' (m)\,,$
 it~follows from [5, statement~2.10.1] that there is an $\F'\-$morphism 
 $\xi'\,\colon N_{P'}(\hat R')\to N_{P'}(\hat Q')$ extending the composition
  $\hat\varphi'\,\colon \hat R'\cong \hat Q'$ of the inverse of the isomorphism $\hat\frak r' (\delta (n))\cong \hat R'$ induced by $\hat\frak r'(\delta (n)\bullet m)$ with 
the isomorphism $\hat\mu'_n\,\colon \hat\frak r' (\delta (n))\cong \hat Q'\,;$ then, always from  [5,~Proposition~18.16], we know that, for some element $x'$ of  $\L'(\hat \frak q')$ such that $\lambda'_{\hat\mu'}\big(N_{P'}(\hat\frak r')\big)\i
 N_{P'}(\hat \frak q')^{x'}$ and any $u'\in  N_{P'}(\hat\frak r')\,,$ we may assume that we have
 $\lambda'_{\hat\mu'}(u') = \xi'(u')^{x'}\,.$ Moreover, since we have (cf.~£4.2.3)
$$(\aut_{\frak f^{^{\frak X'}}_\alpha})_{(\hat\frak q',\Delta_n)}\big(\F'_{\!P'}(\hat\frak q')\big)\i
\F_{\!P} (\frak f^{^{\frak X'}}_\alpha\circ \hat\frak q')
\eqno £4.4.20,$$ 
setting $x =  \lambda_{\zeta_*}\big(\lambda_{\hat\frak q'}(x')\big)\,,$ we get 
(cf.~Proposition~£4.3)
 $$\eqalign{(\lambda_{\zeta_*}\circ \lambda_{\hat\frak q'}\circ\lambda'_{\hat\mu'})
 \big(N_{P'}(\frak r')\big)&\i  \lambda_{\zeta_*}\Big(\lambda_{\hat\frak q'}
 \big( N_{P'}(\hat\frak q')\big)\Big)^x\cr
& \i  \lambda_{\zeta_*}\big(N_P(\frak f^{^{\frak X'}}_\alpha\circ \hat\frak q')\big)^x
\i N_P (\hat Q)\cr}
\eqno £4.4.21.$$

 \smallskip
 Similarly, setting $\hat R = \alpha (\hat R')\,,$ it follows from [5, statement~2.10.1] that there is an $\F\-$morphism $\xi\,\colon N_{P}(\hat R)\to N_{P}(\hat Q)$ extending the composition $\hat\varphi\,\colon \hat R\cong \hat Q$ of the inverse of the isomorphism $ \alpha\big(\hat\frak r' (\delta (n))\big)\cong \hat R$ determined  by $\frak f^{^{\frak X'}}_\alpha \big(\hat\frak r'(\delta (n)\bullet m)\big)$  with the isomorphism
 $$\zeta_*\circ (\frak f^{^{\frak X'}}_\alpha * \hat\mu')_n :
 \alpha\big(\hat\frak r' (\delta (n))\big) \cong \hat Q
 \eqno £4.4.22;$$
by [5,~Proposition~18.16], since $\lambda_{\zeta_*}\circ\lambda_{\frak f^{^{\frak X'}}_\alpha * \hat\mu'}$ is a representative of the composition
$\loc_{\F^{^{\rm sc}}}\big((\zeta_*,\iota_n)\circ (\frak f^{^{\frak X'}}_\alpha * \hat\mu',\delta)\big)\,,$ up to modifying our choice of $\zeta$ according to [5, Remark~18.17], we may assume that 
$$(\lambda_{\zeta_*}\circ\lambda_{\frak f^{^{\frak X'}}_\alpha * \hat\mu'})(u) 
= \xi (u)^x
\eqno £4.4.23\phantom{.}$$
for any $u\in N_P (\frak f^{^{\frak X'}}_\alpha\circ \hat\frak r')\,.$

\smallskip
Moreover, since we have $\hat\varphi = \zeta_*\circ\frak f^{^{\frak X'}}_\alpha 
 (\hat\varphi')\,,$ it is easily checked that,
denoting by $\alpha_{\hat Q'}$ and $\alpha_{\hat R'}$ the respective restrictions
 of $\alpha$ from $N_{P'}(\hat Q')$ to $N_P (\hat Q)\,,$ and from 
 $N_{P'}(\hat R')$ to $N_P(\hat R)\,,$  the compositions 
$$\zeta\circ \alpha_{\hat Q'}\circ \xi'\qq  \xi\circ \alpha_{\hat R'}
\eqno £4.4.24\phantom{.}$$
 restricted to $\hat R'$ coincide with each other.  But, $\zeta\circ\alpha_{\hat Q'}
 \circ \xi'$ maps $N_{P'}(\hat R')$ on a subgroup of~$N_P (\hat Q)\,;$ hence, 
 since $\hat R$ is $\F\-$selfcentralizing, $\xi\circ \alpha_{\hat R'}$ maps $N_{P'}(\hat R')$ on the same subgroup 
 of~$N_P  (\hat Q)$ and it follows from [5, Proposition~4.6] that a suitable modification of our choice of~$\xi'$ suffices to guarantee that we get
 $$\zeta\circ \alpha_{\hat Q'}\circ \xi' = \xi\circ \alpha_{\hat R'}
 \eqno £4.4.25.$$
 At this point, for any $u'\in N_{P'}(\hat\frak r')$ we have (cf.~£4.4.23 and~£4.4.25)
 $$\eqalign{(\lambda_{\zeta_*}\circ \lambda_{\hat\frak q'} \circ
 \lambda'_{\hat\mu'})(u') &= (\lambda_{\zeta_*}\circ \lambda_{\hat\frak q'})\big(\xi' (u')^{x'}\big)\cr
 &= \lambda_{\zeta_*}\Big(\alpha \big(\xi' (u')\big)^{\lambda_{\hat\frak q'}(x')}\Big)
= (\zeta\circ \alpha_{\hat Q'}\circ \xi' )(u')^x \cr
&= (\xi\circ \alpha_{\hat R'})(u')^x = (\lambda_{\zeta_*} \circ
\lambda_{\frak f^{^{\frak X'}}_\alpha * \hat\mu'})\big(\alpha (u')\big)\cr
&= (\lambda_{\zeta_*} \circ
\lambda_{\frak f^{^{\frak X'}}_\alpha * \hat\mu'}\circ \lambda_{\hat\frak r'})(u)\cr}
\eqno £4.4.26;$$
in particular, the compositions 
$$ \lambda_{\hat\frak q'} \circ \lambda'_{\hat\mu'}\qq
 \lambda_{\frak f^{^{\frak X'}}_\alpha * \hat\mu'}\circ \lambda_{\hat\frak r'}
 \eqno £4.4.27\phantom{.}$$
  restricted to $N_{P'}(\hat\frak r')$ coincide with each other. Now, as announced above, 
 according to this statement and to equality~£4.4.17, it suffices to apply the uniqueness part of [5,~Lemma~18.8] and [5,~Remark~18.9]. We are done.
 \eject
 
 \medskip
£4.5\phantom{.} Now, it follows from~£3.3 that the natural map in~£4.1.2
$$\widehat\aut_{\frak f^{^{\frak X'}}_\alpha} : \widehat\aut_{\F'^{^{\frak X'}}} \too \widehat\aut_{\F^{^{\rm sc}}} \circ \ch^*(\frak f^{^{\frak X'}}_\alpha)
\eqno £4.5.1\phantom{.}$$
and the natural map in theorem~£4.4 above
$$\loc_{\frak f^{^{\frak X'}}_\alpha} : \loc_{\F'^{^{\frak X'}}}
\too \loc_{\F^{^{\rm sc}}}\circ \ch^*(\frak f^{^{\frak X'}}_\alpha)
\eqno £4.5.2\phantom{.}$$
determine a new natural map
$$\widehat\loc_{\frak f^{^{\frak X'}}_\alpha} : \widehat\loc_{\F'^{^{\frak X'}}}
\too \widehat\loc_{\F^{^{\rm sc}}}\circ \ch^*(\frak f^{^{\frak X'}}_\alpha)
\eqno £4.5.3;$$
then, composing this natural map with the {\it contravariant\/} functors $\frak g_\K$ 
and $\frak g_k$ in~£3.4 above, we get the  $\O\-\mod\-$valued natural maps
$$\eqalign{\frak g_\K *\widehat\loc_{\frak f^{^{\frak X'}}_\alpha} : \frak g_\K\circ
\widehat\loc_{\F^{^{\rm sc}}}\circ \ch^*(\frak f^{^{\frak X'}}_\alpha) &\too
 \frak g_\K\circ \widehat\loc_{\F'^{^{\frak X'}}}\cr
 \frak g_k *\widehat\loc_{\frak f^{^{\frak X'}}_\alpha} : \frak g_k\circ
\widehat\loc_{\F^{^{\rm sc}}}\circ \ch^*(\frak f^{^{\frak X'}}_\alpha) &\too
 \frak g_k\circ \widehat\loc_{\F'^{^{\frak X'}}}\cr}
 \eqno £4.5.4\phantom{.}$$
 which, together with the functor $\ch^*(\frak f^{^{\frak X'}}_\alpha)\,,$ determine 
 the $\O\-$module homomorphisms
 $$\matrix{\G_\K (\F,\widehat\aut_{\F^{^{\rm sc}}})\cr
\big\downarrow\cr
{\displaystyle \lim_{\longleftarrow}}\big(\frak g_\K\circ
\widehat\loc_{\F^{^{\rm sc}}}\circ \ch^*(\frak f^{^{\frak X'}}_\alpha)\big) \cr
\big\downarrow\cr
{\displaystyle \lim_{\longleftarrow}}\big(\frak g_\K\circ \widehat\loc_{\F'^{^{\frak X'}}}\big) \cr}\qq
\matrix{\G_k (\F,\widehat\aut_{\F^{^{\rm sc}}})\cr
\big\downarrow\cr
{\displaystyle \lim_{\longleftarrow}}\big(\frak g_k\circ
\widehat\loc_{\F^{^{\rm sc}}}\circ \ch^*(\frak f^{^{\frak X'}}_\alpha)\big) \cr
\big\downarrow\cr
{\displaystyle \lim_{\longleftarrow}}\big(\frak g_k\circ \widehat\loc_{\F'^{^{\frak X'}}}\big) \cr}
 \eqno £4.5.5.$$

 \medskip
£4.6\phantom{.}  But, from Proposition~£3.7 above and from [5, Proposition~£14.6] together with Proposition~£2.7 above, we know that
$$\eqalign{{\displaystyle \lim_{\longleftarrow}}\big(\frak g_\K \circ 
\widehat\loc_{\F'^{^{\frak X'}}}\big)\cong \G_\K 
(\F',\widehat\aut_{\F'^{^{\rm sc}}})\cr
{\displaystyle \lim_{\longleftarrow}}\big(\frak g_k \circ 
\widehat\loc_{\F'^{^{\frak X'}}}\big)\cong 
\G_k (\F',\widehat\aut_{\F'^{^{\rm sc}}})\cr}
 \eqno £4.6.1.$$
 Consequently, we get $\O\-$module homomorphisms
$$\eqalign{{\rm Res}_{\frak f_\alpha} : \G_\K (\F,\widehat\aut_{\F^{^{\rm sc}}})  &\too \G_\K (\F',\widehat\aut_{\F'^{^{\rm sc}}}) \cr
{\rm res}_{\frak f_\alpha} : \G_k (\F,\widehat\aut_{\F^{^{\rm sc}}})  &\too 
\G_k (\F',\widehat\aut_{\F'^{^{\rm sc}}}) \cr}
\eqno  £4.6.2\phantom{.}$$
which are clearly independent of our choice of $\frak X'\,;$
moreover, the natural map $\partial\,\colon \G_\K\to \G_k$ induces a commutative
diagram
$$\matrix{\G_\K (\F,\widehat\aut_{\F^{^{\rm sc}}})  &\too &\G_\K (\F',\widehat\aut_{\F'^{^{\rm sc}}}) \cr
\hskip-50pt {\scriptstyle\partial_{(\F,\widehat\aut_{\F^{^{\rm sc}}})}} \big\downarrow
&\phantom{\bigg\downarrow}&\big\downarrow
{\scriptstyle \partial_{(\F',\widehat\aut_{\F'^{^{\rm sc}}})}}\hskip-50pt\cr
\G_k (\F,\widehat\aut_{\F^{^{\rm sc}}})  &\too 
&\G_k (\F',\widehat\aut_{\F'^{^{\rm sc}}}) \cr}
\eqno  £4.6.3.$$

\medskip
£4.7\phantom{.}  Let $P''$ be a third finite $p\-$group, $\F''$ a Frobenius
$P''\-$category and $\alpha'\,\colon P''\to P'$ an $(\F'',\F')\-$functorial group homorphism, so that $\alpha\circ\alpha'$ is an $(\F'',\F)\-$functorial group homomorphism. Assume that  $\frak f_{\alpha'}$ and $\frak f_{\alpha\circ\alpha'}$
 map any $\F''\-$radical subgroup of $P''$ on an $\F'\-$ and an 
$\F\-$selfcentralizing subgroups of~$P'$ and $P$ respectively, and  let 
$\frak X''$ be a set of subgroups $Q''$ of~$P''$ such that $\alpha' (Q'')$ belongs to 
$\frak X'\,,$ which contains the $\F''\-$radical subgroups of $P''$ and is stable by 
$\F''\-$isomorphisms. Then, it easily follows from Theorem~£4.4 that we have natural maps
$$\eqalign{\loc_{\frak f^{^{\frak X''}}_{\alpha'}} : \loc_{\F''^{^{\frak X''}}}
&\too \loc_{\F'^{^{\frak X'}}}\circ \ch^*(\frak f^{^{\frak X''}}_{\alpha'})\cr
\loc_{\frak f^{^{\frak X''}}_{\alpha\circ \alpha'}} : \loc_{\F''^{^{\frak X''}}}
&\too \loc_{\F^{^{\rm sc}}}\circ \ch^*(\frak f^{^{\frak X''}}_{\alpha\circ\alpha'})\cr}
\eqno £4.7.1\phantom{.}$$
and therefore from~£4.5 we get the  $\O\-$module homomorphisms
$$\eqalign{{\rm Res}_{\frak f_{\alpha'}} :  \G_\K (\F',\widehat\aut_{\F'^{^{\rm sc}}})  &\too \G_\K (\F'',\widehat\aut_{\F''^{^{\rm sc}}})\cr
{\rm res}_{\frak f_{\alpha'}} :  \G_k (\F',\widehat\aut_{\F'^{^{\rm sc}}})  
&\too \G_k (\F'',\widehat\aut_{\F''^{^{\rm sc}}})\cr
{\rm Res}_{\frak f_{\alpha\circ\alpha'}} :  \G_\K (\F,\widehat\aut_{\F^{^{\rm sc}}})  &\too \G_\K (\F'',\widehat\aut_{\F''^{^{\rm sc}}})\cr
{\rm res}_{\frak f_{\alpha\circ\alpha'}} :  \G_k (\F,\widehat\aut_{\F^{^{\rm sc}}})  
&\too \G_k (\F'',\widehat\aut_{\F''^{^{\rm sc}}})\cr}
\eqno  £4.7.2.$$

\bigskip
\noindent
{\bf Proposition~£4.8}\phantom{.} {\it With the notation and the hypothesis above,
we have
$$\eqalign{{\rm Res}_{\frak f_{\alpha\circ \alpha'}} &= {\rm Res}_{\frak f_{\alpha'}}\circ {\rm Res}_{\frak f_{\alpha}}\cr
{\rm res}_{\frak f_{\alpha\circ \alpha'}} &= {\rm res}_{\frak f_{\alpha'}}\circ {\rm res}_{\frak f_{\alpha}}\cr}
\eqno  £4.8.1.$$\/}
\medskip
\noindent
{\bf Proof:} It is quite clear that 
$$\frak f_{\alpha\circ\alpha'}^{^{\frak X''}} = \frak f_{\alpha}^{^{\frak X'}}\circ 
\frak f_{\alpha'}^{^{\frak X''}}
\eqno £4.8.2,$$
so that $\,\ch^* (\frak f_{\alpha\circ\alpha'}^{^{\frak X''}}) = 
\ch^* (\frak f_{\alpha}^{^{\frak X'}}) \circ \ch^* (\frak f_{\alpha'}^{^{\frak X''}})\,;$
in particular, from~£4.4.1 we get a natural map
$$\loc_{\frak f^{^{\frak X'}}_{\alpha}} * \ch^* (\frak f^{^{\frak X''}}_{\alpha'}) : \loc_{\F'^{^{\frak X'}}}\circ \ch^*(\frak f^{^{\frak X''}}_{\alpha'})\too 
\loc_{\F^{^{\rm sc}}}\circ \ch^*(\frak f^{^{\frak X''}}_{\alpha\circ\alpha'})
\eqno £4.8.3\phantom{.}$$
and we claim that
$$\loc_{\frak f^{^{\frak X''}}_{\alpha\circ \alpha'}} = 
\big(\loc_{\frak f^{^{\frak X'}}_{\alpha}} * \ch^* (\frak f^{^{\frak X''}}_{\alpha'})\big) \circ \loc_{\frak f^{^{\frak X''}}_{\alpha'}}
\eqno £4.8.4.$$
\eject
\noindent

\smallskip
Indeed, by the uniqueness part of Theorem~£4.4, it suffices to prove that, for any 
 $\F''^{^{\frak X''}}\!\-$chain $\frak q''\,\colon \Delta_n\to \F''^{^{\frak X''}}$
 such that $\frak q'' (n)$ is fully normalized in~$\F''$ and that $\F''_{\!P''}\big(\frak q''(n)\big)$ contains a Sylow $p\-$subgroup of $\F'' (\frak q'')\,,$ both members of equality~£4.8.4 map $(\frak q'',\Delta_n)$ on the same $Z\big((\alpha\circ \alpha') (\frak q''(n))\big)\-$conjugacy class of group homomorphisms from $\L''(\frak q'') = \loc_{\F''^{^{\frak X''}}}(\frak q'',\Delta_n)$ to  $\L (\frak f^{^{\frak X''}}_{\alpha\circ\alpha'}\circ \frak q'')\,;$ that is to say, it suffices to prove that 
 $$\lambda_{\frak q''} : \L''(\frak q'')\too \L (\frak f^{^{\frak X''}}_{\alpha\circ\alpha'}\circ \frak q'')
 \eqno £4.8.5\phantom{.}$$
 is $Z\big((\alpha\circ \alpha') (\frak q''(n))\big)\-$conjugate to the composition of a representative $\lambda_{\frak f^{^{\frak X''}}_{\alpha'}\circ \frak q''}$ of
 $$\big(\loc_{\frak f^{^{\frak X'}}_{\alpha}} * \ch^* (\frak f^{^{\frak X''}}_{\alpha'})\big)_{(\frak q'',\Delta_n)}  = (\loc_{\frak f^{^{\frak X'}}_{\alpha}})_{
 (\frak f^{^{\frak X''}}_{\alpha'}\circ \frak q'',\Delta_n)}
  \eqno £4.8.6\phantom{.}$$
  with the corresponding group homomorphism (cf.~Proposition~£4.3)
  $$\lambda'_{\frak q''} : \L''(\frak q'')\too \L (\frak f^{^{\frak X''}}_{\alpha'}
  \circ \frak q'')
 \eqno £4.8.7.$$

 \smallskip
 But, we already know that $\lambda_{\frak q''}\,,$ $\lambda_{\frak f^{^{\frak X''}}_{\alpha'}\circ \frak q''}$ and $\lambda'_{\frak q''}$ respectively lift 
$$(\aut_{\frak f^{^{\frak X''}}_{\alpha\circ \alpha'}})_{(\frak q'',\Delta_n)}\quad ,
\quad (\aut_{\frak f^{^{\frak X'}}_{\alpha}})_{(\frak f^{^{\frak X''}}_{\alpha'}
\circ \frak q'',\Delta_n)}\qq (\aut_{\frak f^{^{\frak X''}}_{\alpha'}
  \circ \frak q''})_{(\frak q'',\Delta_n)}
  \eqno £4.8.8\phantom{.}$$
  and therefore the composition $\lambda_{\frak f^{^{\frak X''}}_{\alpha'}\circ \frak q''}\circ \lambda'_{\frak q''}$ also lifts $(\aut_{\frak f^{^{\frak X''}}_{\alpha\circ \alpha'}})_{(\frak q'',\Delta_n)}\,;$ moreover it follows from Theorem~£4.4 that, for any $u''\in N_{P''}(\frak q'')\,,$ we have 
  $$\lambda_{\frak q''}(u'')   = (\alpha\circ\alpha')(u'')\qq
  \lambda'_{\frak q''}(u'') = \alpha'(u'')
  \eqno £4.8.9,$$
  and that we may assume that $\lambda_{\frak f^{^{\frak X''}}_{\alpha'}\circ \frak q''}\big( \alpha'(u'')\big)  = \alpha\big(\alpha'(u'')\big)\,.$ Now, our claim follows from Proposition~£4.3.

\smallskip
  From the surjectivity of the {\it Brauer decomposition natural 
map\/} $\partial$ (cf.~3.4) and from the commutativity of diagram~£4.6.3, it follows that in~£4.8.1 above it suffices to prove the top equality; but, from our claim and from~£4.5 we get the commutative diagram of natural maps
$$\matrix{\frak g_\K\circ \widehat\loc_{\F^{^{\rm sc}}}\circ \ch^*(\frak f^{^{\frak X''}}_{\alpha\circ\alpha'})
&\buildrel \frak g *\widehat\loc_{\frak f^{^{\frak X''}}_{\alpha\circ\alpha'}} \over{\hbox to 40pt{\rightarrowfill} }
&\frak g_\K\circ \widehat\loc_{\F''^{^{\frak X''}}}\cr
\hskip-60pt{\scriptstyle \frak g *\widehat\loc_{\frak f^{^{\frak X'}}_{\alpha}} * \ch^* (\frak f^{^{\frak X''}}_{\alpha'})}\hskip4pt\big\downarrow&\nearrow\hskip4pt
{\scriptstyle \frak g *\widehat\loc_{\frak f^{^{\frak X''}}_{\alpha'}}}\hskip-55pt&\phantom{\Bigg\downarrow}\cr
\frak g_\K\circ\widehat\loc_{\F'^{^{\frak X'}}}\circ 
\ch^*(\frak f^{^{\frak X''}}_{\alpha'})\cr}
\eqno £4.8.10\phantom{.}$$
which forces the commutative diagram of the corresponding inverse limits;\break
\eject
\noindent
moreover, we have the obvious commutative $\O\-\mod\-$diagram
$$\matrix{{\displaystyle \lim_{\longleftarrow}}\big(\frak g_\K\circ \widehat\loc_{\F^{^{\rm sc}}}\circ \ch^*(\frak f^{^{\frak X'}}_{\alpha})\big)&\too &{\displaystyle \lim_{\longleftarrow}}\big(\frak g_\K\circ \widehat\loc_{\F^{^{\rm sc}}}\circ \ch^*(\frak f^{^{\frak X''}}_{\alpha\circ\alpha'})\big)\cr
\big\downarrow&&\big\downarrow\cr
{\displaystyle \lim_{\longleftarrow}}(\frak g_\K\circ\widehat\loc_{\F'^{^{\frak X'}}}) &\too &{\displaystyle \lim_{\longleftarrow}}\big(\frak g_\K \circ
\widehat\loc_{\F'^{^{\frak X'}}}\circ  \ch^*(\frak f^{^{\frak X''}}_{\alpha'})\big)\cr}
\eqno £4.8.11\phantom{.}$$
and the canonical $\O\-$module homomorphisms (cf.~Proposition~£3.7)
$$\eqalign{\G_\K (\F,\widehat\aut_{\F^{^{\rm sc}}})&\too {\displaystyle \lim_{\longleftarrow}}\big(\frak g_\K\circ \widehat\loc_{\F^{^{\rm sc}}}\circ \ch^*(\frak f^{^{\frak X'}}_{\alpha})\big)\cr
\G_\K (\F',\widehat\aut_{\F'^{^{\rm sc}}})&\cong {\displaystyle \lim_{\longleftarrow}}(\frak g_\K\circ\widehat\loc_{\F'^{^{\frak X'}}}) \cr
\G_\K (\F'',\widehat\aut_{\F''^{^{\rm sc}}})&\cong {\displaystyle \lim_{\longleftarrow}}\big(\frak g_\K\circ \widehat\loc_{\F''^{^{\frak X''}}}\big)\cr}
\eqno £4.8.12.$$

\smallskip
Finally, by the very definition of the ordinary Grothendieck group in~£4.6, we get the
 commutative $\O\-\mod\-$diagram
$$\matrix{\G_\K (\F,\widehat\aut_{\F^{^{\rm sc}}})&\buildrel {\rm Res}_{\frak f_{\alpha\circ \alpha'}}\over{\hbox to 40pt{\rightarrowfill} }
&\G_\K (\F'',\widehat\aut_{\F''^{^{\rm sc}}})\cr
\searrow\hskip-70pt&&\hskip-80pt\nearrow\cr
&{\displaystyle \lim_{\longleftarrow}}\big(\frak g_\K\circ \widehat\loc_{\F^{^{\rm sc}}}\circ \ch^*(\frak f^{^{\frak X''}}_{\alpha\circ\alpha'})\big)\cr
\hskip-20pt{\scriptstyle {\rm Res}_{\frak f_{\alpha}}}\Big\downarrow&\big\downarrow&\Vert\cr
&{\displaystyle \lim_{\longleftarrow}}\big(\frak g_\K \circ
\widehat\loc_{\F'^{^{\frak X'}}}\circ  \ch^*(\frak f^{^{\frak X''}}_{\alpha'})\big)\cr
\nearrow\hskip-70pt&&\hskip-80pt\searrow\cr
\G_\K (\F',\widehat\aut_{\F'^{^{\rm sc}}})&\buildrel {\rm Res}_{\frak f_{\alpha'}}\over{\hbox to 40pt{\rightarrowfill} }&\G_\K (\F'',\widehat\aut_{\F''^{^{\rm sc}}})\cr}
\eqno £4.8.13.$$
We are done.

\medskip
£4.9\phantom{.} We apply these results to the normalizers
and the centralizers of the subgroups of $P\,.$ Let $Q$ be a subgroup of $P$ and 
$K$ a subgroup of ${\rm Aut}(Q)$ containing ${\rm Int}(Q)\,,$ assume that $Q$ 
is fully $K\-$normalized in $\F$ and consider the Frobenius $N_P^K(Q)\-$category 
$N_\F^K(Q)\,;$ it follows from Lemma~£2.5 that any $N_\F^K(Q)\-$radical subgroup $R$ of  $N_P^K(Q)$  contains $Q\,;$ consequently, since the inclusion $\iota_Q^K\,\colon N_P^K(Q) \to P$ is $(N_\F^K(Q),\F)\-$functorial, we have $\O\-$module homomorphisms
$$\eqalign{{\rm Res}_{\iota_Q^K} : \G_\K (\F,\widehat\aut_{\F^{^{\rm sc}}}) &\too 
 \G_\K \big(N_\F^K (Q),\widehat{\aut}_{N_{\F}^K (Q)^{^{\rm sc}}}\big)\cr
 {\rm res}_{\iota_Q^K} : \G_k (\F,\widehat\aut_{\F^{^{\rm sc}}}) &\too 
 \G_k \big(N_\F^K (Q),\widehat{\aut}_{N_{\F}^K (Q)^{^{\rm sc}}}\big)\cr}
 \eqno £4.9.1\phantom{.}$$
and a commutative diagram 
$$\matrix{\G_\K (\F,\widehat\aut_{\F^{^{\rm sc}}}) &\too 
 &\G_\K \big(N_\F^K (Q),\widehat{\aut}_{N_{\F}^K (Q)^{^{\rm sc}}}\big)\cr
 \hskip-50pt {\scriptstyle\partial_{(\F,\widehat\aut_{\F^{^{\rm sc}}})}}\big\downarrow&\phantom{\Big\downarrow}&\big\downarrow  
 {\scriptstyle\partial_{(N_\F^K (Q),\widehat\aut_{N_\F^K (Q)^{^{\rm sc}}})}} \hskip-70pt\cr
\G_k (\F,\widehat\aut_{\F^{^{\rm sc}}}) &\too 
 &\G_k \big(N_\F^K (Q),\widehat{\aut}_{N_{\F}^K (Q)^{^{\rm sc}}}\big)\cr}
 \eqno £4.9.2.$$

 \medskip
£4.10\phantom{.} Moreover,  let $R$ and $J$ be respective subgroups of $N_P^K (Q)$ and ${\rm Aut} (R)\,,$ denote by $I$ the subgroup of automorphisms
of~$Q\.R$ which stabilize $Q$ and~$R\,,$ and act on them {\it via\/}
elements of $K$ and $J$ respectively, and assume that $R$ is fully $J\-$normalized in $N_\F^K(Q)\,;$ thus, $Q\.R$ is fully$I\-$normalized in~$\F$ [5, Lemma~2.17] and
we have
$$N_P^I( Q\.R) = N_{N_P^K(Q)}^J (R)\qq N_\F^I( Q\.R) = N_{N_\F^K(Q)}^J (R)
\eqno £4.10.1;$$
then, according to Proposition~£4.8 above, we have  the obvious commutative diagram
$$\matrix{\G_\K (\F,\widehat\aut_{\F^{^{\rm sc}}})\hskip-30pt &\too 
 &\hskip-30pt\G_\K \big(N_\F^I (Q\.R),\widehat{\aut}_{N_{\F}^I (Q\.R)^{^{\rm sc}}}\big)\cr
 \searrow\hskip-50pt&\phantom{\Big\downarrow}&\hskip-90pt\nearrow\cr
\big\downarrow&\G_\K \big(N_\F^K (Q),\widehat{\aut}_{N_{\F}^K (Q)^{^{\rm sc}}}\big)&\big\downarrow\cr
&\big\downarrow\cr
\G_k (\F,\widehat\aut_{\F^{^{\rm sc}}})\hskip-30pt &\too 
 &\hskip-30pt\G_k \big(N_\F^I (Q\.R),\widehat{\aut}_{N_{\F}^I (Q\.R)^{^{\rm sc}}}\big)\cr
 \searrow\hskip-50pt&\phantom{\Big\downarrow}&\hskip-90pt\nearrow\cr
&\G_k \big(N_\F^K (Q),\widehat{\aut}_{N_{\F}^K (Q)^{^{\rm sc}}}\big)\cr}
\eqno £4.10.2 .$$

\bigskip
\bigskip
\noindent
{\bf £5\phantom{.} Character decomposition of the functor $\frak g_\K$  }
\bigskip

£5.1\phantom{.} In order to determine the $\O\-$rank of~$\G_\K
(\F,\widehat\aut_{\F^{^{\rm sc}}})\,,$ we need  a suitable decomposition of the functor  $\frak g_\K$ --- analogous to the decomposition of~$\frak g_k$ in [5, 14.22] --- that we develop below. First of all, for any $h\in \Bbb N -\{0\}\,,$ fix a primitive $h\-$th root of 
unity~$\xi_h$ in $\K$ and set
$U_h = \langle \,\xi_h\rangle $  and $\hat U_h = k^*\times U_h\,;$ 
 note that the kernel of the canonical group homomorphism
$${\rm Aut}_{k^*}(\hat U_h)\too {\rm Aut} (U_h)
\eqno  £5.1.1\phantom{.}$$
can be identified with ${\rm Hom} (U_h,k^*)$ and that, denoting by $U'_h$ the subgroup of $p'\-$elements of $U_h\,,$ $\xi_h$ determines a group isomorphism 
$${\rm Hom} (U_h,k^*)\cong U'_h
\eqno £5.1.2.$$

\medskip
£5.2\phantom{.} Let $\hat G$ be a $k^*\-$group with finite $k^*\-$quotient $G\,;$ 
 for any injective $k^*\-$group homomorphism  $\hat\eta\,\colon \hat U_h \to 
\hat G\,,$ the so-called {\it ordinary characters\/}
determine an $\O\-$module homomorphism $\G_\K (\hat G)\to \O$ mapping the
class in $\G_\K (\hat G)$ of a  $\K_*\hat G\-$module~$M$ on the linear trace 
${\rm tr}_M \big(\hat\eta (1,\xi_h)\big)\,;$ actually, this $\O\-$module 
homomorphism only depends on the $G\-$conjugacy class of $\hat\eta\,.$ 
Moreover, if we have $\hat\sigma (1,\xi_h) = (\lambda,\xi_h)$ for some $k^*\-$automorphism $\hat\sigma$ 
of $\hat U_h$ inducing the identity on $U_h$ then we still have 
$${\rm tr}_M\big((\hat\eta\circ\hat\sigma) (1,\xi_h)\big) = \hat\lambda\.{\rm tr}_M
\big(\hat\eta (1,\xi_h)\big)
\eqno £5.2.1\phantom{.}$$
where $\hat\lambda$ denotes the corresponding lifting of $\lambda\in k^*$ to  $\O^*$ (cf.~£1.7). That is to say, denoting by~${\rm Mon}_{k^*}(\hat U_h,\hat G)$ the set of injective $k^*\-$homomorphisms from  $\hat U_h$ 
to $\hat G\,,$ $G$ acts by conjugation on this set and  $U'_h$ acts on it {\it via\/} isomorphism~£5.1.2 
centralizing the action of $G\,,$ and on $\O$ {\it via\/} the inclusion~$U'_h\i \O^*\,;$ moreover, equality~£5.2.1 shows that these actions are preserved by our correspondence. In conclusion, 
denoting by $\F\!ct_{U'_h}\big({\rm Mon}_{k^*} (\hat U_h,\hat G), \O\big)$ the set of $\O\-$valued functions which preserve the corresponding $U'_h\-$actions, we have  obtained an $\O\-$module homomorphism
$$\G_\K (\hat G)\too \F\!ct_{U'_h}\big({\rm Mon}_{k^*}(\hat U_h,\hat G),\O\big)^G
\eqno £5.2.2.$$

\medskip
£5.3\phantom{.} On the other hand, since we have $\hat G = k^*\. G'$ for a suitable 
finite subgroup $G'$ of $\hat G$ and then, setting
$$Z' = k^*\cap G'\qq e' = {1\over \vert Z'\vert}\.\sum_{z'\in Z'} z'
\eqno £5.3.1,$$
  we have $\G_\K (\hat G)\cong \G_\K (\K G'e')$ [4, Proposition~5.15],
it is not difficult to prove that homomorphisms~£5.2.2 when $h$ runs over $\Bbb N -\{0\}$ determine a $\K\-$module isomorphism
$${}^\K \G_\K (\hat G)\cong \prod_{h\in \Bbb N -\{0\}} \F\!ct_{U'_h}
\big({\rm Mon}_{k^*}(\hat U_h,\hat G),\K\big)^G
\eqno  £5.3.2\phantom{.}$$
where  $\F\!ct_{U'_h}\big({\rm Mon}_{k^*} (\hat U_h,\hat G), \K\big)$ is the 
set of $\K\-$valued $U'_h\-$invariant  functions.
Note that any function $\chi\in \F\!ct_{U_h}\big({\rm Mon}_{k^*}
(\hat U_h,\hat G), \K\big)$ {\it vanish\/} over the $k^*\-$group monomorphisms  $\hat\eta\,\colon  \hat U_h\to \hat G$ mapping $(1,\xi_h)$ on an element $\hat x\in \hat G$ such that, setting $x = k^*\.\hat x\,,$ the image of~$C_{\hat G}  (\hat x)$ in $G$ is a proper subgroup of~$C_G (x)\,.$

\medskip
£5.4\phantom{.} Pushing it further, denote by~${\rm Mon}(U_h,G)$ the set of injective group homomorphisms from $U_h$ to~$G$ and by 
$$\varpi_{h,\hat G} :  {\rm Mon}_{k^*}(\hat U_h,\hat G)\too {\rm Mon}(U_h,G) 
\eqno £5.4.1\phantom{.}$$
the canonical map, so that $U'_h$ acts {\it regularly\/} on the {\it fibers\/} (cf.~£5.1.2); thus, we still have the obvious decomposition
$$\F\!ct_{U'_h}\big({\rm Mon}_{k^*}(\hat U_h,\hat G),\O\big)\cong 
\prod_{\eta\in {\rm Mon}(U_h,G) }\F\!ct_{U'_h} 
\big((\varpi_{h,\hat G})^{-1}(\eta),\O\big)
\eqno £5.4.2\phantom{.}$$
and therefore we get
$${}^\K \G_\K (\hat G)\cong \prod_{h\in \Bbb N -\{0\}}\,\Big(\prod_{\eta\in {\rm Mon}(U_h,G)} \F\!ct_{U'_h} \big((\varpi_{h,\hat G})^{-1}(\eta),\K\big)\Big)^G
\eqno £5.4.3.$$
Note that, for any $\eta\in {\rm Mon}(U_h,G)\,,$ the term $\F\!ct_{U'_h} 
\big((\varpi_{h,\hat G})^{-1}(\eta),\O\big)$ is a {\it free $\O\-$module of rank one\/}.

\medskip
£5.5\phantom{.} As in [5, 14.16], isomorphisms~£5.3.2 and~£5.4.3 are actually 
{\it natural\/} and, in order to show this {\it naturality\/}, we have to develop a suitable functorial framework. Moreover,
in our present situation, we have to extend our construction to  the respective subcategories
 $k^*\-\frak i\widetilde\Loc$ and $\frak i\widetilde\Loc$ of~$k^*\-\widetilde\Loc$ (cf.~£3.3) and 
 $\widetilde\Loc$ (cf.~£1.8) formed by the same objects and by the classes of {\it injective\/} homomorphisms. As in [5,~14.16], denote 
by~${}^{U'_h}\aleph$ the category  of finite sets endowed with a $U'_h\-$action and by 
$$\res^{U'_h}_1 : {}^{U'_h}\aleph\too \aleph
\eqno £5.5.1\phantom{.}$$
 the corresponding {\it forgetful\/} functor; then, we consider the evident functors
$$\hat\frak u_h : k^*\-\frak i\widetilde\Loc\too {}^{U'_h}\aleph \qq \frak u_h : 
\frak i\widetilde\Loc\too \aleph
\eqno £5.5.2\phantom{.}$$
mapping any $k^*\-\frak i\widetilde\Loc\-$object $(\hat L,Z)$ on the $U'_h\-$set
$\overline{\rm Mon}_{k^*}(\hat U_h,\hat L)$ of $Z\-$conju-gacy classes in ${\rm Mon}_{k^*}(\hat U_h,\hat L)\,,$ 
and any $\frak i\widetilde\Loc\-$object $(L,Z)$ on the corresponding set
$\overline{\rm Mon} (U_h,G)$ respectively; moreover, denoting by $\qt\,\colon k^*\-\frak i\widetilde\Loc\to \frak i\widetilde\Loc$ the obvious {\it $k^*\-$quotient\/} functor, we clearly have a natural map $$\varpi_h : \res^{U'_h}_1 \circ \hat\frak u_h\too \frak u_h \circ\qt
\eqno £5.5.3.$$
 sending~$\hat G$ to~$\varpi_{h,\hat G}$ (cf.~£5.4.1).

\medskip
£5.6\phantom{.} Furthermore, identifying the category of sets $\aleph$ with the full
subcategory~of the {\it category of small categories\/} $\SC$ [5, A1.6] over the small
categories which have no other morphisms than the corresponding {\it identity\/} 
morphisms, we can consider the functor 
$$\frak u_h\circ\qt : k^*\-\frak i\widetilde\Loc\too \aleph\i \CC
\eqno £5.6.1\phantom{.}$$
 as a so-called  {\it representation\/} of~$k^*\-\frak i\widetilde\Loc$ [5, A2.2] and then we can consider the corresponding {\it semidirect product\/} [5, A2.7] 
$$ k^*\-({\frak u_h}\!\rtimes\frak i\widetilde\Loc)= (\frak u_h \circ \qt) \rtimes ( k^*\-\frak i\widetilde\Loc)
\eqno £5.6.2\phantom{.}$$
  where the objects are the pairs $(\bar\eta,\hat \L)$ formed by a $k^*\-\frak i\widetilde\Loc\-$object 
  $\hat\L = (\hat L,Z)$ and, setting $(L,Z)= \qt(\hat\L)\,,$ by the $Z\-$conjugacy class of an injective group homomorphism $\eta\,\colon U_h\to L\,,$ and the morphisms from  $(\bar\eta,\hat\L)$ to a 
  $k^*\-\frak i\widetilde\Loc\-$ob-ject  $(\bar\eta',\hat\L')$ are the $k^*\-\frak i\widetilde\Loc\-$morphisms 
  $\skew3\hat{\bar\varphi}\,\colon \hat\L\to \hat\L'$ fulfilling $\bar\eta' = \qt (\skew3\hat{\bar\varphi})
  \circ \bar\eta\,.$

\medskip
£5.7\phantom{.} Coherently, for any $h'$ dividing~$h\,,$ choosing the identification 
between ${\rm Hom}(U_{h'},k^*)$ and $U_{h'}$ determined by~$(\xi_h)^{h/h'}$ (cf.~£5.1), the
inclusion $U_{h'}\i U_h$ induces a functor and two natural maps
$$\eqalign{\res^{U'_h}_{U'_{h'}} : {}^{U'_h}\aleph&\too {}^{U'_{h'}}\aleph\cr
\hat\rho_{h',h} : \res^{U'_h}_{U'_{h'}}\circ\hat\frak u_h\too \hat\frak u_{h'}
&\qq \rho_{h',h} : \frak u_h\too \frak u_{h'}\cr}
\eqno £5.7.1,$$
and it is easily checked that [5, A1.5.1]
$$\varpi_{h'}\circ (\hat\rho_{h',h} * \res^{U'_h}_{U'_{h'}}) = (\rho_{h',h} * \qt)\circ 
\varpi_h
\eqno £5.7.2.$$
 Moreover, the natural map $\rho_{h',h}$ above determines a functor
[5,~Proposition~A2.17]
$$ (\rho_{h',h} * \qt)\rtimes \id_{k^*\-\frak i \widetilde\Loc} : k^*\-({\frak u_h}\!\rtimes\frak i\widetilde\Loc)\too k^*\-({\frak u_{h'}}\!\rtimes\frak i\widetilde\Loc)
\eqno £5.7.3\phantom{.}$$
mapping $(\bar\eta,\hat\L)$ on $(\bar\eta',\hat\L)$ where $\bar\eta'$ denotes the $Z\-$conjugacy class of the restriction of $\eta\in \bar\eta\,.$

\bigskip
\noindent
{\bf Proposition~£5.8}\phantom{.} {\it With the notation above, we have a functor
$$\hat\frak w_h : k^*\-({\frak u_h}\!\rtimes\frak i\widetilde\Loc) \too {}^{U'_h}\aleph
\eqno £5.8.1\phantom{.}$$
mapping any $k^*\-({\frak u_h}\!\rtimes\frak i\widetilde\Loc)\-$object 
$(\bar\eta,\hat\L)$  on the regular $U'_h\-$set $(\varpi_{h,\hat\L})^{-1}
(\bar\eta)$ and mapping any  $k^*\-({\frak u_h}\!\rtimes\frak i\widetilde\Loc)\-$morphism 
$\skew3\hat{\bar\varphi}\,\colon (\bar\eta,\hat\L)\to (\bar\eta',\hat\L')$ on the bijective map
$$(\varpi_{h,\hat\L})^{-1}(\bar\eta)\cong (\varpi_{h,\hat\L'})^{-1}(\bar\eta')
\eqno £5.8.2\phantom{.}$$
 determined by $\,\hat\frak u_h (\skew3\hat{\bar\varphi})\,.$ Moreover, $\hat\frak u_h$ is the {\it
direct image\/}  of $\hat\frak w_h$ {\it via\/} the structural functor 
$$\hat\frak p_h : k^*\-({\frak u_h}\!\rtimes\frak i\widetilde\Loc)\too
k^*\-\frak i\widetilde\Loc
\eqno £5.8.3.$$\/}

\par
\noindent
{\bf Proof:} As we mention above, the $\,k^*\-({\frak u_h}\!\rtimes\frak i
\widetilde\Loc)\-$objects  are the pairs formed by a $k^*\-\frak i\widetilde\Loc\-$ob-ject $\hat\L = (\hat L,Z)$ and, setting $(L,Z) = \qt (\hat\L)\,,$ by an object of the ``category'' 
$\overline{\rm Mon}(U_h,L)\,,$ namely by the $Z\-$conjugacy class $\bar\eta$
of an injective group homomorphism $\eta\,\colon U_h\to L\,;$ then, in the map (cf.~£5.4.1)
$$\varpi_{h,\hat\L} : \overline{\rm Mon}_{k^*}(\hat U_h,\hat L)\too 
\overline{\rm Mon}(U_h,L)
\eqno £5.8.4,$$
\eject
\noindent
$U'_h$ acts on the left end, stabilizing and acting regularly on the fibers, so that
$(\varpi_{h,\hat\L})^{-1}(\bar\eta)$ is indeed a {\it $U'_h\-$set\/}. Analogously, the
$k^*\-({\frak u_h}\!\rtimes\frak i\widetilde\Loc)\-$morphisms between two objects 
$(\bar\eta,\hat\L)$ and $(\bar\eta',\hat\L')$ are the pairs formed by a 
$k^*\-\frak i\widetilde\Loc\-$mor-phism $\skew3\hat{\bar\varphi}\,\colon \hat\L\to \hat\L'$ and by 
a ``$\overline{\rm Mon}(U_h,L')\-$morphism'' from  the image of 
$\bar\eta$ by the functor $(\frak u_h \circ \qt)(\skew3\hat{\bar\varphi})$ 
to~$\bar\eta'\,,$ which actually forces the equality $\qt (\skew3\hat{\bar\varphi})\circ \bar\eta  = \bar\eta'\,.$ Then, it is quite clear that the map
$$\hat\frak u (\skew3\hat{\bar\varphi}) : \overline{\rm Mon}_{k^*}(\hat U_h,\hat L)\too \overline{\rm Mon}_{k^*}(\hat U_h,\hat L')
\eqno £5.8.5\phantom{.}$$
sends $(\varpi_{h,\hat\L})^{-1}(\bar\eta)$ bijectively onto $(\varpi_{h,\hat\L'})^{-1} (\bar\eta')\,,$ determining a {\it $U'_h\-$set map\/}.  The proofs of the
functoriality and of the last statement are straightforward.

\medskip
\noindent
{\bf Remark £5.9}\phantom{.} Note that, for any $\xi\in U_h\,,$ the inner $k^*\-$group automorphism of~$\hat L$ determined by an element of $\hat L$ lifting $\eta (\xi)$ acts trivially on $(\varpi_{h,\hat\L})^{-1}(\bar\eta)\,.$

\bigskip
\noindent
{\bf Proposition~£5.10}\phantom{.} {\it With the notation above, for any $h'$ 
dividing $h$ we have a natural map 
$$\hat\tau_{h',h} : {\res}^{U'_h}_{U'_{h'}}\circ \hat\frak w_h\too \hat\frak w_{h'} \circ
\big((\rho_{h',h} * \qt)\rtimes \id_{k^*\-\frak i \widetilde\Loc}\big)
\eqno £5.10.1\phantom{.}$$
which sends  any $\,k^*\-({\frak u_h}\!\rtimes\frak i \widetilde\Loc)\-$object  $(\bar\eta,\hat\L)$ to 
the $U'_{h'}\-$morphism
$$\res^{U'_h}_{U'_{h'}}\big((\varpi_{h,\hat\L})^{-1}(\bar\eta)\big)\too 
(\varpi_{h',\hat\L})^{-1}\big({\rm Res}^{U_h}_{U_{h'}}(\bar\eta)\big)
\eqno £5.10.2\phantom{.}$$
mapping any $\skew2\hat{\bar\eta}\in (\varpi_{h,\hat\L})^{-1}(\bar\eta)\i  \overline{\rm Mon}_{k^*}
(\hat U_h,\hat L)$ on its restriction to $\hat U_{h'}\,.$\/}

\medskip
\noindent
{\bf Proof:} For any $k^*\-\frak i \widetilde\Loc\-$morphism $\skew3\hat{\bar\varphi}\,\colon 
\hat\L\to  \hat\L'\,,$ setting  $\bar\varphi = \qt (\skew3\hat{\bar\varphi})\,,$ the functor 
${\res}^{U'_h}_{U'_{h'}}\circ \hat\frak w_h$ maps the
$\,k^*\-({\frak u_h}\!\rtimes\frak i \widetilde\Loc)\-$morphism $(\bar\eta,\hat\L)\to 
(\bar\varphi\circ \bar\eta,\hat\L')$ on the $U'_{h'}\-$set map
$$\res^{U'_h}_{U'_{h'}}\big((\varpi_{h,\hat\L})^{-1}(\bar\eta)\big)\too
\res^{U'_h}_{U'_{h'}}\big((\varpi_{h,\hat\L'})^{-1}(\bar\varphi\circ \bar\eta)\big)
\eqno £5.10.3\phantom{.}$$
sending $\skew2\hat{\bar\eta}\in (\varpi_{h,\hat\L})^{-1}(\bar\eta)$ to
$ \skew3\hat{\bar\varphi}\circ \skew2\hat{\bar\eta}\,;$ whereas $\hat\frak w_{h'} \circ \big((\rho_{h',h} * \qt)\rtimes \id_{k^*\-\frak i \widetilde\Loc}\big)$ maps this morphism on the analogous
$U'_{h'}\-$set map
$$(\varpi_{h',\hat\L})^{-1}(\bar\eta')\too (\varpi_{h',\hat\L'})^{-1}( \bar\varphi \circ \bar\eta')
\eqno £5.10.4\phantom{.}$$
where we are setting $\bar\eta' = {\rm Res}^{U_h}_{U_{h'}}(\bar\eta)\,;$ thus, the corresponding diagram is indeed commutative since the restriction to $U_{h'}$ 
is compatible with the composition with $\skew3\hat{\bar\varphi}$ on the left. We are done.

\medskip
£5.11\phantom{.} We are ready to discuss the {\it naturality\/} of isomorphisms~£5.3.2 and £5.4.3. As in [5,~14.21], consider the evident  {\it contravariant\/} functor 
$$\F\!ct_{U'_h} : {}^{U'_h}\aleph\too \O\-\mod
\eqno £5.11.1\phantom{.}$$
mapping any finite $U'_h\-$set~$X$ on the $\O\-$module $\F\!ct_{U'_h} (X,\O)$ of the $\O\-$valued functions 
over~$X$ preserving the $U'_h\-$actions; note that if  $\xi\. x = x$ for some $x\in X$ and some $\xi\in U_h -\{1\}$ then we have $f(x) = 0$ for any $f\in \F\!ct_{U'_h} (X,\O)\,.$
 On the other hand, note that if we have a {\it contravariant\/}
functor 
$$\frak m : k^*\-\frak i\widetilde\Loc\too \K\-\mod
\eqno £5.11.2,$$
 for any $k^*\-\frak i\widetilde\Loc\-$object $\hat\L =(\hat L,Z)\,,$ setting 
 $(L,Z) = \qt(\hat\L)\,,$ $\frak m (\hat\L)$ has an obvious $\K L\-$module structure, so that it makes sense to
consider
$$\Bbb H^0 \big(L,\frak m(\hat\L)\big) = \frak m (\hat\L )^L
\eqno £5.11.3;$$ 
further, if $\skew3\hat{\bar\varphi}\,\colon \hat\L\to \hat\L' = (\hat L',Z')$ is a $k^*\-\frak i\widetilde\Loc\-$morphism, it is easily checked that $\frak m (\skew3\hat{\bar\varphi})$ maps $\frak m (\hat\L')^{L'}$ on 
an $\K\-$submodule of $\frak m (\hat\L)^L\,;$ that is to say, we get a new   {\it contravariant\/} functor 
from $k^*\-\frak i\widetilde\Loc$ to $\K\-\mod$
--- noted
$\frak h^0(\frak m)$ --- mapping $\hat\L$ on $\frak m (\hat\L )^L$ and $\skew3\hat{\bar\varphi}$ on the map from $\frak m (\hat\L')^{L'}$ to~$\frak m (\hat\L)^L$ induced by~$\frak m (\skew3\hat{\bar\varphi})\,.$

\medskip
£5.12\phantom{.} Finally, still denote by  $\frak g_\K\,\colon  k^*\-\frak i \widetilde\Loc\to \O\-\mod$ the obvious functor mapping any  $ k^*\-\frak i \widetilde\Loc\-$object $(\hat L,Z)$ on $\G_\K (\hat L)\,.$
With all this notation, it is now quite clear that isomorphism~£5.3.2 actually defines a natural isomorphism
$${}^\K \frak  g_\K\cong \prod_{h\in \Bbb N -\{0\}} \frak h^0({}^\K \F\!ct_{U'_h}\circ  \hat\frak u_h)
\eqno £5.12.1\phantom{.}$$
where ${}^\K g_\K$ and ${}^\K \F\!ct_{U'_h}$ denote the respective compositions of $\frak g_\K$ and 
$\F\!ct_{U'_h}$ with the scalar extension from $\O$ to $\K\,.$
Consequently, considering the composition ${}^\K \frak  g_\K\circ \widehat\loc_{\F^{^{\rm sc}}}$ in~£3.5.1 above, from definition~£3.5.2 we have
$${}^\K \G_\K (\F,\widehat\aut_{\F^{^{\rm sc}}}) = \lim_{\longleftarrow} \big(\prod_{h\in \Bbb N -\{0\}}\,
\frak h^0 ({}^\K \F\!ct_{U'_h}\circ \hat\frak u_h)\circ \widehat\loc_{\F^{^{\rm sc}}}\big)
\eqno £5.12.2;$$
that is to say, we still have 
$${}^\K \G_\K (\F,\widehat\aut_{\F^{^{\rm sc}}}) =
\prod_{h\in \Bbb N -\{0\}} {}^\K \G_\K (\F,\widehat\aut_{\F^{^{\rm nc}}})_h
\eqno £5.12.3\phantom{.}$$
where, for any $h\in \Bbb N -\{0\} \,,$ we set
$$\eqalign{{}^\K \G_\K (\F,\widehat\aut_{\F^{^{\rm sc}}})_h &= {\displaystyle\lim_{\longleftarrow}}
\,\big(\frak h^0 ({}^\K \F\!ct_{U'_h}\circ \hat\frak u_h)\circ  \widehat\loc_{\tilde\F^{^{\rm sc}}}\big)\cr
&\cong {\displaystyle\lim_{\longleftarrow}}
\, ({}^\K \F\!ct_{U'_h}\circ \hat\frak u_h \circ  \widehat\loc_{\tilde\F^{^{\rm sc}}})\cr}
\eqno £5.12.4,$$
the last isomorphism being obvious.

\medskip
£5.13\phantom{.} But, according to Proposition~£5.8 above, the functor
$$\hat\frak u_h : k^*\-\frak i\widetilde\Loc\too {}^{U'_h}\aleph
\eqno £5.13.1\phantom{.}$$
is the {\it direct image\/} of the functor 
$$\hat\frak w_h : k^*\-(\frak u_h\!\rtimes\frak i\widetilde\Loc) =  \too {}^{U'_h}\aleph
\eqno £5.13.2\phantom{.}$$
throughout the structural functor 
$$\hat\frak p_h : k^*\-(\frak u_h\!\rtimes\frak i\widetilde\Loc) =(\frak u_h \circ \qt) \rtimes ( k^*\-\frak i\widetilde\Loc)\too 
k^*\-\frak i\widetilde\Loc
\eqno £5.13.3.$$
Moreover, considering the semidirect product
$${}^{\frak u_h}\frak l\ch^* (\F^{^{\rm sc}}) = (\frak u_h\circ 
\loc_{\F^{^{\rm sc}}}) \!\rtimes  \ch^* (\F^{^{\rm sc}}) 
\eqno £5.13.4,$$
we have the evident commutative diagram of functors
$$\matrix{k^*\-(\frak u_h\!\rtimes\frak i\widetilde\Loc) &
\buildrel \hat\frak p_h\over\too & k^*\-\frak i\widetilde\Loc\cr
\hskip-70pt{\scriptstyle {\rm id}_{\frak u_h \circ \loc_{\F^{^{\rm sc}}}} \!\rtimes \widehat\loc_{\F^{^{\rm sc}}}}\big\uparrow &\phantom{\Big\uparrow}
&\big\uparrow {\scriptstyle\widehat\loc_{\F^{^{\rm sc}}}}\hskip-30pt\cr 
{}^{\frak u_h}\frak l\ch^* (\F^{^{\rm sc}}) &\too & \ch^* (\F^{^{\rm sc}})\cr}
\eqno £5.13.5.$$
Then, it is not difficult to check that the composition $\hat\frak u_h
\circ \widehat\loc_{\F^{^{\rm sc}}}$ is also the {\it direct image\/} of the
composition
$$\hat\frak w\frak l_h = \hat\frak w_h\circ ({\rm id}_{\frak u_h \circ 
\loc_{\F^{^{\rm sc}}}} \!\!\rtimes \widehat\loc_{\F^{^{\rm sc}}}) :
{}^{\frak u_h}\frak l\ch^* (\F^{^{\rm sc}})\too {}^{U'_h}\aleph
\eqno £5.13.6\phantom{.}$$
throughout the bottom functor in diagram~£5.13.5; further, the
 {\it direct image\/} is clearly compatible with the functor
$\F\!ct_{U'_h}\,\colon {}^{U'_h}\aleph\too \O\-\mod\,.$ Thus, we finally get [5,~1.6]
$${\displaystyle\lim_{\longleftarrow}} \,(\F\!ct_{U'_h}\circ \hat\frak u_h\circ \widehat\loc_{\F^{^{\rm sc}}}) \cong {\displaystyle\lim_{\longleftarrow}}
\,(\F\!ct_{U'_h}\circ \hat\frak w\frak l_h)
\eqno £5.13.7.$$
\eject

\bigskip
\bigskip
\noindent
{\bf £6\phantom{.}  An equivalence of categories }
\bigskip

£6.1\phantom{.} With the notation above, the point is that the semidirect product 
$${}^{\frak u_h}\frak l\ch^* (\F^{^{\rm sc}}) = (\frak u_h\circ 
\loc_{\F^{^{\rm sc}}}) \!\rtimes  \ch^* (\F^{^{\rm sc}}) 
\eqno £6.1.1\phantom{.}$$
 admits another description in terms of the following category 
 ${}^{^h}\!(\F^{^{\rm sc}})\,.$ The ${}^{^h}\!(\F^{^{\rm sc}})\-$objects are the pairs $Q^{\bar\rho}$ formed by an $\F\-$selfcentralizing subgroup~$Q$ of~$P$ and by a $Z(Q)\-$conjugacy class of {\it injective\/} group homomorphisms $\bar\rho\,\colon U_h\to  \L (Q)$ or, equivalently, an $\frak i\widetilde\Loc\-$morphism 
$$\bar\rho : (U_h,1)\too \big(\L (Q), Z(Q)\big)
\eqno £6.1.2,$$
 whereas the  ${}^{^h}\!(\F^{^{\rm sc}})\-$morphisms from 
another ${}^{^h}\!(\F^{^{\rm sc}})\-$object $R^{\bar\sigma}$ to  $Q^{\bar\rho}$ are the 
$\F^{^{\rm sc}}\-$morphisms $\varphi\,\colon R\to Q$ such that, denoting by 
$\L (\varphi)$ the {\it localizer\/} of the $\F^{^{\rm sc}}\-$chain $\Delta_1\to \F^{^{\rm sc}}$ determined by $\varphi\,,$ there is an  $\frak i\widetilde\Loc\-$morphism 
$$\bar\alpha :  (U_h,1)\too \big(\L (\varphi),Z(Q)\big)
\eqno £6.1.3$$
such that we have the following commutative $\frak i\widetilde\Loc\-$diagram
 $$\matrix{\big(\L (R), Z(R)\big)&\buildrel \bar\sigma \over\longleftarrow &(U_h,1)&\buildrel \bar\rho \over\too &\big(\L (Q), Z(Q)\big)\cr
{\scriptstyle \loc_{\tilde\F^{^{\rm sc}}}({\rm id}_Q,\delta^0_0)}\nwarrow\hskip-30pt&\phantom{\Big\uparrow} &\hskip-10pt{\scriptstyle \bar\alpha}\downarrow&&\hskip-30pt\nearrow{\scriptstyle \loc_{\tilde\F^{^{\rm sc}}}
({\rm id}_R,\delta^0_1)}\cr
 &&\big(\L (\varphi),Z(Q)\big)\cr}
 \eqno £6.1.4;$$
note that such an $\frak i\widetilde\Loc\-$morphism $\bar\alpha$ is unique
and that $\L (\varphi)$ determines well-defined subgroups of $\L (Q)$ and $\L (R)\,.$

\medskip
£6.2\phantom{.} The composition of the $\F\-$morphisms induces a composition in 
${}^{^h}\!(\F^{^{\rm sc}})$ since, for a third ${}^{^h}\!(\F^{^{\rm sc}})\-$object $T^{\bar\tau}$ 
and an ${}^{^h}\!(\F^{^{\rm sc}})\-$morphism $\psi\,\colon T^{\bar\tau}\to R^{\bar\sigma}\,,$ 
denoting by $\L (\psi,\varphi)$ the localizer of the $\F^{^{\rm sc}}\-$chain $\Delta_2\to \F^{^{\rm sc}}$ determined by $\varphi$ and $\psi\,,$ we have the {\it $\frak i\widetilde\Loc\-$pull-back\/}
$$\matrix{&&\big(\L (R),Z(R)\big)\cr
&{\loc_{\F^{^{\rm sc}}}({\rm id}_T,\delta^0_1)\atop}\hskip-4pt\nearrow\hskip-20pt&
&\hskip-20pt\nwarrow\hskip-4pt{\loc_{\F^{^{\rm sc}}}({\rm id}_Q,\delta^0_0)\atop}\cr
\big(\L (\psi),Z(R)\big) \hskip-60pt&&&&\hskip-60pt \big(\L(\varphi),Z(Q)\big)\cr
&{\atop\loc_{\F^{^{\rm sc}}}({\rm id}_\psi,\delta^1_2)}\hskip-4pt\nwarrow\hskip-20pt&
&\hskip-20pt\nearrow\hskip-4pt{\atop\loc_{\F^{^{\rm sc}}}({\rm id}_\varphi,\delta^1_0)}\cr
&&\big(\L(\psi,\varphi),Z(Q)\big)\cr}
\eqno £6.2.1\phantom{.}$$
and therefore, denoting by $\bar\beta\,\colon (U_h,1)\to \big(\L (\psi),Z(R)\big)$ the corresponding 
$\frak i\widetilde\Loc\-$morphism, the equalities (cf.~diagram~£6.1.4)
$$\loc_{\F^{^{\rm sc}}}({\rm id}_Q,\delta^0_0)\circ \bar\alpha = \bar\sigma = 
\loc_{\F^{^{\rm sc}}}({\rm id}_T,\delta^0_1)\circ \bar\beta 
\eqno £6.2.2\phantom{.}$$
force the existence of an $\frak i\widetilde\Loc\-$morphism $\bar\varepsilon\,\colon  (U_h,1)
\to \big(\L (\psi,\varphi),Z(Q)\big)$ ful-filling
$$\loc_{\F^{^{\rm sc}}}({\rm id}_\varphi,\delta^1_0)\circ \bar\varepsilon = \bar\alpha
\qq \loc_{\F^{^{\rm sc}}}({\rm id}_\psi,\delta^1_2)\circ \bar\varepsilon = \bar\beta
\eqno £6.2.3\,,$$
 so that, setting $\bar\gamma = \loc_{\F^{^{\rm sc}}} ({\rm id}_{\varphi\circ\psi},\delta^1_1)\,,$ we finally get
 $$\eqalign{\loc_{\F^{^{\rm sc}}}({\rm id}_T,\delta^0_0)\circ (\bar\gamma\circ \bar\varepsilon) = 
 \loc_{\F^{^{\rm sc}}}({\rm id}_R,\delta^0_1)\circ \loc_{\F^{^{\rm sc}}}({\rm id}_\varphi,\delta^1_0)\circ \bar\varepsilon = \bar\rho\cr
 \loc_{\F^{^{\rm sc}}}({\rm id}_Q,\delta^0_1)\circ (\bar\gamma\circ \bar\varepsilon) = 
 \loc_{\F^{^{\rm sc}}}({\rm id}_R,\delta^0_0)\circ \loc_{\F^{^{\rm sc}}}({\rm id}_\psi,\delta^1_2)\circ \bar\varepsilon = \bar\tau\cr}
 \eqno £6.2.4.$$

\medskip
£6.3\phantom{.}  Note that we have a {\it faithfully forgetful\/} functor ${}^{^h}\!(\F^{^{\rm sc}})
\to \F^{^{\rm sc}}\,;$ for short, we denote by
$$\frak v_h = \frak v_{{}^{^h}\!(\F^{^{\rm sc}})} :
\ch^*\big({}^{^h}\!(\F^{^{\rm sc}})\big)\too {}^{^h}\!(\F^{^{\rm sc}})
\eqno £6.3.1\phantom{.}$$
the corresponding  evaluation functor.
Moreover, for any~$h'$ dividing~$h\,,$ it~is~clear that the inclusion $U_{h'}\i U_h$
induces a faithful functor
$$\frak r_{h',h} : {}^{^h}\!(\F^{^{\rm sc}})\too {}^{^{h'}}\!\!(\F^{^{\rm sc}})
\eqno £6.3.2.$$
On the other hand, since the category of chains $\ch (\F^{^{\rm sc}})$ is already 
a {\it semidirect product\/} [5,~£A2.8], the ${}^{\frak u_h}\frak l\ch^*(\F^{^{\rm sc}})\-$objects 
can be identified~with the triples  $(\bar\eta,\frak q,\Delta_n)$ formed by an $\F^{^{\rm sc}}\-$chain
 $\frak q\,\colon \Delta_n\to \F^{^{\rm sc}}$ and by an $\frak i\widetilde\Loc\-$morphism 
$$\bar\eta : (U_h,1)\too \big(\L (\frak q),{\rm Ker}
(\pi_\frak q)\big)
\eqno £6.3.3;$$
 for short, for any $i\in
\Delta_n\,,$ denote by  
$$\bar\iota_i^{\frak q} : \big(\L (\frak q),{\rm Ker} (\pi_\frak q)\big)\too 
\Big(\L\big(\frak q (i)\big),Z\big(\frak q(i)\big)\Big)
\eqno £6.3.4\phantom{.} $$
 the image by the functor $\loc_{\F^{^{\rm sc}}}$ of the 
$\ch^* (\F^{^{\rm sc}})\-$morphism from~$(\frak q,\Delta_n)$ to $(\frak q
(i),\Delta_0)$ determined by the identity map of $\frak q (i)\,.$

\bigskip
\noindent
{\bf Proposition~£6.4}\phantom{.} {\it For any $h\in \Bbb N -\{0\}\,,$ we have an
equivalence of categories
$$\frak j_h : {}^{\frak u_h}\frak l\ch^* (\F^{^{\rm sc}}) \cong 
\ch^*\big(\,{}^{^h}\!(\F^{^{\rm sc}})\big)
\eqno £6.4.1\phantom{.}$$
which maps any ${}^{\frak u_h}\frak l\ch^* (\F^{^{\rm sc}})\-$object
$(\bar\eta,\frak q,\Delta_n)$  on the chain $\frak q^{\bar\eta}\, \,\colon {\Delta}_n\to
\! {}^{^h}\!(\F^{^{\rm sc}})$ 
sending any  $i\in \Delta_n$ to the ${}^{^h}\!(\F^{^{\rm sc}})\-$object 
$\frak q(i)^{\bar\iota_i^{\frak q} \circ \bar\eta}$ and any $\Delta_n\-$morphism
$(j \bullet i$) to~$\frak q (j\bullet i)\,.$ Moreover, for any $h'$ dividing~$h\,,$ we have the commutative diagram
$$\matrix{{}^{\frak u_{h'}}\frak l\ch^* (\F^{^{\rm sc}}) &\buildrel \frak j_{h'}\over\cong 
& \ch^* \big(\,{}^{^{h'}}\!\!(\F^{^{\rm sc}}) \big)\cr 
\hskip-90pt{\scriptstyle (\rho_{h',h} *  \loc_{\F^{^{\rm sc}}}) \rtimes \id_{\ch^*
(\F^{^{\rm sc}})}}\big\uparrow &\phantom{\Big\uparrow}
&\big\uparrow {\scriptstyle \ch^*(\frak r_{h',h}) }\hskip-30pt\cr 
{}^{\frak u_h}\frak l\ch^* (\F^{^{\rm sc}})
&\buildrel \frak j_h\over\cong  & \ch^* \big(\,{}^{^h}\!(\F^{^{\rm sc}})\big)\cr}
\eqno£6.4.2.$$\/}

\par
\noindent
{\bf Proof:} Considering the $\F\-$chain $\Delta_1\cong \{j,i\}\to \F^{^{\rm sc}}$
obtained from the restriction of $\frak q\,,$ and the 
$\ch^* (\F^{^{\rm sc}})\-$morphism from~$(\frak q,\Delta_n)$ to the
$\ch^* (\F^{^{\rm sc}})\-$object determined by this $\F\-$chain, it is easily checked that $\frak q (j\!\bullet \!i)$ is indeed an ${}^{^h}\!(\F^{^{\rm sc}})\-$morphism from $\frak q(j)^{\bar\iota_j^{\frak q} \circ\bar\eta}$ to 
$\frak q(i)^{\bar\iota_i^{\frak q} \circ\bar\eta}\,.$ Moreover, a ${}^{\frak u_h}
\frak l\ch^* (\F^{^{\rm sc}})\-$morphism to $(\bar\eta,\frak q,\Delta_n)$ from a 
${}^{\frak u_h}\frak l\ch^* (\F^{^{\rm sc}})\-$object 
$(\bar\theta,\frak r,\Delta_m)$  is defined by a $\ch^*(\F^{^{\rm sc}})\-$mor-phism
$$(\mu,\delta) : (\frak r,\Delta_m)\too (\frak q,\Delta_n)
\eqno £6.4.3\phantom{.}$$
 fulfilling $\,\loc_{\F^{^{\rm sc}}}(\mu,\delta)\circ \bar\theta = \bar\eta\,$ [5,~condition~A2.6.2], and therefore $(\mu,\delta)$ is also a $\ch^*\big(\,{}^{^h}\!(\F^{^{\rm sc}})\big)\-$morphism from $(\frak r^{\bar\theta},\Delta_m)$ to  $(\frak q^{\bar\eta},\Delta_n)\,.$
Thus, we have obtained a functor
$$\frak j_h : {}^{\frak u_h}\frak l\ch^* (\F^{^{\rm
sc}})
\too \ch^* \big({}^{^h}\!(\F^{^{\rm sc}})\big)
\eqno £6.4.4.$$

\smallskip
On the other hand, any ${}^{^h}\!(\F^{^{\rm sc}})\-$chain $\hat\frak q\,
\,\colon {\Delta}_n\to {}^{^h}\!(\F^{^{\rm sc}})$ clearly determines a $\F^{^{\rm sc}}\-$chain
$\frak q\,\colon {\Delta}_n\to \F^{^{\rm sc}}\,;$ consequently, by the very
definition of ${}^{^h}\!(\F^{^{\rm sc}})\,,$ we have
$\hat\frak q(i) = \frak q (i)^{\bar\eta_i}$ where 
$$\bar\eta_i : (U_h,1)\too \big(\L\big(\frak q (i)\big),Z\big(\frak q(i)\big)\big)
\eqno £6.4.5$$
 is an $\frak i\widetilde\Loc\-$morphism for any $i\in \Delta_n$ and, since 
$\hat\frak q (j\!\bullet\! i)$ is a morphism from  $\frak q (j)^{\bar\eta_j}$ to $\frak q (i)^{\bar\eta_i}$ 
for any $0\le j\le i\le n\,,$ arguing by induction on $n$ it is not difficult to prove that
there is a unique $\frak i\widetilde\Loc\-$morphism 
$$\bar\eta : (U_h,1)\to \big(\L (\frak q), {\rm Ker}(\pi_\frak q)\big)
\eqno £6.4.6\phantom{.}$$
 fulfilling $\bar\iota_i^{\frak q} \circ \bar\eta = \bar\eta_i$ for any $i\in \Delta_n\,.$

\smallskip 
Similarly, if $\,\hat\frak r \,\colon {\Delta}_m\to {}^{^h}\!(\F^{^{\rm sc}})$ is 
a chain and, for any $j\in \Delta_m\,,$ we have 
$\hat\frak r(j) = \frak r (j)^{\bar\iota_j^{\frak r} \circ \bar\theta}$ for~a suitable $\frak i\widetilde\Loc\-$morphism 
$$\bar\theta : (U_h,1)\too \big(\L (\frak r),{\rm Ker}(\pi_\frak r)\big)
\eqno £6.4.7,$$
 then any $\ch^*\big({}^{^h}\!(\F^{^{\rm sc}})\big)\-$morphism
$(\hat\mu,\delta)\, \,\colon (\hat\frak r, \Delta_m)\to
(\hat\frak q, \Delta_n)$ induces a $\ch^*(\F^{^{\rm sc}})\-$ morphism 
$(\mu,\delta)\,\colon (\frak r,\Delta_m)\to (\frak q,\Delta_n)$
 such that 
$$\loc_{\F^{^{\rm sc}}}(\mu_i, {\rm id}_{\Delta_0}) \circ
\bar\iota_{\delta(i)}^{\frak r} \circ \bar\theta = \bar\eta_i = \bar\iota_i^{\frak q} \circ \bar\eta
\eqno£6.4.8\phantom{.}$$ 
for any $i\in \Delta_n\,,$ and therefore, since we have
$$\loc_{\F^{^{\rm sc}}}(\mu_i, {\rm id}_{\Delta_0}) \circ \bar\iota_{\delta
(i)}^{\frak r} = \bar\iota_i^{\frak q} \circ \loc_{\F^{^{\rm sc}}}
(\mu,\delta)
\eqno£6.4.9,$$
 we get $\,\loc_{\F^{^{\rm sc}}}(\mu,\delta)\circ \bar\theta = \bar\eta$ by the
uniqueness of $\bar\eta\,.$ Hence, the functor~$\frak j_h$ is an equivalence of
categories. The commutativity of diagram~£6.4.2 is easily checked. We are done.

\bigskip
\noindent
{\bf Proposition~£6.5}\phantom{.} {\it  For any $h\in \Bbb N -\{0\}\,,$ we have
a factorization of  $\hat\frak w\frak l_h$
$$\matrix{{}^{\frak u_h}\frak l\ch^*(\F^{^{\rm sc}})
&\buildrel \hat\frak w\frak l_h\over\too&{}^{U'_h}\aleph\cr
\hskip-6pt {\scriptstyle \frak j_h}\wr\hskip-3pt\Vert&
&\big\uparrow {\scriptstyle \frak t_h}\cr
\ch^*\big(\,{}^{^h}\!(\F^{^{\rm sc}})\big)& \buildrel \frak v_h \over \too
&{}^{^h}\!(\F^{^{\rm sc}})\cr}
\eqno £6.5.1\phantom{.}$$
throughout a functor $\frak t_h\,\colon {}^{^h}\!(\F^{^{\rm sc}})\to {}^{U'_h}\aleph$ mapping any 
${}^{^h}\!(\F^{^{\rm sc}})\-$object $Q^{\bar\rho}$ on the regular
$U'_h\-$set $(\varpi_{h,(\hat\L (Q),Z(Q))})^{-1}(\bar\rho)$ and any
${}^{^h}\!(\F^{^{\rm sc}})\-$morphism $\varphi\,\colon R^{\bar\sigma}\to Q^{\bar\rho}$
on a $U'_h\-$set bijection
$$(\varpi_{h,(\hat\L (R),Z(R))})^{-1}(\bar\sigma)\cong (\varpi_{h,(\hat\L (Q),Z(Q))})^{-1}(\bar\rho)
\eqno £6.5.2.$$
In particular, we have
$${}^\K \G_\K (\F,\widehat\aut_{\F^{^{\rm sc}}})_h\cong \Bbb H^0 \big(\,
{}^{^h}\!(\F^{^{\rm sc}}),{}^\K \F\!ct_{U'_h}\circ \frak t_h\big)
\eqno £6.5.3.$$
Moreover, for any $h'$ dividing $h$ we have a natural map 
$$\theta_{h',h} : \res^{U'_h}_{U'_{h'}}\circ\frak t_h\too \frak t_{h'}\circ\frak r_{h',h}
\eqno £6.5.4\phantom{.}$$ 
sending $Q^{\bar\rho}$ to the $U'_{h'}\-$set map
$$\res^{U'_h}_{U'_{h'}}\big((\varpi_{h,(\hat\L(Q),Z(Q))})^{-1}(\bar\rho)\big)\too
(\varpi_{h',(\hat\L(Q),Z(Q))})^{-1}\big({\rm Res}_{U_{h'}}^{U_h}(\bar\rho)\big)
\eqno £6.5.5.$$
 induced by the restriction throughout the inclusion $U_{h'}\i U_h\,.$\/}

\medskip
\noindent
{\bf Proof:} Let us denote by $R$ and $Q$ the obvious $\F^{^{\rm sc}}\-$chains  and by
$\varphi$ the $\F^{^{\rm sc}}\-$chain mapping $0$ on~$R\,,$ $1$ on $Q$ and the
$\Delta_1\-$morphism $(0\bullet\! 1)$ on $\varphi\,;$ thus, we have evident $\ch^*(\F^{^{\rm sc}})\-$morphisms
$$(R,\Delta_0)\longleftarrow (\varphi,\Delta_1)\too (Q,\Delta_0)
\eqno £6.5.6\phantom{.}$$
and then the natural map $\widehat\loc_{\F^{^{\rm sc}}}\to
\loc_{\F^{^{\rm sc}}}$ (cf.~£3.3) sends these $\ch^*(\F^{^{\rm sc}})\-$mor-phisms to a commutative diagram
$$\matrix{\big(\L (R),Z(R)\big)&\longleftarrow &\big(\L (\varphi),Z(Q)\big)&\too
&\big(\L (Q),Z(Q)\big)\cr
\uparrow &\phantom{\big\uparrow}&\uparrow&&\uparrow\cr
\big(\hat\L(R),Z(R)\big)&\longleftarrow &\big(\hat\L(\varphi),Z(Q)\big)&\too
&\big(\hat\L(Q),Z(Q)\big)\cr}
\eqno £6.5.7;$$
hence, the bottom $k^*\-\frak i\widetilde\Loc\-$morphisms induce a $U'_h\-$set bijection (cf.~diagram~£6.1.4)
$$\frak t_h (\varphi) : (\varpi_{h,(\hat\L (R),Z(R))})^{-1}(\bar\sigma)\cong
(\varpi_{h,(\hat\L (Q),Z(Q))})^{-1}(\bar\rho)
\eqno £6.5.8.$$
\eject

\smallskip
Now, we claim that the correspondence sending the  ${}^{^h}\!(\F^{^{\rm sc}})\-$morphism $\varphi$ 
to the $U'_h\-$set bijection  $\frak t_h (\varphi)$ defines a functor; indeed, for a third
${}^{^h}\!(\F^{^{\rm sc}})\-$object $T^{\bar\tau}$ and a ${}^{^h}\!(\F^{^{\rm sc}})\-$morphism 
$\psi\,\colon T^{\bar\tau}\to R^{\bar\sigma}\,,$ we have the following evident commutative 
$\ch^*(\F^{^{\rm sc}})\-$diagram 
$$\matrix{(Q,\Delta_0)\hskip-10pt&\hskip-10pt =\hskip-12pt &(Q,\Delta_0)\cr
&&\uparrow\cr
&&(\varphi,\Delta_1)&\to & (R,\Delta_0)\cr
\uparrow&&\uparrow&&\uparrow\cr
&&(\frak c,\Delta_2)&\to &(\psi,\Delta_1)&\to & (T,\Delta_0)\cr
&\hskip-10pt\swarrow \hskip-15pt&&&&&\Vert\cr
(\varphi\circ\psi,\Delta_1)\hskip-10pt&&&&\to&&(T,\Delta_0)\cr}
\eqno £6.5.9\phantom{.}$$
where the $\F^{^{\rm sc}}\-$chain $\frak c\,\colon \Delta_2\to \F^{^{\rm sc}}$ maps $0$ on $T\,,$
$1$ on $R\,,$ $2$ on $Q\,,$ $(0\bullet1)$ on~$\psi$ and $(1\!\bullet 2)$ on
$\varphi\,;$ once again, the functor $\loc_{\F^{^{\rm sc}}}$ maps this $\ch^*(\F^{^{\rm sc}})\-$diagram  on the
commutative $\frak i\widetilde\Loc\-$diagram
$$\matrix{\big(\L(Q),Z(Q)\big)&\leftarrow&\big(\L(\varphi),Z(Q)\big)
&\to & \big(\L(R),Z(R)\big)\cr
&&\uparrow&&\uparrow&\phantom{\big\uparrow}\cr
\uparrow&&\big(\L(\frak c),Z(Q)\big)&\to &\big(\L(\psi),Z(R)\big)\cr 
&\swarrow&\phantom{\big\uparrow}&&\downarrow\cr
\big(\L(\varphi\circ\psi),Z(Q)\big)&&\to&&\big(\L(T),Z(T)\big)\cr}
\eqno £6.5.10;$$
at this point, considering  the corresponding commutative diagrams £6.5.7, it is easily checked that 
 $$\frak t_h (\varphi\circ \psi) = \frak t_h (\varphi)\circ\frak t_h (\psi)
\eqno £6.5.11.$$

\smallskip
Moreover, it follows from Proposition~£5.8 and from definition~£5.13.6 above that the functor
$\hat\frak w\frak l_h$ maps the ${}^{\frak u_h} \frak l\ch^*(\F^{^{\rm sc}})\-$object $(\bar\eta,\frak q, \Delta_n)$ on the $U'_h\-$set
$(\varpi_{h,(\hat\L (\frak q), {\rm Ker}(\pi_\frak q))})^{-1}(\bar\eta)\,;$ but, the image by
$\widehat\loc_{\F^{^{\rm sc}}}$ of the $\ch^* (\F^{^{\rm sc}})\-$mor-phism
$(\frak q,\Delta_n) \to (\frak q (i),\Delta_0)$ determines a lifting of
$\bar\iota_i^{\frak q}$
$$\skew2\hat{\bar\iota}_i^{\,\frak q} : \hat\L (\frak q)\too \hat\L\big(\frak q (i)\big)
\eqno £6.5.12;$$
then, it is quite clear that
$$\skew2\hat{\bar\iota}_i^{\,\frak q}\circ (\varpi_{h,(\hat\L (\frak q), {\rm Ker}(\pi_\frak q))})^{-1}(\bar\eta) 
=  (\varpi_{h,(\hat\L (\frak q(i)),Z(\frak q(i)))})^{-1} (\bar\iota_i^{\,\frak q}\circ \bar\eta)
\eqno £6.5.13,$$
where the left member denotes the set of compositions of $\skew2\hat{\bar\iota}_i^{\,\frak q}$ 
with all the elements of the $U'_h\-$set $(\varpi_{h,(\hat\L(\frak q), {\rm Ker}(\pi_\frak q))})^{-1}(\bar\eta)\,.$
On the other hand, by the very definition of $\frak t_h\,,$ we actually have 
$$(\varpi_{h,(\hat\L (\frak q(0)),Z(\frak q (0)))})^{-1} (\bar\iota_0^{\frak q}\circ \bar\eta) 
= (\frak t_h\circ \frak v_h \circ \frak j_h)(\bar\eta,\frak q,\Delta_n)
\eqno £6.5.14.$$
From these equalities it is esaily checked that we have a natural isomorphism
$$\hat\frak w\frak l_h \cong\frak t_h \circ \frak v_h\circ
\frak j_h
\eqno £6.5.15.$$

\smallskip
Finally, the $\K\-$module isomorphism~£6.5.3 follows from the $\O\-$module
isomorphism~£5.13.7, from Proposition~£6.4 and from this natural isomorphism 
[5,~£A3.9]. The proof of the last statement is straightforward.

\medskip
\noindent
{\bf Remark~£6.6}\phantom{.} By Remark~£5.9, for any $\xi\in U_h$ the functor
$\frak t_h$ maps the ${}^{^h}\!(\F^{^{\rm sc}})\-$automorphism $\pi_Q\big(\rho (\xi)\big)$ of 
$Q^{\bar\rho}\,,$ where $\rho\in \bar\rho\,,$ on the identity map 
of~$(\varpi_{h,(\hat\L (Q),Z(Q))})^{-1}(\bar\rho)\,;$  in
particular, $\frak t_h$ factorizes {\it via\/} the {\it exterior quotient\/}
$\,\tilde{{}^{^{h}}}\!(\F^{^{\rm sc}})$ of $\,{}^{^h}\!(\F^{^{\rm sc}})$
determined by the correspondence mapping any $\,{}^{^h}\!(\F^{^{\rm sc}})\-$ob-ject
$Q^{\bar\rho}$ on the group of ${}^{^h}\!(\F^{^{\rm sc}})\-$automorphism of 
$Q^{\bar\rho}$ induced by $\rho (U_h)$ [5,~6.3].

\medskip
£6.7\phantom{.} Note that the category ${}^h (\F^{^{\rm sc}})$ also admits an {\it interior  structure\/} [5,~1.3]  which maps any ${}^h (\F^{^{\rm sc}})\-$object $Q^{\bar\rho}$ on the group $\F_{Q^{\rho (U'_h)}}(Q)$
 of ${}^h (\F^{^{\rm sc}})\-$ automorphisms of  $Q^{\bar\rho}$ where we are choosing  $\rho\in \bar\rho\,;$ then, denoting by $\widetilde{{}^h (\F^{^{\rm sc}}})$ the corresponding {\it exterior quotient\/}, it is easily checked 
 that the functor $\frak t_h$ admits a factorization
$$\tilde\frak t_h : \widetilde{{}^h (\F^{^{\rm sc}}})\too {}^{U'_h}\aleph
\eqno £6.7.1\phantom{.}$$
 and, in particular, we still have (cf.~£6.5.3)
$${}^\K \G_\K (\F,\widehat\aut_{\F^{^{\rm sc}}})_h\cong \Bbb H^0 \big(\,
\widetilde{{}^h (\F^{^{\rm sc}}}),{}^\K \F\!ct_{U'_h}\circ \tilde\frak t_h\big)
\eqno £6.7.2.$$

\medskip
£6.8\phantom{.} On the other hand, for any $h'\in \Bbb N -p\Bbb N \,,$ recall that in [5, 6.3 and~14.25] we have defined an analogous
category ${}^{h'} (\tilde\F^{^{\rm sc}})$ where  the ${}^{^{h'}}\!(\tilde\F^{^{\rm sc}})\-$ objects are the pairs $Q^{\rho}$ formed by an $\F\-$selfcentralizing subgroup~$Q$ of~$P$ and by an {\it injective\/} group homomorphism $\rho\,\colon U_{h'}\to  \tilde\F (Q)$ and, for a second 
${}^{^{h'}}\!(\tilde\F^{^{\rm sc}})\-$object $R^{\sigma}\,,$ the ${}^{^{h'}}\!(\tilde\F^{^{\rm sc}})\-$morphisms
from $R^{\sigma}$ to $Q^{\rho}$ are the~$\tilde\F^{^{\rm sc}}\-$mor-phisms $\tilde\varphi\,\colon R\to Q$ fulfilling
$$\rho (\xi)\circ \tilde\varphi = \tilde\varphi\circ \rho (\xi)
\eqno £6.8.1\phantom{.}$$
for any $\xi\in U_{h'}\,,$ and that in [5,.~Proposition~14.28]  we have obtained a factorization analogous to £6.5.1, {\it via\/} a functor $\frak s_{h'}\,\colon {}^{h'} (\tilde\F^{^{\rm sc}})\to {}^{U_{h'}}\aleph$ which maps any ${}^{h'} (\tilde\F^{^{\rm sc}})\-$object $Q^{\rho}$ on the $U_{h'}\-$set
$(\varpi_{h',\skew4\hat{\tilde\F}(Q)})^{-1}(\rho)$ [5, 14.18.4] and any ${}^{h'} (\tilde\F^{^{\rm sc}})\-$morphism 
$R^{\sigma}\to Q^{\rho}$ on a suitable $U_{h'}\-$set bijection
$$(\varpi_{h',\skew4\hat{\tilde\F}(R)})^{-1}(\sigma)\cong 
(\varpi_{h',\skew4\hat{\tilde\F}(Q)})^{-1}(\rho)
\eqno £6.8.2.$$

\medskip
£6.9\phantom{.} Moreover, respectively  denoting by $h^{p'}$ and $h^p$ the $p'\-$ and the $p\-$part of any 
$h\in \Bbb N-\{0\}\,,$ let us consider the functor
$$\frak x_{h} : {}^{h^{p'}}\! (\tilde\F^{^{\rm sc}})\too \aleph
\eqno £6.9.1\phantom{.}$$
which maps any ${}^{^{h^{p'}}}\!(\tilde\F^{^{\rm sc}})\-$ object $Q^{\rho}$ on the set $\widetilde{\rm Mon}(U_{h^p},Q^{\widehat{\rho (U_{h^{p'}})}})$ of  $Q^{\widehat{\rho (U_{h^{p'}})}}\-$ conjugacy classes of injective group homomorphisms  $U_{h^p}\to Q^{\widehat{\rho (U_{h^{p'}})}}\,,$ where we choose a lifting $\widehat{\rho (U_{h^{p'}})}$ of $\rho (U_{h^{p'}})$ to $\F (Q)$ and denote by $ Q^{\widehat{\rho (U_{h^{p'}})}}$ the subgroup of 
$\widehat{\rho (U_{h^{p'}})}\-$fixed elements of $Q\,,$ and sends any  ${}^{^{h'}}\!
(\tilde\F^{^{\rm sc}})\-$morphism from $R^{\sigma}$ to $Q^{\rho}\,,$ determined by 
a~$\tilde\F^{^{\rm sc}}\-$morphism $\tilde\varphi\,\colon R\to Q\,,$ to the map
$$\widetilde{\rm Mon}(U_{h^p},R^{\widehat{\sigma (U_{h^{p'}})}})\too \widetilde{\rm Mon}
(U_{h^p},Q^{\widehat{\rho (U_{h^{p'}})}})
\eqno £6.9.2\phantom{.}$$
determined by  the group homomorphism $R^{\widehat{\sigma (U_{h^{p'}}})}\to  Q^{\widehat{\rho (U_{h^{p'}}})}$ induced by a suitable representative of~$\tilde\varphi\,.$
Finally, the map ${}^{U_{h^{p'}}}\aleph\times \aleph\to {}^{U_{h^{p'}}}\aleph$ induced by the direct product determines a new functor
$$\frak s_{h^{p'}}\times \frak x_h  : 
{}^{h^{p'}}\! (\tilde\F^{^{\rm sc}})\too {}^{U_{h^{p'}}}\aleph
\eqno £6.9.3.$$

\bigskip
\noindent
{\bf Proposition~£6.10}\phantom{.} {\it  For any $h\in \Bbb N -\{0\}$ we have a functor
$$\frak d_\F^h : \widetilde{{}^h {(\F^{^{\rm sc}}})}\too   {}^{h^{p'}}\!  (\tilde\F^{^{\rm sc}})
\eqno £6.10.1\phantom{.}$$
mapping any $ {}^h {(\F^{^{\rm sc}}})\-$object $Q^{\bar\rho}$ such that  $Q$ is fully normalized in $\F$ and that, choosing $\rho\in \bar\rho\,,$ $\rho (U_{h^p})$ is contained in $N_P (Q)\,,$
on the ${}^{h^{p'}}\! (\tilde\F^{^{\rm sc}})\-$object formed by $\hat Q = Q\.\rho (U_{h^p})$
endowed  with the $U_{h^{p'}}\-$action $\hat\rho'\,\colon U_{h^{p'}}\to 
\tilde\F(\hat Q)$ induced by~$\rho\,,$ and mapping the class of any 
$ {}^h {(\F^{^{\rm sc}}})\-$morphism $\varphi\,\colon R^{\bar\sigma}  \to Q^{\bar\rho}$ on the unique
${}^{h^{p'}}\!  (\tilde\F^{^{\rm sc}})\-$morphism $\hat R^{\hat\sigma'} 
\to \hat Q^{\hat\rho'}$ extending $\tilde\varphi\,\colon R\to Q\,.$ Moreover, we have
 $$(\frak d_\F^h)_*(\F\!ct_{U_{h^{p'}}}\circ \tilde\frak t_h) = \F\!ct_{U_{h^{p'}}}
 \circ  (\frak s_{h^{p'}} \times \frak x_h)
 \eqno £6.10.2\phantom{.}$$
 and, in particular, for any $n\in \Bbb N$ we get a group isomorphism
 $$\Bbb H^n \big(\, \widetilde{{}^h {(\F^{^{\rm sc}}})},\F\!ct_{U_{h^{p'}}}\circ \tilde\frak t_h\big)\cong \Bbb H^n \big({}^{h^{p'}}\!  (\tilde\F^{^{\rm sc}}), \F\!ct_{U_{h^{p'}}} \circ  (\frak s_{h^{p'}} \times \frak x_h)\big)
\eqno £6.10.3.$$\/}

\par
\noindent
{\bf Proof:} First of all, consider the analogous functor
$$\frak y_{h} : \widetilde{{}^{^{h^{p'}}}\!(\F^{^{\rm sc}}})\too \aleph
\eqno £6.10.4$$
mapping any  ${}^{^{h^{p'}}}\!(\F^{^{\rm sc}})\-$object $Q^{\bar\rho'}$ on the set 
$\widetilde{\rm Mon}(U_{h^p},Q^{\rho' (U_{h^{p'}})})$ of 
$Q^{\rho' (U_{h^{p'}})}\-$ conjugacy classes of injective group homomorphisms 
$U_{h^p}\to Q^{\rho' (U_{h^{p'}})}\,,$ where\break
\eject
\noindent
 $\rho'\in \bar\rho'$ and $Q^{\rho' (U_{h^{p'}})}$ denotes the 
subgroup of fixed points of $\rho' (U_{h^{p'}})\i \L (Q)$ in $Q\i \L (Q)\,,$  and mapping the class of any 
${}^{^{h^{p'}}}\!(\F^{^{\rm sc}})\-$morphism $\varphi\,\colon R^{\bar\sigma'} \to Q^{\bar\rho'}$ to the map
$$\widetilde{\rm Mon}(U_{h^p},R^{\sigma' (U_{h^{p'}})})\too \widetilde{\rm Mon}(U_{h^p},Q^{\rho' (U_{h^{p'}})})
\eqno £6.10.5\phantom{.}$$
determined by  the group homomorphism $R^{\sigma' (U_{h^{p'}})}\to  Q^{\rho' (U_{h^{p'}})}$ induced by 
$\varphi\,,$ where $\sigma'\in \bar\sigma'\,.$ As above, the functor $\frak y_h$ can be viewed as a {\it representation\/} of the category $ \widetilde{\,{}^{^{h^{p'}}}\!(\F^{^{\rm sc}}})$ [5,~A2.2] and we consider the corresponding 
{\it semidirect product\/} $\frak y_h\!\rtimes \widetilde{{}^{^{h^{p'}}}\!(\F^{^{\rm sc}}})$ [5,~A2.7];
then, we define an {\it adjoint pair\/} of functors
$$\frak f_h : \frak y_h\!\rtimes \widetilde{{}^{^{h^{p'}}}\!(\F^{^{\rm sc}}}) \too 
\widetilde{{}^{^h}\!(\F^{^{\rm sc}}}) \qq \frak g_h : \widetilde{{}^{^h}\!(\F^{^{\rm sc}}}) \too \frak y_h\!\rtimes \widetilde{{}^{^{h^{p'}}}\!(\F^{^{\rm sc}}})
\eqno £6.10.6\phantom{.}$$
as follows.

\smallskip
 An $\frak y_h\!\rtimes \widetilde{{}^{^{h^{p'}}}\!(\F^{^{\rm sc}}})\-$object is a pair formed by a  
 ${}^{^{h^{p'}}}\!(\F^{^{\rm sc}})\-$object $Q^{\bar\rho'}$ and by a $Q^{\rho' (U_{h^{p'}})}\-$ conjugacy class
 $\tilde\rho''$ of injective group homomorphisms
$$\rho'' : U_{h^{p}}\too Q^{\rho' (U_{h^{p'}})}
\eqno £6.10.7;$$
it is quite clear that $\rho'$ and $\rho''$ define an injective group homomorphism $\rho\,\colon U_h\to \L (Q)$
and therefore they determine a $\widetilde\Loc\-$morphism
$$\bar\rho : (U_h,1)\too \big(\L (Q),Z(Q)\big)
\eqno £6.10.8;$$  
hence, we get a ${}^{^{h}}\! (\F^{^{\rm sc}})\-$object $Q^{\bar\rho}$ and  then, we  define $\frak f_h 
(\tilde\rho'',Q^{\bar\rho'}) = Q^{\bar\rho}$ for a choice of $\rho''\in \tilde\rho''\,.$ Similarly, a  
$\frak y_h\!\rtimes \widetilde{{}^{^{h^{p'}}}\!(\F^{^{\rm sc}}})\-$morphism to 
$ (\tilde\rho'',Q^{\bar\rho'})$ from another $\frak y_h\!\rtimes \widetilde{{}^{^{h^{p'}}}\!(\F^{^{\rm sc}}})\-$object 
$(\tilde\sigma'',R^{\bar\sigma'})$ is the class of an ${}^{^{h^{p'}}}\!\! (\F^{^{\rm sc}})\-$morphism 
$\varphi\,\colon R^{\bar\sigma'}\to Q^{\bar\rho'}$ such that, denoting by $\bar\varphi$ the restriction of 
$\varphi$ to $R^{\sigma'(U_{h^{p'}})}$ for a choice $\sigma'\in \bar\sigma'\,,$ we have $\skew4\tilde{\bar\varphi}\circ \tilde\sigma''  = \tilde\rho''$ and we may assume that $\bar\varphi\circ \sigma'' = \rho''\,;$ it is then clear 
that~$\varphi$ still determines an 
$\widetilde{{}^{^h}\!  (\F^{^{\rm sc}}})\-$morphism from $R^{\bar\sigma} = 
\frak f_h (\tilde\sigma'',R^{\bar\sigma'})$ to $Q^{\bar\rho}$ and, coherently, we define $\frak f_h$ mapping the class
of $\varphi$ in $\widetilde{{}^{^{h^{p'}}}\!(\F^{^{\rm sc}}})$ on the class
of $\varphi$ in $\widetilde{{}^{^{h}}\!(\F^{^{\rm sc}}})\,.$ The functoriality of this correspondence is clear.

 \smallskip
 Conversely,  note that any $\widetilde{{}^h (\F^{^{\rm sc}}})\-$object admits 
an ${}^h (\F^{^{\rm sc}})\-$isomorphic one $Q^{\bar\rho}$ such that $Q$ is fully
normalized in $\F$ and,  identifying $N_P (Q)$ with its structural image in $\L (Q)\,,$
$\rho (U_{h^p})$ is contained in $N_P (Q)$ for $\rho\in \bar\rho\,;$ thus, in order to define the functor $\frak g_h\,,$ we may replace  $\widetilde{{}^h (\F^{^{\rm sc}}})$ by its full subcategory 
 over those objects.

 \smallskip
If $Q^{\bar\rho}$ is such an ${}^{^{h}}\!(\F^{^{\rm sc}})\-$object,
 the product  $\hat Q = Q\.\rho (U_{h^p})$ is a subgroup of $N_P(Q)$ and  it is not difficult to check that the {\it localizer\/} $\L (\iota_Q^{\hat Q})$ of the $\F\-$chain  determined by the inclusion map from $Q$ to~$\hat Q$ can be identified with the normalizer $N_{\L (Q)}(\hat Q)$ which clearly contains~$\rho (U_{h^{p'}})\,;$ thus, choosing a representative $\mu$ 
 of the $\widetilde\Loc\-$morphism
$$\loc_{\F^{^{\rm sc}}}({\rm id}_{\hat Q},\delta^0_0) : \big(\L (\iota_Q^{\hat Q}),Z(\hat Q)\big)
\too \big(\L (\hat Q),Z(\hat Q)\big)
\eqno £6.10.9,$$
the composition of $\mu$ with the restriction $\rho'$ of $\rho$ to $U_{h^{p'}}$ determines a 
$\widetilde\Loc\-$mor-phism
$$\bar\rho' : ( U_{h^{p'}},1)\too \big(\L (\hat Q), Z(\hat Q)\big)
\eqno £6.10.10\phantom{.}$$
and therefore we get an   ${}^{^{h^{p'}}}\!\! (\F^{^{\rm sc}})\-$object 
 $\hat Q^{\bar\rho'}\,;$  now, since $\rho (U_{h^p})$ is contained in~$\hat Q\,,$ the restriction $\rho''$  of $\rho$ to $U_{h^{p}}$ determines an injective group homomorphism
$\rho''\,\colon U_{h^{p}}\to \hat Q^{\,\rho(U_{h^{p'}})}\,;$
then, we define 
$$\frak g_h (Q^{\bar\rho}) = (\tilde\rho'', \hat Q^{\bar\rho'})
\eqno £6.10.11\phantom{.}$$
where $\tilde\rho''$ denotes the $\hat Q^{\,\rho(U_{h^{p'}})}\-$conjugacy class of $\rho''\,.$

\smallskip
If $R^{\bar\sigma}$ is also such an  ${}^{^{h}}\! (\F^{^{\rm sc}})\-$object and
 $\varphi\,\colon R^{\bar\sigma}\to Q^{\bar\rho}$ is an ${}^{^{h}}\! (\F^{^{\rm sc}})\-$mor-phism, 
we know that there is a $\widetilde\Loc\-$morphism (cf.~£6.1)
$$\bar\alpha :  (U_h,1)\too  \big(\L (\varphi),Z(Q)\big)
\eqno £6.10.12\phantom{.}$$
fulfilling (cf.~£6.1.4) 
$$\loc_{\F^{^{\rm sc}}}({\rm id}_R,\delta^0_1)\circ \bar\alpha = \bar\sigma
\qq \loc_{\F^{^{\rm sc}}}({\rm id}_Q,\delta^0_0)\circ \bar\alpha = \bar\rho
\eqno £6.10.13;$$
then, choosing  representatives $\lambda_0$ of $\loc_{\F^{^{\rm sc}}}({\rm id}_Q,\delta^0_0)$ and 
$\lambda_1$ of $\loc_{\F^{^{\rm sc}}}({\rm id}_R,\delta^0_1)$ 
such that, identifying $R$ with its structural image in $\L (\varphi)\,,$ we have
$\lambda_0(v) = v$ and $\lambda_1(v) =\varphi(v)$ for any $v\in R\,,$
 there are representatives $\rho\in \bar\rho\,,$ $\sigma\in \bar\sigma$ and 
$\alpha\in \bar\alpha$ fulfilling
$$\lambda_0\circ \alpha = \sigma
\qq \lambda_1\circ \alpha = \rho
\eqno £6.10.14\phantom{.}$$
and, in particular, inducing group isomorphisms
$$\sigma(U_{h^p})\cong \alpha(U_{h^p})\cong \rho(U_{h^p})
\eqno £6.10.15;$$
since clearly $\varphi^{-1}\big(Z(Q)\big)\i Z(R)\,,$ it is easily checked that $\varphi$
and these isomorphisms determine an injective group homomorphism
$$\hat\varphi : \hat R =  R\. \sigma (U_{h^p}) \too Q\. \rho (U_{h^p}) = \hat Q
\eqno £6.10.16\phantom{.}$$
which agree with $\lambda_0$ and $\lambda_1\,;$ then, it follows from
[5,~Proposition~18.16] that $\hat\varphi$ is actually an $\F\-$morphism.

\smallskip
Now, $\L(\varphi)$ contains $\hat R$ and it is clear that
$$\sigma (U_{h^{p'}})\i  \lambda_1\big(N_{\L (\varphi)}(\hat R)\big) \i N_{\L (R)}(\hat R)
\cong \L (\iota_R^{\hat R})
\eqno £6.10.17;$$
moreover, choosing a representative $\nu$ of the $\widetilde\Loc\-$morphism
$$\loc_{\F^{^{\rm sc}}}({\rm id}_{\hat R},\delta^0_0) : \big(\L (\iota_R^{\hat R}),Z(\hat R)\big)
\too \big(\L (\hat R),Z(\hat R)\big)
\eqno £6.10.18,$$ 
it is easily checked that $(\nu\circ \lambda_1)\big(N_{\L (\varphi)}(\hat R)\big)$ is contained in the image 
of $\L(\hat\varphi)$ in $\L (\hat R)$ {\it via\/} any representative  $\hat\lambda_1$ of the  
$\widetilde\Loc\-$morphism
$$\loc_{\F^{^{\rm sc}}}({\rm id}_{\hat R},\delta^0_1) : \big(\L (\hat\varphi),Z (\hat Q)\big)
\too \big(\L (\hat R),Z (\hat R)\big)
\eqno £6.10.19\phantom{.}$$
and therefore the restriction of $\alpha$ to $U_{h^{p'}}$ and the composition $\nu\circ \lambda_1$
determine a $\widetilde\Loc\-$morphism 
$$\bar\alpha' :  (U_{h^{p'}},1)\too  \big(\L (\hat\varphi),Z(\hat Q)\big)
\eqno £6.10.20\phantom{.}$$ 
fulfilling the corresponding equalities
$$\loc_{\F^{^{\rm sc}}}({\rm id}_{\hat R},\delta^0_1)\circ \bar\alpha' = \bar\sigma'
\qq \loc_{\F^{^{\rm sc}}}({\rm id}_{\hat Q},\delta^0_0)\circ \bar\alpha'= \bar\rho'
\eqno £6.10.21,$$
so that $\hat\varphi$ is also an  ${}^{^{h^{p'}}}\!\!(\F^{^{\rm sc}})\-$morphism
from $\hat R^{\bar\sigma'}$ to $\hat Q^{\bar\rho'}\,;$ it is easy to check that the  
$\widetilde{{}^{^{h^{p'}}}\!(\F^{^{\rm sc}}})\-$morphism determined by $\hat\varphi$ does not depend on our choice.

\smallskip
Finally, respectively denoting by 
$$\rho'' : U_{h^{p}}\too \hat Q^{\,\rho(U_{h^{p'}})}\qq \sigma'' : U_{h^{p}}\too \hat R^{\,\sigma(U_{h^{p'}})}
\eqno £6.10.22\phantom{.}$$
 the restrictions to $U_{h^p}$ of $\rho$  and $\sigma\,,$ by the   very definition of $\hat\varphi$ 
 we obtain $\hat\varphi\big(\sigma'' (\xi)\big) = \rho''(\xi)$ 
 for any~$\xi\in U_{h^{p'}}\,;$  in conclusion, denoting by $\tilde\rho''$ and $\tilde\sigma''$ the respective
 $\hat Q^{\,\rho(U_{h^{p'}})}\-$ and $\hat R^{\,\sigma(U_{h^{p'}})}\-$conjugacy classes of $\rho''$ and 
 $\sigma''\,,$ the class of $\hat\varphi$ still defines a 
 $\frak y_h\!\rtimes \widetilde{{}^{^{h^{p'}}}\!(\F^{^{\rm sc}}})\-$morphism from 
 $(\tilde\sigma'', \hat R^{\bar\sigma'})$ to $(\tilde\rho'', \hat Q^{\bar\rho'})$
and we define $\frak g_h$ mapping the class of $\varphi$ in $ \widetilde{{}^{^{h}}\!(\F^{^{\rm sc}}})$ on 
the class of $\hat\varphi$ in $ \widetilde{{}^{^{h^{p'}}}\!(\F^{^{\rm sc}}})\,.$ Once again, the proof of the
functoriality of $\frak g_h$ is straightforward.

\smallskip
 At this point, it is quite clear that 
 $$\frak g_h\circ \frak f_h\cong \id_{\frak y_h\!\rtimes \widetilde{{}^{^{h^{p'}}}\!(\F^{^{\rm sc}}})}
 \eqno £6.10.23;$$
  moreover, it is not difficult to verify that the correspondence sending the  $\widetilde{{}^{^h}\!
(\F^{^{\rm sc}}})\-$object $Q^{\bar\rho}$ above to the $\widetilde{{}^{^h}\!
(\F^{^{\rm sc}}})\-$morphism from  $Q^{\bar\rho}$ to $\hat Q^{\,\bar\rho'}$ induced by
the inclusion $Q\i \hat Q$ defines a {\it natural map\/} 
$$\eta : \id_{\widetilde{{}^{^h}\!(\F^{^{\rm sc}}})} \too \frak f_h\circ\frak g_h
\eqno £6.10.24\phantom{.}$$
which fulfills $\eta *\frak f_h = {\rm id}_{\frak f_h}$ and $\frak g_h * \eta = 
{\rm id}_{\frak g_h}\,;$ this proves that $\frak f_h$ and $\frak g_h$ form an {\it adjoint pair\/} and,
in particular, we get [5, 1.6]
$$(\frak g_h)_*(\F\!ct_{U_{h^{p'}}} \circ \tilde\frak t_h) = \F\!ct_{U_{h^{p'}}} \circ 
\tilde\frak t_h\circ \frak f_h
\eqno £6.10.25.$$

\smallskip
On the other hand, for any $h'\in \Bbb N  - p\Bbb N$ note that the categories
$\widetilde{{}^{^{h'}}\!(\F^{^{\rm sc}}})$ and ${}^{^{h'}}\!(\tilde\F^{^{\rm sc}})$
are equivalent. Indeed, for any ${}^{^{h'}}\!(\F^{^{\rm sc}})\-$object 
$Q^{\bar\rho'}\,,$ it is quite clear that the $Q\-$conjugacy class of the group homomorphism 
$\rho'\,\colon U_{h'}\to \L (Q)\,,$ where we choose $\rho'\in \bar\rho'\,,$ is determined by
the composition 
$$\tilde\pi_Q\circ \rho' :  U_{h'}\too \tilde\F (Q)\cong \L (Q)/Q
\eqno £6.10.26,$$
namely by the ${}^{^{h'}}\!(\tilde\F^{^{\rm sc}})\-$object $Q^{\tilde\pi_Q\circ \rho'}\,.$ 
Similarly, any $\widetilde{{}^{^{h'}}\!(\F^{^{\rm sc}}})\-$morphism from an ${}^{^{h'}}\!\!(\F^{^{\rm sc}})\-$object 
$R^{\bar\sigma'}$  to $Q^{\bar\rho'}$ is the  $Q^{\rho'(U_{h'})}\-$conjugacy class of an $\F\-$morphism 
$\varphi\,\colon R\to Q$ which admits a group homomorphism $\alpha'\,\colon U_{h'}\to \L (\varphi)$ fulfilling
$$\lambda_0\circ \alpha' = \sigma' \qq \lambda_1\circ \alpha' = \rho'
\eqno £6.10.27\phantom{.}$$
for suitable   representatives $\lambda_0$ of $\loc_{\F^{^{\rm sc}}}({\rm id}_Q,\delta^0_0)$ and $\lambda_1$ 
of $\loc_{\F^{^{\rm sc}}}({\rm id}_R,\delta^0_1)\,;$ hence, according to [5, Proposition~A2.10], for any $\xi'\in U_{h'}$ we have 
$$\eqalign{\pi_Q\big(\rho' (\xi')\big)\circ \varphi 
&= (\pi_Q\circ \lambda_1)\big(\alpha' (\xi')\big)\circ \varphi\cr
& = \big(\aut_{\F^{^{\rm sc}}}  ({\rm id}_R,\delta^0_1)\circ \pi_\varphi\big)
\big(\alpha' (\xi')\big)\circ \varphi\cr
& = \big(\aut_{\F^{^{\rm sc}}} ({\rm id}_R,\delta^0_1)\big) 
\big( (\pi_\varphi\circ \alpha') (\xi')\big)\circ \varphi\cr
& = \varphi\circ \big(\aut_{\F^{^{\rm sc}}} ({\rm id}_R,\delta^0_0)\big) 
\big( (\pi_\varphi\circ \alpha') (\xi')\big)\cr
& = \varphi\circ \big(\aut_{\F^{^{\rm sc}}} ({\rm id}_R,\delta^0_0)\circ 
\pi_\varphi\big) \big(\alpha' (\xi')\big)\cr
& = \varphi\circ (\pi_R\circ\lambda_0) \big(\alpha' (\xi')\big) = 
\varphi\circ \pi_R\big(\sigma' (\xi')\big)\cr}
\eqno £6.10.28\phantom{.}$$
which proves that $\varphi$ also induces an ${}^{^{h'}}\!
(\tilde\F^{^{\rm sc}})\-$morphism from $R^{\tilde\pi_R\circ \sigma'}$ to $Q^{\tilde\pi_Q\circ \rho'}\,.$

 \smallskip
  Conversely, any representative $\psi$ of an ${}^{^{h'}}\! (\tilde\F^{^{\rm sc}})\-$morphism 
$$\tilde\psi :  R^{\tilde\pi_R\circ \sigma'}\too Q^{\tilde\pi_Q\circ \rho'}
\eqno £6.10.29\phantom{.}$$
determines a group homomorphism 
$$ U_{h'}\too \tilde\F(\psi)\cong \L (\psi)/R
\eqno £6.10.30\phantom{.}$$
 which clearly can by lifted to a group homomorphism $\beta\,\colon U_{h'}\to \L (\psi)\,,$ and it is easily checked 
 that this group homomorphism fulfills equalities 6.10.27 for a suitable choice 
 of the representatives $\lambda_0$ and~$\lambda_1\,;$ hence, the $\F\-$morphism $\psi\,\colon R\to Q$ 
 also determines an $\widetilde{{}^{^{h'}}\!(\F^{^{\rm sc}}})\-$morphism from $R^{\bar\sigma'}$ to 
 $Q^{\bar\rho'}\,,$ completing the proof of the equivalence
 $$\widetilde{{}^{^{h'}}\!(\F^{^{\rm sc}}})\cong {}^{^{h'}}\!(\tilde\F^{^{\rm sc}})
  \eqno £6.10.31.$$
  \eject
\noindent
Moreover, the functors $\frak y_{h}$ and $\frak x_{h}$ agree with this equivalence; consequently, for any $h\in \Bbb N -\{0\}\,,$ we also get an equivalence of categories
 $$\frak e_h : \frak y_h\!\rtimes \widetilde{{}^{^{h^{p'}}}\!(\F^{^{\rm sc}}})\cong \frak x_{h} \!\rtimes
 {}^{^{h^{p'}}}\!\! (\tilde\F^{^{\rm sc}})
 \eqno £6.10.32.$$

 \smallskip
 Finally, denoting by 
 $$\frak p^{h}_\F :  \frak x_{h}\!\rtimes  {}^{^{h^{p'}}}\! \!(\tilde\F^{^{\rm sc}})\too 
  {}^{^{h^{p'}}}\! (\tilde\F^{^{\rm sc}})
  \eqno £6.10.33\phantom{.}$$
the structural functor, we set
$$\frak d_\F^h = \frak p^{h}_\F\circ \frak e_h \circ \frak g_h
\eqno £6.10.34\phantom{.}$$
and we claim that the {\it direct image\/} of $\F\!ct_{U_{h^{p'}}}\circ \tilde\frak t_h$
throughout this composition coincides with $\F\!ct_{U_{h^{p'}}}\circ (\frak s_{h^{p'}}\times \frak x_{h})\,;$
 indeed, according to~equality~£6.10.25,  for any $ {}^{^{h^{p'}}}\!\! (\tilde\F^{^{\rm sc}})\-$object $Q^{\rho'}\,,$  we have
 $$\eqalign{\big((\frak p^{h}_\F\circ \frak e_h)_*(\F\!ct_{U_{h^{p'}}}&\circ 
 \tilde\frak t_h \circ \frak f_h)\big)(Q^{\rho'})\cr
&= \prod_{\tilde\rho''\in \widetilde{\rm Mon} (U_{h^{p}},Q^{\widehat{\rho'(U_{h^{p'}})}})}
\F\!ct_{U_{h^{p'}}}\big(\frak s_{h^{p'}} (Q^{\rho'})\big)\cr
&= \F\!ct_{U_{h^{p'}}}\big(\bigsqcup_{\tilde\rho''\in \widetilde{\rm Mon} 
(U_{h^{p}},Q^{\widehat{\rho'(U_{h^{p'}}))}})}\frak s_{h^{p'}} (Q^{\rho'})\big)\cr
&= \F\!ct_{U_{h^{p'}}}\big( (\frak s_{h^{p'}}\times \frak x_{h}) (Q^{\rho'})\Big)\cr}
\eqno £6.10.35.$$
Then,  isomorphism~£6.10.3 follows from this fact and from Lemma~£6.11 below.

\bigskip
\noindent
{\bf Lemma~£6.11} \phantom{.} {\it Let $\frak C$ be a small category, $\frak m\,\colon \frak C\to \CC$ a representation such that, for any $\frak C\-$object $C\,,$
the category ${}^\frak m C$ only has the identity morphisms, and $\frak a\,\colon
\frak m\rtimes \frak C\to \Ab$ a contravariant functor. Denote by 
$\frak a^\frak m\,\colon \frak C\to \Ab$ the contravariant functor mapping any $\frak C\-$object $C$ 
on~$\prod_X \frak a (X,C)\,,$ where $X$ runs over the set of 
${}^\frak m C\-$objects, and any $\frak C\-$morphism $f\,\colon C\to C'$ on~the
group homomorphism
$$\frak a^\frak m (C') = \prod_{X'} \frak a (X',C')\too \prod_X \frak a (X,C) = 
\frak a^\frak m (C)
\eqno £6.11.1\phantom{.}$$
sending $\big(a_{X'}\big)_{X'}\in\frak a^\frak m (C')$ to $\Big(\big(\frak a ({\rm id}_{{}^\frak m\! f (X)},f)\big) 
(a_{\,{}^\frak m\! f (X)})\Big)_X \in \frak a^\frak m (C)\,.$ Then, for any $n\in \Bbb N$ the structural functor 
$\frak m\rtimes\frak C\to \frak C$ induces a group isomorphism
$$\Bbb H^n(\frak m\rtimes\frak C,\frak a)\cong \Bbb H^n (\frak C,\frak a^\frak m)
\eqno £6.11.2.$$\/}
\eject

\par
\noindent
{\bf Proof:} By its very definition in [5, A3.8], we know that $\Bbb H^n(\frak m\rtimes\frak C,\frak a)$ is the $n\-$th
homology group of the functor 
$$(\frak a\circ \frak v_{\frak m\rtimes\frak C}^\frak o)^{\Fct(\st,\frak m\rtimes
\frak C)^\frak o} : \bf\Delta\too \Ab
\eqno £6.11.3\phantom{.}$$
where $\bf\Delta$ denotes the {\it simplicial $2\-$category\/} [5, A1.7], $\st\,\colon
{\bf\Delta}\to \CC$ the {\it standard representation\/} of $\bf\Delta$ [5, A2.2],
$\Fct(\st,\frak m\rtimes\frak C)^\frak o\,\colon {\bf\Delta^{\!\circ}}\to \CC$ the 
{\it naive $\frak m\rtimes\frak C\-$dual representation} of~$\st$
mapping $\Delta_n$ on $\Fct (\Delta_n,\frak m\rtimes \frak C)^\frak o$ [5, A2.5], and 
$$\frak v_{\frak m\rtimes\frak C }^\frak o : \Fct(\st,\frak m\rtimes\frak C)^\frak o\rtimes {\bf\Delta}^{\!\circ} = \ch^\frak o(\frak m\rtimes\frak C)\too 
\frak m\rtimes\frak C
\eqno £6.11.4\phantom{.}$$
the corresponding {\it evaluation functor\/} [5,~A3.7].

\smallskip
But, any functor $\hat\frak q\,\colon \Delta_n \to\frak m\rtimes \frak C$
determines a functor $\frak q\,\colon \Delta_n\to \frak C$ and then $\hat\frak q$
is determined~by $\frak q$ and by $\hat\frak q (0) = \big(X_0,\frak q (0)\big)$ since,
according to our hypothesis on $\frak m$ and setting $X_i = {}^\frak m 
\frak q(0\bullet i)(X_0)\,,$ for any $i\in \Delta_n -\{0\}$ we have 
$$\hat\frak q (i) = \big(X_i,\frak q(i)\big)\qq
\hat\frak q (i\!-\!1\bullet i) = \big({\rm id}_{X_i},\frak q  (i\!-\!1\bullet i)\big)
\eqno £6.11.5.$$
Consequently, it is quite clear that
$$\eqalign{\big((\frak a\circ \frak v_{\frak m\rtimes\frak C}^\frak o)^{\Fct(\st,\frak m\rtimes\frak C)^\frak o}\big)(\Delta_n) &= \prod_{\hat\frak q\in \Fct (\Delta_n,\frak m\rtimes\frak C)}\frak a\big(\hat\frak q(0)\big)\cr
& = \prod_{\frak q\in \Fct (\Delta_n,\frak C)}
\frak a^\frak m\big(\frak q(0)\big)\cr
&= \big((\frak a^\frak m\circ \frak v_\frak C^\frak o)^{\Fct (\st,\frak C)^\frak o}
\big)(\Delta_n)\cr}
\eqno £6.11.6\phantom{.}$$
and it is easily checked the coincidence of the functors $(\frak a\circ \frak v_{\frak m\rtimes\frak C}^\frak o)^{\Fct(\st,\frak m\rtimes\frak C)^\frak o}$ and 
$(\frak a^\frak m\circ \frak v_\frak C^\frak o)^{\Fct (\st,\frak C)^\frak o}\,.$ We are done.

\bigskip
\bigskip
\noindent
{\bf £7\phantom{.} \bf A vanishing cohomological result }
\bigskip

£7.1\phantom{.} As in [5, Corollary~14.32] for  the determination of the $\O\-$rank of the modular
Grothendieck group $\G_k (\F,\widehat\aut_{})\,,$ the determination of the $\O\-$rank of the ordinary
Grothendieck group $\G_\K (\F,\widehat\aut_{})$ ultimately depends on a  vanishing cohomological result.
However, the general result [5, Theorem~6.26] we employ there it is not powerful enough, as it stands, to discuss
 our present situation; but, as a matter of fact, essentially the same arguments prove a sufficiently general result.
Nevertheless, even if our proof below mainly repeats the proof of [5, Theorem~6.26], we write it completely in 
order to clarify some arguments.

\medskip
£7.2\phantom{.} In this section, our setting is just a finite $p\-$group $P$ and a Frobenius $P\-$category $\F\,.$
  Let $K$ be a finite $p'\-$group and,  as in [5,~6.3], consider the category ${}^K\!\ad(\tilde\F^{^{\rm sc}})$ of
  the $K\-$objects of $\ad (\tilde\F^{^{\rm sc}})$ [5,~6.2]; recall that this category admits {\it direct products\/}
   and {\it pull-backs\/} [5,~Propositions~6.14 and~6.21]. As in [5,~6.25], consider the object $\bigoplus_{x\in K} P$
of the category $\ad (\tilde\F^{^{\rm sc}})$ [5,~6.2] endowed with the $K\-$action $\pi$
defined by the {\it regular\/} action of~$K$ on itself and by the identity on $P$
between the corresponding terms, so that  $\big(\bigoplus_{x\in K} P\big)^\pi$ is an
indecomposable $K\-$object of~$\ad (\tilde\F^{^{\rm sc}})$  [5,~6.3].

\medskip
£7.3\phantom{.}  If $\frak F$ is a full subcategory  of~${}^K\!\ad(\tilde\F^{^{\rm sc}})$ and 
$\frak m\,\colon \frak F\to \O\-\mod$ is a {\it contravariant\/} functor, we say that a functor
$\frak m^\circ\,\colon \frak F\to \O\-\mod$ is a {\it right-hand sectional functor\/}
of $\frak m$ if it coincides with $\frak m$ over the $\frak F\-$objects and, for any 
$\frak F\-$morphism $\tilde\varphi\,\colon R^\sigma\to Q^\rho\,,$
we have
$$\frak m (\tilde\varphi)\circ \frak m^\circ (\tilde\varphi) =  {\rm id}_{\frak m (R^\sigma)}
\eqno £7.3.1.$$
Note that $\ad (\frak F)$ still can be identified to a subcategory of ${}^K\!\ad(\tilde\F^{^{\rm sc}})\,.$
Moreover, recall that the center $Z(K)$ defines an  {\it exterior quotient\/} 
$\widetilde{{}^K}\!\ad(\tilde\F^{^{\rm sc}})$ of ${}^K\!\ad(\tilde\F^{^{\rm sc}})$ [5,~6.3]; then,
since~$\vert Z(K)\vert$ is invertible in~$\O\,,$ if $\frak m$  factorizes {\it via\/} the image
$\tilde\frak F$ of~$\frak F$ in $\widetilde{{}^K}\!\ad(\tilde\F^{^{\rm sc}})$ throughout a 
{\it
contravariant\/} functor $\tilde\frak m\,\colon 
\tilde\frak F\to \O\-\mod\,,$ it follows from [5, Proposition~A4.13] that, for any~$n\in
\Bbb N\,,$ we have
$$\Bbb H^n (\tilde\frak F,\tilde\frak m) = \Bbb H^n (\frak F,\frak m)
\eqno £7.3.2\phantom{.}$$
since, considering the subcategory $\frak Z$ of $\frak F$ formed by the same objects and by the automorphisms of the $\frak F\-$objects induced by $Z(K)\,,$ $\Bbb H^n (\tilde\frak F,\tilde\frak m)$ clearly coincides with the {\it $\frak Z\-$stable $n\-$cohomology group\/} of~$\frak F$ over $\frak m$ [5,~A3.18].

\bigskip
\noindent
{\bf Theorem £7.4}\phantom{.} {\it  With the notation above, let  $\frak F$ be a full subcategory of~${}^K\!\ad(\tilde\F^{^{\rm sc}})$ over indecomposable $K\-$objects of $\ad(\tilde\F^{^{\rm sc}})$ including $\big(\bigoplus_{x\in K} P\big)^\pi\,,$ such that the subcategory~$\ad (\frak F)$ is closed by direct products and pull-backs. For any
 contravariant functor  $\frak m\,\colon \frak F\to \O\-\mod$ admitting a right-hand
sectional functor $\frak m^\circ\,\colon \frak F\to \O\-\mod\,,$ we have $\,\Bbb H^n (\frak F,\frak m) = \{0\}$ 
for any $n\ge 1\,.$\/}

\medskip
\noindent
{\bf Proof:} First of all, we prove the statement assuming that $p\.\frak m = 0$ or,
equi-valently, that $\frak m$ is a {\it contravariant\/} functor from $\frak F$ to
$k\-\mod\,;$ coherently, $\frak m^\circ$ is a functor from $\frak F$ to
$k\-\mod\,.$ Moreover, it follows from [5, Proposition~A4.11] that for any
$n\ge 1$ we have
$$\Bbb H^n (\frak F,\frak m)\cong \Bbb H^n \big(\ad
(\frak F), \ad (\frak m)\big)
\eqno £7.4.1,$$
so that it suffices to prove that $\Bbb H^n \big(\ad (\frak F), \ad (\frak m)\big) =
\{0\}$ for any $n\ge 1\,.$

\smallskip
Set  $S = \bigoplus_{x\in K} P\,;$ it follows from [5, Proposition~6.14] that the {\it direct product  by\/}  $S^{\pi}\,,$
 or the {\it exterior intersection with\/} $S^{\pi}$ [5,~definition~6.13.3], defines a functor 
$$\int_{S^{\pi}} : \ad  (\frak F)\too  \ad (\frak F)
\eqno £7.4.2;$$
then, the existence of the structural $\ad  (\frak F)\-$morphism
$Q^{\rho}\,\widetilde\cap\,S^{\pi}\to S^{\pi}$ for any $\ad  (\frak F)\-$object
$Q^{\rho}$ shows that $\int_{S^{\pi}}$ factorizes throughout the evident
{\it forgetful\/} functor [5,~1.7] 
$$\ft_{S^{\pi}}\, \,\colon \ad (\frak F)_{S^{\pi}}\too \ad  
(\frak F)
\eqno £7.4.3;$$
explicitly, it suffices to consider the functor $\ad (\frak F)
\to \ad  (\frak F)_{S^{\pi}}$ mapping any $\ad  (\frak F)\-$object
$Q^{\rho}$ on the structural $\ad  (\frak F)\-$morphism above and any $\ad  (\frak F)\-$mor-phism 
$\tilde\alpha\,\colon R^{\sigma}\to Q^{\rho}$ on $\tilde\alpha\,\tilde\cap\,\widetilde{\rm id}_{S^{\pi}}$ [5,~Proposition~6.14]. But, since the category $\ad  (\frak F)_{S^{\pi}}$ has the final object 
$\widetilde{\rm id}_{S^{\pi}} \,\colon {S^{\pi}}\to {S^{\pi}}\,,$ it follows from [5, Corollary~A4.8] that 
for any $n\ge 1$ we have
$$ \Bbb H^n \big(\ad  (\frak F)_{S^{\pi}}, \ad (\frak m) \circ
\ft_{S^{\pi}}\big) = \{0\}
\eqno £7.4.4\phantom{.}$$
and therefore, we still have [5,~A3.10.4]
$$\Bbb H^n \big(\ad  (\frak F),\ad (\frak m)\circ \int_{S^{\pi}}\big) = \{0\}
\eqno £7.4.5.$$

\smallskip
 Moreover, the existence of the structural morphism $\tilde\omega_{Q^{\rho}} \,\colon
Q^{\rho}\,\widetilde\cap \,S^{\pi}\to Q^{\rho}$ for any $\ad  (\frak F)\-$object
$Q^{\rho}$ shows the existence of a natural map
$$\omega : \int_{S^{\pi}}\too \id_{\ad (\frak F)}
\eqno £7.4.6\phantom{.}$$
 sending $Q^{\rho}$ to $\tilde\omega_{Q^{\rho}}\,;$
thus, in order to prove that $\Bbb H^n \big(\ad (\frak F),\ad (\frak m)\big) =
\{0\}\,,$ it suffices to prove that the natural map $\ad (\frak m) *\omega$ admits a
{\it natural section\/}
$$\theta : \ad (\frak m)\circ \int_{S^{\pi}}\too \ad (\frak m)
\eqno £7.4.7,$$
so that $\ad (\frak m)$ becomes a direct summand of $\ad (\frak m)\circ
\int_{S^{\pi}}\,.$

\smallskip 
Explicitly, for any $\frak F\-$object $Q^{\rho} = (\bigoplus_{i\in I}
Q_i)^\rho\,,$ we have [5,~6.13]
$$Q^{\rho}\,\widetilde\cap\,S^{\pi} = \Big(\bigoplus_{(i,x)\in I\times K}\,
\bigoplus_{(\tilde\tau,T,\tilde\iota_{T}^P) \in \check\frak T_{Q_i,P}}
T_i\Big)^{\hat\rho}
\eqno £7.4.8\phantom{.}$$ 
for a set of representatives $\check\frak T_{Q_i,P}$ of $\tilde\frak T_{Q_i,P}$ in
$\frak T_{Q_i,P}$ [5,~6.9] and a suitable action $\hat\rho$ of $K$ on the 
$\ad (\tilde\F^{^{\rm sc}})\-$object $Q\,\widetilde\cap\,S\,;$ in particular, $K$ acts
freely on the disjoint union
$$\hat I = \bigsqcup_{(i,x)\in I\times K}  \check\frak T_{Q_i,P}
\eqno £7.4.9\phantom{.}$$
and  let us denote by $\hat I/K$ the set of $K\-$orbits on $\hat I$ and, for any
$O\in \hat I/K\,,$ by~$(T_{_O})^{\hat\rho_{_O}}$ the corresponding
indecomposable ``direct summand'' of $Q^{\rho}\,\widetilde\cap\,S^{\pi}$  and by $\tilde
\tau_{_O}$ the composition
$$\tilde\tau_{_O} : (T_{_O})^{\hat\rho_{_O}}\too Q^{\rho}\,\widetilde\cap\,S^{\pi} \too
Q^\rho
\eqno £7.4.10\phantom{.}$$
of the structural $\ad (\frak F)\-$morphism with $\tilde \omega_{Q^\rho}\,.$ Moreover,
we denote by $\hat I^\circ/K$ the set of {\it special\/} orbits $O\in \hat I/K$ where
the $\tilde\F^{^{\rm sc}}\-$morphisms determining $\tilde\tau_{_O}$ are isomorphisms;
note that, according to [5, Proposition~6.14], we have a canonical bijection
$$\hat I^\circ/K\cong \bigsqcup_{i\in I}\tilde\F (P,Q_i)
\eqno £7.4.11.$$

\smallskip
Then, we consider the homomorphism
$$\matrix{\theta_{Q^\rho} :\hskip-5pt &\big(\ad (\frak m)\big)
(Q^{\rho}\,\widetilde\cap\,S^{\pi})&\hskip-5pt\too&
\frak m (Q^\rho)\cr
&\Vert\cr
&\prod_{O\in \hat I/K} \frak m \big((T_{_O})^{\hat\rho_{_O}}\big)\cr}
\eqno £7.4.12\phantom{.}$$
sending an element $m = (m_{_O})_{O\in \hat I/K}$ of this product to
$$\theta_{Q^\rho} (m) = \vert\hat I^\circ/K\vert^{-1}\.\sum_{O\in \hat I^\circ/K} 
\frak m^\circ(\tilde \tau_{_O})(m_{_O})
\eqno £7.4.13;$$
since for any {\it special\/} orbit $O$ the composition $\tilde \tau_{_O}$ above is an isomorphism and
since  $\frak m(\tilde\tau_{_O})\circ \frak m^\circ (\tilde\tau_{_O}) 
= {\rm id}_{(T_{_O})^{\hat\rho_{_O}}}\,,$ we actually have $\frak m^\circ (\tilde\tau_{_O}) =\frak m (\tilde\tau_{_O})^{-1}$ and therefore we clearly get 
$$\theta_{Q^\rho}\circ (\ad(\frak m) * \omega)_{Q^\rho} = 
{\rm id}_{\frak m (Q^\rho)}
\eqno £7.4.14.$$
 By the {\it distributivity\/} of the {\it exterior intersection\/} [5,~6.13], we easily can extend this 
 correspondence to all the $\ad(\frak F)\-$objects and then we claim that the extended correspondence 
 is a natural map from $\ad (\frak m)\circ \int_{S^{\pi}}$ to
$\ad (\frak m)\,;$ actually, it suffices to consider an $\frak F\-$morphism
$\tilde\alpha\,\colon R^{\sigma}\to Q^{\rho}$ and to prove the commutativity of the
following diagram
$$\matrix{\big(\ad (\frak m)\big)(Q^{\rho}\,\widetilde\cap\,S^{\pi})&\buildrel
\theta_{Q^\rho}\over \too& \frak m (Q^\rho)\cr
\hskip-70pt{\scriptstyle 
(\ad (\frak m))(\tilde\alpha\,\widetilde\cap\,{\rm id}_{S^{\pi}})} \downarrow&
\phantom{\Big\downarrow}&\downarrow{\scriptstyle
\frak m(\tilde \alpha)}\hskip-10pt\cr
\big(\ad (\frak m)\big) (R^{\sigma}\,\widetilde\cap\,S^{\pi})&\buildrel
\theta_{R^\sigma}\over \too&
\frak m (R^\sigma)\cr}
\eqno £7.4.15.$$

\smallskip
Explicitly, if $R = \bigoplus_{j\in J} R_j$ is the structural decomposition of $R\,,$ then $\tilde\alpha$ is
given by a $K\-$compatible map $f\,\colon J\to I$ and by a $K\-$compatible family of
$\tilde\F\-$morphisms $\tilde\alpha_j\,\colon R_j\to Q_{f(j)}$ where $j$ runs 
over $J\,,$ and as above we have
$$R^{\sigma}\,\widetilde\cap\,S^{\pi} = \Big(\bigoplus_{(j,x)\in J\times K}\,
\bigoplus_{(\tilde\upsilon,U,\tilde\iota_{U}^P) \in \check\frak T_{R_j,P}} U
\Big)^{\hat\sigma}
\eqno £7.4.16\phantom{.}$$ 
for a set of representatives $\check\frak T_{R_j,P}$ of $\tilde\frak T_{R_j,P}$ in
$\frak T_{R_j,P}$ [5,~6.9] and a suitable action $\hat\sigma$ of $K$ on the 
$\ad (\tilde\F^{^{\rm sc}})\-$object $R\,\widetilde\cap\,S\,;$ again, we set
$$\hat J = \bigsqcup_{(j,x)\in J\times K}  \check\frak T_{R_j,P}
\eqno £7.4.17\phantom{.}$$
\eject
\noindent
and denote by $\hat J/K$ the set of $K\-$orbits on $\hat J\,,$ by $\hat J^\circ/K$ the
set of {\it special\/} $K\-$orbits on $\hat J$ and, for any
$O\in \hat J/K\,,$ by $(U_{_O})^{\hat\sigma_{_O}}$ the
corresponding indecomposable ``direct summand'' of $R^{\sigma}\,\widetilde\cap\,S^{\pi}$
 and by $\tilde \upsilon_{_O}\,\colon (U_{_O})^{\hat\sigma_{_O}}\to
R^{\sigma}$ the ana-logous composition~£7.4.10;
moreover, it is clear that the map $f$ and the family $\{\tilde\alpha_j\}_{j\in J}$
determine a
$K\-$compatible map $\hat f
\,\colon \hat J\to \hat I$ and, for any $O\in \hat J/K\,,$ an $\frak F\-$morphism
$$\tilde\alpha_{_O} : (U_{_O})^{\hat\sigma_{_O}}\too
(T_{_{\hat f(O)}})^{\hat\rho_{_{\hat f(O)}}}
\eqno £7.4.18.$$

\smallskip
It is easily checked from [5, Propositions~6.14 and~6.21] that  [5,~6.18.2]
$$R^{\sigma}\,{}_{\tilde\alpha}\widetilde\cap_{\tilde\omega_{Q^{\rho}}}\,
(Q^{\rho}\,\widetilde\cap\, S^{\pi})
\cong R^{\sigma}\,\widetilde\cap\, S^{\pi}
\eqno £7.4.19\phantom{.}$$
and, by the {\it distributivity\/} property,  we may assume that the exterior intersection
$R^{\sigma}\,\widetilde\cap\, S^{\pi}$ coincides with
$$\Big(\bigoplus_{(j,x)\in J\times K}\, \bigoplus_{(\tilde\tau,T,
\tilde\iota_{T}^P) \in \check\frak T_{Q_{f(j)},P}}  R_j\,{}_{\tilde\alpha_j}
\widetilde\cap_{\tilde\tau}\, T\Big)^{\hat\sigma}
\eqno £7.4.20.$$
Then, for any $(j,x)\in J\times K$ and any $\frak t = (\tilde\tau,T,
\tilde\iota_{T}^P) \in \check\frak T_{Q_{f(j)},P}\,,$ we choose representatives $\tau\in \tilde\tau$
and $\alpha_{f(j)}\in \tilde\alpha_{f(j)}\,,$ and a set of representatives $W_{(j,x,\frak t)}$ in~$Q_{f(j)}$ for the set of double classes $\tau(T)\big\backslash Q_{f(j)}\big/ \alpha_j(R_j)\,;$  we denote
by $W^{^{\rm sc}}_{(j,x,\frak t)}$ the set of $w\in W_{(j,x,\frak t)}$
such that the subgroup
$$U_w = \big(\kappa_{_{Q_{f(j)}}}(w)\circ\alpha_j\big)^{-1}\big(\tau (T)\big)
\eqno £7.4.21\phantom{.}$$
remains $\F\-$selfcentralizing, and by  $\tilde\beta_{(j,x,\frak t,w)}
\,\colon U_w\to T$  the $\tilde\F^{^{\rm sc}}\-$morphism determined by the compositions
$\kappa_{_{Q_{f(j)}}}(w)\circ \alpha_j\,,$ where $\kappa_{_{Q_{f(j)}}}(w)$ denotes the corresponding 
conjugation by $w\,.$

\smallskip
Now, with all this notation, it follows from [5,~Proposition~6.19] that, for any $(j,x)\in J\times K$ and any triple 
$\frak t = (\tilde\tau_\frak t,T_\frak t, \tilde\iota_{T_\frak t}^P) \in \check\frak T_{Q_{f(j)},P}\,,$ 
we have
$$R_j\,{}_{\tilde\alpha_j}\widetilde\cap_{\tilde\tau}\, T_\frak t
= \bigoplus_{w\in W^{^{\rm sc}}_{(j,x,\frak t)}} U_w
\eqno £7.4.22\phantom{.}$$
and isomorphism~7.4.19 determines a ``graded'' bijection between the disjoint union
$$\bigsqcup_{(j,x)\in J\times K}\,\bigsqcup_{\frak t\in \check\frak T_{Q_{f(j)},P}}
 W^{^{\rm sc}}_{(j,x,\frak t)}
 \eqno £7.4.23\phantom{.}$$
 and $\hat J\,;$ moreover, the $\ad(\tilde\F^{^{\rm sc}})\-$morphism
$$\tilde\alpha\,\widetilde\cap\,{\rm id}_{S^{\pi}} :
R^{\sigma}\,\widetilde\cap\,S^{\pi}\too Q^{\rho}\,\widetilde\cap\,S^{\pi}
\eqno £7.4.24\phantom{.}$$
is the ``direct sum'' over the set $\bigsqcup_{(j,x)\in J\times K}
\check\frak T_{Q_{f(j)},P}$ of the $\ad(\tilde\F^{^{\rm sc}})\-$morphisms
$$\tilde\beta_{(j,x,\frak t)} : \bigoplus_{w\in W^{^{\rm sc}}_{(j,x,\frak t)}}  U_w\too T_\frak t
\eqno £7.4.25\phantom{.}$$
defined by the $\tilde\F^{^{\rm sc}}\-$morphisms $\tilde\beta_{(j,x,\frak t,w)}
\,\colon U_w\to T_\frak t$ above.

\smallskip
Furthermore, $K$ acts on all this situation, and let us denote by $O_{(j,x,\frak t,w)}$ the $K\-$orbit --- which is
 actually {\it regular\/} --- of the element of $\hat J$ determined by
 $(j,x,\frak t,w)\,,$ by $O_{(f(j),x,\frak t)}$ the image in $\hat I$ {\it via\/} $\hat f$ of 
$O_{(j,x,\frak t,w)}$ --- which actually does not depend on $w$ ---
and by 
$$\eqalign{\tilde\tau_{_{O_{(f(j),x,\frak t)}}} : (T_{_{O_{(f(j),x,\frak t)}}})^{\hat\rho_{_{O_{(f(j),x,\frak t)}}}}&\too Q^\rho\cr
\tilde\beta_{_{O_{(j,x,\frak t_,w)}}} : 
(U_{_{O_{(j,x,\frak t,w)}}})^{\hat\sigma_{_{O_{(j,x,\frak t,w)}}}} &\too (T_{_{O_{(f(j),x,\frak t)}}})^{\hat\rho_{_{O_{(f(j),x,\frak t)}}}}\cr 
\tilde\upsilon_{_{O_{(j,x,\frak t,w)}}} : 
(U_{_{O_{(j,x,\frak t,w)}}})^{\hat\sigma_{_{O_{(j,x,\frak t_,w)}}}}  &\too R^\sigma\cr}
\eqno £7.4.26\phantom{.}$$
 the $\frak F\-$morphisms respectively determined by the $K\-$orbits of $\tilde\tau_\frak t\,,$ $\tilde\beta_{(j,x,\frak t,w)}$ and~$\tilde\iota_{U_w}^{R_j}\,;$ then, the naturality of
$\ad (\frak m) *\omega$ forces
$$\eqalign{\frak m(\tilde\upsilon_{_{O_{(j,x,\frak t,w)}}})&\circ \frak m (\tilde\alpha) = \frak m (\tilde\beta_{_{O_{(j,x,\frak t,w)}}})\circ \frak m
(\tilde\tau_{_{O_{(f(j),x,\frak t)}}} )\cr}
\eqno £7.4.27,$$
so that we still have
$$\eqalign{\frak m(\tilde\upsilon_{_{O_{(j,x,\frak t,w)}}})\circ\frak m (\tilde \alpha)&\circ \frak m^\circ (\tilde\tau_{_{O_{(f(j),x,
\frak t)}}}) =  \frak m (\tilde\beta_{_{O_{(j,x,\frak t),w)}}})\cr}
\eqno £7.4.28\,.$$

\smallskip
Now, we are ready to prove the commutativity of the diagram~£7.4.15;
according to our definition, the composition $\frak m(\tilde\alpha) \circ
\theta_{Q^\rho}$ sends the element $m = (m_{_O})_{O\in \hat I/K}$
where $m_{_O}\in \frak m \big( (T_{_O})^{\hat\rho_{_O}}\big)\,,$ to the sum
$$\vert \hat I^\circ/K\vert^{-1}\.\sum_{O\in \hat I^\circ/K} 
(\frak m (\tilde\alpha))\big(\frak m^\circ (\tilde \tau_{_O}) (m_{_O})\big)
\eqno £7.4.29;$$
on the other hand, we have
$$\eqalign{\big((\ad (\frak m))&(\tilde\alpha\,\widetilde\cap\,{\rm id}_{S^{\pi}})
\big)(m)\cr
& = \sum_{j\in J}\,\sum_{\frak t\in \check\frak T_{Q_{f(j)},P}}\, 
\sum_{w\in W^{^{\rm sc}}_{(j,1,\frak t)}} \big(\frak m
(\tilde\beta_{_{O_{(j,1,\frak t,w)}}})\big) (m_{_{O_{(f(j),1,\frak t)}}})\cr} 
\eqno £7.4.30\phantom{.}$$
and therefore, denoting by $ W^{R_j}_{(j,1,\frak t)}$ the set of 
$w\in  W^{^{\rm sc}}_{(j,1,\frak t)}$ such that $U_w = R_j\,,$ so that then we have
$\frak m^\circ (\tilde\upsilon_{_{O_{(j,x,\frak t,w)}}}) = \frak m(\tilde\upsilon_{_{O_{(j,x,\frak t,w)}}})^{-1}\,,$
 it follows from our definition of
$\theta_{R^\sigma}$ and from equality~£7.4.28 that
$$\eqalign{\big(\theta_{R^\sigma} &\circ (\ad (\frak m)) (\tilde\alpha\, \widetilde\cap\,
{\rm id}_{S^{\pi}})\big)(m)\cr  
&= \sum_{j\in J}\,\sum_{\frak t\in \check\frak T_{Q_{f(j)},P}} {\vert W^{R_j}_{(j,1,\frak t)} \vert\over 
\vert \hat J^\circ/K\vert}\.  \frak m (\tilde\alpha)\big(\frak m^\circ (\tilde\tau_{_{O_{(f(j),1,\frak t)}}})  (m_{_{O_{(f(j),1,\frak t)}}})\big)\cr}
\eqno £7.4.31.$$

\smallskip
But, note that if $O_{(f(j),1,\frak t)}$ belongs to $\hat I^\circ/K$ then we have 
$\vert W^{R_j}_{(j,1,\frak t)}\vert = 1\,;$ moreover, it follows from [5,~Corollary~£4.9] that
$\tilde\alpha_j$ induces an injective map from $\tilde\F(P,Q_{f(j)})$ to
$\tilde\F(P,R_j)$ and it is clear that $\vert f^{-1}(i)\vert
= \vert J\vert/\vert I\vert$ for any $i\in I\,;$ furthermore, according to [5,~6.7.2] and to
bijection~£7.4.11 above, we have $$\vert \hat I^\circ/K\vert \equiv \vert I\vert\vert\tilde\F (P)\vert \qq
\vert \hat J^\circ/K\vert \equiv \vert J\vert\vert\tilde\F (P)\vert \pmod p
\eqno £7.4.32;$$
hence, the sum of all the corresponding terms in the second member of equality~£7.4.31 coincides with the sum~£7.4.29 above.

\smallskip
Consequently, in order to show the commutativity of diagram~£7.4.15, it suffices
to prove that, for any $j\in J$ and any $\frak t\in \check\frak T_{Q_{f(j)},P}$
such that $\tilde\tau_\frak t\,\colon T_\frak t\to Q_{f(j)}$ is {\it not\/} an
isomorphism, $p$ divides~$\vert W^{R_j}_{(j,1,\frak t)}\vert\,;$ but, it is clear that $W^{R_j}_{(j,1,\frak t)}$ is a set of representatives for the quotient set
$$\tau_\frak t (T_\frak t)
\big\backslash \T_{Q_{f(j)}} \big(\tau_\frak t (T_\frak t),\alpha_j (R_j)\big)
\eqno £7.4.33\phantom{.}$$
 and that the nontrivial $p\-$group $\bar N_{Q_{f(j)}} \big(\tau_\frak t (T_\frak t)\big)$
acts freely on this set. This completes the proof of the naturality of
$\theta$ and therefore the proof of the theorem in the case where $p\.\frak m =0\,.$
\smallskip

In the general case, there is a subfunctor $\frak m^{\rm tor}\,\colon \frak F\to
\O\-\mod$ mapping any $\frak F\-$object $Q^\rho$  on the {\it torsion\/}
$\O\-$submodule of $\frak m (Q^\rho)$ and then we have the quotient functor 
$$\frak m/\frak m^{\rm tor} : \frak F\too \O\-\mod
\eqno £7.4.34\phantom{.}$$
which maps any object on a {\it free\/} $\O\-$module; consequently, since we have an
exact sequence [5,~A3.11.4]
$$\Bbb H^n (\frak F,\frak m^{\rm tor})\too \Bbb H^n (\frak F,\frak m)\too \Bbb
H^n (\frak F,\frak m/\frak m^{\rm tor})
\eqno £7.4.35,$$
for any $n\in \Bbb N \,,$ we already may assume that either $\frak m = \frak m^{\rm tor}$ or $\frak m^{\rm tor}
= 0\,.$
\eject

\smallskip
In the first case we have $p^\ell\.\frak m = 0$ for some $\ell\in \Bbb N -\{0\}$ and, 
considering the exact sequence [5,~A3.11.4]
$$\Bbb H^n (\frak F,p\.\frak m)\too \Bbb H^n (\frak F,\frak m)\too \Bbb
H^n (\frak F,\frak m/p\.\frak m)= \{0\}
\eqno £7.4.36,$$
for any $n\in \Bbb N\,,$ it suffices to argue by induction on $\ell\,.$
In the second case,  if $c_0$ is an {\it $n\-$cocycle\/} for $n\ge 1\,,$ we already have proved that
$$c_0 \equiv d_{n-\!1}(a_0)\pmod p
\eqno £7.4.37\phantom{.}$$ 
for a suitable {\it $(n\!-\!1)\-$cochain\/} $a_0\,,$ so that we have $c_0
-d_{n-\!1}(a_0) = p\.c_1$ for a suitable {\it $n\-$cocycle\/}~$c_1$ since we are dealing
with free $\O\-$modules; thus, inductively, we can define {\it $n\-$cocycles\/} $c_i$ and
{\it $(n\!-\!1)\-$cochains\/}
$a_i$ fulfilling 
$$c_i\equiv d_{n-\!1}(a_i)\pmod p\qq c_i - d_{n-\!1} (a_i) = p\. c_{i+\!1}
\eqno £7.4.38\phantom{.}$$
and then, according to the completeness of $\O\,,$ it is quite clear that
$$c_0 = d_{n-\!1} \big(\sum_{i\in \Bbb N -\{0\}} p^i\.a_i\big)
\eqno £7.4.39.$$
We are done.

\bigskip
\bigskip
\noindent
{\bf £8\phantom{.} \bf The $\O\-$rank of the Grothendieck groups of a folded Frobenius $P\-$category }
\bigskip

£8.1\phantom{.} As in \S3, let $P$ be a finite $p\-$group, $\F$ a Frobenius $P\-$category and
$$\widehat\aut_{\F^{^{\rm sc}}} : \ch^*(\F^{^{\rm sc}})\too k^*\-\Gr
\eqno £8.1.1\phantom{.}$$
a functor lifting $\aut_{\F^{^{\rm sc}}}$ (cf.~£2.8); we are ready to determine the $\O\-$rank 
of~$\G_\K(\F,\widehat\aut_{\F^{^{\rm sc}}})$. As in [5,~Corollary~14.32]  for the determination of the 
$\O\-$rank of~$\G_k(\F,\widehat\aut_{\F^{^{\rm sc}}})\,,$ our argument is an easy consequence of 
the charac-ter decomposition of the functor $\frak g_\K$ obtained in \S5, and
of the following vanishing cohomological result which, setting  $h' = h^{p'}$ for any $h\in \Bbb N - \{0\}\,,$
involves  the quotient category  
$\widetilde{{}^h (\F^{^{\rm sc}}})$
and the functor
$$\tilde\frak t_h : \widetilde{{}^h (\F^{^{\rm sc}}})\too {}^{U_{h'}}\aleph
\eqno £8.1.1\phantom{.}$$
introduced in~£6.7 from the factorization in Proposition~£6.5.

\bigskip
\noindent
{\bf Theorem~£8.2}\phantom{.} {\it For any $h\in \Bbb N -\{0\}$ and any $n\ge 1$ we have 
$$\Bbb H^n \big(\, \widetilde{{}^h (\F^{^{\rm sc}}}),\F\!ct_{U_{h'}}\circ \tilde\frak t_h\big)  = \{0\}
\eqno £8.2.1.$$\/}
\par
\noindent
{\bf Proof:} According to Proposition~£6.10, for  any $n\ge 1$ it suffices to prove that we have
$$\Bbb H^n \big({}^{h'}\!  (\tilde\F^{^{\rm sc}}), \F\!ct_{U_{h'}}
 \circ  (\frak s_{h'} \times  \frak x_{h})\big) = \{0\}
 \eqno £8.2.2.$$
In order to apply Theorem~£7.4, let us consider the full subcategory~${}^{h'}\frak F$  of the category  ${}^{U_{h'}} \ad (\tilde\F^{^{\rm sc}})$ of $U_{h'}\-$objects of~$\ad (\tilde\F^{^{\rm sc}})$ 
[5,~6.2] over the set of {\it faithful indecomposable\/} $U_{h'}\-$objects, 
namely over the {\it indecomposable $U_{h'}\-$objects $Q^\rho$\/}
of~$\ad (\tilde\F^{^{\rm sc}})$ (cf.~£6.24) such that the group homomorphism 
$\rho\,\colon U_{h'}\to \tilde\F (Q)$ is injective. Note that the indecomposable 
$U_{h'}\-$object $\big(\bigoplus_{u\in U_{h'}} P\big)^\pi\,,$ defined by the {\it regular\/} action of $U_{h'}$ on itself, is faithful.

\smallskip
Moreover, if $Q^\rho = \big(\bigoplus_{i\in I} Q_i)^\rho$ and
$R^\sigma = \big(\bigoplus_{j\in J} R_j)^\sigma$ are faithful indecomposable 
${}^{U_{h'}}\ad (\tilde\F^{^{\rm sc}})\-$objects then,
according to £6.11 and~£6.13, the exterior intersection of $Q = \bigoplus_{i\in I} Q_i$
and $R =\bigoplus_{j\in J} R_j$ in
$\ad (\tilde\F^{^{\rm sc}})$ yields
$$Q\,\tilde{\cap}\, R = \bigoplus_{(i,j)\in I\times J}
\,\bigoplus_{(\tilde\alpha,T,\tilde\beta)\in \check\frak T_{Q_i,R_j}} T
\eqno £8.2.3\phantom{.}$$
and, for any $\xi\in U_{h'}\,,$ $\rho (\xi)$ and $\sigma (\xi)$ induce an automorphism of this intersection.
Thus, they induce a permutation of the disjoint union $\bigsqcup_{(i,j)\in I\times J} \check\frak T_{Q_i,R_j}$ and if  $\rho(\xi)$ and $\sigma (\xi)$ respectively fix $i$ and $j\,,$ and $(\tilde\alpha,T,\tilde\beta)\in 
\check\frak T_{Q_i,R_j}$ is a fixed element, then  $\rho (\xi)$ and $\sigma (\xi)$ induce $\tilde\F\-$automorphisms 
$\rho_i (\xi) $ of $Q_i\,,$ $\sigma_j (\xi)$ of~$R_j$ and $\tau (\xi)$ of $T$ fulfilling
$$\rho_i (\xi)\circ \tilde\alpha = \tilde\alpha\circ \tau (\xi)\qq \sigma_j (\xi)\circ
\tilde\beta = \tilde\beta\circ \tau (\xi)
\eqno £8.2.4.$$

\smallskip
In particular, if $\tau (\xi)$ is trivial then it follows from Corollary~£4.9 that
$\rho_i (\xi)$ and $\sigma_j (\xi)$ are trivial too and therefore we get $\xi = 1\,;$
thus, the exterior intersection $Q^\rho\,\tilde{\cap}\, R^\sigma$ in the category 
${}^{U_h}\ad (\tilde\F^{^{\rm sc}})$ is a {\it direct sum\/} of {\it faithful
indecomposable $U_{h'}\-$objects\/}. Consequently, since an indecomposable direct
summand of a {\it pull-back\/}  is also a direct summand of some exterior intersection  [5,~6.18],
{\it the subcategory  $\ad ({}^{h'}\frak F)$ of $\,{}^{U_{h'}}\ad (\tilde\F^{^{\rm sc}})$ is closed by  direct products and pull-backs.\/}

\smallskip
On the other hand, it is clear that ${}^{^{h'}}\!(\tilde\F^{^{\rm sc}})$ is a full 
subcategory of ${}^{h'}\frak F$ and we claim that the {\it contravariant\/} functor
(cf.~£6.9)
$$\frak n_h = \F\!ct_{U_{h'}}\circ (\frak s_{h'} \times \frak x_{h}) : {}^{^{h'}}\!(\tilde\F^{^{\rm sc}})
\too \O\-\mod
\eqno £8.2.5\phantom{.}$$ 
can be extended to a {\it contravariant\/} functor $\frak m_{h}\,\colon {}^{h'}\frak F\to \O\-\mod$ admitting a {\it right-hand sectional functor\/} $\frak m_h^\circ\,.$ Indeed, for any faithful indecomposable ${}^{U_{h'}}\ad(\tilde\F^{^{\rm sc}})\-$object $Q^\rho = \big(\bigoplus_{i\in I} Q_i\big)^\rho$ choose $i\in I\,;$ it is clear that $\rho$ induces a group homomorphism $\rho_i\,\colon U_{h_i} \to \tilde\F (Q_i)$ where $U_{h_i}$ denotes the stabilizer of $i$ in $U_{h'}\,,$ and then
we define
$$\frak m_{h} (Q^\rho) = \F\!ct_{U_{h_i}}\big((\varpi_{h_i,\skew4\hat{\tilde\F}(Q_i)})^{-1} (\rho_i)
\times \widetilde{\rm Mon}  (U_{h^p},(Q_i)^{\widehat{\rho_i(U_{h_i})}}),\O\big)
\eqno £8.2.6\phantom{.}$$
for the lifting $\widehat{\rho_i(U_{h_i})}\i \F (Q_i)$ of $\rho_i(U_{h_i})$ chosen in~£6.9.

\smallskip
For any faithful indecomposable ${}^{U_{h'}}\ad (\tilde\F^{^{\rm sc}})\-$object
$R^\sigma = \big(\bigoplus_{j\in J} R_j)^\sigma$ and any ${}^{h'}\frak F\-$morphism
$\tilde\varphi\,\colon R^\sigma\to Q^\rho\,,$ $\tilde\varphi$ determines a necessarily
surjective $U_{h'}\-$set map $f\,\colon J\to I$ and  an $\tilde\F^{^{\rm sc}}\-$morphism $\tilde\varphi_j\,\colon R_j\to Q_{i'}\,,$ where $j\in J$ is the chosen element and we set $i' = f(j)\,;$ in particular, it is clear that $U_{h_j}\i U_{h_{i'}}$ or, equivalently, that $h_j$
divides $h_{i'}$ and therefore, denoting by $\rho_j$ the restriction to $U_{h_j}$ of
the group homomorphism
$$\rho_{i'} : U_{h_{i'}}\too \tilde\F (Q_{i'})
\eqno £8.2.7,$$
$\tilde\varphi_j$ becomes an ${}^{h_j}(\tilde\F^{^{\rm sc}})\-$morphism from
$(R_j)^{\sigma_j}$ to $(Q_{i'})^{\rho_{j}}\,,$ so that from [5,~Pro-position~14.28] we get a $U_{h_j}\-$set bijection$$\tilde\frak s_{h_j}(\tilde\varphi_j) :
(\varpi_{h_j,\skew4\hat{\tilde\F}(R_j)})^{-1}(\sigma_j)\cong
(\varpi_{h_j,\skew4\hat{\tilde\F}(Q_{i'})})^{-1}(\rho_{j})
\eqno £8.2.8.$$
On the other hand, it is clear that the inclusion $U_{h_j}\i U_{h_{i'}}$ determines
the commutative diagram [5,~14.16.3]
$$\matrix{{\rm Mon}_{k^*} \big(\hat U_{h_j},\skew4\hat{\tilde\F}(Q_{i'})\big) 
&\buildrel \varpi_{h_j,\skew4\hat{\tilde\F}(Q_{i'})}\over 
{\hbox to 40pt{\rightarrowfill}}&{\rm Mon} \big(U_{h_j},\tilde\F(Q_{i'})\big)\cr
\big\uparrow&\phantom{\Big\uparrow}&\big\uparrow\cr
{\rm Mon}_{k^*} \big(\hat U_{h_{i'}},\skew4\hat{\tilde\F}(Q_{i'})\big)
&\buildrel \varpi_{h_{i'},\skew4\hat{\tilde\F}(Q_{i'})}
\over{\hbox to 40pt{\rightarrowfill}}
&{\rm Mon} \big(U_{h_{i'}},\tilde\F(Q_{i'})\big)\cr} 
\eqno £8.2.9.$$

\smallskip
Similarly, for the  chosen lifting $\widehat{\sigma_j (U_{h_j}}) \i \F (R_j)$ of $\sigma_j (U_{h_j})\,,$ since we assume that $\tilde\varphi$  is an ${}^{h'}\frak F\-$morphism,   a suitable representative $\varphi_j$ of  
$\tilde\varphi_j$ induces a group homomorphism from 
$(R_j)^{\widehat{\sigma_j (U_{h_j})}}$ to $(Q_{f(j)})^{\widehat{\rho_{i'} (U_{h_{i'}}})}\,,$ which defines a map
$$\widetilde{\rm Mon} \big (U_{h^p},(R_j)^{\widehat{\sigma_j (U_{h_j}})}\big)\too 
\widetilde{\rm Mon} \big(U_{h^p},(Q_{i'})^{\widehat{\rho_{i'} (U_{h_{i'}}})}\big)
\eqno £8.2.10.$$
Then, applying the functor $\F\!ct_{U_{h_j}}$ and  the inclusion
$\F\!ct_{U_{h_{i'}}}\i \F\!ct_{U_{h_j}}\circ \res^{U_{h_{i'}}}_{U_{h_j}}$ to
diagram~£8.2.9 and to its bottom map respectively, the bijection~£8.2.8
and the map~£8.2.10 determine an $\O\-$module homomorphism
$$\matrix{ \F\!ct_{U_{h_j}}\Big((\varpi_{h_j,\skew4\hat{\tilde\F}(R_j)})^{-1}(\sigma_j)\times \widetilde{\rm Mon} \big (U_{h^p},(R_j)^{\widehat{\sigma_j (U_{h_j}})}\big),\O\Big)\cr 
\hskip-20pt{\scriptstyle {\rm res}_{\tilde\varphi_j}}\hskip2pt \big\uparrow\phantom{\Big\uparrow}\cr 
\F\!ct_{U_{h_{i'}}}\Big((\varpi_{h_{i'},\skew4\hat{\tilde\F} (Q_{i'})})^{-1}
(\rho_{i'})\times \widetilde{\rm Mon} \big(U_{h^p},(Q_{i'})^{\widehat{\rho_{i'} (U_{h_{i'}}})}\big),\O\Big)\cr}
\eqno £8.2.11\phantom{.}$$
which admits a section ${\rm pro}_{\tilde\varphi_j}$ extending the $\O\-$valued functions by zero.
\eject

\smallskip
Moreover, there is $\xi\in U_{h'}$ such that $\rho (\xi)$ maps $i' = f(j)$ on the chosen
element $i\in I$ and therefore we have $h_{i'} = h_i$ and $\xi$ induces an
${}^{h_i}(\tilde\F^{^{\rm sc}})\-$morphism 
$$\rho (\xi)^{i'}_i : (Q_{i'})^{\rho_{i'}}\too (Q_i)^{\rho_i}
\eqno £8.2.12,$$
 so that it follows again from   [5,~Proposition~14.28] that we get a
$U_{h_i}\-$set bijection
$$\matrix{ \F\!ct_{U_{h_{i}}}\Big((\varpi_{h_{i},\skew4\hat{\tilde\F} (Q_{i})})^{-1}
(\rho_{i})\times \widetilde{\rm Mon} \big(U_{h^p},(Q_{i})^{\widehat{\rho_{i} (U_{h_{i}}})}\big),\O\Big)\cr 
\hskip-50pt {\scriptstyle \frak n_{ h_i h^p} \big(\rho (\xi)^{i'}_i \big)}\hskip4pt \wr\!\Vert\phantom{\bigg\uparrow}\cr 
\F\!ct_{U_{h_{i'}}}\Big((\varpi_{h_{i'},\skew4\hat{\tilde\F} (Q_{i'})})^{-1}
(\rho_{i'})\times \widetilde{\rm Mon} \big(U_{h^p},(Q_{i'})^{\widehat{\rho_{i'} (U_{h_{i'}}})}\big),\O\Big)\cr}
\eqno £8.2.13\phantom{.}$$
which clearly does not depend on the choice of $\xi\,.$ Finally, we consider the
compositions
$$\eqalign{\frak m_h (\tilde\varphi) =  {\rm res}_{\tilde\varphi_j}\circ \frak n_{ h_i h^p}
\big(\rho (\xi)^{i'}_i\big) &: \frak m_{h} (Q^\rho)\too \frak m_{h} (R^\sigma)\cr 
\frak m_h^\circ (\tilde\varphi)  =   \frak n_{ h_ih^p} \big(\rho (\xi)^{i'}_i\big)^{-1}
\circ {\rm pro}_{\tilde\varphi_j} &: \frak m_{h} (R^\sigma)\too \frak m_{h} (Q^\rho)\cr}
\eqno £8.2.14;$$
 we claim that the correspondence $\frak m_{h}$ is the announced {\it contravariant\/} functor and that $\frak m_h^\circ$ is a {\it right-hand sectional functor\/} of  $\frak m_h\,.$

\smallskip
 Indeed, it is clear that $\frak m_h$ extends $\frak n_h$ and that we have
 $$\frak m_h (\tilde\varphi)\circ \frak m_h^\circ (\tilde\varphi) = {\rm id}_{\frak m_{h} (R^\sigma)}
\eqno £8.2.15; $$
further, for any faithful indecomposable ${}^{U_{h'}}\ad (\tilde\F^{^{\rm sc}})\-$object $T^\tau = \big(\bigoplus_{\ell\in L} T_\ell)^\tau$ and any ${}^{h'}\frak F\-$morphism
$\tilde\psi\,\colon T^\tau\to R^\sigma\,,$ as above we have a surjective
$U_{h'}\-$set map $g\,\colon L\to J\,,$ a chosen element~$\ell$ in~$L$ and, setting
$j' = g(\ell)$ and denoting by $\sigma_\ell$ the restriction of $\sigma_{j'}$ to 
$U_{h_\ell}\,,$ an ${}^{h_\ell}(\tilde\F^{^{\rm sc}})\-$morphism 
$\tilde\psi_\ell\,\colon (T_\ell)^{\tau_\ell}\to (R_{j'})^{\sigma_\ell}\,,$ together with a
$U_{h_\ell}\-$set bijection and a map
$$\eqalign{\tilde\frak s_{h_\ell}(\tilde\psi_\ell) :
(\varpi_{h_\ell,\skew4\hat{\tilde\F}(T_\ell)})^{-1}(\tau_\ell) &\cong
(\varpi_{h_\ell,\skew4\hat{\tilde\F}(R_{g(\ell)})})^{-1}(\sigma_\ell)\cr
\widetilde{\rm Mon} \big (U_{h^p},(T_\ell)^{\widehat{\tau_{\ell} (U_{h_{\ell}}})}\big) &\too 
\widetilde{\rm Mon} \big(U_{h^p},(R_{j'})^{\widehat{\sigma_{j'} (U_{h_{j'}}})}\big)\cr}
\eqno £8.2.16\phantom{.}$$
for the chosen lifting $\widehat{\tau_\ell (U_{h_\ell})}\i \F (T_\ell)$ of $\tau_\ell (U_{h_\ell})\,.$
Analogously, we have an $\O\-$module isomorphism
$$\matrix{ \F\!ct_{U_{h_\ell}}\Big((\varpi_{h_\ell,\skew4\hat{\tilde\F}(T_\ell)})^{-1}(\tau_\ell)\times \widetilde{\rm Mon} \big (U_{h^p},(T_\ell)^{\widehat{\tau_{\ell} (U_{h_{\ell}}})}\big),\O\Big)\cr 
\hskip-20pt{\scriptstyle {\rm res}_{\tilde\psi_\ell}}\hskip2pt \big\uparrow\phantom{\Big\uparrow}\cr 
\F\!ct_{U_{h_{j'}}}\Big((\varpi_{h_{j'},\skew4\hat{\tilde\F} (R_{j'})})^{-1}
(\sigma_{j'})\times \widetilde{\rm Mon} \big(U_{h^p},(R_{j'})^{\widehat{\sigma_{j'} (U_{h_{j'}}})}\big),\O\Big)\cr}
\eqno £8.2.17\,.$$

\smallskip
Moreover, chosing $\zeta\in U_{h'}$ such that $\sigma (\zeta)$ maps 
$g(\ell)$ on $j\,,$  as above we have $h_{j'} = h_j$ and $\zeta$ induces an
${}^{h_j}(\tilde\F^{^{\rm sc}})\-$isomorphism 
$$\sigma (\zeta)^{j'}_j : (R_{j'})^{\sigma_{j'}}\cong (R_j)^{\sigma_j}
\eqno £8.2.18,$$ 
and therefore we get a $U_{h_j}\-$set bijection (cf.~Proposition~£14'.27)
$$\matrix{ \F\!ct_{U_{h_{j}}}\Big((\varpi_{h_{j},\skew4\hat{\tilde\F} (R_{j})})^{-1}
(\sigma_{j})\times \widetilde{\rm Mon} \big(U_{h^p},(R_{j})^{\widehat{\sigma_{j} (U_{h_{j}}})}\big),\O\Big)\cr 
\hskip-60pt {\scriptstyle \frak n_{ h_j h^p} \big(\sigma (\zeta)^{j'}_j \big)}\hskip4pt \wr\!\Vert\phantom{\bigg\uparrow}\cr 
\F\!ct_{U_{h_{j'}}}\Big((\varpi_{h_{j'},\skew4\hat{\tilde\F} (R_{j'})})^{-1}
(\sigma_{j'})\times \widetilde{\rm Mon} \big(U_{h^p},(R_{j'})^{\widehat{\sigma_{j'} (U_{h_{j}}})}\big),\O\Big)\cr}
\eqno £8.2.19.$$
Finally, we also consider
$$\eqalign{\frak m_h (\tilde\psi) =  {\rm res}_{\tilde\psi_\ell}\circ \frak n_{ h_j h^p}
\big(\sigma (\zeta)^{j'}_j\big) &: \frak m_{h} (R^\sigma)\too \frak m_{h} (T^\tau)\cr 
\frak m_h^\circ (\tilde\psi)  =   \frak n_{h_j h^p } \big(\sigma (\zeta)^{j'}_j\big)^{-1}
\circ {\rm pro}_{\tilde\psi_\ell} &: \frak m_{h} (T^\tau)\too \frak m_{h} (R^\sigma)\cr}
\eqno £8.2.20.$$

\smallskip
Consequently, we get
$$\eqalign{\frak m_h(\tilde\psi)&\circ \frak m_h (\tilde\varphi) = {\rm res}_{\tilde\psi_\ell}\circ \frak n_{ h_j h^p}\big(\sigma (\zeta)^{j'}_j\big) \circ   {\rm res}_{\tilde\varphi_j}\circ \frak n_{h_i h^p } \big(\rho (\xi)^{i'}_i\big)\cr
\frak m_h^\circ (\tilde\varphi)&\circ \frak m_h^\circ (\tilde\psi) = \frak n_{ h_i h^p} \big(\rho (\xi)^{i'}_i \big)^{-1} \circ  {\rm pro}_{\tilde\varphi_j}\circ  \frak n_{h_j h^p } \big(\sigma (\zeta)^{j'}_j\big)^{-1} \circ 
{\rm pro}_{\tilde\psi_\ell}\cr }
\eqno £8.2.21\phantom{.}$$
and in order to prove our claim it suffices to prove that
$$\eqalign{\frak n_{h_j h^p }\big(\sigma (\zeta)^{j'}_j\big) \circ {\rm res}_{\tilde\varphi_j} 
&= {\rm res}_{\tilde\varphi_{j'}}\circ \frak n_{ h_{i'} h^p}\big(\rho (\zeta)^{f(j')}_{i'}\big) \cr
 {\rm pro}_{\tilde\varphi_j}\circ  \frak n_{h_j h^p } \big(\sigma (\zeta)^{j'}_j\big)^{-1} &= \frak n_{ h_j h^p} \big(\rho (\zeta)^{f(j')}_{i'}\big)^{-1} \circ  
 {\rm pro}_{\tilde\varphi_{j'}}\cr}
 \eqno £8.2.22\phantom{.}$$
 since it is easily checked that 
 $${\rm res}_{\tilde\psi_\ell}\circ {\rm res}_{\tilde\varphi_{j'}} = 
 {\rm res}_{(\tilde\varphi\circ\tilde\psi)_\ell}\qq {\rm pro}_{\tilde\varphi_{j'}}
 \circ {\rm pro}_{\tilde\psi_\ell} = {\rm pro}_{(\tilde\varphi\circ\tilde\psi)_\ell}
  \eqno £8.2.23\phantom{.}$$

\smallskip
In order to prove equalities~£8.2.22 note that,
 since $f$ is a $U_{h'}\-$set map,  $\rho (\zeta)$ maps $i'' = f(j')$ on $f(j) = i'$ 
 and we have $h_{i''} = h_{i'} = h_i\,;$ thus, $\zeta$ induces~an
${}^{h_{i'}}(\tilde\F^{^{\rm sc}})\-$isomorphism 
$$\rho (\zeta)^{i''}_{i'} : (Q_{i''})^{\rho_{i''}}\cong
(Q_{i'})^{\rho_{i'}}
\eqno £8.2.24\phantom{.}$$
and denoting by $\rho_{j'}$ the restriction to $U_{h_{j'}}$  of the
group homomorphism $\rho_{i''}\,,$ the ${}^{h'}\frak F\-$morphism
$\tilde\varphi\,\colon R^\sigma\to Q^\rho$ forces the following commutative 
${}^{h_j}\frak F\-$diagram
$$\matrix{(R_{j})^{\sigma_{j}}&\buildrel \tilde\varphi_{j}\over\too 
& (Q_{i'})^{\rho_{j}}\cr
\wr\Vert&\phantom{\Big\uparrow}&\wr\Vert\cr
(R_{j'})^{\sigma_{j'}}&\buildrel \tilde\varphi_{j'}\over\too &(Q_{i''})^{\rho_{j'}} \cr}
\eqno £8.2.25;$$
hence,  we get the commutative ${}^{U_{h_j}} \aleph\-$diagram
$$\matrix{(\varpi_{h_{j},\skew4\hat{\tilde\F}(R_{j})})^{-1}(\sigma_{j})
&\buildrel \over\cong
&(\varpi_{h_{j},\skew4\hat{\tilde\F}(Q_{i'})})^{-1}(\rho_{j})\cr
\wr\Vert&\phantom{\Big\uparrow}&\wr\Vert\cr
(\varpi_{h_{j'},\skew4\hat{\tilde\F}(R_{j'})})^{-1}(\sigma_{j'})
&\buildrel \over\cong
&(\varpi_{h_{j'},\skew4\hat{\tilde\F}(Q_{i''})})^{-1}(\rho_{j'})\cr} 
\eqno £8.2.26.$$
On the other hand, the natural map $\theta_{h_\ell,h_j}$ in [5,~Proposition~14.28] applied to the ${}^{h_j}(\tilde\F^{^{\rm sc}})\-$morphism $\tilde\varphi_{j'}\,\colon
(R_{j'})^{\sigma_{j'}}\to (Q_{i''})^{\rho_{j'}}$ yields the following
commutative ${}^{U_{h_j}} \aleph\-$diagram
$$\matrix{(\varpi_{h_{j},\skew4\hat{\tilde\F}(R_{j'})})^{-1}(\sigma_{j'})
&\buildrel \over\cong
&(\varpi_{h_{j},\skew4\hat{\tilde\F}(Q_{i''})})^{-1}(\rho_{j'})\cr
\big\downarrow&\phantom{\Big\uparrow}&\big\downarrow\cr
{\res}^{U_{h_\ell}}_{U_{h_j}}\big((\varpi_{h_\ell,\skew4\hat{\tilde\F}(R_{j'})})^{-1}(\sigma_{\ell})\big)
&\buildrel \over\cong
&{\res}^{U_{h_\ell}}_{U_{h_j}}\big((\varpi_{h_\ell,\skew4\hat{\tilde\F}(Q_{i''})})^{-1}(\rho_{\ell})\big)\cr} 
\eqno £8.2.27\phantom{.}$$
 where $\rho_\ell$ is the restriction of the group homomorphism $\rho_{i''}$ 
 to~$U_{h_\ell}\,.$

\smallskip
Similarly, from the commutative $\widetilde\Gr\-$diagram
$$\matrix{(R_{j})^{\widehat{\sigma_{j} (U_{h_{j}}})}&\buildrel \over\too 
& (Q_{i'})^{\widehat{\rho_{i'} (U_{h_{i'}}})}\cr
\hskip-15pt\wr\Vert&\phantom{\Big\uparrow}&\hskip-15pt\wr\Vert\cr
(R_{j'})^{\widehat{\sigma_{j'} (U_{h_{j'}}})}&\buildrel \over\too &(Q_{i''})^{\widehat{\rho_{i''} (U_{h_{i''}}})} \cr}
\eqno £8.2.28\phantom{.}$$
we get the commutative diagram
$$\matrix{\widetilde{\rm Mon} \big (U_{h^p},(R_j)^{\widehat{\sigma_{j} (U_{h_{j}}})}\big)
&\too &\widetilde{\rm Mon} \big(U_{h^p},(Q_{i'})^{\widehat{\rho_{i'} (U_{h_{i'}}})}\big)\cr
\hskip-15pt\wr\Vert&\phantom{\Big\uparrow}&\hskip-15pt\wr\Vert\cr
\widetilde{\rm Mon} \big (U_{h^p},(R_{j'})^{\widehat{\sigma_{j'} (U_{h_{j'}}})}\big)&\too &
\widetilde{\rm Mon} \big(U_{h^p},(Q_{i''})^{\widehat{\rho_{i''} (U_{h_{i''}}})}\big)\cr}
\eqno £8.2.29.$$
At this point, considering the direct product of this diagram with the square diagram 
obtained by ``composing'' diagrams~£8.2.26 and~£8.2.27, and applying the functor
$\F\!ct_{U_{h_{j}}}\,,$ the equalities~£8.2.22 follow easily, proving our claim.
 Hence, it follows from Theorem~£7.4 that, for any $n\ge 1\,,$ we have
$$\Bbb H^n ({}^{h'}\!\frak F,\frak m_h) = \{0\}
\eqno £8.2.30.$$

\smallskip
If $h' =1$ then ${}^1\frak F = {}^1(\tilde\F^{^{\rm sc}}) = \tilde\F^{^{\rm sc}}$ and we are done. Otherwise, we consider the full subcategory ${}^{h'}\frak E$ of 
${}^{h'}\frak F$ over the faithful indecomposable ${}^{U_{h'}}
\ad (\tilde\F^{^{\rm sc}})\-$objects $Q^\rho =\big(\bigoplus_{i\in I} Q_i\big)^\rho$ such that $\vert I\vert > 1\,,$ namely over all the  ${}^{h'}\frak F\-$objects which are {\it not\/} ${}^{^{h'}}\!(\tilde\F^{^{\rm sc}}) \-$objects; note that there is {\it no\/}  ${}^{h'}\frak F\-$morphism from an ${}^{^{h'}}\!(\tilde\F^{^{\rm sc}})\-$object to an
${}^{h'}\frak E\-$object.  Denoting by $\frak l_h\,\colon {}^{h'}\frak E\to
\O\-\mod$ the restriction of~$\,\frak m_h$ to ${}^{h'}\frak E\,,$ it is quite clear that 
Theorem~£7.4 applies again to ${}^{h'}\frak E$ and~$\frak l_h\,,$ so that for any 
$n\ge 1$ we get
$$\Bbb H^n ({}^{h'}\frak E,\frak l_h) = \{0\}
\eqno £8.2.31.$$

\smallskip
Recall that, denoting by $\bf \Delta$ the {\it simplicial $2\-$category\/} [5,~A1.7],
the cohomology groups $\Bbb H^n ({}^{h'}\frak F,\frak m_h)\,,$ 
$\Bbb H^n ({}^{h'}\frak E,\frak l_h)$ and $\Bbb H^n\big( {}^{^{h'}}\!
(\tilde\F^{^{\rm sc}}), \frak n_h\big)$ are nothing but the homology groups of the respective evident functors  
$\frak c_{\frak m_h}\,,$  $\frak c_{\frak l_h}$ and $\frak c_{\frak n_h}$ from~$\bf
\Delta$ to~$\O\-\mod\,,$ mapping $\Delta_n$ on [5,~A3.8]
$$\eqalign{\Bbb C^n ({}^{h'}\frak F,\frak m_h) 
&= \prod_{\tilde\frak q\in \Fct (\Delta_n,{}^{h'}\frak F) }\frak m_h \big(\tilde\frak q
(0)\big)\cr 
\Bbb C^n ({}^{h'}\frak E,\frak l_h) 
&= \prod_{\tilde\frak q\in \Fct(\Delta_n,{}^{h'}\frak E)}
\frak m_h\big(\tilde\frak q(0)\big)\cr 
\Bbb C^n \big({}^{^{h'}}\!(\tilde\F^{^{\rm sc}}),\frak n_h\big) 
&= \prod_{\tilde\frak q\in \Fct (\Delta_n,{}^{^{h'}}\!(\tilde\F^{^{\rm sc}})) } \frak m_h
\big(\tilde\frak q (0)\big)\cr}
\eqno £8.2.32;$$
then, the inclusion ${}^{^{h'}}\!(\tilde\F^{^{\rm sc}}) \i {}^{h'}\frak F$ clearly 
determines a surjective natural map 
$$\mu_h : \frak c_{\frak m_h}\too \frak c_{\frak n_h}
\eqno £8.2.33,$$
so that we obtain a fourth functor 
$${\Ker}(\mu_h) : {\bf \Delta}\too \O\-\mod
\eqno £8.2.34\phantom{.}$$

\smallskip
Moreover, since there is {\it no\/} ${}^{h'}\frak F\-$morphism from an ${}^{^{h'}}
\!(\tilde\F^{^{\rm sc}})\-$object to any ${}^{h'}\frak E\-$object,  we have an $\O\-$module isomorphism
$$\big({\Ker}(\mu_h)\big)(\Delta_n)\cong  \prod_{\tilde\frak q} \frak m_h 
\big(\tilde\frak q (0)\big)
\eqno £8.2.35,$$
where $\tilde\frak q$ runs over the set $\E^{h'}(\Delta_n)$ of functors from $\Delta_n$ to  ${}^{h'}\frak F$ such that  $\tilde\frak q (0)$ is a ${}^{h'}\frak E\-$object  and 
therefore the inclusion ${}^{h'}\frak E\i {}^{h'}\frak F$ also determines a surjective natural map 
$$\lambda_h : {\Ker}(\mu_h)\too \frak c_{\frak l_h}
\eqno £8.2.36.$$
But on the one hand, we already know that $\Bbb H_{n} (\frak c_{\frak m_h}) = \{0\} 
=\Bbb H_{n} (\frak c_{\frak l_h})$ for any $n\ge 1$ (cf.~equalities~£8.2.30
and~£8.2.31) and on the other hand, setting 
$$\Bbb H_{-n} (\frak c_{\frak n_h}) =\Bbb H_{-n} (\frak c_{\frak m_h}) = \Bbb H_{-n}
\big({\Ker}(\mu_h)\big) = \{0\}
\eqno £8.2.37\phantom{.}$$
 for any $n>0\,,$ there is a $1\-${\it graded connecting homomorphism\/} [5,~A3.3.4]
$$\delta : \bigoplus_{n\in \Bbb Z} \Bbb H_{n} (\frak c_{\frak n_h})\too 
\bigoplus_{n\in \Bbb Z} \Bbb H_{n} \big({\Ker}(\mu_h)\big)
\eqno £8.2.38\phantom{.}$$
such that we have the following {\it exact triangle\/} [5,~A3.3.5] 
$$\matrix{\bigoplus_{n\in \Bbb Z} \Bbb H_{n} (\frak c_{\frak n_h})&\buildrel{\delta}
\over\too &\bigoplus_{n\in \Bbb Z} \Bbb H_{n} \big({\Ker}(\mu_h)\big)\cr
{\scriptstyle{\oplus_{n\in \Bbb Z} \Bbb H_{n}
(\mu_h)}}\quad\nwarrow&\phantom{\bigg\downarrow}&\hskip-50pt\swarrow\quad\cr 
&\bigoplus_{n\in \Bbb Z} \Bbb H_{n}(\frak c_{\frak m_h})&\cr}
\eqno £8.2.39.$$

\smallskip
Consequently, in order to prove that $ \Bbb H_{n} (\frak c_{\frak n_h}) =\{0\}$ for any
$n\ge 1\,,$ it suffices to show that $\Bbb H_n (\lambda_h)$ is injective. Actually, since
the functors $\frak s_{h'}$ and $\frak x_{h}$ factorize throughout the exterior quotient
$\widetilde{{}^{^{h'}}}\!(\tilde\F^{^{\rm sc}})$ of ${}^{^{h'}}\!(\tilde\F^{^{\rm sc}})$
[5,~Remark~14.29], it is easily checked that the {\it contravariant\/} functors 
$\frak n_h\,,$ $\frak m_h$ and~$\frak l_h$ respectively determine {\it contravariant\/}
functors
$$\eqalign{\tilde\frak n_h : \widetilde{{}^{^{h'}}}\!(\tilde\F^{^{\rm sc}})&\too
\O\-\mod\cr
\tilde\frak m_h :  \widetilde{{}^{h'}}\frak F\too \O\-\mod&\qq
\tilde\frak l_h : \widetilde{{}^{h'}}\frak E\too \O\-\mod\cr}
\eqno £8.2.40,$$
where $\widetilde{{}^{h'}}\frak F$ and $\widetilde{{}^{h'}}\frak E$ denote the corresponding {\it exterior quotients\/} [5,~6.3]. Coherently, we get the corresponding
functors $\frak c_{\tilde\frak n_h}\,,$ $\frak c_{\tilde\frak m_h}$ and 
$\frak c_{\tilde\frak l_h}$ (cf.~£8.2.32), and the corresponding natural maps
(cf.~£8.2.33 and~£8.2.36)
$$\tilde\mu_h : \frak c_{\tilde\frak m_h}\too \frak c_{\tilde\frak n_h}
\qq \tilde\lambda_h : {\Ker}(\tilde\mu_h)\too \frak c_{\tilde\frak l_h}
\eqno £8.2.41.$$

\smallskip
Then, it follows from equality~£7.3.2 that, for any $n\ge 1\,,$ it suffices to prove
that $\Bbb H_n (\tilde\lambda_h)$ is injective. First of all note that, up to
isomorphisms, any  $\widetilde{{}^{h'}} \frak F\-$object has the {\it canonical\/} form $Q^\rho =\big(\bigoplus_{\bar\xi\in U_{h'}/U} Q_{\bar\xi}\big)^\rho$ for a suitable subgroup $U$ of $U_{h'}\,,$ and, denoting by $\tilde\frak i_{h'}\,\colon \widetilde{{}^{h'}}\frak E \to \widetilde{{}^{h'}}\frak F$ the {\it inclusion\/} functor, we claim that we have a functor and a natural map
$$\tilde\frak e_{h'} : \widetilde{{}^{h'}}\frak F\too \widetilde{{}^{h'}}\frak E \qq \tilde\varepsilon_{h'}
: \tilde\frak i_{h'}\circ \tilde\frak e_{h'}\too \id_{\tilde{{}^{h'}}\frak F}
\eqno £8.2.42\phantom{.}$$
which respectively map $Q^\rho$ on the ${}^{U_{h'}}
\ad (\tilde\F^{^{\rm sc}})\-$object $\hat Q^{\hat\rho}$ formed by 
$$\hat Q = {\textstyle\bigoplus_{\xi\in U_{h'}}} Q_{\bar\xi}
\eqno £8.2.43\phantom{.}$$
and by the group homomorphism $\hat\rho\,\colon U_{h'}\to \widetilde{{}^{h'}}\frak F (\hat Q)$ defined by the regular action of $U_{h'}$ on itself, together with the 
$\tilde\F^{^{\rm sc}}\-$isomorphisms from the $\xi\-$summand to the 
$\xi'\-$summand induced by $\rho(\xi'\xi^{-1})\,,$ and  on the 
$\widetilde{{}^{h'}}\frak F\-$morphism
$$(\tilde\varepsilon_{h'})_{Q^\rho} : \hat Q^{\hat\rho}\too Q^\rho
\eqno £8.2.44\phantom{.}$$
defined by the canonical map $U_{h'}\to  U_{h'}/U$ and by the identity automorphism
of~$Q_{\bar\xi}$ for any~$\xi\in U_h\,.$
\eject

\smallskip
Indeed, if $R^\sigma = \big(\bigoplus_{\tilde\xi\in U_{h'}/V} R_{\tilde\xi}\big)^\sigma$ 
is another $\widetilde{{}^{h'}}\frak F\-$object, an $\widetilde{{}^{h'}}\frak F\-$morphism from $R^\sigma$ to $Q^\rho$ forces the inclusion $V\i U$ and admits a {\it canonical\/} representative  formed by the canonical map $U_{h'}/V\to U_{h'}/U$ and by an $\tilde\F^{^{\rm sc}}\-$morphism  $\tilde\varphi_{\tilde\xi}\,\colon R_{\tilde\xi}\to 
Q_{\bar\xi}$ for any $\tilde\xi\in U_{h'}/V\,;$ then, the functor  $\tilde\frak e_{h'}$ above maps this $\widetilde{{}^{h'}}\frak F\-$morphism on the $\widetilde{{}^{h'}}\frak E\-$morphism
$$\hat R^{\hat\sigma} = \big({\textstyle\bigoplus_{\xi\in U_{h'}}}
R_{\tilde\xi}\big)^{\hat\sigma}Ê\too \hat Q^{\hat\rho} = 
\big({\textstyle\bigoplus_{\xi\in U_{h'}}} Q_{\bar\xi}\big)^{\hat\rho}
\eqno £8.2.45\phantom{.}$$
admitting a representative formed by the identity map of $U_{h'}$  and by the
$\tilde\F^{^{\rm sc}}\-$morphism  $\tilde\varphi_{\tilde\xi}\,\colon R_{\tilde\xi}\to
Q_{\bar\xi}$ for any $\xi\in U_{h'}\,;$ it is quite clear that this correspondence preserves
the composition of $\widetilde{{}^{h'}}\frak F\-$morphisms and is compatible with  the 
$\widetilde{{}^{h'}}\frak F\-$morphisms~£8.2.44.

\smallskip
Secondly, we consider the {\it contravariant\/} functor  $\tilde\frak m_h\circ
\tilde\frak i_{h'}\circ\tilde\frak e_{h'} = \tilde\frak l_h\circ\tilde\frak e_{h'}$ which, as above, determines a functor 
$$\frak c_{\tilde\frak l_h\circ\tilde\frak e_{h'}} : {\bf \Delta}\too \O\-\mod
\eqno £8.2.46\phantom{.}$$ 
mapping $\Delta_n$ on 
$$\Bbb C^n ({}^{\widetilde {h'}}\frak F,\tilde\frak l_h\circ\tilde\frak e_{h'}) =
\prod_{\tilde\frak q\in \Fct(\Delta_n,{}^{\widetilde {h'}}\frak F)}\tilde\frak m_h
\Big(\hat\frak e_{h'}\big(\tilde\frak q(0)\big)\Big)
\eqno £8.2.47,$$
and then we get a natural map $\kappa_{\tilde\frak e_{h'}}\,\colon
\frak c_{\tilde\frak m_h}\to \frak c_{\tilde\frak l_h\circ\tilde\frak e_{h'}}$ sending
$\Delta_n$ to the $\O\-$mo-dule homomorphism
$$(\kappa_{\tilde\frak e_{h'}})_n :\!\!\prod_{\tilde\frak q\in\Fct(\Delta_n,
\widetilde{{}^{h'}}\frak F)}
\tilde\frak m_h\big(\tilde\frak q(0)\big)\too
\!\!\prod_{\tilde\frak q\in \Fct(\Delta_n,\widetilde{{}^{h'}}\frak F)} \tilde\frak m_h
\Big(\hat\frak e_{h'}\big(\tilde\frak q(0)\big)\Big)
\eqno £8.2.48\phantom{.}$$
mapping $m = (m_{\tilde\frak q})_{\tilde\frak q\in \Fct(\Delta_n,
\widetilde{{}^{h'}}\frak F)}$ on $(m_{\tilde\frak i_{h'}\circ\tilde\frak e_{h'}\circ
\tilde\frak q})_{\tilde\frak q\in \Fct(\Delta_n,\widetilde{{}^{h'}}\frak F)}\,.$ Note that if $m$ belongs to the
kernel ${\rm Ker}\big((\tilde\lambda_h)_{\Delta_n}\big)$ (cf.~£8.2.41) then, for any  $\tilde\frak q\in \Fct(\Delta_n,\widetilde{{}^{h'}}\frak E)\,,$ we have $m_{\tilde\frak q} = 0$ and therefore we get $(\kappa_{\tilde\frak e_{h'}})_n ( m) 
= 0\,.$ Moreover, according to the very definition of $\tilde\frak m_h$ (cf.~£8.2.14 and~£8.2.40), the natural map (cf.~£8.2.42)
 $$\tilde\frak m_h *\tilde\varepsilon_{h'} : \tilde\frak m_h \too \tilde\frak l_h\circ
\tilde\frak e_{h'}
\eqno £8.2.49\phantom{.}$$ 
determines a natural isomorphism
$$ \kappa_{\tilde\varepsilon_{h'}} : \frak c_{\tilde\frak m_h}\cong \frak c_{\tilde\frak
l_h\circ\tilde\frak e_{h'}}
\eqno £8.2.50\phantom{.}$$
sending $\Delta_n$ to the $\O\-$module isomorphism
$$(\kappa_{\tilde\varepsilon_{h'}})_n : \prod_{\tilde\frak q\in \Fct(\Delta_n,
\widetilde{{}^{h'}}\frak F)} \tilde\frak m_h\big(\tilde\frak q(0)\big)\cong
\prod_{\tilde\frak q\in \Fct(\Delta_n,\widetilde{{}^{h'}}\frak F)} \tilde\frak m_h
\Big(\hat\frak e_{h'}\big(\tilde\frak q(0)\big)\Big)
\eqno £8.2.51\phantom{.}$$
mapping $m = (m_{\tilde\frak q})_{\tilde\frak q\in \Fct(\Delta_n,
\widetilde{{}^{h'}}\frak F)}$
on $\bigg(\Big(\tilde\frak m_h \big((\tilde\varepsilon_{h'})_{\tilde\frak q (0)}\big)\Big)
(m_{\tilde\frak q})\bigg)_{\tilde\frak q\in \Fct(\Delta_n,\widetilde{{}^{h'}}\frak F)}\,.$
\eject

\smallskip
 At this point, for any $n\ge 1\,,$ following the notation introduced in [5,~Lemma~A4.2],
 we consider the $\O\-$module  homomorphism
$$h_{n-1} : \frak c_{\tilde\frak m_h}(\Delta_n)\too
\frak c_{\tilde\frak l_h \circ \tilde\frak e_{h'}}(\Delta_{n-1})
\eqno £8.2.52\phantom{.}$$
mapping $m = (m_{\tilde\frak q})_{\tilde\frak q\in \Fct(\Delta_n,
\widetilde{{}^{h'}} \frak F)} \in  \frak c_{\tilde\frak m_h}(\Delta_n)$ on
$$h_{n-1} (m) = \Big( \sum_{i=0}^{n-1}(-1)^i \,m_{\frak h_i^{n-1} 
(\tilde\varepsilon_{h'} * \tilde\frak r)}\Big)_{\tilde\frak r\in
\Fct(\Delta_n,\widetilde{{}^{h'}}\frak F)}
\eqno £8.2.53\phantom{.}$$
\smallskip
 Then, respectively denoting by $d_n$ and $\hat d_n$ the differential maps for the functors $\frak c_{\tilde\frak m_h}$ and~$\frak c_{\tilde\frak l_h \circ 
 \tilde\frak e_{h'}}$ [5,~A3.2], we claim that
$$(\kappa_{\tilde\varepsilon_{h'}})_n (m) - (\kappa_{\tilde\frak e_{h'}})_n  (m) 
= (\hat d_{n-1} \circ h_{n-1} + h_{n}\circ d_n)(m)
\eqno £8.2.54.$$
Indeed, setting $\theta = \frak m_h \Big(\tilde\frak e_{h'}
\big(\tilde\frak q (0\!\bullet\! 1)\big)\Big)\,,$ for any $\tilde\frak q\in
\Fct(\Delta_n,\widetilde{{}^{h'}}\frak F)$ we have
$$\eqalign{&\hat d_{n-1} \big(h_{n-1} (m)\big)_{\tilde\frak q} = \sum_{i=0}^{n-1}(-1)^i \hat d_{n-1}
\big(( m_{\frak h_i^{n-1}  (\tilde\varepsilon_{h'} *  \tilde\frak r)})_{\tilde\frak r} \big)_{\tilde\frak q}\cr        
& = \sum_{i = 0}^{n-1}(- 1)^i \big(\theta( m_{\frak h_i^{n-\!1} 
(\tilde\varepsilon_{h'} * (\tilde\frak q\circ\delta_0^{n-\!1}))}) + \sum_{j=1}^{n} (-
1)^j m_{\frak h_i^{n-\!1} (\tilde\varepsilon_{h'} * (\tilde\frak q\circ\delta_j^{n-\!1}))}\big)\cr}
\eqno £8.2.55\phantom{.}$$
where $\frak r$ runs over $\Fct(\Delta_{n-1},
\widetilde{{}^{h'}}\frak F)\,;$ analogously,   we still have
$$\eqalign{h_{n}\big( d_n&(m)\big)_{\tilde\frak q} =  h_{n} \Big(\big(\theta ( m_{\tilde\frak r\circ \delta_0^n})\big)_{\frak r }\Big)_{\tilde\frak q} + \sum_{j=1}^{n+1} (-1)^j  h_{n} 
\big((m_{\frak r\circ \delta_j^n})_{\frak r }\big)_{\tilde\frak q}\cr     
&= \sum_{i= 0}^{n}(-1)^i\big(\theta ( m_{\frak h_i^n (\tilde\varepsilon_{h'} * 
\tilde\frak q) \circ \delta_0^n}) + \sum_{j=1}^{n+1}(-1)^j\, m_{\frak h_i^n
(\tilde\varepsilon_{h'} * \tilde\frak q)\circ \delta_j^n}\big) \cr  }
\eqno £8.2.56\phantom{.}$$ 
where $\frak r$ runs over $\Fct(\Delta_{n+1}, \widetilde{{}^{h'}}\frak F)\,.$ But from  [5,~Lemma~A4.2] we know that
$$\eqalign{\frak h_{i+1}^n(\tilde\varepsilon_{h'}* \tilde\frak q) \circ\delta_{i+1}^n 
&= \frak h_i^n (\tilde\varepsilon_{h'} * \tilde\frak q)\circ\delta_{i+1}^n\cr 
\frak h_i^{n-1} \big(\tilde\varepsilon_{h'} * (\tilde\frak q\circ\delta_j^{n-1})\big) 
&= \cases{  \frak h_{i+1}^{n} (\tilde\varepsilon_{h'} * \tilde\frak q) \circ \delta_j^{n}
&{if $j\le i$}\cr  
\frak h_i^{n} (\tilde\varepsilon_{h'} * \tilde\frak q)\circ \delta_{j+1}^n
&{if $i < j\,.$}\cr}\cr}
\eqno £8.2.57.$$
Consequently, in equality~£8.2.56 the terms where $j = i$ and 
$j = i+1$ cancel with each other for any $1\le j\le n\,;$ moreover, the term $(i,j)$ in
equality~£8.2.55 cancel either with the term $(i+1,j)$ if $1\le j\le i\le n-1\,,$ or
with the term $(i,j+1)$ if $0\le i <j\le n$ in equality~£8.2.56.
Finally, the term $(i,0)$ in equality~£8.2.52 cancel with the
term $(i+1,0)$ in equality~£8.2.56 for any $0\le i\le n-1\,,$
whereas the terms $(0,0)$ and $(n,n+1)$ respectively coincide with
$(\kappa_{\tilde\varepsilon_{h'}})_n (m)$ and 
$- (\kappa_{\tilde\frak e_{h'}})_n  (m)\,,$ which proves the claim.
\eject

\smallskip
In conclusion, if $m\in {\rm Ker}(d_n)$ lifts an element of $\Bbb H_n
(\tilde\lambda_n)$ then we may assume that $m$ belongs to ${\rm Ker} \big((\lambda_h)_{\Delta_n}\big)\,,$ so that  we get $(\kappa_{\tilde\frak e_{h'}})_n ( m) = 0\,,$ and therefore it follows from isomorphism~£8.2.50 and from 
equality~£8.2.54 that 
$$\eqalign{(\kappa_{\tilde\varepsilon_{h'}})_n (m)    
&= \hat d_{n-1} \big( h_{n-1} (m)\big)\cr
& = \hat d_{n-1}\Big( \big((\kappa_{\tilde\varepsilon_{h'}})_{n-1}\circ
\big((\kappa_{\tilde\varepsilon_{h'}})_{n-1}\big)^{-1}\circ h_{n-1}\big) (m)\Big)\cr
&= \big((\kappa_{\tilde\varepsilon_{h'}})_{n}\circ d_{n-1}\big)\Big( 
\big((\kappa_{\tilde\varepsilon_{h'}})_{n-1}\big)^{-1}\circ h_{n-1}\big) (m)\Big)\cr}
\eqno £8.2.58\phantom{.}$$
which proves that $m$ belongs to ${\rm Im}(d_{n-1})\,.$ We are done.

\medskip
£8.3\phantom{.} Let us say that a $\ch^*(\F^{^{\rm sc}})\-$object $(\frak q,\Delta_n)$ is
{\it regular\/} if $\frak q (i-1,i)$ is {\it not\/} an isomorphism for any $1\le i\le n$
[5,~A5.2]; note that there is a canonical bijection between  a set of representatives 
for the set of isomorphism classes of regular $\ch^*(\F^{^{\rm sc}})\-$objects and a set of
representatives for the set of $\F\-$isomorphism classes of nonempty sets of
$\F\-$selfcentralizing subgroups of~$P$, which are totally ordered by the inclusion.

\bigskip
\noindent
{\bf Corollary~£8.4}\phantom{.} {\it With the notation above, we have 
$${\rm rank}_{\O}\big(\G_\K (\F,\widehat\aut_{\F^{^{\rm sc}}})\big) = \sum_{(\frak
q,\Delta_n)} (-1)^{n}\,{\rm rank}_{\O} \Big(\G_\K \big(\hat\L
(\frak q)\big)\Big)
\eqno £8.4.1\phantom{.}$$
where $(\frak q,\Delta_n)$ runs over a set of representatives for the set of isomorphism
classes of regular $\ch^*(\F^{^{\rm sc}})\-$objects.\/}

\medskip
\noindent
{\bf Proof:} It follows from the decomposition~£5.12.3 and from Proposition~£6.5 that we have
$${}^\K \G_\K (\F,\widehat\aut_{\F^{^{\rm sc}}}) \cong
\prod_{h\in \Bbb N -\{0\}} \Bbb H^0 \big(\,
{}^{^h}\!(\F^{^{\rm sc}}),{}^\K \F\!ct_{U_{h'}}\circ \frak t_h\big)
\eqno £8.4.2\phantom{.}$$
where we set $h' = h^{p'}\,.$ On the other hand,  it is more or less well-known that the cohomology groups over functors
to $\K\-$mod coincide with the corresponding cohomology groups computed from the {\it regular chains\/}; more precisely, in our situation let us denote  by ${}^\K\frak n_h$ the extension of~$\frak n_h = 
 \F\!ct_{U_{h'}}\circ \tilde\frak t_h$ from $\O$ to $\K$ for any  $h\in \Bbb N -\{0\}\,;$ then, it follows from Propositions~A4.13 and~A5.7 in [5] that, with the notation there, for any $n\in \Bbb N $ we have
$$\Bbb H^n \big( \widetilde{{}^{^h}\! (\F^{^{\rm sc}}}),{}^\K\frak n_h \big)
\cong \Bbb H^n_* \big( \widetilde{{}^{^h}\! (\F^{^{\rm sc}}}),{}^\K\frak n_h\big)\cong 
\Bbb H^n_{\rm r} \big( \widetilde{{}^{^h}\! (\F^{^{\rm sc}}}),{}^\K\frak n_h\big)
\eqno £8.4.3.$$

\smallskip
In particular, it follows from Theorem~£8.2 that for any $n\ge 1$ we get
$$\Bbb H^n_{\rm r} \big(\widetilde{ {}^{^h}\! (\F^{^{\rm sc}}}),{}^\K\frak n_h\big) = \{0\}
\eqno £8.4.4\phantom{.}$$
which amounts to saying that we have an infinite exact sequence
$$0\too \Bbb H^0_{\rm r} \big(\widetilde{{}^{^h}\! (\F^{^{\rm sc}}}),{}^\K\frak n_h\big)
\too \dots\too C^n\to C^{n+1}\too \dots
\eqno £8.4.5\phantom{.}$$
where $C^n = \Bbb C^n_{\rm r}\big(\widetilde{ {}^{^h}\! (\F^{^{\rm sc}}}),{}^\K\frak n_h\big)$
is the set of elements (cf.~Proposition~£6.4)
$$(\chi_{\widetilde{\frak q^{\bar\eta}}})_{\widetilde{\frak q^{\bar\eta}}\in \Fct_{\rm r} 
(\Delta_n, \widetilde{{}^{^h}\! (\F^{^{\rm sc}}}))} \in \prod_{\widetilde{\frak q^{\bar\eta}}\in \Fct_{\rm r} 
(\Delta_n, \widetilde{{}^{^h}\! (\F^{^{\rm sc}}}))} {}^\K\frak n_h\big(\widetilde{\frak q^{\bar\eta}} (0)\big)
\eqno £8.4.6\phantom{.}$$
such that, for any natural isomorphism $\tilde\nu\,\colon \widetilde{\frak q^{\bar\eta}}\cong   
\widetilde{\frak r^{\bar\theta}}$ between  {\it regular $\widetilde{{}^{^h}\! (\F^{^{\rm sc}}})\-$va-lued $n\-$chains\/} $\widetilde{\frak q^{\bar\eta}}$  and~$\widetilde{\frak r^{\bar\theta}}\,,$ ${}^\K\frak n_h (\tilde\nu_0)$ maps  
$\chi_{\widetilde{\frak r^{\bar\theta}}}$ on $\chi_{\widetilde{\frak q^{\bar\eta}}}\,.$
That is to say,  since  we have a bijection between the sets of isomorphism classes of   ${}^{^h}\! (\F^{^{\rm sc}})\-$
and  $\widetilde{ {}^{^h}\! (\F^{^{\rm sc}}})\-$objects, we actually have
$$\Bbb C^n_{\rm r}\big(\widetilde{ {}^{^h}\! (\F^{^{\rm sc}}}),{}^\K\frak n_h\big)\cong 
\prod_{\frak q^{\bar\eta}} {}^\K \F\!ct_{U_{h'}} 
\Big(\frak t_h \big(\frak q^{\bar\eta} (0)\big)\Big)^{\F(\frak q)_{\bar\eta}}
\eqno £8.4.7\phantom{.}$$
where $\frak q^{\bar\eta}$ runs over a set of representatives for the set of isomorphism
classes in~$\Fct_{\rm r} \big(\Delta_n, {}^{^h}\! (\F^{^{\rm sc}})\big)$
[5,~A5.3] and then $\F(\frak q)_{\bar\eta}$ denotes the stabilizer of~$\bar\eta$ in~$\F(\frak q)$.

\smallskip
On the other hand, it is clear that for $n$ big enough there are no {\it regular 
$\widetilde{{}^{^h}\! (\F^{^{\rm sc}}})\-$valued $n\-$chains\/} and therefore, in the exact sequence above only finitely many terms are not zero; thus, we still get
$$\eqalign{{\rm dim}_\K \Big(\Bbb H^0 &\big( \widetilde{{}^{^h}\! (\F^{^{\rm sc}}}),
{}^\K\frak n_h\big)\Big)\cr
&= \sum_{(\frak q^{\bar\eta},\Delta_n)} (-1)^{n}\,{\rm dim\,}_\K \bigg({}^\K \F\!ct_{U_{h'}} 
\Big(\frak t_h \big(\frak q^{\bar\eta} (0)\big)\Big)^{\F(\frak q)_{\bar\eta}}\bigg)\cr}
\eqno £8.4.8\phantom{.}$$
where $(\frak q^{\bar\eta},\Delta_n)$ runs over a set of representatives for the set of isomorphism classes 
of {\it regular\/} $\ch^*\big({}^{^h}\! (\F^{^{\rm sc}})\big)\-$objects [5,~A5.3]. Consequently, since the functor 
$\hat\frak w_h$ in Proposition~£5.8 maps the $k^*\-\frak i\widetilde\Loc\-$morphism
$$\skew2\hat{\bar\iota}_0^{\,\frak q} : \big(\hat\L (\frak q),{\rm Ker} (\pi_\frak q)\big)\too 
\Big(\hat\L\big(\frak q (0)\big),Z\big(\frak q(0)\big)\Big)
\eqno 8.4.9,$$
which lifts the corresponding $\frak i\widetilde\Loc\-$morphism~£6.3.4, on a $U_{h'}\-$set bijection (cf.~£5.8.2)
$$(\varpi_{h,\hat\L (\tilde\frak q (0))})^{-1} (\bar\iota^{\frak q}_0\circ \bar\eta) \cong 
(\varpi_{h,\hat\L (\tilde\frak q)})^{-1}(\bar\eta)
\eqno £8.4.10,$$
where, for short, we write $\hat\L\big(\frak q (0)\big)$ and $\hat\L (\tilde\frak q)$ instead of
$\Big(\hat\L\big(\frak q (0)\big),Z\big(\frak q(0)\big)\Big)$ and $\big(\hat\L (\frak q),{\rm Ker} 
(\pi_\frak q)\big)\,,$  it follows from equality~£8.4.2 and from Propositions~£6.4 and~£6.5 that we actually have
$$\eqalign{{\rm rank}_{\O}&\big(\G_\K (\F,\widehat\aut_{\F^{^{\rm sc}}})\big)\cr
& = \sum_{h} \sum_{(\bar\eta,\frak q,\Delta_n)}\!\! (-1)^{n} {\rm dim\,}_\K \Big({}^\K \F\!ct_{U_{h'}}  \big((\varpi_{h,\hat\L (\frak q)})^{-1} (\bar\eta)\big)^{\F(\frak q)_{\bar\eta}}\Big)\cr}
\eqno £8.4.11\phantom{.}$$
where $h$ runs over $\Bbb N -\{0\}$ and $(\bar\eta,\frak q,\Delta_n)$ over a set of
representatives for the isomorphism classes of ${}^{\frak u_h} 
\frak l\ch^* (\F^{^{\rm sc}})\-$objects such that $(\frak q,\Delta_n)$ is  regular.

\smallskip
On the other hand, for any $\ch^*(\F^{^{\rm sc}})\-$object $(\frak
q,\Delta_n)\,,$ it follows from isomorphism~£5.4.3 that
$$\eqalign{\bigoplus_{h\in \Bbb N -\{0\}} &\bigoplus_{\bar\eta }\,{}^\K \F\!ct_{U_{h'}} 
\big((\varpi_{h,\hat\L (\frak q)})^{-1}(\bar\eta)\big)^{\F(\frak q)_{\bar\eta}}\cr 
&\cong \bigoplus_{h\in \Bbb N -\{0\}}  \Big(\bigoplus_{\eta\in 
{\rm Mon}(U_h,\L(\frak q)) }\,{}^\K \F\!ct_{U_{h'}} 
\big((\varpi_{h,\hat\L (\frak q)})^{-1}(\eta)\big)\Big)^{\F(\frak q)}\cr
&\cong  {}^\K\G_K\big(\hat\L (\frak q)\big)\cr}
\eqno £8.4.12\phantom{.}$$
where $\bar\eta$ runs over a set of representatives for the orbits of 
$\L(\frak q)$ on the set $\overline{\rm Mon}\big(U_h,\L(\frak q)\big)\,.$ We
are done.

 \bigskip
\bigskip
\noindent
{\bf £9\phantom{.} \bf General decomposition maps in a folded Frobenius $P\-$category}
\bigskip

£9.1\phantom{.} With the notation of \S8, let us choose a set of representatives $\P\i P$ for the set of $\F\-$isomorphism classes of the elements of $P$ in such a way that, for any $u\in \P\,,$ the subgroup $\langle u\rangle$ is fully centralized in $\F$ [5,~Proposition~2.7]. For any $u\in \P\,,$  we have the Frobenius $C_P (u)\-$category 
$C_\F (u)$ [5,~Proposition~2.16] and we know that the inclusion $\iota_u\,\colon C_P (u)\to P$ is $(C_\F (u),\F)\-$functorial [5,~12.1]; since a $C_\F (u)\-$selfcentralizing subgroup of~$C_P (u)$ contains~$u\,,$ it is also
 an $\F\-$selfcentralizing subgroup of~$P\,;$ co-herently, we write $C_{\F^{^{\rm sc}}} (u)$
 instead of $C_{\F} (u)^{^{\rm sc}}$ and, according to~£4.9 above, we have the $\O\-$module homomorphism
 $${\rm Res}_{\iota_u} : \G_\K (\F,\widehat\aut_{\F^{^{\rm sc}}}) \too 
 \G_\K \big(C_\F (u),\widehat{\aut}_{C_{\F^{^{\rm sc}}} (u)}\big)
  \eqno £9.1.1.$$

  \medskip
£9.2\phantom{.} Following Brou\'e [2, Appendice], for any finite group $G$ and any central $p\-$element $z$
of $G\,,$ we consider the {\it $z\-$twist\/} 
$$\omega_G^z :  \G_\K (G)\cong \G_\K (G)
\eqno £9.2.1\phantom{.}$$
determined by the {\it $z\-$translation\/} map in the $\O\-$valued functions $\F\! ct(G,\O)$ induced by the multiplication by $z\,;$ explicitly, if $\chi$ is an {\it irreducible ordinary\break
\eject 
\noindent
 character\/} of $G\,,$ the 
{\it $z\-$translated function\/} maps $x\in G$ on $\chi(xz)$ and therefore it coincides with 
${\chi(z)\over \chi(1)}\.\chi$ which still belongs to the image of $\G_\K (G)\,;$ actually,this definition
can be easily extended to the {\it finite\/} $k^*\-$groupes [4,~Proposition~5.15]. Thus,  for any $u\in \P$
and any
$C_{\F^{^{\rm sc}}} (u)\-$chain 
  $\frak q\,\colon \Delta_n\to C_{\F^{^{\rm sc}}} (u)$ we can consider the 
  {\it $u\-$twist\/} 
  $$\omega^u_\frak q : \G_\K\big(\widehat{\loc}_{C_{\F^{^{\rm sc}}} (u)}(\frak q)\big)\cong \G_\K\big(\widehat{\loc}_{C_{\F^{^{\rm sc}}} (u)}(\frak q)\big)
   \eqno £9.2.2;$$
then, it is clear that we get a 
 {\it natural automorphism\/}
 $$\omega^u : \frak g_\K\circ \widehat{\loc}_{C_{\F^{^{\rm sc}}} (u)}\cong
 \frak g_\K\circ \widehat{\loc}_{C_{\F^{^{\rm sc}}} (u)}
 \eqno £9.2.3\phantom{.}$$
and therefore an $\O\-$module automorphism  
$$\Omega^u : \G_\K \big(C_\F (u),\widehat{\aut}_{C_{\F^{^{\rm sc}}} (u)}\big)\cong  \G_\K \big(C_\F (u),\widehat{\aut}_{C_{\F^{^{\rm sc}}} (u)}\big)
 \eqno £9.2.4.$$
Finally, we can define the {\it $u\-$general decomposition map\/}
$$\partial^u_{(\F,\widehat\aut_{\F^{^{\rm sc}}})} : 
\G_\K (\F,\widehat\aut_{\F^{^{\rm sc}}}) \too \G_k \big(C_\F (u),
\widehat{\aut}_{C_{\F^{^{\rm sc}}} (u)}\big)
\eqno £9.2.5\phantom{.}$$
as the composition $\partial_{(C_\F (u),\widehat{\aut}_{C_{\F^{^{\rm sc}}} (u)})}
\circ \Omega^u \circ {\rm Res}_{\iota_u}$ (cf.~£3.4.2).

\bigskip
\noindent
{\bf Theorem~£9.3}\phantom{.} {\it The family of general decomposition maps $\{\partial^u_{(\F,
\widehat\aut_{\F^{^{\rm sc}}})}\}_{u\in \P}$ determines a $\K\-$module isomorphism
$${}^\K \G_\K (\F,\widehat\aut_{\F^{^{\rm sc}}})\cong 
\bigoplus_{u\in \P} {}^\K \G_k \big(C_\F (u), \widehat{\aut}_{C_{\F^{^{\rm sc}}} (u)}\big)
\eqno £9.3.1.$$\/}

\par
\noindent
{\bf Proof:} According to decomposition~£5.12.3, to Proposition~£6.5 and to isomorphism~£6.10.3, setting $h' =
h^{p'}$ for any $h\in \Bbb N-\{0\}\,,$ we have an injective $\O\-$module  homomorphism
$$\G_\K (\F,\widehat\aut_{\F^{^{\rm sc}}})\too 
\prod_{h\in \Bbb N-\{0\}}\,\prod_{Q^{\rho'}} \F\!ct_{U_{h'}}\big( (\frak s_{h'}\times \frak x_{h}) (Q^{\rho'}),\K\big)
\eqno £9.3.2\phantom{.}$$
 where $Q^{\rho'}$ runs over a set of representatives for the set of 
 ${}^{h'}\! (\tilde\F^{^{\rm sc}})\-$isomorphism classes  of ${}^{h'}\! (\tilde\F^{^{\rm sc}})\-$objects; thus, in order to prove the injectivity of the map determined by the family of {\it general decomposition maps\/}, it suffices to prove that, for any 
$X\in \G_\K (\F,\widehat\aut_{\F^{^{\rm sc}}})$ in the kernel of this map,
any  $h\in \Bbb N -\{0\}$ and any ${}^{h'}\! (\tilde\F^{^{\rm sc}})\-$object 
$Q^{\rho'}\,,$ for the chosen  lifting $\widehat{\rho'(U_{h'})}$ to $\F (Q)$ of $\rho'(U_{h'})$ (cf.~£6.9)
the corresponding projection 
$$\chi_{X,h,Q^{\rho'}} : (\varpi_{h',\skew4\hat{\tilde\F}(Q)})^{-1}(\rho')
\times \widetilde{\rm Mon}(U_{h^p},Q^{\widehat{\rho'(U_{h'})}})\too \K
\eqno £9.3.3\phantom{.}$$
 is the zero function, namely that $\chi_{X,h,Q^{\rho'}}$ vanish
 over $(\varpi_{h',\skew4\hat{\tilde\F}(Q)})^{-1}(\rho')\times \{\tilde\rho''\}$
 for any injective group homomorphism $\rho''\,\colon U_{h^p}\to Q^{\widehat{\rho'(U_{h'})}}\,.$
 \eject

 \smallskip
Setting $u = \rho''(\xi_{h^p})\,,$ we may assume that $u$ belongs to $\P$
and then we are actually assuming that $\partial^u_{(\F,\widehat\aut_{\F^{^{\rm sc}}})} 
(X) = 0\,;$ consider the subgroup $C_Q (u)$ which is clearly selfcentralizing  in $Q$
and, for a suitable $n\in \Bbb N$ and any $i\in \Delta_{n-1}\,,$ set $R_n = C_Q (u)$ and $R_{i}= N_Q (R_{i+1})$ in such a way that $R_0 = Q\,;$ it is clear that the lifting
$\sigma\in \widehat{\rho'(U_{h'})}$ of $\rho' (\xi_{h'})$ stabilizes the family 
$\{R_i\}_{i\in \Delta_n}$ and then it follows from Lemma~£9.4 below that
 there are a family $\{Q_i\}_{i\in \Delta_n}$ of $\F\-$selfcentralizing subgroups of $P$
such that $Q_0 = Q\,,$ a family of $\F\-$morphisms $ \theta_{i}\,\colon Q_{i}
\to P\,,$  where $i$ runs over~$\Delta_{n-1}\,,$ such that $\theta_{i}(Q_{i})$ normalizes $Q_{i+1}\,,$ and  a family of $p'\-$elements 
 $$\sigma_{i}\in \F \big(\theta_{i}(Q_{i})\.Q_{i+1}\big)
 \eqno £9.3.4,$$
  where $i$ runs over $\Delta_{n-1} \,,$ stabilizing $\theta_{i}(Q_{i})$ and
  $Q_{i+1}\,,$ in such a way that the action of $\sigma_i$ over $Q_{i}$ defined via 
  $\theta_i$ coincides  with  $\sigma$ if $i= 0$ or with the action of~$\sigma_{i-1}$ otherwise, and that, setting $T_{i,0 } = R_i$ for any $i\in \Delta_n$ and 
  then $T_{i,j+1}   = \theta_j (T_{i,j})$ for any $j < i\,,$ we have $Q_i = T_{i,i}\.C_{Q_i} (T_{i,i})$ for any~$i\in \Delta_n\,.$

\smallskip
In particular, for any $i\in \Delta_{n-1}\,,$  respectively denoting by 
$$\tilde\F \big(\theta_{i}(Q_{i})\.Q_{i+1}\big)_{\theta_{i} (Q_{i})} \qq 
\tilde\F \big(\theta_{i}(Q_{i})\.Q_{i+1}\big)_{Q_{i+1}}
 \eqno £9.3.5,$$
 the stabilizers of  $\theta_{i}(Q_{i})$ and  $Q_{i+1}$ in $\tilde\F \big(\theta_{i}(Q_{i})\.Q_{i+1}\big)\,,$ we already know that the restriction {\it via\/}~$\theta_i$ and the ordinary restriction respectively induce group homomorphisms
$$\eqalign{\alpha_i : &\,\tilde\F \big(\theta_{i}(Q_{i})\.Q_{i+1}\big)_{\theta_{i} (Q_{i})}\too \tilde\F (Q_{i}) \cr 
\beta_i : &\,\tilde\F \big(\theta_{i}(Q_{i})\.Q_{i+1}\big)_{Q_{i+1}}\too 
\tilde\F (Q_{i+1})\cr}
\eqno £9.3.6\phantom{.}$$
such that their kernels are $p\-$groups [5,~Corollary~4.7]; hence, arguing by induction, the order of~$\sigma$ coincides with the order of $\sigma_i$ for any $i\in\Delta_{n-1}\,,$ and 
therefore we have an injective group homomorphism
$$\rho'_i : U_{h'}\too \tilde\F \big(\theta_{i}(Q_{i})\.Q_{i+1}\big)
\eqno £9.3.7\phantom{.}$$
mapping $\xi_{h'}$ on $\tilde\sigma_i\,,$ so that we get 
$$\alpha_{i}\circ \rho'_{i} = \cases{\beta_{i-1}\circ \rho'_{i-1} & if $i\not= 0$\cr
\rho'& if $i = 0$\cr}
\eqno £9.3.8.$$

\smallskip
That is to say, setting $\eta'_0 = \rho'$ and $\eta'_{i+1}= \beta_i\circ \rho'_i $ for any $i\in \Delta_{n-1}\,,$ and  considering the ${}^{h'}\! 
(\tilde\F^{^{\rm sc}})\-$objects  $(Q_{i})^{\eta'_{i}}$ and $\big(\theta_{i}(Q_{i})\.Q_{i+1}\big)^{\rho'_i}\,,$ the homomorphism $\theta_i$ and the inclusion map
$Q_{i+1}\to \theta_{i}(Q_{i})\.Q_{i+1}$ respectively determine   
${}^{h'}\! (\tilde\F^{^{\rm sc}})\-$morphisms
$$(Q_i)^{\eta'_i}\too \big(\theta_{i}(Q_{i})\.Q_{i+1}\big)^{\rho'_i}\longleftarrow
(Q_{i+1})^{\eta'_{i+1}}
\eqno £9.3.9\phantom{.}$$
and therefore, for any $i\in \Delta_{n-1}\,,$ we get the $U_{h'}\-$set bijections (cf.~£6.8.2)
$$\matrix{(\varpi_{h',\skew4\hat{\tilde\F} (Q_i)})^{-1}(\eta'_i)\cr 
\wr\Vert\phantom{\big\uparrow}\cr
(\varpi_{h',\skew4\hat{\tilde\F} (\theta_{i}(Q_{i})\.Q_{i+1})})^{-1}(\rho'_i)\cr
\wr\Vert\phantom{\big\uparrow}\cr
(\varpi_{h',\skew4\hat{\tilde\F} (Q_{i+1})})^{-1}(\eta'_{i+1})\cr}
\eqno £9.3.10.$$
Analogously, we set $u_0 = u$ and inductively define $u_{i+1} = \theta_i (u_i)$ 
for any $i\in \Delta_{n-1}\,,$ which determines injective group homomorphisms
$$\eqalign{&\eta''_i : U_{h^p}\too (Q_i)^{\widehat{\eta'_i (U_{h'})}}\cr
&\rho''_i : U_{h^p}\too \big(\theta_{i}(Q_{i})\.Q_{i+1}\big)^{\widehat{\rho'_i (U_{h'})}}\cr}
\eqno £9.3.11\phantom{.}$$
for the corresponding chosen liftings $\widehat{\eta'_i (U_{h'})}$ and $\widehat{\rho'_i (U_{h'})}$
(cf.~£6.9).

\smallskip
More precisely, the functor $\frak s_{h'}\times \frak x_{h}$ maps the ${}^{h'}\! (\tilde\F^{^{\rm sc}})\-$morphisms £9.3.9 above on  $U_{h'}\-$set maps (cf.~£6.8~and~£6.9)
$$\matrix{(\varpi_{h',\skew4\hat{\tilde\F} (Q_i)})^{-1}(\eta'_i)
\times \widetilde{\rm Mon}\big(U_{h^p}, (Q_i)^{\widehat{\eta'_i (U_{h'})}}\big)\cr 
\downarrow\phantom{\big\uparrow}\cr
(\varpi_{h',\skew4\hat{\tilde\F} (\theta_{i}(Q_{i})\.Q_{i+1})})^{-1}(\rho'_i)
\times \widetilde{\rm Mon}\big(U_{h^p},\big(\theta_{i}(Q_{i})\.Q_{i+1}\big)^{\widehat{\rho'_i (U_{h'})}}\big)\cr
\uparrow\phantom{\big\uparrow}\cr
(\varpi_{h',\skew4\hat{\tilde\F} (Q_{i+1})})^{-1}(\eta'_{i+1})
\times \widetilde{\rm Mon}\big(U_{h^p}, (Q_{i+1})^{\widehat{\eta'_{i+1} (U_{h'})}}\big)\cr}
\eqno £9.3.12\phantom{.}$$
which send the $\K\-$valued function $\chi_{X,h,(\theta_{i}(Q_{i})\.Q_{i+1})^{\rho'_i}}$ to the $\K\-$valued functions $\chi_{X,h,(Q_i)^{\eta'_i }}$   and~$\chi_{X,h,(Q_{i+1})^{\eta'_{i+1} }}\,,$ and the functor 
${}^\K\F \!ct_{U_{h'}}$ actually determines the  $\K\-$module  isomorphisms
$$\matrix{\F \!ct_{U_{h'}} \big((\varpi_{h',\skew4\hat{\tilde\F} (Q_i)})^{-1}(\eta'_i)\times \{\tilde\eta''_i\},\K\big)\cr 
\wr\Vert\phantom{\big\uparrow}\cr
\F \!ct_{U_{h'}} \big((\varpi_{h',\skew4\hat{\tilde\F} (\theta_{i}(Q_{i})\.Q_{i+1})})^{-1}(\rho'_i)\times \{\tilde\rho''_i\},\K\big)\cr
\wr\Vert\phantom{\big\uparrow}\cr
\F \!ct_{U_{h'}}\big((\varpi_{h',\skew4\hat{\tilde\F} (Q_{i+1})})^{-1}
(\eta'_{i+1})\times \{\tilde\eta''_{i+1}\},\K\big)\cr}
\eqno £9.3.13.$$
In conclusion,  $\chi_{X,h,Q^{\rho'}}$ vanish
 over $(\varpi_{h',\skew4\hat{\tilde\F}(Q)})^{-1}(\rho')\times \{\tilde\rho''\}$
 if and only if $\chi_{X,h,(Q_n)^{\eta'_n}}$ vanish
 over $(\varpi_{h',\skew4\hat{\tilde\F}(Q_n)})^{-1}(\eta'_n)\times 
 \{\tilde\eta''_n\}\,.$

 \smallskip
 Finally, since $T_{n,0} = C_Q (u)$ and $Q_n = T_{n,n}\.C_{Q_n}(T_{n,n})\,,$
the element $u_n$ belongs to $Z(Q_n)\,;$ but, since $\langle u\rangle$ is fully centralized in $\F\,,$ there is an $\F\-$morphism $\theta_n\,\colon C_P(u_n)\to P$ fulfilling  $\theta_n(u_n) = u$ [5,~Proposition~2.7]; then, it is easily checked that $Q_{n+1} = \theta_n (Q_n)$ is a $C_\F (u)\-$selfcentralizing subgroup of~$C_P (u)$
and therefore, denoting by 
$$\eta'_{n+1} = {}^{\theta_n} \eta'_n : U_{h'}\too C_{\tilde\F} \big(\theta_n (Q_n)\big)
\eqno £9.3.14\phantom{.}$$
the corresponding action of $ U_{h'}$ on the group $Q_{n+1}\,,$ we have the 
$C_{\tilde\F^{^{\rm sc}}} (u)\-$ob-ject $(Q_{n+1})^{\eta'_{n+1}}$
and, denoting by $\frak s_{h'}^u$ and $\frak x_{h}^u$ the functors defined in~£6.8 and~£6.9 for the Frobenius
$C_P(u)\-$category $C_{\F} (u)\,,$ the corresponding projection map
$$\G_\K\big(C_\F (u),\widehat{\aut}_{C_{\F^{^{\rm sc}}} (u)}\big)\too
{}^\K\Fct_{U_{h'}}\big((\frak s_{h'}^u\times \frak x_{h}^u) \big((Q_{n+1})^{\eta'_{n+1}}\big)\Big)
\eqno £9.3.15,$$
sends ${\rm Res}_{\iota_u}(X)$ to following the $U_{h'}\-$set map 
$\chi_{_h} = \chi_{{\rm Res}_{\iota_u}(X),h,(Q_{n+1})^{\eta'_{n+1}}}$
$$(\varpi_{h',(C_{\skew4\hat{\tilde\F}}(u))(Q_{n+1})})^{-1}
(\eta'_{n+1}) \times \widetilde{\rm Mon}\big(U_{h^p},(Q_{n+1})^{\eta'_{n+1}(U_{h'})}\big)
\buildrel \chi_{_h}\over\too \K
\eqno £9.3.16\,.$$

\smallskip
Moreover, the element $u = \theta_n (u_n)$ determines an injective group homomorphism
$$\eta''_{n+1} : U_{h^p}\too (Q_{n+1})^{\eta'_{n+1}(U_{h'})}
\eqno £9.3.17\phantom{.}$$
and the condition $\partial^u_{(\F,\widehat\aut_{\F^{^{\rm sc}}})} (X) = 0$
clearly implies that the $\K\-$valued  function $ \chi_{_h}$ vanish over the $U_{h'}\-$set
$$(\varpi_{h',(C_{\skew4\hat{\tilde\F}}(u))(Q_{n+1})})^{-1}
(\eta'_{n+1}) \times \{\eta''_{n+1}\}
\eqno £9.3.18.$$
But, from the inclusion functor $C_\F (u)\to\F$ and from  the isomorphism $Q_{n+1} \cong Q_n$ determined by 
$\theta_n\,,$ we get an injective $k^*\-$group homomorphism
$$\big(C_{\skew4\hat{\tilde\F}}(u)\big)(Q_{n+1})\too \skew4\hat{\tilde\F} (Q_n)
\eqno £9.3.19\phantom{.}$$
inducing a $U_{h'}\-$set bijection [5,~Proposition~14.18]
$$(\varpi_{h',(C_{\skew4\hat{\tilde\F}}(u))(Q_{n+1})})^{-1}(\eta'_{n+1}) \cong (\varpi_{h',\skew4\hat{\tilde\F}(Q_n)})^{-1}(\eta'_n)
\eqno £9.3.20\phantom{.}$$
and it is clear that the above isomorphism sends $\eta''_{n+1}$ to $\eta''_{n}\,.$
At last, it is easily checked that the corresponding $U_{h'}\-$set bijection
$$(\varpi_{h',(C_{\skew4\hat{\tilde\F}}(u))(Q_{n+1})})^{-1}(\eta'_{n+1})
\times \{\eta''_{n+1}\}\cong (\varpi_{h',\skew4\hat{\tilde\F}(Q_n)})^{-1}(\eta'_n)\times \{\eta''_{n}\}
\eqno £9.3.21\phantom{.}\phantom{.}$$
sends the restriction of $\chi_{X,h,(Q_n)^{\eta'_n}}$ to the restriction of
$\chi_{_h}\,,$
so that   $\chi_{X,h,Q^{\rho'}}$ vanish
 over $(\varpi_{h',\skew4\hat{\tilde\F}(Q)})^{-1}(\rho')\times \{\tilde\rho''\}\,,$
 proving the injectivity in~£9.3.1.
 \eject

 \smallskip
 At this point, it suffices to prove that both members of isomorphism £9.3.1 have the same dimension, namely that the following equality holds
 $${\rm rank}_\O \big(\G_\K (\F,\widehat{\aut}_{\F^{^{\rm sc}}})\big) = 
 \sum_{u\in \P} {\rm rank}_\O \Big(\G_k \big(C_\F(u),\widehat{\aut}_{C_{\F^{^{\rm sc}}}(u)}\big)\Big)
 \eqno £9.3.22;$$
we already know that we have (cf.~Corollary~£8.4)
 $${\rm rank}_\O \big(\G_\K (\F,\widehat{\aut}_{\F^{^{\rm sc}}})\big) = 
 \sum_{(\frak q,\Delta_n)} (-1)^n\, {\rm rank}_\O\Big(\G_\K \big(\hat\L (\frak q)\big)\Big)
 \eqno £9.3.23\phantom{.}$$
 where $(\frak q,\Delta_n)$ runs over a set of representatives for the set of $\F\-$isomorphism classes of
$\ch^*_{\rm r}(\F^{^{\rm sc}})\-$ objects [5,~A5.3] which are {\it fully normalized\/} in~$\F$ [5,~2.18]. But, denoting by
$\P_{\!\frak q}$  a set of representatives  for the set of $N_\F(\frak q)\-$isomorphism
classes of elements of $N_P (\frak q)$ in such a way that, for any $u\in \P_\frak q\,,$ the subgroup 
$\langle u\rangle$ is fully centralized in $N_\F(\frak q)\,,$ and identifying $N_P (\frak q)$ with its structural image in 
$\hat\L (\frak q)\,,$ it easily follows from [5,~Proposition~19.5] that $N_\F(\frak q)$ coincides with the Frobenius category associated with $\hat\L (\frak q)$ [5,~1.8] and then, it is well-known that we have (cf. isomorphism~1.6.1)
$${\rm rank}_\O\Big(\G_\K \big(\hat\L (\frak q)\big)\Big) = \sum_{u\in \P_\frak q}
{\rm rank}_\O \Big(\G_k \big(C_{\hat\L (\frak q)}(u)\big)\Big)
\eqno £9.3.24.$$
On the other hand, for any $u\in \P\,,$ we also have [5,~Corollary~14.32]
$${\rm rank}_\O \Big(\G_k \big(C_\F(u),\widehat{\aut}_{C_{\F^{^{\rm sc}}}(u)}\big)\Big)= \sum_{(\frak q_u,\Delta_n)}  (-1)^n\,{\rm rank}_\O \bigg(\G_k \Big(\big(C_{\hat\F}(u)\big) (\frak q_u)\Big)\bigg)
\eqno £9.3.25\phantom{.}$$
where $(\frak q_u,\Delta_n)$ runs over a set of representatives for the set of $C_\F (u)\-$isomor-phism classes of
$\ch^*_{\rm r} \big(C_{\F^{^{\rm sc}}}(u)\big)\-$objects  [5,~A5.3]  which are {\it fully normalized\/} in~$C_\F (u)$ [5,~2.18] and $\big(C_{\hat\F}(u)\big)(\frak q_u)$ is the converse image of $\big(C_{\F}(u)\big)(\frak q_u)$
in~$\hat\F (\frak q_u)\,.$

\smallskip 
 Consequently, the left-hand member in £9.3.22 is equal to
$$\sum_{(\frak q,\Delta_n)}\,\sum_{u\in \P_\frak q} (-1)^n\, {\rm rank}_\O 
\Big(\G_k \big(C_{\hat\L (\frak q)}(u)\big)\Big)
\eqno £9.3.26\phantom{.}$$
where $(\frak q,\Delta_n)$ runs over a set of representatives for the set of $\F\-$isomorphism classes of
$\ch^*_{\rm r}(\F^{^{\rm sc}})\-$objects   which are {\it fully normalized\/} in~$\F$  [5,~2.18], whereas the right-hand  member is equal to
$$\sum_{u\in\P}\, \sum_{(\frak q_u,\Delta_n)}  (-1)^n\,{\rm rank}_\O \bigg(\G_k \Big(\big(C_{\hat\F}(u)\big) 
(\frak q_u)\Big)\bigg)
\eqno £9.3.27\phantom{.}$$
where, for any $u\in \P\,,$  $(\frak q_u,\Delta_n)$ runs over a set of representatives for the set of 
$C_\F (u)\-$isomorphism classes of $\ch^*_{\rm r} \big(C_{\F^{^{\rm sc}}}(u)\big)\-$objects which are 
{\it fully normalized\/} in~$C_\F (u)$ [5,~2.18]. But, it is clear that the element $u$ belongs to~$Z\big(\frak q_u (n)\big)$\break
\eject
\noindent
and then it is easily checked that  $\big(C_{\hat\F}(u)\big) (\frak q_u) $ coincides with the stabilizer of $u$
in $\hat\F (\frak q_u)\,,$ so that we have
$$\G_k \Big(\big(C_{\hat\F}(u)\big) 
(\frak q_u)\Big) = \G_k \big(C_{\hat\L (\frak q_u)}(u)\big)
\eqno £9.3.28\,.$$
Moreover, if $(\frak q,\Delta_n)$ is a $\ch^*_{\rm r}(\F^{^{\rm sc}})\-$object 
{\it fully normalized\/}  in $\F$ such that all the group homomorphisms $\frak q(j\bullet i)$ are {\it inclusion maps\/} then, for any
 element $u$ of~$Z\big(\frak q (n)\big)$ such that $\langle u\rangle$ is fully centralized in $\F\,,$  it is quite clear that $(\frak q,\Delta_n)$ remains a  $\ch^*_{\rm r} \big(C_{\F^{^{\rm sc}}}(u)\big)\-$object which is {\it fully normalized\/} in~$C_\F (u)\,.$

 \smallskip
 Hence, the sum~£9.3.27 above coincides with
$$\sum_{(\frak q,\Delta_n)}\,\sum_{u\in \Z_\frak q} (-1)^n \,{\rm rank}_\O 
\Big(\G_k \big(C_{\hat\L (\frak q)}(u)\big)\Big)
\eqno £9.3.29\phantom{.}$$
where $(\frak q,\Delta_n)$ runs over a set of representatives for the set of 
$\ch^*_{\rm r}(\F^{^{\rm sc}})\-$isomor-phism classes of
$\ch^*_{\rm r}(\F^{^{\rm sc}})\-$objects which are {\it fully normalized\/} in~$\F$ and, for such a $\ch^*_{\rm r}(\F^{^{\rm sc}})\-$object $(\frak q,\Delta_n)\,,$ $\Z_\frak q$ is a set of representatives for the set of
orbits of~$\F (\frak q)$ in $Z\big(\frak q(n)\big)\,.$ Finally, we may assume that $\P_\frak q$ 
contains $\Z_\frak q$ and then equality~£9.3.22 above is equivalent to the following one
$$0 =\sum_{(\frak q,\Delta_n)}\,\sum_{u\in \P_\frak q - \Z_\frak q} (-1)^n\, {\rm rank}_\O 
\Big(\G_k \big(C_{\hat\L (\frak q)}(u)\big)\Big)
\eqno £9.3.30\phantom{.}$$
where $(\frak q,\Delta_n)$ runs over a set of representatives for the set of 
$\ch^*_{\rm r}(\F^{^{\rm sc}})\-$isomor-phism classes of
$\ch^*_{\rm r}(\F^{^{\rm sc}})\-$objects which are {\it fully normalized\/} in~$\F\,.$

 \smallskip
Since any $\ch^*_{\rm r}(\F^{^{\rm sc}})\-$object is $\ch^*_{\rm r}(\F^{^{\rm sc}})\-$isomorphic to one which is {\it fully normalized\/} in~$\F\,,$ actually we are considering a set of representatives $\C$ for the~set of isomorphism classes of pairs formed by a  $\ch^*_{\rm r}(\F^{^{\rm sc}})\-$object $(\frak q,\Delta_n)$ such that
 all the group homomorphisms $\frak q(j\bullet i)$ are {\it inclusions\/},
 and an ele-ment $u\in P$ which normalizes but does not centralize $\frak q\,.$
Thus, it suffices to exhibit an {\it involutive permutation\/} $t$ of $\C$ such that, setting 
$$t\big((\frak q,\Delta_n),u\big) = \big((\frak q',\Delta_{n'}),u'\big)
\eqno £9.3.31$$
 and assuming that $(\frak q,\Delta_n)$ and $(\frak q',\Delta_n')$ are both {\it fully normalized\/} 
in~$\F\,,$ we have
$$\G_k\big(C_{\hat\L (\frak q')}(u')\big)\cong \G_k \big(C_{\hat\L (\frak q)}(u)\big)
\qq n'\not \equiv n \pmod 2
\eqno £9.3.32.$$
First of all, we consider the set $\C'$ of  pairs $\big((\frak q,\Delta_n),u\big)\in \C$ such that $u$ does not belong to $\frak q (0)\,;$ in this case, let $i$ be the last element of $\Delta_n$ such that $u$ does not belong to $\frak q (i)\,.$ If $i = n$ or  $ \frak q (i)\.\langle u\rangle\not= \frak q (i +1)$ then we set $n' = n+1$ and  consider the functor 
$\frak q'\,\colon \Delta_{n'}\to \F^{^{\rm sc}}$ mapping any $0\le \ell\le i$ on $\frak q (\ell)\,,$ $i +1$~on~$\frak q (i)\.\langle u\rangle\,,$ any $i +2\le \ell\le n'$ 
 on $\frak q (\ell -1)\,,$ and all  the $\Delta_{n'}\-$morphisms on the corresponding inclusions.
 \eject

\smallskip
 In this case, up to replacing the  $\ch^*_{\rm r}(\F^{^{\rm sc}})\-$object  $(\frak q',\Delta_{n'})$ 
 by a $\ch^*_{\rm r}(\F^{^{\rm sc}})\-$ isomorphic one {\it fully
normalized\/} in $\F\,,$ we get the $\ch^*_{\rm r}(\F^{^{\rm sc}})\-$morphism
$$(\nu,\delta^n_{i +1}) : (\frak q',\Delta_{n'})\too (\frak q,\Delta_{n})
\eqno £9.3.33\phantom{.}$$
for a suitable {\it natural isomorphism\/} $\nu\,\colon \frak q'\circ\delta^n_{i +1}
\cong \frak q\,,$ and thus we still get the $k^*\-\frak i\widetilde\Loc\-$morphism
$$\widehat{\loc}_{\F^{^{\rm sc}}}(\nu,\delta^n_{i +1}) : \hat\L (\frak q')\too 
\hat\L(\frak q)
\eqno £9.3.34.$$
Moreover, since we have [5,~14.8]
$$\big(\aut_{\tilde\F^{^{\rm sc}}}(\nu,\delta^n_{i +1})\big)
\big(\tilde\F (\frak q')\big) = \tilde\F (\frak q)_{\frak q'(i +1)}
\eqno £9.3.35\phantom{.} $$
where $\tilde\F (\frak q)_{\frak q'(i +1)}$ denotes the set of $\tilde\sigma\in 
\tilde\F (\frak q)$ fulfilling
$$\tilde\iota_{\frak q (0)}^{\,\frak q' (i +1)}\circ\tilde\sigma_0 \in \tilde\F \big(\frak q'
(i +1)\big)\circ \tilde\iota_{\frak q (0)}^{\,\frak q' (i +1)}
\eqno £9.3.36,$$
which clearly contains the image of the stabilizer $\F (\frak q)_u$ of $u$ in 
$\F (\frak q)\,,$  a suitable representative of the exomorphism $\widehat{\loc}_{\F^{^{\rm sc}}}(\nu,\delta^n_{i +1})$ induces a $k^*\-$group isomorphism
$$C_{\hat\L (\frak q')}(u')\cong C_{\hat\L (\frak q)}(u)
\eqno £9.3.37\phantom{.}$$
where we are  setting $u' = (\nu_0)^{-1}(u)\,.$

\smallskip
If $i  +1\le  n$ and  $\frak q (i)\.\langle u\rangle = \frak q (i +1)\,,$  we~set $n' = n-1$ and consider  a 
chain $\frak q'\,\colon \Delta_{n'}\to \F^{^{\rm sc}}$ mapping any $0\le \ell \le i$ on  $\frak q (\ell)\,,$ any $i +1 \le \ell\le n'$  on~$\frak q (\ell +1)$ and, as before, all the  $\Delta_{n'}\-$morphisms on the corresponding inclusions.
 In this case, up to replacing the 
$\ch^*_{\rm r}(\F^{^{\rm sc}})\-$object  $(\frak q',\Delta_{n'})$ by a $\ch^*_{\rm r}(\F^{^{\rm sc}})\-$isomorphic one {\it fully
normalized\/} in $\F\,,$ we get the $\ch^*_{\rm r}(\F^{^{\rm sc}})\-$morphism
$$(\nu',\delta^{n'}_{i +1}) : (\frak q,\Delta_{n})\too (\frak q',\Delta_{n'})
\eqno £9.3.38\phantom{.}$$
for a suitable {\it natural isomorphism\/} $\nu'\,\colon \frak q\circ\delta^{n'}_{i +1}
\cong \frak q' \,,$ and thus we still get the $k^*\-\frak i\widetilde\Loc\-$morphism
$$\widehat{\loc}_{\F^{^{\rm sc}}}(\nu',\delta^{n'}_{i +1}) : \hat\L (\frak q)\too 
\hat\L(\frak q')
\eqno £9.3.39.$$
Moreover, since we have [5,~14.8]
$$\big(\aut_{\tilde\F^{^{\rm sc}}}(\nu',\delta^{n'}_{i +1})\big)
\big(\tilde\F (\frak q)\big) = \tilde\F (\frak q')_{\frak q(i +1)}
\eqno £9.3.40\phantom{.} $$
where as above $\tilde\F (\frak q')_{\frak q(i +1)}$ denotes the set of $\tilde\sigma\in 
\tilde\F (\frak q')$ fulfilling
$$\tilde\iota_{\frak q (0)}^{\,\frak q (i +1)}\circ\tilde\sigma_0 \in \tilde\F \big(\frak q
(i +1)\big)\circ \tilde\iota_{\frak q (0)}^{\,\frak q (i +1)}
\eqno £9.3.41,$$
which clearly contains the image of the stabilizer $\F (\frak q')_u$ of $u$ in 
$\F (\frak q')\,,$ a suitable representative of the $k^*\-\frak i\widetilde\Loc\-$morphism
 $\widehat{\loc}_{\F^{^{\rm sc}}}(\nu',\delta^{n'}_{i +1})$ induces a $k^*\-$group isomorphism
$$C_{\hat\L (\frak q)}(u)\cong C_{\hat\L (\frak q')}(u')
\eqno £9.3.42\phantom{.}$$
where we are   setting $u' = \nu'_0 (u)\,.$ Finally, since in both cases $\frak q' (0) 
= \frak q (0)\,,$ we may assume that the pair $\big((\frak q',\Delta_{n'}),u\big)$ still belongs to $\C'$ and, defining $t\big((\frak q,\Delta_n),u\big) = \big((\frak q',\Delta_{n'}),u\big)\,,$ {\it mutatis mutandis\/} it is easily checked that
$$t\big((\frak q',\Delta_{n'}),u\big)  = \big((\frak q,\Delta_{n}),u\big)
\eqno £9.3.43.$$

\smallskip
 From now on, we consider the set $\C''$ of  pairs $\big((\frak q,\Delta_n),u\big)\in \C$ such that  $u$ belongs to $\frak q(0)\,;$ 
note that the product 
$$R' = C_{\frak q (0)} (u)\.\big[\frak q(n),C_{\frak q (0)} (u)\big]
\eqno £9.3.44\phantom{.}$$
 is a normal subgroup of $\frak q (n)$ and therefore it follows from [5, Proposition~2.7] that, up to replacing
the pair $\big((\frak q,\Delta_{n}),u\big)$ by its image throughout a suitable $\F\-$morphism 
$\varphi \,\colon q (n)\to P\,,$ we may assume that $R'$ is fully normalized and fully centralized in $\F\,;$ then,  
$Q' = R'\.C_P (R')$ is $\F\-$selfcentralizing [5,~4.10]. Since we may assume that $Q'\not= \{1\}\,,$ setting 
$\frak q (-1) = \{1\}\,,$ let~$i$ be the last element in $\Delta_n\cup \{-1\}$ such that $Q'\not\i \frak q (i)\,,$ so that we have $\frak q (i)\not= Q'\.\frak q (i)\,.$ If $i = n$ or 
$Q'\.\frak q (i)\not= \frak q (i +1)$ then we set $n' = n+1$ and consider
the $\F^{^{\rm sc}}\-$chain $\frak q'\,\colon \Delta_{n'}\to \F^{^{\rm sc}}$ mapping any $0\le \ell\le i$ 
on~$\frak q (\ell)\,,$ $i +1$ on  $Q'\.\frak q (i)\,,$ any $i +2\le \ell\le n'$ 
 on $\frak q (\ell -1)$ and all  the $\Delta_{n'}\-$morphisms on the corresponding inclusions.

\smallskip
As above, up to replacing the $\ch^*_{\rm r}(\F^{^{\rm sc}})\-$object  $(\frak q',\Delta_{n'})$ by a $\ch^*_{\rm r}(\F^{^{\rm sc}})\-$iso-morphic one {\it fully
normalized\/} in $\F\,,$ we get the $\ch^*_{\rm r}(\F^{^{\rm sc}})\-$morphism
$$(\nu,\delta^n_{i +1}) : (\frak q',\Delta_{n'})\too (\frak q,\Delta_{n})
\eqno £9.3.45\phantom{.}$$
for a suitable {\it natural isomorphism\/} $\nu\,\colon \frak q'\circ\delta^n_{i +1}
\cong \frak q \,,$ and thus we still get the $k^*\-\frak i\widetilde\Loc\-$morphism
$$\widehat{\loc}_{\F^{^{\rm sc}}}(\nu,\delta^n_{i +1}) : \hat\L (\frak q')\too 
\hat\L(\frak q)
\eqno £9.3.46.$$
Moreover, as above we have [5,~14.8]
$$\big(\aut_{\F^{^{\rm sc}}}(\nu,\delta^n_{i +1})\big)
\big(\F (\frak q')\big) = \F (\frak q)_{Q'\.\frak q (i)}
\eqno £9.3.47\phantom{.} $$
where $\F (\frak q)_{Q'\.\frak q (i)}$ denotes the set of $\sigma\in \F (\frak q)$ fulfilling
$$\tilde\iota_{\frak q (0)}^{\,Q'\.\frak q (i)}\circ\tilde\sigma_0 \in \tilde\F \big(Q'\.\frak q (i)\big)\circ \tilde\iota_{\frak q (0)}^{\,Q'\.\frak q (i)}
\eqno £9.3.48;$$
actually, we may assume that $\F_P (\frak q)$ contains a Sylow $p\-$subgroup of 
$\F (\frak q)_u$ and then, denoting by  $C_{\F (\frak q)_u} (R')$ the kernel of the action of~$\F (\frak q)_u$ on~$R'\,,$ $\F_{Q'}(\frak q)$ is a Sylow $p\-$subgroup of the product $\F_{Q'}(\frak q)\.C_{\F (\frak q)_u} (R')\,;$ thus, by the {\it Frattini argument\/} we get
 $$\F (\frak q)_u \i C_{\F (\frak q)_u} (R')\.\F (\frak q)_{Q'\.\frak q (i)}
 \eqno £9.3.49.$$
 \eject
 \noindent
On the other hand, since we have $C_{\frak q (0)}(R')\i R'\,,$ a $p'\-$subgroup of $\F (\frak q)_u$ which acts trivially on $R'$ is necessarily trivial [3, Ch. 5, Theorem~3.4], so that $C_{\F (\frak q)_u} (R')$ is
a $p\-$group. Consequently, setting $u' = (\nu_0)^{-1}(u)\,,$ the $k^*\-\frak i\widetilde\Loc\-$morphism
$\widehat{\loc}_{\F^{^{\rm sc}}}(\nu,\delta^n_{i +1})$ induces a group isomorphism
$$\G_k \big(C_{\hat\L (\frak q')}(u')\big)\cong \G_k \big(C_{\hat\L (\frak q)}(u)\big)
\eqno £9.3.50.$$

 \smallskip
Finally, if $i +1\le  n$ and  $Q'\.\frak q (i) = \frak q (i +1)\,,$  we set $n' = n-1$ and consider the 
$\F^{^{\rm sc}}\-$chain $\frak q'\,\colon \Delta_{n'}\to \F^{^{\rm sc}}$ mapping any $0\le \ell \le i$ on $\frak q (\ell)\,,$ any $i +1 \le \ell\le n'$  on 
$ \frak q (\ell +1)$ and, as before, all the $\Delta_{n'}\-$morphisms on the corresponding inclusions. Once again, up to replacing the  $\ch^*_{\rm r}(\F^{^{\rm sc}})\-$object  
 $(\frak q',\Delta_{n'})$ by a $\ch^*_{\rm r}(\F^{^{\rm sc}})\-$isomorphic one {\it fully
normalized\/} in $\F\,,$ we get the $\ch^*_{\rm r}(\F^{^{\rm sc}})\-$morphism
$$(\nu',\delta^{n'}_{i +1}) : (\frak q,\Delta_{n})\too (\frak q',\Delta_{n'})
\eqno £9.3.51\phantom{.}$$
for a suitable {\it natural isomorphism\/} $\nu'\,\colon \frak q\circ\delta^{n'}_{i +1}
\cong \frak q' \,,$ and thus we still get the $k^*\-\frak i\widetilde\Loc\-$morphism
$$\widehat{\loc}_{\F^{^{\rm sc}}}(\nu',\delta^{n'}_{i +1}) : \hat\L (\frak q)\too 
\hat\L(\frak q')
\eqno £9.3.52;$$
as above, we have [5,~14.8]
$$\big(\aut_{\F^{^{\rm sc}}}(\nu',\delta^{n'}_{i +1})\big)
\big(\F (\frak q)\big) = \F (\frak q')_{Q'\.\frak q (i)}
\eqno £9.3.53\phantom{.} $$
where $\F (\frak q')_{Q'\.\frak q (i)}$ denotes the set of $\sigma\in \F (\frak q')$ fulfilling
$$\tilde\iota_{\frak q (0)}^{\,Q'\.\frak q (i)}\circ (\nu_0)^{-1}\circ\tilde\sigma_0 \in \tilde\F \big(Q'\.\frak q (i)\big)\circ \tilde\iota_{\frak q (0)}^{\,Q'\.\frak q (i)}\circ (\nu_0)^{-1}
\eqno £9.3.54.$$

\smallskip
Moreover, setting $u' = \nu_0 (u)\,,$ we may assume that $\F_P (\frak q')$ contains a Sylow $p\-$subgroup 
of~$\F (\frak q')_{u'}$ and then, denoting by  $C_{\F (\frak q')_{u'}} (R')$ the kernel of the action of~$\F (\frak q')_{u'}$ on~$R'$ {\it via\/} 
the group isomorphism $\nu'_{n'}\,\colon
\frak q (n)\cong \frak q' (n')\,,$ $\F_{\nu'_{n'}(Q')}(\frak q')$ is a Sylow $p\-$subgroup of the product $\F_{\nu'_{n'}(Q')}(\frak q')\.C_{\F (\frak q')_{u'}} (R')\,;$ once again, by the {\it Frattini argument\/} we get
 $$\F (\frak q')_{u'} \i C_{\F (\frak q')_{u'}} (R')\.\F (\frak q')_{Q'\.\frak q (i)}
 \eqno £9.3.55.$$
On the other hand, since we have $C_{\frak q (0)}(R')\i R'\,,$ a $p'\-$subgroup of 
$\F (\frak q')_{u'}$ which acts trivially on $R'$ is necessarily trivial [3, Ch. 5, Theorem~3.4], so that 
$C_{\F (\frak q')_{u'}} (R')$ is
a $p\-$group. Consequently, the $k^*\-\frak i\widetilde\Loc\-$morphism
 $\widehat{\loc}_{\F^{^{\rm sc}}}(\nu,\delta^n_{i +1})$ induces a group isomorphism
$$\G_k \big(C_{\hat\L (\frak q')}(u')\big)\cong \G_k \big(C_{\hat\L (\frak q)}(u)\big)
\eqno £9.3.56.$$

\smallskip
Since $u$ belongs to $R'\,,$ in both cases $u$ belongs to   $\frak q' (0)$ and we may assume that the pair 
$\big((\frak q',\Delta_{n'}),u\big)$ still belongs to $\C''\,;$ in this situation, we define $t\big((\frak q,\Delta_n),u\big) = \big((\frak q',\Delta_{n'}),u\big)$ and claim that 
 $$t\big((\frak q',\Delta_{n'}),u\big) = \big((\frak q,\Delta_n),u\big)
 \eqno £9.3.57.$$
 \eject
 \noindent
Indeed, set $t\big((\frak q',\Delta_{n'}),u\big) = \big((\frak q'',\Delta_{n''}),u\big)\,;$
if $\frak q' (n')  = \frak q (n)$ then {\it mutatis mutandis\/} we consider
$$R'' = C_{\frak q' (0)} (u)\.\big[\frak q'(n'),C_{\frak q' (0)} (u)\big]
\eqno £9.3.58;$$
since $C_{\frak q (0)}(u)\i R'\i Q'\,,$ we also  have $C_{\frak q' (0)}(u) = 
C_{\frak q (0)}(u)$ and therefore we get 
$$R'' = R'\qq Q'' = R''\.C_P (R'') = Q'
 \eqno £9.3.59;$$
moreover, once again setting $\frak q' (-1) = \{1\}\,,$ $i$ is also the last element in $\Delta_{n'}\cup \{-1\}$ such that 
$Q'$ is not contained in $\frak q' (i) = \frak q (i)\,;$ in this situation, the product $Q'\.\frak q' (i)$ is different from $\frak q' (n')  = \frak q (n)$ and equality~£9.3.57 is easily checked.

\smallskip
 If  $\frak q' (n')  \not= \frak q (n)$ then we have either $\frak q' (n') = Q'\.\frak q (n)$ or
$$\frak q (n) = Q'\.\frak q (n - 1) \qq \frak q' (n') = \frak q (n -1)
\eqno £9.3.60;$$
note that in both cases we have $\frak q' (0) = \frak q (0)\,;$ in the first case, we  have $\frak q' (n')  = C_P (R')\.\frak q (n)$ and, since $R'$ contains $C_{\frak q (0)}(u)
= C_{\frak q' (0)}(u)\,,$ we get
$$[\frak q' (n),C_{\frak q' (0)}(u)] = [\frak q (n),C_{\frak q (0)}(u)]
\eqno £9.3.61;$$
consequently, equalities~£9.3.59 still hold and therefore equality~£9.3.57 is easily checked;
in the second case, since $[\frak q (n),C_{\frak q (0)}(u)]$ is contained in the Frattini
subgroup of $\frak q (n)\,,$ we similarly obtain
$$\frak q (n) = C_P(R')\.C_{\frak q (0)}(u)\.\frak q (n -1) = C_P(R')\.\frak q' (n')
\eqno £9.3.62;$$
once again, we get
$$[\frak q (n),C_{\frak q (0)}(u)] = [\frak q' (n),C_{\frak q' (0)}(u)]
\eqno £9.3.63,$$
 equalities~£9.3.59 still hold and  equality~£9.3.57 is easily checked. We are done.

\bigskip
\noindent
{\bf Lemma~£9.4}\phantom{.} {\it Let $Q$ be an $\F\-$selfcentralizing subgroup
of $P$ and $\{R_i\}_{i\in \Delta_n}$ a family of selfcentralizing subgroups of $Q$
such that $R_0 = Q$ and $R_{i+1}\triangleleft R_{i}$ for~any $i\in \Delta_{n-1}\,.$
Then, for any $p'\-$element $\sigma$ in $\F (Q)$ stabilizing this family, there are
a family $\{Q_i\}_{i\in \Delta_n}$ of $\F\-$selfcentralizing subgroups of $P$
such that $Q_0 = Q\,,$ a family of $\F\-$morphisms $ \theta_{i}\,\colon Q_{i}\to P$  where $i$ runs over~$\Delta_{n-1}\,,$ such that $\theta_{i}(Q_{i})$ normalizes $Q_{i+1}\,,$ and  a family of $p'\-$elements 
 $$\sigma_{i}\in \F \big(\theta_{i}(Q_{i})\.Q_{i+1}\big)
 \eqno £9.4.1\phantom{.}$$
  where $i$ runs over $\Delta_{n-1} \,,$ stabilizing $\theta_{i}(Q_{i})$ and
  $Q_{i+1}\,,$ in such a way that the action of $\sigma_i$ over $Q_{i}$ defined via 
  $\theta_i$ coincides  with  $\sigma$ if $i= 0$ or with the action of~$\sigma_{i-1}$ otherwise, and that, 
  setting $T_{i,0 } = R_i$ for any $i\in \Delta_n$ and then $T_{i,j+1}   = \theta_j (T_{i,j})$ for any $j < i\,,$ we have $Q_i = T_{i,i}\.C_{Q_i} (T_{i,i})$ for any~$i\in \Delta_n\,.$\/}
  \eject
\medskip
\noindent
{\bf Proof:} We argue by induction on $n$  and may assume that $n\not= 0\,;$
thus, we assume the existence of a family $\{Q_i\}_{i\in \Delta_{n-1}}$ of 
$\F\-$selfcentralizing subgroups of $P$ such that $Q_0 = Q\,,$ a family of 
$\F\-$morphisms $ \theta_{i}\,\colon Q_{i} \to P\,,$  where $i$ runs 
over~$\Delta_{n-2}\,,$ such that $\theta_{i}(Q_{i})$ normalizes $Q_{i+1}\,,$ and  a family of $p'\-$elements 
 $\sigma_{i}\in \F \big(\theta_{i}(Q_{i})\.Q_{i+1}\big)$   where $i$ runs over $\Delta_{n-2} \,,$ stabilizing 
 $\theta_{i}(Q_{i})$ and    $Q_{i+1}\,,$ which fulfill the corresponding conditions above; thus,  setting $T_{i,0 } = R_i$ for any $i\in \Delta_{n-1}$  and then $T_{i,j+1}   = \theta_j (T_{i,j})$ for any $j < i\,,$ the
 following diagram summarizes our situation
 $$\matrix{&\hskip-2pt Q &\hskip-4pt\buildrel \theta_0\over\too &\hskip-6pt{}^{\theta_0}Q\.Q_1\cr
&\hskip-27pt{\sigma\atop}\hskip-5pt\nearrow&\hskip-4pt{\sigma_0\atop}\hskip-5pt\nearrow\hskip-40pt&\hskip-6pt\triangledown \cr
Q&\hskip-6pt\buildrel \theta_0\over\too &\hskip-4pt{}^{\theta_0}Q\.Q_1&\hskip-6ptQ_1&\hskip-8pt\buildrel \theta_1\over\too &\hskip-10pt\dots\cr
&&\hskip-4pt\triangledown &\hskip-36pt\nearrow&&\cr
\triangledown &&\hskip-4pt Q_1&\hskip-6pt\buildrel \theta_1\over\too &\hskip-8pt\dots&&\hskip-12pt\hskip-15pt\dots\cr
&&\hskip-4pt\triangledown &&&&\hskip-27pt\triangledown\cr
R_1&\hskip-6pt\cong &\hskip-4ptT_{1,1}&&\dots\hskip-4pt&&\hskip-20ptQ_{n-2}&\hskip-34pt\buildrel \theta_{n-2}\over\too &\hskip-40pt{}^{\theta_{n-2}}Q_{n-2}\.Q_{n-1}\cr
\dots&&\hskip-4pt\dots&&\triangledown\hskip-6pt&&\hskip-62pt\nearrow&\hskip-10pt{\sigma_{n-2}\atop}\hskip-5pt\nearrow\hskip-20pt&\hskip-20pt\triangledown\cr
&&&&&\hskip-30pt Q_{n-2}&\hskip-20pt\buildrel \theta_{n-2}\over\too &\hskip-10pt{}^{\theta_{n-2}}Q_{n-2}\.Q_{n-1}&\hskip-10pt Q_{n-1}\cr
\triangledown &&\hskip-4pt\triangledown &&&\hskip-30pt\triangledown &&\hskip-10pt\triangledown &\hskip-55pt\nearrow\cr
R_{n-2}&\hskip-6pt\cong &\hskip-4ptT_{n-2,1}&\hskip-10pt\cong&\hskip-30pt\dots&\hskip-20ptT_{n-2,n-2}&&\hskip-6ptQ_{n-1}\cr
\triangledown&&\hskip-4pt\triangledown&&&\hskip-30pt\triangledown &&\hskip-12pt\triangledown\cr
R_{n-1}&\hskip-6pt\cong &\hskip-4ptT_{n-1,1}&\hskip-10pt\cong&\hskip-30pt\dots&\hskip-20ptT_{n-1,n-2}
&\cong &\hskip-20pt T_{n-1,n-1}\cr}
\eqno £9.4.2\phantom{.}$$
where  we set ${}^{\theta_i} Q_i = \theta_i (Q_i)$ and we have $Q_i = T_{i,i}\.C_{Q_i} (T_{i,i})$ for any~$i\in \Delta_{n-1}\,.$

  \smallskip
  Setting $T_{n,0} = R_n\,,$ inductively define $T_{n,j+1}   = \theta_j (T_{n,j})$
  for any $j < n-1\,,$  which makes sense since inductively we get $T_{n,j}\i T_{j,j}
  \i Q_j\,;$ moreover, since we have 
  $$Q_{n-1} = T_{n-1,n-1}\.C_{Q_{n-1}}(T_{n-1,n-1})\qq T_{n,n-1}\triangleleft T_{n-1,n-1}
\eqno  £9.4.3,$$
 the group $Q_{n-1,0} = Q_{n-1}$ contains and normalizes 
 $$Q_{n,0} = T_{n,n-1}\.C_{Q_{n-1,0}}(T_{n,n-1})
 \eqno  £9.4.4,$$
 and, in particular, $Q_{n,0}$ is selfcentralizing in $Q_{n-1,0}\,;$
  then, according to [5,~Corollary~2.21], there is an $\F\-$morphism
 $\theta_{n-1,0}\,\colon Q_{n-1,0}\to P$ such that $\theta_{n-1,0} (Q_{n-1,0})$ and  $\theta_{n-1,0} (Q_{n,0})$ are both fully centralized in~$\F\,,$ and we set
 $$\eqalign{Q_{n,1} &= \theta_{n-1,0} (Q_{n,0})\.N_{C_P(\theta_{n-1,0} 
 (Q_{n,0}))}\big(\theta_{n-1,0} (Q_{n-1,0})\big)\cr
 Q_{n-1,1} &= \theta_{n-1,0} (Q_{n-1,0})\.N_{C_P(\theta_{n-1,0} 
 (Q_{n,0}))}\big(\theta_{n-1,0} (Q_{n-1,0})\big)\cr}
 \eqno  £9.4.5;$$
once again,  $Q_{n,1}$ is selfcentralizing in $Q_{n-1,1}\,.$

\smallskip
On the other hand, we denote by $\tau_{i,0}$ the image of $\sigma$ in $\F (R_i)$ and,  for any $i\in \Delta_{n-1}$ and any $j < i\,,$ we inductively denote by   $\tau_{i,j+1}$ the image of $\tau_{i,j}$ in $\F (T_{i,j+1})$ {\it via\/} $\theta_j\,;$ from our conditions,  it is easily checked that  the action $\bar\sigma_{n-2}$ of $\sigma_{n-2}$ over $Q_{n-1}$ stabilizes $T_{n-1,n-1}$ and induces $\tau_{n-1,n-1}$ over this group, so that it stabilizes $T_{n,n-1}$ and  then  it stabilizes~$Q_{n,0}\,;$ thus, the action of $\sigma_{n-2}$ over $\theta_{n-1,0} (Q_{n-1,0})$ defined {\it via\/} $\theta_{n-1,0}$
stabilizes $\theta_{n-1,0} (Q_{n,0})$ and it follows from [5,~statement~2.10]
that this action can be extended to an $\F\-$morphism $Q_{n-1,1}\to P\,;$
but, the image of the group
$$N_{C_P(\theta_{n-1,0}  (Q_{n,0}))}\big(\theta_{n-1,0} (Q_{n-1,0})\big)
 \eqno  £9.4.6\phantom{.}$$ 
 clearly normalizes $\theta_{n-1,0} (Q_{n-1,0})$ and centralizes $\theta_{n-1,0} (Q_{n,0})\,;$ hence, this  $\F\-$morphism determines an $\F\-$automorphism 
 $\sigma_{n-1,0}$ of $Q_{n-1,1}$ which stabilizes $\theta_{n-1,0} (Q_{n-1,0})$ and $Q_{n,1}\,,$ and, since $\sigma_{n-2}$ is a $p'\-$element, we may assume that  
$\sigma_{n-1,0}$ is also a $p'\-$element.

\smallskip
Now, arguing by induction on $j\in \Bbb N\,,$ assume that we have two families 
$\{Q_{n-1,j'}\}_{j'\le j}$ and $\{Q_{n,j'}\}_{j'\le j}$ of subgroups of $P$ such that $Q_{n,j'}$ is a normal and a selfcentralizing subgroup of~$Q_{n-1,j'}\,,$ a family of $\F\-$morphisms 
$$\theta_{n-1,j'} : Q_{n-1,j'}\too P
\eqno  £9.4.7\phantom{.}$$
where $j'$ runs over $\Delta_{j-1}\,,$ such that 
$$\eqalign{\theta_{n-1,j'} (Q_{n-1,j'})\, &\triangleleft\, Q_{n-1,j'+1}\cr
Q_{n-1,j'+1}  &= \theta_{n-1,j'} (Q_{n-1,j'})\.C_{Q_{n-1,j'+1} }\big(\theta_{n,j'} (Q_{n,j'})\big)\cr 
Q_{n,j'+1} &= \theta_{n,j'} (Q_{n,j'})\.C_{Q_{n,j'+1}}\big(\theta_{n,j'} (Q_{n,j'})\big)\cr}
\eqno  £9.4.8,$$
and a family of $p'\-$elements $\sigma_{n-1,j'}\in \F (Q_{n-1,j'+1})\,,$ where 
$j'$ runs over $\Delta_{j-1}\,,$ stabilizing $\theta_{n-1,j'} (Q_{n-1,j'})$ and 
$Q_{n,j'}\,,$ and, if $j'\not= 0\,,$ inducing $\sigma_{n-1,j'-1}$ on~$Q_{n-1,j'}$
{\it via\/} $\theta_{n-1,j'}\,;$ once again, the following diagram summarizes our situation
$$\matrix{&Q_{\bar n}&&\hskip-20pt{}^{\theta_{\bar n,0}}Q_{\bar n,0}\.N_{n,0}& &&\hskip-20pt{}^{\theta_{\bar n,\bar\jmath}}Q_{\bar n,\bar\jmath}\. N_{n,\bar\jmath}\cr
&\Vert&&\hskip-20pt \Vert&&&\hskip-20pt\Vert\cr
&Q_{\bar n,0}&\hskip-20pt\buildrel \theta_{\bar n,0}\over\too &\hskip-20pt Q_{\bar n,1}&\hskip-20pt\too 
&\dots &\hskip-10pt Q_{\bar n,j}\cr
{\bar\sigma_{n-2}\atop}\hskip-4pt\nearrow\hskip-20pt&\hskip-10pt\triangledown &{\sigma_{\bar n,1}\atop}\hskip-4pt\nearrow\hskip-10pt&\hskip-30pt\triangledown &&\hskip-20pt{\sigma_{\bar n,j}\atop}\hskip-4pt\nearrow\hskip-20pt&\hskip-20pt\triangledown \cr
Q_{\bar n,0}\hskip-20pt&\buildrel \theta_{\bar n,0}\over\too &Q_{\bar n,1}
&\hskip-30pt\too &\dots &Q_{\bar n,j}\cr
&Q_{n,0}&\hskip-20pt\too &\hskip-20ptQ_{n,1}& \hskip-20pt\too&\dots&\hskip-10ptQ_{n,j}\cr
\hskip-25pt\triangledown \hskip10pt\nearrow\hskip-60pt&&\hskip-5pt\triangledown \hskip10pt\nearrow\hskip-20pt&&&\triangledown \hskip4pt\nearrow\hskip-10pt\cr
Q_{n,0}\hskip-20pt&\too &Q_{n,1}&\hskip-40pt\too&\dots&Q_{n,j}\cr
\Vert\hskip-20pt&&\Vert&&&\Vert\cr
T_{n,\bar n}\.C_{Q_{\bar n}}(T_{n,\bar n})\hskip-30pt&&{}^{\theta_{\bar n,0}}Q_{n,0}\.N_{n,0}\hskip-20pt&&&\hskip-20pt {}^{\theta_{\bar n,\bar\jmath}}Q_{n,\bar\jmath}\.N_{n,\bar\jmath}\hskip-30pt\cr}
\eqno £9.4.9\phantom{.}$$
where we are setting $\bar n = n-1\,,$ $\bar\jmath = j-1$ and 
$$N_{n,j'} = N_{C_P({}^{\theta_{n-1,j'-1}}Q_{n,j'-1})}({}^{\theta_{n-1,0}}Q_{n-1,j'-1})
\eqno £9.4.10.$$
\eject

\smallskip
As above, by [5,~Corollary~2.21], there is an $\F\-$morphism
 $$\theta_{n-1,j} : Q_{n-1,j}\too P
 \eqno  £9.4.11\phantom{.}$$
  such that $\theta_{n-1,j} (Q_{n-1,j})$ and  $\theta_{n-1,j} (Q_{n,j})$ are both fully centralized in~$\F\,;$ then, we set
 $$\eqalign{Q_{n,j+1} &= \theta_{n-1,j} (Q_{n,j})\.N_{C_P(\theta_{n-1,j} 
 (Q_{n,j}))}\big(\theta_{n-1,j} (Q_{n-1,j})\big)\cr
 Q_{n-1,j+1} &= \theta_{n-1,j} (Q_{n-1,j})\.N_{C_P(\theta_{n-1,j} 
 (Q_{n,j}))}\big(\theta_{n-1,j} (Q_{n-1,j})\big)\cr}
 \eqno  £9.4.12\phantom{.}$$
 and it is clear that $Q_{n,j+1}$ is again a normal and a selfcentralizing subgroup 
 of~$Q_{n-1,j+1}\,;$ similarly, if $j\not= 0\,,$ the action of $\sigma_{n-1,j-1}$ over $\theta_{n-1,j} (Q_{n-1,j})$ defined {\it via\/} $\theta_{n-1,j}$
stabilizes $\theta_{n-1,j} (Q_{n,j})$ and, again from [5,~statement~2.10.1],
this action can be extended to an $\F\-$mor-phism $Q_{n-1,j+1}\to P\,;$
but, the image of the group
$$N_{C_P(\theta_{n-1,j}  (Q_{n,j}))}\big(\theta_{n-1,j} (Q_{n-1,j})\big)
 \eqno  £9.4.13\phantom{.}$$ 
 clearly normalizes $\theta_{n-1,j} (Q_{n-1,j})$ and centralizes $\theta_{n-1,j} 
 (Q_{n,j})\,;$ hence, this  $\F\-$morphism determines an $\F\-$automorphism 
 $\sigma_{n-1,j}$ of $Q_{n-1,j+1}$ which stabilizes $\theta_{n-1,j} (Q_{n-1,j})$ and $Q_{n,j+1}\,,$ and, since $\sigma_{n-1,j-1}$ is a $p'\-$element, we may assume that  
$\sigma_{n-1,j}$ is also a $p'\-$element.

\smallskip
Finally, if we have $Q_{n,j+1} = \theta_{n-1,j} (Q_{n,j})$ then we still have
$$\eqalign{ Q_{n-1,j+1} &= \theta_{n-1,j} (Q_{n-1,j})\cr
& = N_{\theta_{n-1,j} (Q_{n-1,j})\.C_P(\theta_{n-1,j}  (Q_{n,j}))} \big(\theta_{n-1,j} (Q_{n-1,j})\big)\cr}
 \eqno  £9.4.14\phantom{.}$$
 and therefore $\theta_{n-1,j} (Q_{n-1,j})$ contains $C_P 
 \big(\theta_{n-1,j}  (Q_{n,j})\big) \,;$ but, since  $Q_{n,j}$ is selfcentralizing in $Q_{n-1,j}\,,$ $\theta_{n-1,j}  (Q_{n,j})$ is selfcentralizing in $\theta_{n-1,j}  (Q_{n-1,j})$ and thus  $\theta_{n-1,j} (Q_{n,j})$ contains $C_P 
 \big(\theta_{n-1,j}  (Q_{n,j})\big) \,,$ so that  $Q_{n,j+1}$ is $\F\-$self-centralizing;
 in this situation, it suffices to consider $Q_n = Q_{n,j+1}\,,$ to define 
 $$\theta_{n-1} : Q_{n-1}\too P
\eqno  £9.4.15\phantom{.}$$
mapping $u = u_{n-1,0}\in Q_{n-1}$ on $u_{n-1,j+1}$ where we inductively set
$$u_{n-1,j'+1} = \theta_{n-1,0}(u_{n-1,j'})
\eqno  £9.4.16\phantom{.}$$
for any $j'\in \Delta_j\,,$ and to denote by $\sigma_n$ the restriction of 
$\sigma_{n-1,j+1}$ to the product $\theta_{n-1}(Q_{n-1})\.Q_n\,.$

\smallskip
Indeed, since $\theta_{n-1}(Q_{n-1})$ is contained in $Q_{n-1,j+1}\,,$ this group
normalizes $Q_n\,;$ moreover, arguing by induction on $j'\in \Delta_j\,,$ it is easily checked that $\sigma_{n-1,j+1}$ stabilizes $\theta_{n-1}(Q_{n-1})$ and that the 
action over $Q_{n-1}$ defined {\it via\/} $\theta_{n-1}$ coincides with the action
of $\sigma_{n-1}\,;$ similarly, it is clear that for any $j'\in \Delta_j$ we have
$$Q_{n,j'+1} = \theta_{n-1,j'}(Q_{n,j'})\.C_{Q_{n,j'+1}}\big(\theta_{n-1,j'}(Q_{n,j'})\big)
\eqno £9.4.17\phantom{.}$$
and by induction we get $Q_n = \theta_{n-1}(Q_{n,0})\.C_{Q_n}\big(\theta_{n-1}
(Q_{n,0}\big)\,;$ thus, setting $T_{n,n} = \theta_{n-1}(T_{n,n-1})\,,$  it follows from equality~£9.4.4 that
 $$Q_n = T_{n,n}\.C_{Q_n} (T_{n,n})
 \eqno £9.4.18.$$
 We are done.
\bigskip
\bigskip
\noindent
{\bf References}
\medskip
\noindent
{\cds[1]  Richard Brauer, {\cdt On blocks and sections in finite groups, I and II\/}, Amer. J. Math. 89(1967), 1115-1136,
90(1968), 895-925
\smallskip
\noindent
[2] Michel Brou\'e, {\cdt Radical, hauteurs, p-sections et blocs\/}, Ann. of Math. 107(1978), 89-107
\smallskip
\noindent
[3] Daniel Gorenstein, {\cdt ``Finite Groups''\/}, Harper's Series, 1968, Harper and Row
\smallskip
\noindent
[4] Llu\'\i s Puig, {\cdt Pointed groups and  construction of modules}, Journal of Algebra, 116(1988), 7-129
\smallskip
\noindent
[5] Llu\'\i s Puig, {\cdt ``Frobenius Categories versus Brauer Blocks''\/}, Progress in Math. 274, 2009, Birkh\"auser, Basel
\smallskip
\noindent
[6] Llu\'\i s Puig, {\cdt Block source algebras in p-solvable groups\/}, Michigan Math. J. 58(2009), 323-328
\bigskip
\noindent
{\bf Author's address

\medskip
\noindent
{\cds CNRS, Institut de Math\'ematiques de Jussieu}
\smallskip
\noindent
{\cds 6 Av Bizet, 94340 Joinville-le-Pont, France}
\smallskip
\noindent
{\cds puig@math.jussieu.fr}

\end